С. В. Курапов
М. В. Давидовский

# АЛГОРИТМИЧЕСКИЕ МЕТОДЫ КОНЕЧНЫХ ДИСКРЕТНЫХ СТРУКТУР

# ТОПОЛОГИЧЕСКИЙ РИСУНОК ГРАФА

## часть 3

(на правах рукописи)

**Курапов С.В.**

**Давидовский М.В.**

В работе рассматриваются математические модели для создания топологического рисунка графа, основанные на методах теории вращения вершина Г. Рингеля. Представлен алгоритм генерации топологического рисунка плоской части графа на основе выделения базиса подпространства циклов $C(G)$ методом Монте-Карло. Описан метод наискорейшего спуска для построения топологического рисунка плоского суграфа. Построение топологического рисунка графа осуществляется с помощью комбинации разработанных Л. И. Раппортом методов векторной алгебры пересечений векторов. Представлены три этапа построения плоского суграфа несепарабельного графа.

Рассматриваются вопросы построения гамильтонова цикла на базе построения плоского суграфа. Описан способ построения гамильтонова цикла графа на основе графа циклов плоского суграфа.

Для научных работников, студентов и аспирантов, специализирующихся на применении методов прикладной математики.

# ОГЛАВЛЕНИЕ

# Введение

Теория вращения вершин [31] и факт существования множества изометрических циклов в несепарабельном графе, позволяют построить алгебраические методы теории графов для описания топологического рисунка графа, не прибегая к промежуточным метрическим представлениям на плоскости.

Анализ результатов генерации плоских конфигураций, полученных с помощью программной реализации алгебраических методов построения топологического рисунка графа, показали, что применения одного критерия Маклейна для построения плоского топологического рисунка графа недостаточно. Нужно осуществлять проверку суграфа на связность и замкнутость циклического порядка вращения вершин.

Для оптимизации процесса вычисления базиса подпространства циклов графа C(G), состоящего из изометрических циклов, представлен метод наискорейшего спуска, имеющий полиномиальную вычислительную сложность.

На основе построения топологического рисунка плоского суграфа представлен метод построения гамильтонова цикла графа. Построение гамильтонова цикла графа основано на выделении подмножества простых циклов, удовлетворяющего нулевому значению кубического функционала Маклейна [16,20,24].

# Глава 11. ГЕНЕРАЦИЯ БАЗИСОВ ПОДПРОСТРАНСТВА ЦИКЛОВ

## 11.1. Основные понятия и определения

**Определение 11.1.** *Несепарабельный граф* – связный неориентированный граф без мостов и точек сопряжения, без петель и кратных ребер, без вершин с валентностью меньшей или равной двум.

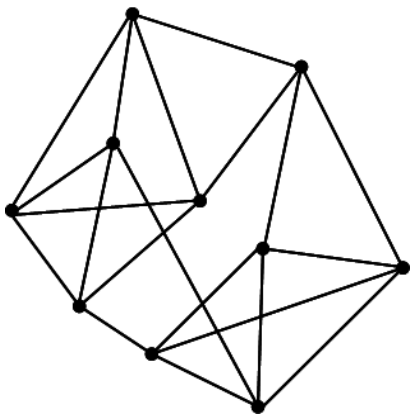

Рис. 11.1. Несепарабельный граф.

**Определение 11.2.** *Изометрическим циклом* в графе называется простой цикл, для которого кратчайший путь между любыми двумя его вершинами состоит из рёбер этого цикла. Изометрический цикл – частный случай изометрического подграфа [18,19].

**Определение 11.3.** *Множество изометрических циклов* графа является инвариатом графа и обозначается $C_\tau$ и представляет собой кортеж

$C_\tau = <c_1, c_2, c_3, \ldots, c_k>$, где $k$ – количество циклов в $C_\tau$.

Рассмотрим несепарабельный граф $G_1$.

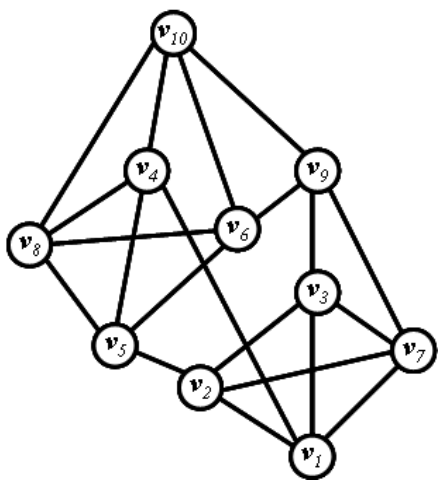

Рис. 11.2. Несепарабельный граф $G_1$.

Количество вершин графа = 10.
Количество ребер графа = 20.
Количество изометрических циклов = 16.

Смежность графа:

вершина $v_1$: $v_2$ $v_3$ $v_4$ $v_7$
вершина $v_2$: $v_1$ $v_3$ $v_5$ $v_7$
вершина $v_3$: $v_1$ $v_2$ $v_7$ $v_9$
вершина $v_4$: $v_1$ $v_5$ $v_8$ $v_{10}$
вершина $v_5$: $v_2$ $v_4$ $v_6$ $v_8$
вершина $v_6$: $v_5$ $v_8$ $v_9$ $v_{10}$
вершина $v_7$: $v_1$ $v_2$ $v_3$ $v_9$
вершина $v_8$: $v_4$ $v_5$ $v_6$ $v_{10}$
вершина $v_9$: $v_3$ $v_6$ $v_7$ $v_{10}$
вершина $v_{10}$: $v_4$ $v_6$ $v_8$ $v_9$

Инцидентность графа:

ребро $e_1$: $(v_1,v_2)$ или $(v_2,v_1)$;   ребро $e_2$: $(v_1,v_3)$ или $(v_3,v_1)$;
ребро $e_3$: $(v_1,v_4)$ или $(v_4,v_1)$;   ребро $e_4$: $(v_1,v_7)$ или $(v_7,v_1)$;
ребро $e_5$: $(v_2,v_3)$ или $(v_3,v_2)$;   ребро $e_6$: $(v_2,v_5)$ или $(v_5,v_2)$;
ребро $e_7$: $(v_2,v_7)$ или $(v_7,v_2)$;   ребро $e_8$: $(v_3,v_7)$ или $(v_7,v_3)$;
ребро $e_9$: $(v_3,v_9)$ или $(v_9,v_3)$;   ребро $e_{10}$: $(v_4,v_5)$ или $(v_5,v_4)$;
ребро $e_{11}$: $(v_4,v_8)$ или $(v_8,v_4)$;   ребро $e_{12}$: $(v_4,v_{10})$ или $(v_{10},v_4)$;
ребро $e_{13}$: $(v_5,v_6)$ или $(v_6,v_5)$;   ребро $e_{14}$: $(v_5,v_8)$ или $(v_8,v_5)$;
ребро $e_{15}$: $(v_6,v_8)$ или $(v_8,v_6)$;   ребро $e_{16}$: $(v_6,v_9)$ или $(v_9,v_6)$;
ребро $e_{17}$: $(v_6,v_{10})$ или $(v_{10},v_6)$;   ребро $e_{18}$: $(v_7,v_9)$ или $(v_9,v_7)$;
ребро $e_{19}$: $(v_8,v_{10})$ или $(v_{10},v_8)$;   ребро $e_{20}$: $(v_9,v_{10})$ или $(v_{10},v_9)$.

Множество изометрических циклов графа $G_1$:

цикл $c_1 = \{e_1,e_2,e_5\} \leftrightarrow \{v_1,v_2,v_3\}$;
цикл $c_2 = \{e_1,e_3,e_6,e_{10}\} \leftrightarrow \{v_1,v_2,v_4,v_5\}$;
цикл $c_3 = \{e_1,e_4,e_7\} \leftrightarrow \{v_1,v_2,v_7\}$;
цикл $c_4 = \{e_2,e_4,e_8\} \leftrightarrow \{v_1,v_3,v_7\}$;
цикл $c_5 = \{e_2,e_3,e_9,e_{12},e_{20}\} \leftrightarrow \{v_1,v_3,v_4,v_9,v_{10}\}$;
цикл $c_6 = \{e_3,e_4,e_{12},e_{18},e_{20}\} \leftrightarrow \{v_1,v_4,v_7,v_{10},v_9\}$;
цикл $c_7 = \{e_5,e_7,e_8\} \leftrightarrow \{v_2,v_3,v_7\}$;
цикл $c_8 = \{e_5,e_6,e_9,e_{13},e_{16}\} \leftrightarrow \{v_2,v_3,v_5,v_9,v_6\}$;
цикл $c_9 = \{e_6,e_7,e_{13},e_{16},e_{18}\} \leftrightarrow \{v_2,v_5,v_7,v_6,v_9\}$;   (11.1)
цикл $c_{10} = \{e_8,e_9,e_{18}\} \leftrightarrow \{v_3,v_7,v_9\}$;
цикл $c_{11} = \{e_{10},e_{11},e_{14}\} \leftrightarrow \{v_4,v_5,v_8\}$;
цикл $c_{12} = \{e_{10},e_{12},e_{13},e_{17}\} \leftrightarrow \{v_4,v_5,v_{10},v_6\}$;
цикл $c_{13} = \{e_{11},e_{12},e_{19}\} \leftrightarrow \{v_4,v_8,v_{10}\}$;
цикл $c_{14} = \{e_{13},e_{14},e_{15}\} \leftrightarrow \{v_5,v_6,v_8\}$;
цикл $c_{15} = \{e_{15},e_{17},e_{19}\} \leftrightarrow \{v_6,v_8,v_{10}\}$;
цикл $c_{16} = \{e_{16},e_{17},e_{20}\} \leftrightarrow \{v_6,v_9,v_{10}\}$.

**Определение 11.4.** *Вектором циклов по ребру* называется кортеж количества циклов, проходящих по этому ребру, для данного подмножества циклов $Y(G)$ (суграфа) [20].

Множество изометрических циклов графа $G_1$ характеризуется следующим вектором $P_e$:

$P_e = <p_1.p_2.p_3,p_4,p_5,p_6,p_7,p_8,p_9,p_{10},p_{11},p_{12},p_{13},p_{14},p_{15},p_{16},p_{17},p_{18},p_{19}.p_{20}> =$

= <3,3,3,3,3,3,3,3,3,3,2,4,4,2,2,3,3,3,2,3>

Подмножество циклов может характеризоваться значением квадратичного функционала Маклейна [20]:

$$F_1(Y) = \sum_{i=1}^{m} p_i^2 - 3\sum_{i=1}^{m} p_i + 2m \qquad (11.2)$$

или значением кубического функционала Маклейна (также называемым функционалом Понтрягина-Куратовского)

$$F_2(Y) = \sum_{i=1}^{m} p_i^3 - 3\sum_{i=1}^{m} p_i^2 + 2\sum_{i=1}^{m} p_i, \qquad (11.3)$$

$p_i$ – количество циклов, проходящих по ребру $e_i$.

Для множества изометрических циклов графа $G_1$ значение квадратичного функционала Маклейна равно 40, а значение кубического функционала Маклейна равно 132.

**Определение 11.5.** *Вращение вершины* – это циклический порядок обхода инцидентных ребер (представляется диаграммой вращения) [31].

**Определение 11.6.** *Топологический рисунок* плоского несепарабельного графа – это система вращений вершин, индуцирующая (порождающая) систему простых циклов, или система простых циклов, индуцирующая вращение вершин [20].

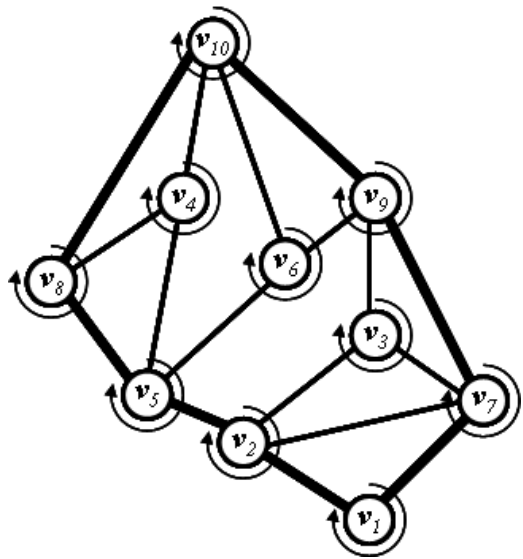

Рис. 11.3. Топологический рисунок плоской части графа.

Вращение вершин в топологическом рисунке:

$\sigma(v_1)$:  $v_2$  $v_7$
$\sigma(v_2)$:  $v_1$  $v_5$  $v_3$  $v_7$
$\sigma(v_3)$:  $v_2$  $v_9$  $v_7$
$\sigma(v_4)$:  $v_8$  $v_{10}$  $v_5$
$\sigma(v_5)$:  $v_8$  $v_4$  $v_6$  $v_2$
$\sigma(v_6)$:  $v_9$  $v_5$  $v_{10}$

$\sigma(v_7)$:    $v_1$   $v_2$   $v_3$   $v_9$

$\sigma(v_8)$:    $v_4$   $v_5$   $v_{10}$

$\sigma(v_9)$:    $v_7$   $v_3$   $v_6$   $v_{10}$

$\sigma(v_{10})$:    $v_9$   $v_6$   $v_4$   $v_8$

Индуцируемая врашением вершин система циклов имеет следующий вид:

$c_3 = \{e_1,e_4,e_7\} \leftrightarrow <v_1,v_7,v_2> \leftrightarrow (v_1,v_7)+(v_7,v_2)+(v_2,v_1)$;

$c_7 = \{e_5,e_7,e_8\} \leftrightarrow <v_2,v_7,v_3> \leftrightarrow (v_2,v_7)+(v_7,v_3)+(v_3,v_2)$;

$c_8 = \{e_5,e_6,e_9,e_{13},e_{16}\} \leftrightarrow <v_2,v_3,v_9,v_6,v_5> \leftrightarrow$

$\leftrightarrow (v_2,v_3)+(v_3,v_9)+(v_9,v_6)+(v_6,v_5)+(v_5,v_2)$;

$c_{10} = \{e_8,e_9,e_{18}\} \leftrightarrow <v_3,v_7,v_9> \leftrightarrow (v_3,v_7)+(v_7,v_9)+(v_9,v_3)$;           (11.4)

$c_{11} = \{e_{10},e_{11},e_{14}\} \leftrightarrow <v_4,v_8,v_5> \leftrightarrow (v_4,v_8)+(v_8,v_5)+(v_5,v_4)$;

$c_{12} = \{e_{10},e_{12},e_{13},e_{17}\} \leftrightarrow <v_4,v_5,v_6,v_{10}> \leftrightarrow (v_4,v_5)+(v_5,v_6)+(v_6,v_{10})+(v_{10},v_4)$;

$c_{13} = \{e_{11},e_{12},e_{19}\} \leftrightarrow \{v_4,v_{10},v_8\} \leftrightarrow (v_4,v_{10})+(v_{10},v_8)+(v_8,v_4)$;

$c_{16} = \{e_{16},e_{17},e_{20}\} \leftrightarrow <v_6,v_9,v_{10}> \leftrightarrow (v_6,v_9)+(v_9,v_{10})+(v_{10},v_6)$

В топологическом рисунке цикл может быть записан не только в виде множества ребер и вершин, но и в виде циклического порядка ориентированных ребер. Такая запись называется векторной записью цикла:

$c_{12} = (v_4,v_5) + (v_5,v_6) + (v_6,v_{10}) + (v_{10},v_4)$.

**Определение 11.7.** *Ободом* называется кольцевая сумма полмножества Y(G) простых циклов графа.

Например, обод плоского суграфа, представленного на рис. 11.3, имеет вид:

$c_0 = c_3 \oplus c_7 \oplus c_8 \oplus c_{10} \oplus c_{11} \oplus c_{12} \oplus c_{13} \oplus c_{16} = \{e_1,e_4,e_6,e_{14},e_{18},e_{19},e_{20}\} \leftrightarrow$

$\leftrightarrow <v_1,v_2,v_5,v_8,v_{10},v_9,v_7,v_1> \leftrightarrow (v_1,v_2)+(v_2,v_5)+(v_5,v_8)+(v_8,v_{10})+(v_{10},v_{90})+(v_9,v_7)+(v_7,v_1)$.

Кольцевая сумма независимых циклов и обод в топологическом рисунке суграфа есть пустое множество [32].

$$\sum_{i=1}^{k} c_i \oplus c_0 = \varnothing \qquad (11.5)$$

**Определение 11.8.** *Топологический рисунок* несепарабельного графа с пересекающимися ребрами – это система вращений основных и мнимых вершин, индуцирующая систему простых циклов.

Здесь мнимая вершина определяет и характеризует местоположение пересечения двух ребер несепарабельного графа. Диаграмма вращения вершин для топологического рисунка графа с пересечениями состоит из двух частей – диаграммы основных вершин (вершин принадлежащих графу) и диаграммы мнимых вершин.

Вращение вершин в топологическом рисунке графа с пересечением ребер (мнимые вершины окрашены алым цветом) [31]:

$\sigma(v_1)$:    $v_2$   $v_{11}$   $v_7$   $v_{13}$

$\sigma(v_2):$   $v_1$   $v_5$   $v_3$   $v_{11}$
$\sigma(v_3):$   $v_2$   $v_9$   $v_7$   $v_{11}$
$\sigma(v_4):$   $v_8$   $v_{13}$   $v_{10}$   $v_6$
$\sigma(v_5):$   $v_8$   $v_{12}$   $v_6$   $v_2$
$\sigma(v_6):$   $v_9$   $v_5$   $v_{12}$   $v_{10}$
$\sigma(v_7):$   $v_1$   $v_{11}$   $v_3$   $v_9$
$\sigma(v_8):$   $v_4$   $v_{12}$   $v_5$   $v_{13}$
$\sigma(v_9):$   $v_7$   $v_3$   $v_6$   $v_{10}$
$\sigma(v_{10}):$   $v_9$   $v_6$   $v_4$   $v_{13}$
$\sigma(v_{11}):$   $v_1$   $v_2$   $v_3$   $v_7$
$\sigma(v_{12}):$   $v_5$   $v_8$   $v_4$   $v_6$
$\sigma(v_{13}):$   $v_8$   $v_1$   $v_{10}$   $v_5$

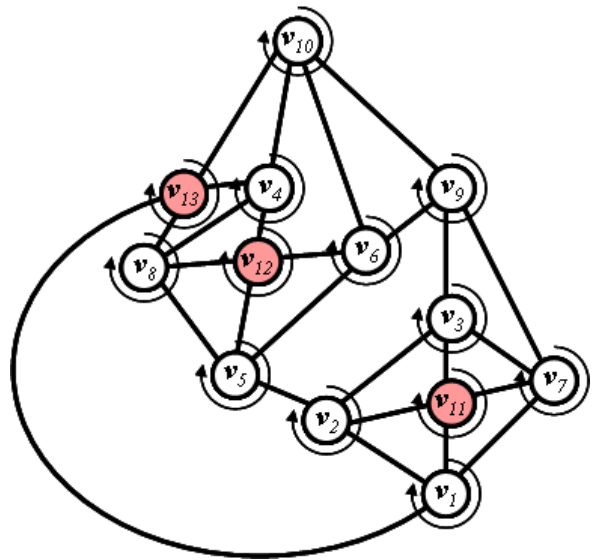

Рис. 11.4. Топологический рисунок суграфа с пересечением ребер.

Индуцированная вращением вершин система циклов имеет вид:

$c_1 = (v_1,v_2)+(v_2,v_5)+(v_5,v_8)+(v_8,v_{13})+(v_{13},v_1);$
$c_2 = (v_8,v_4)+(v_4,v_{13})+(v_{13},v_8);$
$c_3 = (v_8,v_5)+(v_5,v_{12})+(v_{12},v_8);$
$c_4 = (v_5,v_2)+(v_2,v_3)+(v_3,v_9)+(v_9,v_6)+(v_6,v_5);$
$c_5 = (v_2,v_1)+(v_1,v_{11})+(v_{11},v_2);$
$c_6 = (v_4,v_8)+(v_8,v_{12})+(v_{12},v_4);$
$c_7 = (v_5,v_6)+(v_6,v_{12})+(v_{12},v_5);$
$c_8 = (v_7,v_3)+(v_3,v_{11})+(v_{11},v_7);$ (11.6)
$c_9 = (v_1,v_7)+(v_7,v_{11})+(v_{11},v_1);$
$c_{10} = (v_4,v_{10})+(v_{10},v_{13})+(v_{13},v_4);$
$c_{11} = (v_6,v_{10})+(v_{10},v_4)+(v_4,v_{12})+(v_{12},v_6);$
$c_{12} = (v_6,v_9)+(v_9,v_{10})+(v_{10},v_6);$
$c_{13} = (v_3,v_2)+(v_2,v_{11})+(v_{11},v_3);$
$c_{14} = (v_3,v_7)+(v_7,v_9)+(v_9,v_3);$
$c_0 = (v_9,v_7)+(v_7,v_1)+(v_1,v_{13})+(v_{13},v_{10})+(v_{10},v_9)$

**Определение 11.9.** *Базис изометрических циклов* – это линейно независимое подмножество изометрических циклов с мощностью равной цикломатическому числу графа [4-6,13,14,20]:

$$gard(\upsilon(G)) = m - n + 1. \tag{11.7}$$

Например, базис изометричеких циклов графа $G_1$ имеет следубщий вид:

$$\begin{aligned}
&\text{цикл } c_1 = \{e_1, e_2, e_5\} \leftrightarrow \{v_1, v_2, v_3\}; \\
&\text{цикл } c_2 = \{e_1, e_3, e_6, e_{10}\} \leftrightarrow \{v_1, v_2, v_4, v_5\}; \\
&\text{цикл } c_3 = \{e_1, e_4, e_7\} \leftrightarrow \{v_1, v_2, v_7\}; \\
&\text{цикл } c_7 = \{e_5, e_7, e_8\} \leftrightarrow \{v_2, v_3, v_7\}; \\
&\text{цикл } c_8 = \{e_5, e_6, e_9, e_{13}, e_{16}\} \leftrightarrow \{v_2, v_3, v_5, v_9, v_6\}; \\
&\text{цикл } c_{10} = \{e_8, e_9, e_{18}\} \leftrightarrow \{v_3, v_7, v_9\}; \\
&\text{цикл } c_{11} = \{e_{10}, e_{11}, e_{14}\} \leftrightarrow \{v_4, v_5, v_8\}; \\
&\text{цикл } c_{12} = \{e_{10}, e_{12}, e_{13}, e_{17}\} \leftrightarrow \{v_4, v_5, v_{10}, v_6\}; \\
&\text{цикл } c_{13} = \{e_{11}, e_{12}, e_{19}\} \leftrightarrow \{v_4, v_8, v_{10}\}; \\
&\text{цикл } c_{14} = \{e_{13}, e_{14}, e_{15}\} \leftrightarrow \{v_5, v_6, v_8\}; \\
&\text{цикл } c_{16} = \{e_{16}, e_{17}, e_{20}\} \leftrightarrow \{v_6, v_9, v_{10}\}
\end{aligned} \tag{11.8}$$

Вектор циклов по ребрам для базиса:

$F_3(B) = 24, \quad P_e(B) = <3,1,1,1,3,2,2,2,2,3,2,2,3,2,1,2,2,1,1,1>$

**Определение 11.10**. *Вычислимое подмножество циклов* – это подмножество изометрических циклов, состоящее из элементов базиса, позволяющее производить удаление цикла с уменьшением значения функционала Маклейна и одновременным выполнением условия Эйлера (удаление цикла должно приводить к удалению одного, и только одного ребра).

Например, следующий процесс удаления циклов из базиса, приводящий к построению плоского суграфа, описывается следующей последовательностью вычислимых подмножеств:

$$\frac{\partial F_3(B)}{\partial c_1} = 12, \quad P_e = <2,0,1,1,2,2,2,2,2,3,2,2,3,2,1,2,2,1,1,1>$$

$$\frac{\partial F_3^2(B)}{\partial c_1 \partial c_2} = 6, \quad P_e = <1,0,0,1,2,1,2,2,2,2,2,2,3,2,1,2,2,1,1,1>$$

$$\frac{\partial F_3^3(B)}{\partial c_1 \partial c_2 \partial c_{14}} = 0, \quad P_e = <1,0,0,1,2,1,2,2,2,2,2,2,1,0,2,2,1,1,1>$$

Вычислимое подмножество циклов может быть описано как дифференциация структурного числа базиса B с указанием удаленных циклов [2]. Цифровой индекс при функционале Маклейна $F_3(B)$ характеризует кубический функционал. $F_2(B)$ характеризует квадратичный функционал Маклейна.

Подмножество циклов, описываемое вектором циклов по ребрам без единиц, считается *не вычислимым*.

**Определение 11.11.** *Плоский суграф* – это подмножество простых циклов с нулевым значением функционала Маклейна, характеризующее суграф, построенный на множестве всех вершин графа.

**Определение 11.12.** *Максимально плоский суграф* – это плоский суграф с минимальным количеством удаленных ребер.

**Определение 11.13.** *Проская конфигурация* – это подмножество циклов с нулевым значением функционала Маклейна [21].

Плоский суграф характеризуется количеством удаленных ребер графа.

Следует заметить, что может существовать несколько топологических рисунков плоского суграфа с равным количеством удаленных ребер. Например, на рис. 11.3 и рис. 11.5 представлены плоские суграфы с одинаковым количеством удаленных ребер.

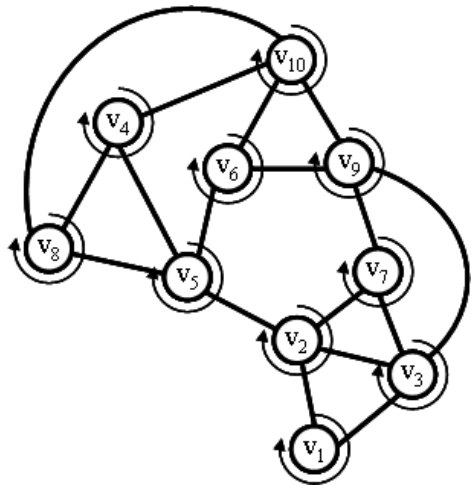

Рис. 11.5. Топологический рисунок суграфа $G_1$.

Топологический рисунок плоского суграфа на рис. 11.5 порожден следующий системой циклов:

$$\begin{aligned}
&\text{цикл } c_1 = \{e_1,e_2,e_5\} \leftrightarrow \{v_1,v_2,v_3\}; \\
&\text{цикл } c_7 = \{e_5,e_7,e_8\} \leftrightarrow \{v_2,v_3,v_7\}; \\
&\text{цикл } c_9 = \{e_6,e_7,e_{13},e_{16},e_{18}\} \leftrightarrow \{v_2,v_5,v_7,v_6,v_9\}; \\
&\text{цикл } c_{10} = \{e_8,e_9,e_{18}\} \leftrightarrow \{v_3,v_7,v_9\}; \\
&\text{цикл } c_{11} = \{e_{10},e_{11},e_{14}\} \leftrightarrow \{v_4,v_5,v_8\}; \\
&\text{цикл } c_{12} = \{e_{10},e_{12},e_{13},e_{17}\} \leftrightarrow \{v_4,v_5,v_{10},v_6\}; \\
&\text{цикл } c_{13} = \{e_{11},e_{12},e_{19}\} \leftrightarrow \{v_4,v_8,v_{10}\}; \\
&\text{цикл } c_{16} = \{e_{16},e_{17},e_{20}\} \leftrightarrow \{v_6,v_9,v_{10}\}
\end{aligned} \quad (11.9)$$

Следующая система изометрических циклов характеризует плоский суграф с одним удаленным ребром (рис. 11.6).

$$\begin{aligned}
&\text{цикл } c_3 = \{e_1,e_4,e_7\} \leftrightarrow \{v_1,v_2,v_7\}; \\
&\text{цикл } c_4 = \{e_2,e_4,e_8\} \leftrightarrow \{v_1,v_3,v_7\}; \\
&\text{цикл } c_5 = \{e_2,e_3,e_9,e_{12},e_{20}\} \leftrightarrow \{v_1,v_3,v_4,v_9,v_{10}\}; \\
&\text{цикл } c_9 = \{e_6,e_7,e_{13},e_{16},e_{18}\} \leftrightarrow \{v_2,v_5,v_7,v_6,v_9\}; \\
&\text{цикл } c_{10} = \{e_8,e_9,e_{18}\} \leftrightarrow \{v_3,v_7,v_9\}; \\
&\text{цикл } c_{11} = \{e_{10},e_{11},e_{14}\} \leftrightarrow \{v_4,v_5,v_8\};
\end{aligned} \quad (11.10)$$

цикл $c_{13} = \{e_{11}, e_{12}, e_{19}\} \leftrightarrow \{v_4, v_8, v_{10}\}$;
цикл $c_{14} = \{e_{13}, e_{14}, e_{15}\} \leftrightarrow \{v_5, v_6, v_8\}$;
цикл $c_{15} = \{e_{15}, e_{17}, e_{19}\} \leftrightarrow \{v_6, v_8, v_{10}\}$;
цикл $c_{16} = \{e_{16}, e_{17}, e_{20}\} \leftrightarrow \{v_6, v_9, v_{10}\}$

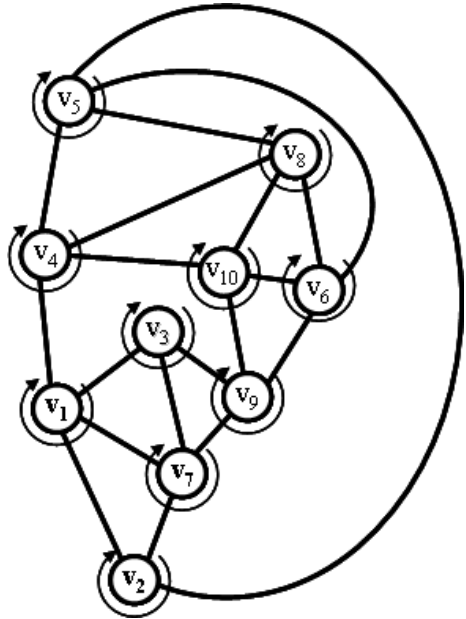

Рис. 11.6. Топологический рисунок плоского суграфа $G_1$ с одним удаленным ребром.

Топологический рисунок, показанный на рис. 11.6, является максимально плоским суграфом (планарность суграфа можно проверить алгоритмом Хопкрофта-Тарьяна) [24].

### 11.2. Генерация базисов изометрических циклов

Мощность множества изометрических циклов всегда больше либо равна мощности независимой системы изометрических циклов равной цикломатическому числу графа $\vartheta(G)$.

Количество таких независимых подмножеств $C_н$ изометрических циклов определяется числом сочетаний

$$|C_н| = |C_\tau|!/|C_\tau - \vartheta(G)|!\,\vartheta(G)! - |C_з| \qquad (11.11)$$

за исключением количества зависимых подмножеств $C_з$.

Для нахождения независимых подмножеств будем применять метод Монте-Карло. Воспользуемся модифицированным алгоритмом Гаусса для случайной перестановки циклов.

Рассмотрим алгоритм [1,3,8-10] выделения базиса изометрических циклов основанный на случайной последовательности расположения изометрических циклов (метод Монте-Карло) [21].

**Шаг 1**. Случайным образом образуем последовательность изометрических циклов графа.

**Шаг 2**. Используя модифицированный алгоритм Гаусса для определения ранга матрицы, выделяем независимую систему изометрических циклов (базис). Мощность базиса определяется цикломатическим числом графа.

**Шаг 3**. Из выделенной независимой системы циклов последовательно удаляются циклы согласно правилу: удаление одного цикла приводит к удалению одного и только одного ребра (для удовлетворения закона Эйлера [30,33]). Удаление циклов прекращается по достижению нулевого значения кубического функционала Маклейна.

Рассмотрим данный алгоритм на примере следующего графа (рис. 11.7).

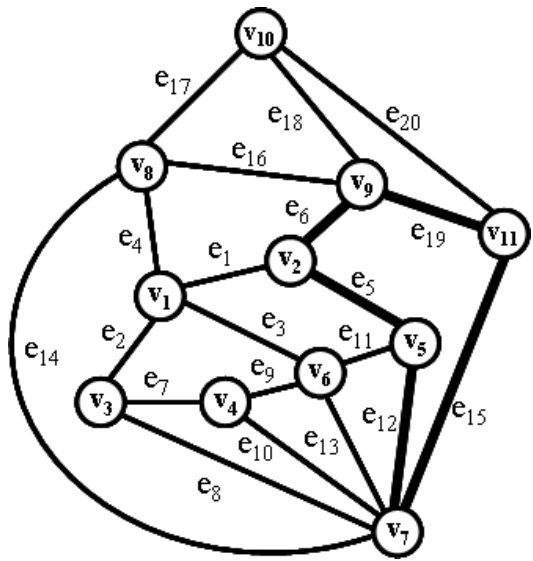

Рис. 11.7. Граф $G_2$ и обод графа.

Количество вершин графа = 11
Количество ребер графа = 20
Количество изометрических циклов = 17

Матрица смежностей графа:
вершина $v_1$: $v_2$ $v_3$ $v_6$ $v_8$
вершина $v_2$: $v_1$ $v_5$ $v_9$
вершина $v_3$: $v_1$ $v_4$ $v_7$
вершина $v_4$: $v_3$ $v_6$ $v_7$
вершина $v_5$: $v_2$ $v_6$ $v_7$
вершина $v_6$: $v_1$ $v_4$ $v_5$ $v_7$
вершина $v_7$: $v_3$ $v_4$ $v_5$ $v_6$ $v_8$ $v_{11}$
вершина $v_8$: $v_1$ $v_7$ $v_9$ $v_{10}$
вершина $v_9$: $v_2$ $v_8$ $v_{10}$ $v_{11}$
вершина $v_{10}$: $v_8$ $v_9$ $v_{11}$
вершина $v_{11}$: $v_7$ $v_9$ $v_{10}$

Элементы матрицы инциденций:
ребро $e_1$: $(v_1,v_2)$ или $(v_2,v_1)$;  ребро $e_2$: $(v_1,v_3)$ или $(v_3,v_1)$;
ребро $e_3$: $(v_1,v_6)$ или $(v_6,v_1)$;  ребро $e_4$: $(v_1,v_8)$ или $(v_8,v_1)$;
ребро $e_5$: $(v_2,v_5)$ или $(v_5,v_2)$;  ребро $e_6$: $(v_2,v_9)$ или $(v_9,v_2)$;
ребро $e_7$: $(v_3,v_4)$ или $(v_4,v_3)$;  ребро $e_8$: $(v_3,v_7)$ или $(v_7,v_3)$;
ребро $e_9$: $(v_4,v_6)$ или $(v_6,v_4)$;  ребро $e_{10}$: $(v_4,v_7)$ или $(v_7,v_4)$;
ребро $e_{11}$: $(v_5,v_6)$ или $(v_6,v_5)$;  ребро $e_{12}$: $(v_5,v_7)$ или $(v_7,v_5)$;
ребро $e_{13}$: $(v_6,v_7)$ или $(v_7,v_6)$;  ребро $e_{14}$: $(v_7,v_8)$ или $(v_8,v_7)$;
ребро $e_{15}$: $(v_7,v_{11})$ или $(v_{11},v_7)$;  ребро $e_{16}$: $(v_8,v_9)$ или $(v_9,v_8)$;
ребро $e_{17}$: $(v_8,v_{10})$ или $(v_{10},v_8)$;  ребро $e_{18}$: $(v_9,v_{10})$ или $(v_{10},v_9)$;

ребро $e_{19}$: $(v_9,v_{11})$ или $(v_{11},v_9)$;   ребро $e_{20}$: $(v_{10},v_{11})$ или $(v_{11},v_{10})$.

Множество изометрических циклов графа:

цикл $c_1 = \{e_1,e_3,e_5,e_{11}\} \leftrightarrow \{v_1,v_2,v_5,v_6\}$;
цикл $c_2 = \{e_1,e_4,e_6,e_{16}\} \leftrightarrow \{v_1,v_2,v_8,v_9\}$;
цикл $c_3 = \{e_2,e_3,e_7,e_9\} \leftrightarrow \{v_1,v_3,v_4,v_6\}$;
цикл $c_4 = \{e_2,e_3,e_8,e_{13}\} \leftrightarrow \{v_1,v_3,v_6,v_7\}$;
цикл $c_5 = \{e_2,e_4,e_8,e_{14}\} \leftrightarrow \{v_1,v_3,v_7,v_8\}$;
цикл $c_6 = \{e_3,e_4,e_{13},e_{14}\} \leftrightarrow \{v_1,v_6,v_7,v_8\}$;
цикл $c_7 = \{e_5,e_6,e_{12},e_{14},e_{16}\} \leftrightarrow \{v_2,v_5,v_7,v_8,v_9\}$;
цикл $c_8 = \{e_5,e_6,e_{12},e_{15},e_{19}\} \leftrightarrow \{v_2,v_5,v_7,v_9,v_{11}\}$;
цикл $c_9 = \{e_7,e_8,e_{10}\} \leftrightarrow \{v_3,v_4,v_7\}$;
цикл $c_{10} = \{e_9,e_{10},e_{13}\} \leftrightarrow \{v_4,v_6,v_7\}$;
цикл $c_{11} = \{e_{11},e_{12},e_{13}\} \leftrightarrow \{v_5,v_6,v_7\}$;
цикл $c_{12} = \{e_1,e_2,e_5,e_8,e_{12}\} \leftrightarrow \{v_1,v_2,v_3,v_5,v_7\}$;
цикл $c_{13} = \{e_1,e_4,e_5,e_{12},e_{14}\} \leftrightarrow \{v_1,v_2,v_5,v_7,v_8\}$;
цикл $c_{14} = \{e_{14},e_{15},e_{16},e_{19}\} \leftrightarrow \{v_7,v_8,v_9,v_{11}\}$;
цикл $c_{15} = \{e_{14},e_{15},e_{17},e_{20}\} \leftrightarrow \{v_7,v_8,v_{10},v_{11}\}$;
цикл $c_{16} = \{e_{16},e_{17},e_{18}\} \leftrightarrow \{v_8,v_9,v_{10}\}$;
цикл $c_{17} = \{e_{18},e_{19},e_{20}\} \leftrightarrow \{v_9,v_{10},v_{11}\}$.

Цикломатическое число графа = 10

Номер $\langle 1,2,3,4,5,6,7,8,9,0,1,2,3,4,5,6,7,8,9,0\rangle$
Вектор количества циклов по ребрам    $P_e = \langle 4,4,4,4,5,3,2,4,2,2,2.5,4,6,3,4,2,2,3,2\rangle$

Номер $\langle 1,2,3,4,5,6,7,8,9,0,1\rangle$
Вектор количества циклов по вершинам $P_v = \langle 8,6,5,3,6,6,12,7,5,3,4\rangle$

Построим следующую последовательность изометрических циклов для графа $G_2$:

$C_\tau = \langle \mathbf{c_{17}}, \mathbf{c_2}, \mathbf{c_7}, \mathbf{c_6}, \mathbf{c_8}, \mathbf{c_{12}}, \mathbf{c_{15}}, c_{14}, \mathbf{c_{11}}, c_5, \mathbf{c_{10}}, c_1, c_{13}, \mathbf{c_3}, c_{16}, c_9, c_4 \rangle$.

Независимая система циклов, определенная модифицированным алгоритмом Гаусса, обозначена в последовательности красным цветом.

Базис изометрических циклов:

цикл $c_2 = \{e_1,e_4,e_6,e_{16}\} \leftrightarrow \{v_1,v_2,v_8,v_9\}$;
цикл $c_3 = \{e_2,e_3,e_7,e_9\} \leftrightarrow \{v_1,v_3,v_4,v_6\}$;
цикл $c_6 = \{e_3,e_4,e_{13},e_{14}\} \leftrightarrow \{v_1,v_6,v_7,v_8\}$;
цикл $c_7 = \{e_5,e_6,e_{12},e_{14},e_{16}\} \leftrightarrow \{v_2,v_5,v_7,v_8,v_9\}$;
цикл $c_8 = \{e_5,e_6,e_{12},e_{15},e_{19}\} \leftrightarrow \{v_2,v_5,v_7,v_9,v_{11}\}$;
цикл $c_{10} = \{e_9,e_{10},e_{13}\} \leftrightarrow \{v_4,v_6,v_7\}$;
цикл $c_{11} = \{e_{11},e_{12},e_{13}\} \leftrightarrow \{v_5,v_6,v_7\}$;
цикл $c_{12} = \{e_1,e_2,e_5,e_8,e_{12}\} \leftrightarrow \{v_1,v_2,v_3,v_5,v_7\}$;
цикл $c_{15} = \{e_{14},e_{15},e_{17},e_{20}\} \leftrightarrow \{v_7,v_8,v_{10},v_{11}\}$;
цикл $c_{17} = \{e_{18},e_{19},e_{20}\} \leftrightarrow \{v_9,v_{10},v_{11}\}$.

Номер $\langle 1,2,3,4,5,6,7,8,9,0,1,2,3,4,5,6,7,8,9,0\rangle$.
Вектор количества циклов по ребрам    $P_e = \langle 2,2,2,2,3,3,1,1,2,1,1,4,3,3,2,2,1,1,2,2\rangle$.

Кубический функционал Маклейна равен 48. Для построения плоского суграфа будем последовательно удалять циклы из базиса, с соблюдением условия Эйлера: удаление цикла порождает удаление одного и только одного ребра. Для записи процесса удаления циклов из базиса будем применять понятие дифференцирования структурного числа

**Определение 11.14**. *Алгебраической производной структурного числа* называется число [2]:

$\partial A / \partial a$, определенное как

$$\frac{\partial A}{\partial a} = A \big| \text{столбцы, не содержащие элемент } a \text{ исключены.} \tag{11.12}$$

Определим для исключения цикл, используя для записи операцию дифференцирования структурного числа. Значение кубического функционала Маклейна для выбранного базиса равно 48.

$$F(\frac{\partial C_b}{\partial c_3}) = 48, \ F(\frac{\partial C_b}{\partial c_{10}}) = 42, \ F(\frac{\partial C_b}{\partial c_{11}}) = 24, \ F(\frac{\partial C_b}{\partial c_{12}}) = 24, \ F(\frac{\partial C_b}{\partial c_{15}}) = 42, \ F(\frac{\partial C_b}{\partial c_{17}}) = 42.$$

Для удаления выбираем цикл $c_{12}$.

Номер  <1,2,3,4,5,6,7,8,9,0,1,2,3,4,5,6,7,8,9,0>.
Вектор количества циклов по ребрам     $P_e$ = <1,1,2,2,1,3,1,0,2,1,1,3,3,3,2,2,1,1,2,2>.

Значение кубического фунционала Маклейна не равно нулю.

Определим цикл для дальнейшего исключения, используя операцию дифференцирования. Удаление цикла $c_3$ приводит к удалению двух ребер, что противоречит условию Эйлера. Для сокращения записи будем рассматривать только те циклы, удаление которых приводит к исключению одного и только одного ребра.

$$F(\frac{\partial^2 C_b}{\partial c_{12} \partial c_{10}}) = 24, \ F(\frac{\partial^2 C_b}{\partial c_{12} \partial c_{11}}) = 12, \ F(\frac{\partial^2 C_b}{\partial c_{12} \partial c_{15}}) = 18, \ F(\frac{\partial^2 C_b}{\partial c_{12} \partial c_{17}}) = 24.$$

Для удаления выбираем цикл $c_{11}$.

Номер  <1,2,3,4,5,6,7,8,9,0,1,2,3,4,5,6,7,8,9,0>.
Вектор количества циклов по ребрам     $P_e$ = <1,1,2,2,1,3,1,0,2,1,0,2,2,3,2,2,1,1,2,2>.

Определим для дальнейшего исключения цикл, используя операцию дифференцирования.

$$F(\frac{\partial^3 C_b}{\partial c_{12} \partial c_{11} \partial c_2}) = 6, \ F(\frac{\partial^3 C_b}{\partial c_{12} \partial c_{11} \partial c_{10}}) = 6, \ F(\frac{\partial^3 C_b}{\partial c_{12} \partial c_{11} \partial c_{15}}) = 6.$$

Для удаления выбираем цикл $c_2$.

Номер  <1,2,3,4,5,6,7,8,9,0,1,2,3,4,5,6,7,8,9,0>.
Вектор количества циклов по ребрам     $P_e$ = <0,1,1,1,2,2,1,0,2,1,0,2,2,3,2,1,1,1,2,1>.

Определим для дальнейшего исключения цикл, используя операцию дифференцирования.

$$F(\frac{\partial^4 C_b}{\partial c_{12} \partial c_{11} \partial c_2 \partial c_7}) = 0,\ F(\frac{\partial^4 C_b}{\partial c_{12} \partial c_{11} \partial c_2 \partial c_{15}}) = 0,\ F(\frac{\partial^4 C_b}{\partial c_{12} \partial c_{11} \partial c_2 \partial c_{17}}) = 6.$$

Для удаления выбираем цикл $c_7$. Значение функционала Маклейна равно нулю. Сформирован плоский суграф.

цикл $c_3 = \{e_2, e_3, e_7, e_9\} \leftrightarrow \{v_1, v_3, v_4, v_6\}$;
цикл $c_6 = \{e_3, e_4, e_{13}, e_{14}\} \leftrightarrow \{v_1, v_6, v_7, v_8\}$;
цикл $c_8 = \{e_5, e_6, e_{12}, e_{15}, e_{19}\} \leftrightarrow \{v_2, v_5, v_7, v_9, v_{11}\}$;
цикл $c_{10} = \{e_9, e_{10}, e_{13}\} \leftrightarrow \{v_4, v_6, v_7\}$;
цикл $c_{15} = \{e_{14}, e_{15}, e_{17}, e_{20}\} \leftrightarrow \{v_7, v_8, v_{10}, v_{11}\}$;
цикл $c_{17} = \{e_{18}, e_{19}, e_{20}\} \leftrightarrow \{v_9, v_{10}, v_{11}\}$.

Номер <1,2,3,4,5,6,7,8,9,0,1,2,3,4,5,6,7,8,9,0>.
Вектор количества циклов по рёбрам $P_e = <0,1,2,1,1,1,1,0,2,1,0,1,2,2,2,0,1,1,2,2>$..

Номер <1,2,3,4,5,6,7,8,9,0,1>.
Вектор количества циклов по вершинам $P_v = <2,1,1,2,1,3,4,2,2,2,3>$.

На основе выделенной системы циклов, формируется топологический рисунок суграфа [11,12] (рис. 11.8). Обод базиса представлен на рис. 11.9.

Рис. 11.8. Плоская конфигурация $G_1^*$ графа $G_2$.

Случайным образом выделим следующую последовательность изометрических циклов:

$C_\tau = <\mathbf{c_{16}, c_4, c_{11}, c_1, c_8, c_{10}, c_5, c_9, c_{15}, c_2}, c_6, c_{12}, c_{17}, c_3, c_{14}, c_7, c_{13}>$.

Базис изометрических циклов графа:

цикл $c_1 = \{e_1, e_3, e_5, e_{11}\} \leftrightarrow \{v_1, v_2, v_5, v_6\}$;
цикл $c_2 = \{e_1, e_4, e_6, e_{16}\} \leftrightarrow \{v_1, v_2, v_8, v_9\}$;
цикл $c_4 = \{e_2, e_3, e_8, e_{13}\} \leftrightarrow \{v_1, v_3, v_6, v_7\}$;
цикл $c_5 = \{e_2, e_4, e_8, e_{14}\} \leftrightarrow \{v_1, v_3, v_7, v_8\}$;
цикл $c_8 = \{e_5, e_6, e_{12}, e_{15}, e_{19}\} \leftrightarrow \{v_2, v_5, v_7, v_9, v_{11}\}$;
цикл $c_9 = \{e_7, e_8, e_{10}\} \leftrightarrow \{v_3, v_4, v_7\}$;
цикл $c_{10} = \{e_9, e_{10}, e_{13}\} \leftrightarrow \{v_4, v_6, v_7\}$;
цикл $c_{11} = \{e_{11}, e_{12}, e_{13}\} \leftrightarrow \{v_5, v_6, v_7\}$;

цикл  $c_{15} = \{e_{14}, e_{15}, e_{17}, e_{20}\} \leftrightarrow \{v_7, v_8, v_{10}, v_{11}\}$;
цикл  $c_{16} = \{e_{16}, e_{17}, e_{18}\} \leftrightarrow \{v_8, v_9, v_{10}\}$.

Номер   <1,2,3,4,5,6,7,8,9,0,1,2,3,4,5,6,7,8,9,0>.
Вектор количества циклов по рёбрам    $P_e$ = <2,2,2,2,2,2,1,3,1,2,2,2,3,2,2,2,2,1,1,1>.

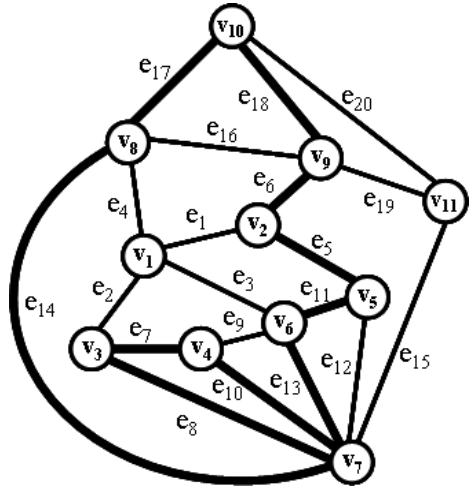

Рис. 11.9. Обод базиса циклов для суграфа $G_1^*$ графа $G_2$.

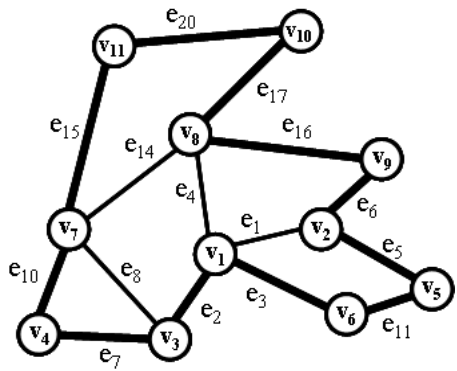

Рис. 11.10. Плоская конфигурация $G_2^*$ графа $G_2$.

Кубический функционал Маклейна равен 12. Удаляем циклы по правилу Эйлера до получения топологического рисунка плоского суграфа $G_2^*$ (рис. 11.10).

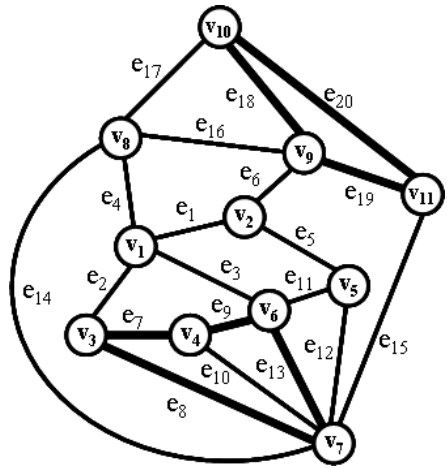

Рис. 11.11. Обод базиса циклов для суграфа $G_2^*$ графа $G_2$.

цикл $c_1 = \{e_1, e_3, e_5, e_{11}\} \leftrightarrow \{v_1, v_2, v_5, v_6\}$;
цикл $c_2 = \{e_1, e_4, e_6, e_{16}\} \leftrightarrow \{v_1, v_2, v_8, v_9\}$;
цикл $c_5 = \{e_2, e_4, e_8, e_{14}\} \leftrightarrow \{v_1, v_3, v_7, v_8\}$;
цикл $c_9 = \{e_7, e_8, e_{10}\} \leftrightarrow \{v_3, v_4, v_7\}$;
цикл $c_{15} = \{e_{14}, e_{15}, e_{17}, e_{20}\} \leftrightarrow \{v_7, v_8, v_{10}, v_{11}\}$;

Номер $<1,2,3,4,5,6,7,8,9,0,1,2,3,4,5,6,7,8,9,0>$.
Вектор количества циклов по рёбрам $P_e = <2,1,1,2,1,1,1,2,0,1,1,0,0,2,1,1,1,0,0,1>$.
Номер $<1,2,3,4,5,6,7,8,9,0,1>$.
Вектор количества циклов по вершинам $P_v = <3,2,2,1,1,1,3,3,1,1,1>$.

Заметим, что обод топологического рисунка представляет собой гамильтонов цикл.

Рассмотрим следующую случайную последовательность циклов графа $G_2$:

$C_\tau = <\mathbf{c_7}, \mathbf{c_2}, c_{13}, \mathbf{c_9}, \mathbf{c_3}, \mathbf{c_{11}}, \mathbf{c_{12}}, \mathbf{c_1}, \mathbf{c_{16}}, c_4, \mathbf{c_{17}}, c_6, \mathbf{c_{14}}, c_{10}, c_8, c_5, c_{15}>$.

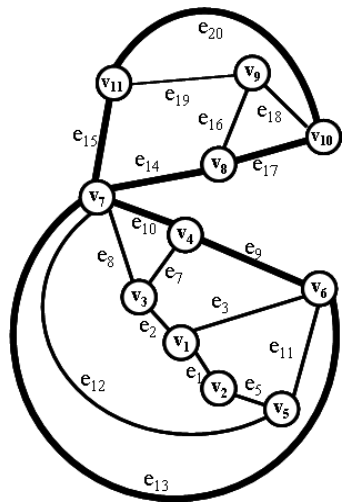

Рис. 11.12. Плоская конфигурация суграфа $G_3^*$ графа $G_2$.

Базис изометрических циклов:

цикл $c_1 = \{e_1, e_3, e_5, e_{11}\} \leftrightarrow \{v_1, v_2, v_5, v_6\}$;
цикл $c_2 = \{e_1, e_4, e_6, e_{16}\} \leftrightarrow \{v_1, v_2, v_8, v_9\}$;
цикл $c_3 = \{e_2, e_3, e_7, e_9\} \leftrightarrow \{v_1, v_3, v_4, v_6\}$;
цикл $c_7 = \{e_5, e_6, e_{12}, e_{14}, e_{16}\} \leftrightarrow \{v_2, v_5, v_7, v_8, v_9\}$;
цикл $c_9 = \{e_7, e_8, e_{10}\} \leftrightarrow \{v_3, v_4, v_7\}$;
цикл $c_{11} = \{e_{11}, e_{12}, e_{13}\} \leftrightarrow \{v_5, v_6, v_7\}$;
цикл $c_{12} = \{e_1, e_2, e_5, e_8, e_{12}\} \leftrightarrow \{v_1, v_2, v_3, v_5, v_7\}$;
цикл $c_{14} = \{e_{14}, e_{15}, e_{16}, e_{19}\} \leftrightarrow \{v_7, v_8, v_9, v_{11}\}$;
цикл $c_{16} = \{e_{16}, e_{17}, e_{18}\} \leftrightarrow \{v_8, v_9, v_{10}\}$;
цикл $c_{17} = \{e_{18}, e_{19}, e_{20}\} \leftrightarrow \{v_9, v_{10}, v_{11}\}$.

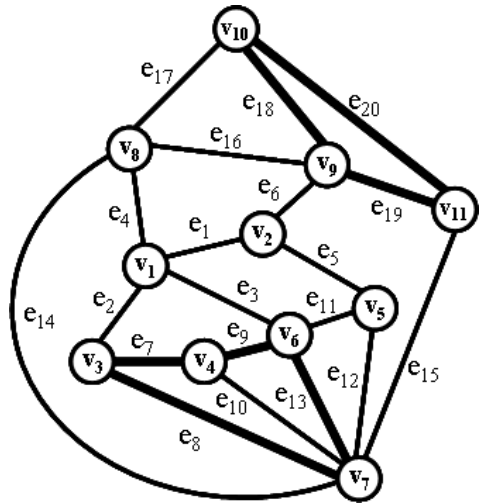

Рис. 11.13. Обод базиса циклов для суграфа $G_3^*$ графа $G_2$.

Номер <1,2,3,4,5,6,7,8,9,0,1,2,3,4,5,6,7,8,9,0>.
Вектор количества циклов по рёбрам $P_e$ = <3,2,2,1,3,2,2,2,1,1,2,3,1,2,1,4,1,2,2,1>.

Кубический функционал Маклейна равен 42. Выделим систему циклов, удовлетворяющую нулевому значению функционала Маклейна, и построим топологический рисунок:

цикл $c_1 = \{e_1, e_3, e_5, e_{11}\} \leftrightarrow \{v_1, v_2, v_5, v_6\}$;
цикл $c_3 = \{e_2, e_3, e_7, e_9\} \leftrightarrow \{v_1, v_3, v_4, v_6\}$;
цикл $c_9 = \{e_7, e_8, e_{10}\} \leftrightarrow \{v_3, v_4, v_7\}$;
цикл $c_{11} = \{e_{11}, e_{12}, e_{13}\} \leftrightarrow \{v_5, v_6, v_7\}$;
цикл $c_{12} = \{e_1, e_2, e_5, e_8, e_{12}\} \leftrightarrow \{v_1, v_2, v_3, v_5, v_7\}$;
цикл $c_{14} = \{e_{14}, e_{15}, e_{16}, e_{19}\} \leftrightarrow \{v_7, v_8, v_9, v_{11}\}$;
цикл $c_{16} = \{e_{16}, e_{17}, e_{18}\} \leftrightarrow \{v_8, v_9, v_{10}\}$;
цикл $c_{17} = \{e_{18}, e_{19}, e_{20}\} \leftrightarrow \{v_9, v_{10}, v_{11}\}$.

Номер <1,2,3,4,5,6,7,8,9,0,1,2,3,4,5,6,7,8,9,0>.
Вектор количества циклов по рёбрам $P_e$ = <2,2,2,0,2,0,2,2,1,1,2,2,1,1,1,2,1,2,2,1>.

Номер <1,2,3,4,5,6,7,8,9,0,1>.
Вектор количества циклов по вершинам $P_v$ = <3,2,3,2,3,3,4,2,3,2,2>.

Данный суграф обладает точкой сопряжения (рис. 11.12).

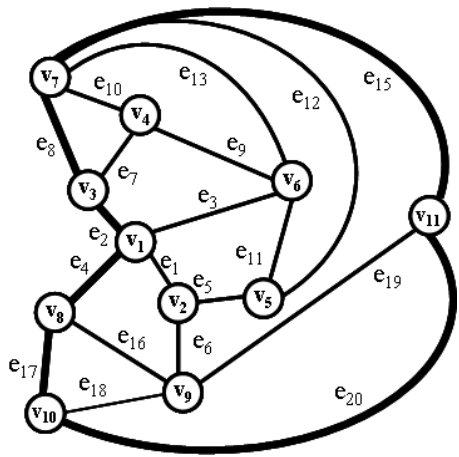

Рис. 11.14. Плоская конфигурация суграфа $G_4^*$.

Рассмотрим следующую случайную последовательность изометрических циклов графа $G_2$:

$C_\tau = <\mathbf{c_3},\mathbf{c_{17}},\mathbf{c_{14}},\mathbf{c_9},\mathbf{c_8},\mathbf{c_2},\mathbf{c_{16}},c_{15},\mathbf{c_{10}},\mathbf{c_1},\mathbf{c_{11}},c_{13},c_5,c_6,c_7,c_4,c_{12}>$.

Базис изометрических циклов:

цикл $c_1 = \{e_1,e_3,e_5,e_{11}\} \leftrightarrow \{v_1,v_2,v_5,v_6\}$;
цикл $c_2 = \{e_1,e_4,e_6,e_{16}\} \leftrightarrow \{v_1,v_2,v_8,v_9\}$;
цикл $c_3 = \{e_2,e_3,e_7,e_9\} \leftrightarrow \{v_1,v_3,v_4,v_6\}$;
цикл $c_8 = \{e_5,e_6,e_{12},e_{15},e_{19}\} \leftrightarrow \{v_2,v_5,v_7,v_9,v_{11}\}$;
цикл $c_9 = \{e_7,e_8,e_{10}\} \leftrightarrow \{v_3,v_4,v_7\}$;
цикл $c_{10} = \{e_9,e_{10},e_{13}\} \leftrightarrow \{v_4,v_6,v_7\}$;
цикл $c_{11} = \{e_{11},e_{12},e_{13}\} \leftrightarrow \{v_5,v_6,v_7\}$;
цикл $c_{14} = \{e_{14},e_{15},e_{16},e_{19}\} \leftrightarrow \{v_7,v_8,v_9,v_{11}\}$;
цикл $c_{16} = \{e_{16},e_{17},e_{18}\} \leftrightarrow \{v_8,v_9,v_{10}\}$;
цикл $c_{17} = \{e_{18},e_{19},e_{20}\} \leftrightarrow \{v_9,v_{10},v_{11}\}$.

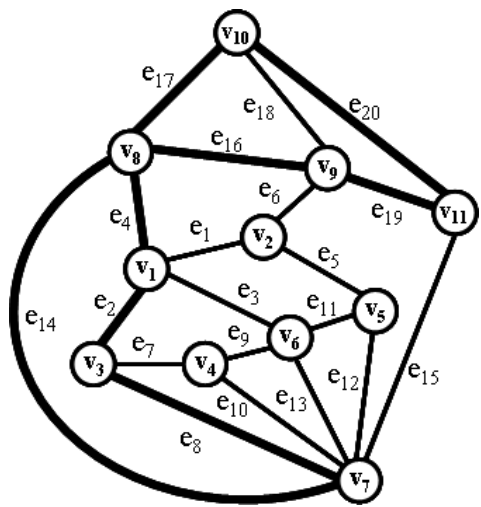

Рис. 11.15. Обод базиса циклов для суграфа $G_4^*$.

Номер $<1,2,3,4,5,6,7,8,9,0,1,2,3,4,5,6,7,8,9,0>$.
Вектор количества циклов по ребрам $P_e = <2,1,2,1,2,2,2,1,2,2,2,2,2,1,2,3,1,2,3,1>$.

Кубический функционал Маклейна равен 12. Выделим систему циклов, удовлетворяющую нулевому значению функционала Маклейна, и построим топологический рисунок (рис. 11.14):

цикл $c_1 = \{e_1,e_3,e_5,e_{11}\} \leftrightarrow \{v_1,v_2,v_5,v_6\}$;
цикл $c_2 = \{e_1,e_4,e_6,e_{16}\} \leftrightarrow \{v_1,v_2,v_8,v_9\}$;
цикл $c_3 = \{e_2,e_3,e_7,e_9\} \leftrightarrow \{v_1,v_3,v_4,v_6\}$;
цикл $c_8 = \{e_5,e_6,e_{12},e_{15},e_{19}\} \leftrightarrow \{v_2,v_5,v_7,v_9,v_{11}\}$;
цикл $c_9 = \{e_7,e_8,e_{10}\} \leftrightarrow \{v_3,v_4,v_7\}$;
цикл $c_{10} = \{e_9,e_{10},e_{13}\} \leftrightarrow \{v_4,v_6,v_7\}$;
цикл $c_{11} = \{e_{11},e_{12},e_{13}\} \leftrightarrow \{v_5,v_6,v_7\}$;
цикл $c_{16} = \{e_{16},e_{17},e_{18}\} \leftrightarrow \{v_8,v_9,v_{10}\}$;
цикл $c_{17} = \{e_{18},e_{19},e_{20}\} \leftrightarrow \{v_9,v_{10},v_{11}\}$.

Номер $<1,2,3,4,5,6,7,8,9,0,1,2,3,4,5,6,7,8,9,0>$.
Вектор количества циклов по ребрам $P_e = <2,1,2,1,2,2,2,1,2,2,2,2,2,0,1,2,1,2,2,1>$.

Номер <1,2,3,4,5,6,7,8,9,0,1>.

Вектор количества циклов по вершинам $P_v = <3,3,2,3,3,4,3,2,4,2,2>$.

Рассмотрим следующую случайную последовательность изометрических циклов графа $G_2$:

$C_\tau = <\mathbf{c_{12}},\mathbf{c_{13}},\mathbf{c_8},\mathbf{c_{17}},\mathbf{c_{15}},\mathbf{c_{11}},\mathbf{c_{16}},c_2,\mathbf{c_1},c_4,c_5,c_7,\mathbf{c_{10}},c_6,\mathbf{c_3},c_9,c_{14}>$

Базис изометрических циклов:

цикл $c_1 = \{e_1,e_3,e_5,e_{11}\} \leftrightarrow \{v_1,v_2,v_5,v_6\}$;
цикл $c_3 = \{e_2,e_3,e_7,e_9\} \leftrightarrow \{v_1,v_3,v_4,v_6\}$;
цикл $c_8 = \{e_5,e_6,e_{12},e_{15},e_{19}\} \leftrightarrow \{v_2,v_5,v_7,v_9,v_{11}\}$;
цикл $c_{10} = \{e_9,e_{10},e_{13}\} \leftrightarrow \{v_4,v_6,v_7\}$;
цикл $c_{11} = \{e_{11},e_{12},e_{13}\} \leftrightarrow \{v_5,v_6,v_7\}$;
цикл $c_{12} = \{e_1,e_2,e_5,e_8,e_{12}\} \leftrightarrow \{v_1,v_2,v_3,v_5,v_7\}$;
цикл $c_{13} = \{e_1,e_4,e_5,e_{12},e_{14}\} \leftrightarrow \{v_1,v_2,v_5,v_7,v_8\}$;
цикл $c_{15} = \{e_{14},e_{15},e_{17},e_{20}\} \leftrightarrow \{v_7,v_8,v_{10},v_{11}\}$;
цикл $c_{16} = \{e_{16},e_{17},e_{18}\} \leftrightarrow \{v_8,v_9,v_{10}\}$;
цикл $c_{17} = \{e_{18},e_{19},e_{20}\} \leftrightarrow \{v_9,v_{10},v_{11}\}$.

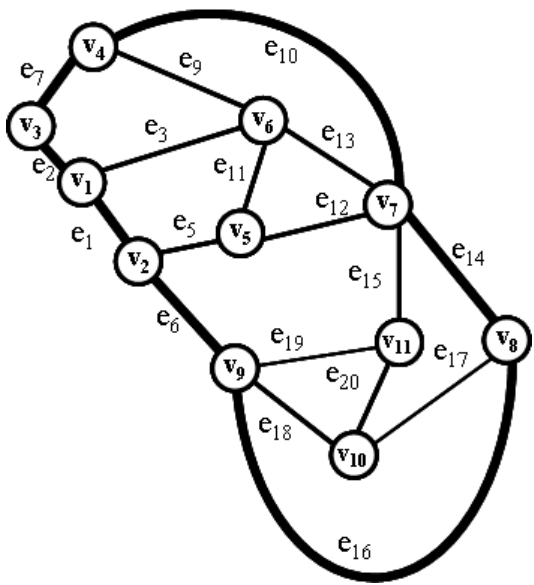

Рис. 11.16. Плоская конфигурация суграфа $G_5^*$.

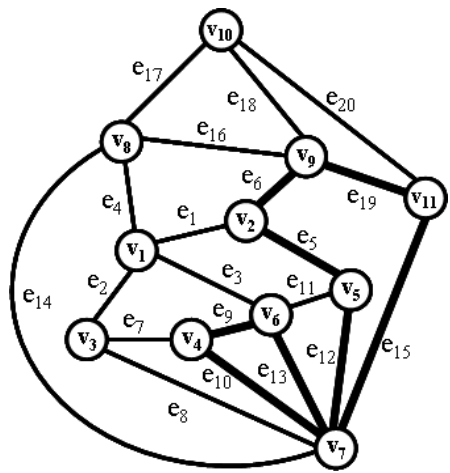

Рис. 11.17. Обод базиса циклов для суграфа $G_5^*$.

Номер <1,2,3,4,5,6,7,8,9,0,1,2,3,4,5,6,7,8,9,0>.
Вектор количества циклов по рёбрам $P_e = <2,2,2,0,3,1,2,2,1,1,2,3,1,2,3,2,2,2,3,2>$.

Кубический функционал Маклейна равен 24 и удалено одно ребро. Выделим систему циклов, удовлетворяющую нулевому значению функционала Маклейна, и построим топологический рисунок (рис. 11.16):

цикл $c_1 = \{e_1,e_3,e_5,e_{11}\} \leftrightarrow \{v_1,v_2,v_5,v_6\}$;
цикл $c_3 = \{e_2,e_3,e_7,e_9\} \leftrightarrow \{v_1,v_3,v_4,v_6\}$;
цикл $c_8 = \{e_5,e_6,e_{12},e_{15},e_{19}\} \leftrightarrow \{v_2,v_5,v_7,v_9,v_{11}\}$;
цикл $c_{10} = \{e_9,e_{10},e_{13}\} \leftrightarrow \{v_4,v_6,v_7\}$;
цикл $c_{11} = \{e_{11},e_{12},e_{13}\} \leftrightarrow \{v_5,v_6,v_7\}$;
цикл $c_{15} = \{e_{14},e_{15},e_{17},e_{20}\} \leftrightarrow \{v_7,v_8,v_{10},v_{11}\}$;
цикл $c_{16} = \{e_{16},e_{17},e_{18}\} \leftrightarrow \{v_8,v_9,v_{10}\}$;
цикл $c_{17} = \{e_{18},e_{19},e_{20}\} \leftrightarrow \{v_9,v_{10},v_{11}\}$.

Номер <1,2,3,4,5,6,7,8,9,0,1,2,3,4,5,6,7,8,9,0>.
Вектор количества циклов по рёбрам $P_e = <1,1,2,0,2,1,1,0,2,1,2,2,2,1,2,1,2,2,2,2>$.

Номер <1,2,3,4,5,6,7,8,9,0,1>.
Вектор количества циклов по вершинам $P_v = <2,2,1,2,3,4,4,2,3,3,3>$.

Рассмотрим следующую случайную последовательность изометрических циклов графа $G_2$:

$C_\tau = <\mathbf{c_{14},c_3,c_4,c_{12},c_1,c_6,c_{16},c_{15},c_8,c_9},c_7,c_{10},c_5,c_2,c_{17},c_{13},c_{11}>$.

Базис изометрических циклов:

цикл $c_1 = \{e_1,e_3,e_5,e_{11}\} \leftrightarrow \{v_1,v_2,v_5,v_6\}$;
цикл $c_3 = \{e_2,e_3,e_7,e_9\} \leftrightarrow \{v_1,v_3,v_4,v_6\}$;
цикл $c_4 = \{e_2,e_3,e_8,e_{13}\} \leftrightarrow \{v_1,v_3,v_6,v_7\}$;
цикл $c_6 = \{e_3,e_4,e_{13},e_{14}\} \leftrightarrow \{v_1,v_6,v_7,v_8\}$;
цикл $c_8 = \{e_5,e_6,e_{12},e_{15},e_{19}\} \leftrightarrow \{v_2,v_5,v_7,v_9,v_{11}\}$;
цикл $c_9 = \{e_7,e_8,e_{10}\} \leftrightarrow \{v_3,v_4,v_7\}$;
цикл $c_{12} = \{e_1,e_2,e_5,e_8,e_{12}\} \leftrightarrow \{v_1,v_2,v_3,v_5,v_7\}$;
цикл $c_{14} = \{e_{14},e_{15},e_{16},e_{19}\} \leftrightarrow \{v_7,v_8,v_9,v_{11}\}$;
цикл $c_{15} = \{e_{14},e_{15},e_{17},e_{20}\} \leftrightarrow \{v_7,v_8,v_{10},v_{11}\}$;
цикл $c_{16} = \{e_{16},e_{17},e_{18}\} \leftrightarrow \{v_8,v_9,v_{10}\}$.

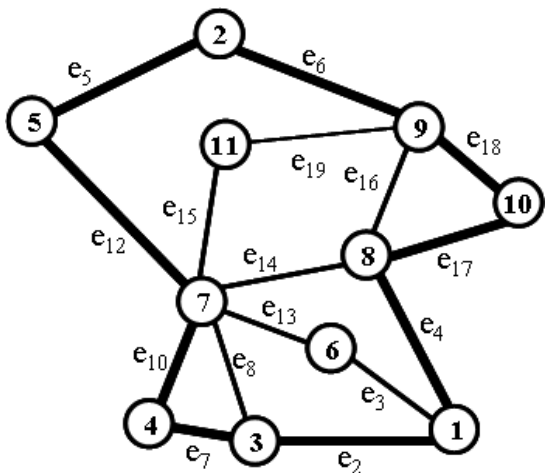

Рис. 11.18. Плоская конфигурация суграфа $G_6^*$.

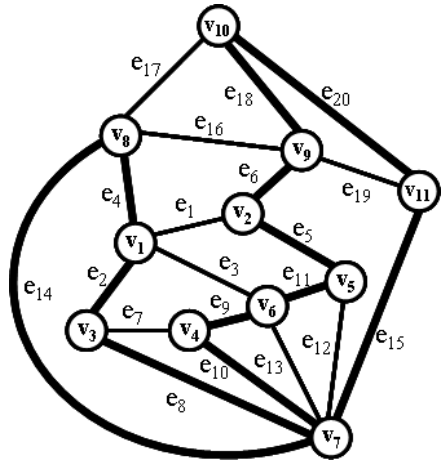

Рис. 11.19. Обод базиса циклов для суграфа $G_6^*$.

Номер <1,2,3,4,5,6,7,8,9,0,1,2,3,4,5,6,7,8,9,0>.
Вектор количества циклов по рёбрам $P_e$ = <2,3,4,1,3,1,2,3,1,1,1,2,2,3,3,2,2,1,2,1>.

Кубический функционал Маклейна равен 54. Выделим систему циклов, удовлетворяющую нулевому значению функционала Маклейна, и построим топологический рисунок (рис. 11.18):

цикл $c_4 = \{e_2, e_3, e_8, e_{13}\} \leftrightarrow \{v_1, v_3, v_6, v_7\}$;
цикл $c_6 = \{e_3, e_4, e_{13}, e_{14}\} \leftrightarrow \{v_1, v_6, v_7, v_8\}$;
цикл $c_8 = \{e_5, e_6, e_{12}, e_{15}, e_{19}\} \leftrightarrow \{v_2, v_5, v_7, v_9, v_{11}\}$;
цикл $c_9 = \{e_7, e_8, e_{10}\} \leftrightarrow \{v_3, v_4, v_7\}$;
цикл $c_{14} = \{e_{14}, e_{15}, e_{16}, e_{19}\} \leftrightarrow \{v_7, v_8, v_9, v_{11}\}$;
цикл $c_{16} = \{e_{16}, e_{17}, e_{18}\} \leftrightarrow \{v_8, v_9, v_{10}\}$.

Номер <1,2,3,4,5,6,7,8,9,0,1,2,3,4,5,6,7,8,9,0>.
Вектор количества циклов по рёбрам $P_e$ = <0,1,2,1,1,1,1,2,0,1,0,1,2,2,1,2,1,1,1,0>.

Номер <1,2,3,4,5,6,7,8,9,0,1>.
Вектор количества циклов по вершинам $P_v$ = <2,1,2,1,1,2,5,3,3,1,2>.

Рассмотрим следующую случайную последовательность изометрических циклов графа $G_2$:

$C_\tau = <\mathbf{c_5}, \mathbf{c_6}, \mathbf{c_{14}}, \mathbf{c_{13}}, c_4, \mathbf{c_{16}}, \mathbf{c_{17}}, \mathbf{c_2}, \mathbf{c_1}, c_{15}, \mathbf{c_9}, c_{11}, \mathbf{c_3}, c_8, c_7, c_{12}, c_{10}>$.

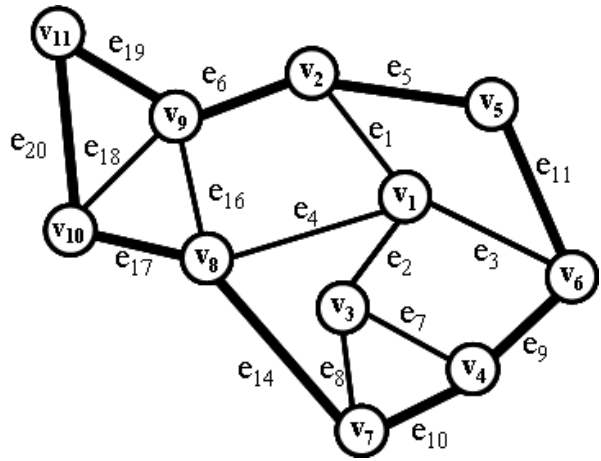

Рис. 11.20. Плоская конфигурация суграфа $G_7^*$.

Базис изометрических циклов:

цикл $c_1 = \{e_1, e_3, e_5, e_{11}\} \leftrightarrow \{v_1, v_2, v_5, v_6\}$;
цикл $c_2 = \{e_1, e_4, e_6, e_{16}\} \leftrightarrow \{v_1, v_2, v_8, v_9\}$;
цикл $c_3 = \{e_2, e_3, e_7, e_9\} \leftrightarrow \{v_1, v_3, v_4, v_6\}$;
цикл $c_5 = \{e_2, e_4, e_8, e_{14}\} \leftrightarrow \{v_1, v_3, v_7, v_8\}$;
цикл $c_6 = \{e_3, e_4, e_{13}, e_{14}\} \leftrightarrow \{v_1, v_6, v_7, v_8\}$;
цикл $c_9 = \{e_7, e_8, e_{10}\} \leftrightarrow \{v_3, v_4, v_7\}$;
цикл $c_{13} = \{e_1, e_4, e_5, e_{12}, e_{14}\} \leftrightarrow \{v_1, v_2, v_5, v_7, v_8\}$;
цикл $c_{14} = \{e_{14}, e_{15}, e_{16}, e_{19}\} \leftrightarrow \{v_7, v_8, v_9, v_{11}\}$;
цикл $c_{16} = \{e_{16}, e_{17}, e_{18}\} \leftrightarrow \{v_8, v_9, v_{10}\}$;
цикл $c_{17} = \{e_{18}, e_{19}, e_{20}\} \leftrightarrow \{v_9, v_{10}, v_{11}\}$.

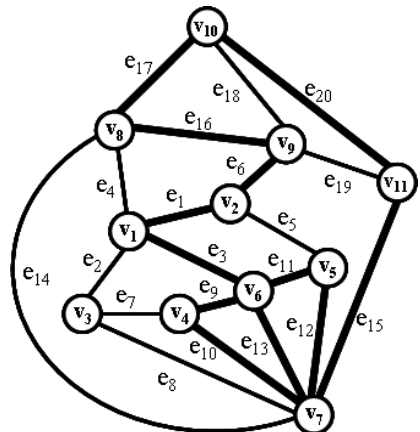

Рис. 11.21. Обод базиса циклов для суграфа $G_7^*$.

Номер <1,2,3,4,5,6,7,8,9,0,1,2,3,4,5,6,7,8,9,0>.
Вектор количества циклов по рёбрам $P_e = $ <3,2,3,4,2,1,2,2,1,1,1,1,1,4,1,3,1,2,2,1>.

Кубический функционал Маклейна равен 66. Выделим систему циклов, удовлетворяющую нулевому значению функционала Маклейна, и построим топологический рисунок (рис. 11.20):

цикл $c_1 = \{e_1, e_3, e_5, e_{11}\} \leftrightarrow \{v_1, v_2, v_5, v_6\}$;
цикл $c_2 = \{e_1, e_4, e_6, e_{16}\} \leftrightarrow \{v_1, v_2, v_8, v_9\}$;
цикл $c_3 = \{e_2, e_3, e_7, e_9\} \leftrightarrow \{v_1, v_3, v_4, v_6\}$;
цикл $c_5 = \{e_2, e_4, e_8, e_{14}\} \leftrightarrow \{v_1, v_3, v_7, v_8\}$;

цикл  $c_9 = \{e_7, e_8, e_{10}\} \leftrightarrow \{v_3, v_4, v_7\}$;
цикл  $c_{16} = \{e_{16}, e_{17}, e_{18}\} \leftrightarrow \{v_8, v_9, v_{10}\}$;
цикл  $c_{17} = \{e_{18}, e_{19}, e_{20}\} \leftrightarrow \{v_9, v_{10}, v_{11}\}$.

Номер <1,2,3,4,5,6,7,8,9,0,1,2,3,4,5,6,7,8,9,0>.
Вектор количества циклов по рёбрам    $P_e = <2,2,2,2,1,1,2,2,1,1,1,0,0,1,0,2,1,2,1,1>$.

Номер <1,2,3,4,5,6,7,8,9,0,1>.
Вектор количества циклов по вершинам $P_v = <4,2,3,2,2,2,2,3,3,2,1>$.

Рассмотрим следующую случайную последовательность изометрических циклов графа $G_2$:

$C_\tau = <\mathbf{c_{10}, c_3, c_8, c_{16}, c_5, c_1, c_9, c_2, c_{17}, c_{11}}, c_6, c_4, c_{15}, c_{13}, c_7, c_{12}, c_{14}>$

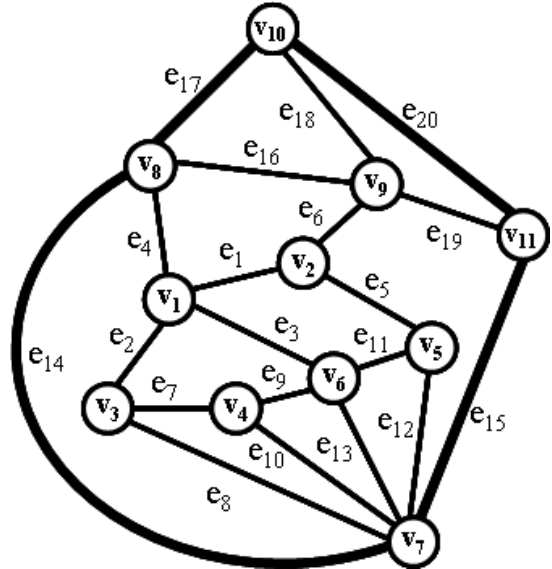

Рис. 11.22. Плоская конфигурация суграфа $G_8^*$ и обод для базиса.

Базис изометрических циклов:

цикл  $c_1 = \{e_1, e_3, e_5, e_{11}\} \leftrightarrow \{v_1, v_2, v_5, v_6\}$;
цикл  $c_2 = \{e_1, e_4, e_6, e_{16}\} \leftrightarrow \{v_1, v_2, v_8, v_9\}$;
цикл  $c_3 = \{e_2, e_3, e_7, e_9\} \leftrightarrow \{v_1, v_3, v_4, v_6\}$;
цикл  $c_5 = \{e_2, e_4, e_8, e_{14}\} \leftrightarrow \{v_1, v_3, v_7, v_8\}$;
цикл  $c_8 = \{e_5, e_6, e_{12}, e_{15}, e_{19}\} \leftrightarrow \{v_2, v_5, v_7, v_9, v_{11}\}$;
цикл  $c_9 = \{e_7, e_8, e_{10}\} \leftrightarrow \{v_3, v_4, v_7\}$;
цикл  $c_{10} = \{e_9, e_{10}, e_{13}\} \leftrightarrow \{v_4, v_6, v_7\}$;
цикл  $c_{11} = \{e_{11}, e_{12}, e_{13}\} \leftrightarrow \{v_5, v_6, v_7\}$;
цикл  $c_{16} = \{e_{16}, e_{17}, e_{18}\} \leftrightarrow \{v_8, v_9, v_{10}\}$;
цикл  $c_{17} = \{e_{18}, e_{19}, e_{20}\} \leftrightarrow \{v_9, v_{10}, v_{11}\}$.

Номер <1,2,3,4,5,6,7,8,9,0,1,2,3,4,5,6,7,8,9,0>.
Вектор количества циклов по рёбрам    $P_e = <2,2,2,2,2,2,2,2,2,2,2,2,2,1,1,2,1,2,2,1>$.

Номер <1,2,3,4,5,6,7,8,9,0,1>.
Вектор количества циклов по вершинам $P_v = <4,3,3,3,3,4,5,3,4,2,2>$.

Кубический функционал Маклейна равен 0 и соответствует топологическому рисунку плоского суграфа (рис. 11.22).

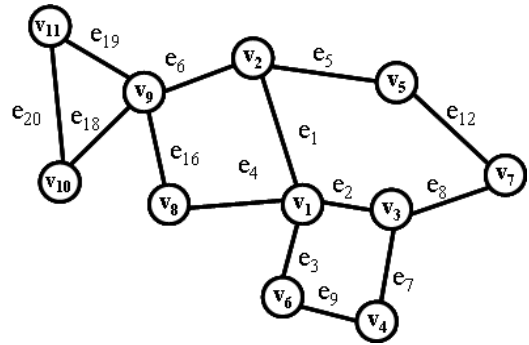

Рис. 11.23. Плоская конфигурация суграфа $G_9^*$.

Рассмотрим следующую случайную последовательность изометрических циклов графа $G_2$:

$C_\tau = <\mathbf{c_{16}},\mathbf{c_{10}},\mathbf{c_{17}},\mathbf{c_{11}},\mathbf{c_6},\mathbf{c_4},\mathbf{c_{14}},\mathbf{c_{12}},\mathbf{c_2},\mathbf{c_5},c_8,c_1,c_7,c_3,c_9,c_{13},c_{15}>$.

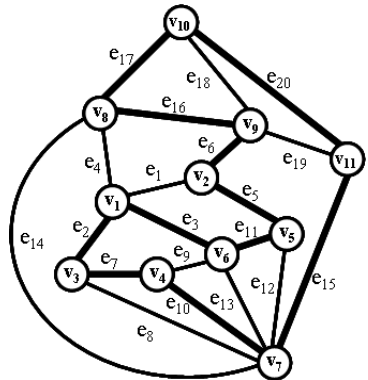

Рис. 11.24. Обод базиса циклов для суграфа $G_9^*$.

Базис изометрических циклов:

цикл $c_2 = \{e_1,e_4,e_6,e_{16}\} \leftrightarrow \{v_1,v_2,v_8,v_9\}$;
цикл $c_3 = \{e_2,e_3,e_7,e_9\} \leftrightarrow \{v_1,v_3,v_4,v_6\}$;
цикл $c_4 = \{e_2,e_3,e_8,e_{13}\} \leftrightarrow \{v_1,v_3,v_6,v_7\}$;
цикл $c_6 = \{e_3,e_4,e_{13},e_{14}\} \leftrightarrow \{v_1,v_6,v_7,v_8\}$;
цикл $c_{10} = \{e_9,e_{10},e_{13}\} \leftrightarrow \{v_4,v_6,v_7\}$;
цикл $c_{11} = \{e_{11},e_{12},e_{13}\} \leftrightarrow \{v_5,v_6,v_7\}$;
цикл $c_{12} = \{e_1,e_2,e_5,e_8,e_{12}\} \leftrightarrow \{v_1,v_2,v_3,v_5,v_7\}$;
цикл $c_{14} = \{e_{14},e_{15},e_{16},e_{19}\} \leftrightarrow \{v_7,v_8,v_9,v_{11}\}$;
цикл $c_{16} = \{e_{16},e_{17},e_{18}\} \leftrightarrow \{v_8,v_9,v_{10}\}$;
цикл $c_{17} = \{e_{18},e_{19},e_{20}\} \leftrightarrow \{v_9,v_{10},v_{11}\}$.

Номер <1,2,3,4,5,6,7,8,9,0,1,2,3,4,5,6,7,8,9,0>.
Вектор количества циклов по ребрам $P_e = <2,3,3,2,1,1,1,2,2,1,1,2,4,2,1,3,1,2,2,1>$

Кубический функционал Маклейна равен 42. Выделим систему циклов, удовлетворяющую нулевому значению функционала Маклейна, и построим топологический рисунок (рис. 11.23):

цикл   $c_2 = \{e_1,e_4,e_6,e_{16}\} \leftrightarrow \{v_1,v_2,v_8,v_9\}$;
цикл   $c_3 = \{e_2,e_3,e_7,e_9\} \leftrightarrow \{v_1,v_3,v_4,v_6\}$;
цикл   $c_{12} = \{e_1,e_2,e_5,e_8,e_{12}\} \leftrightarrow \{v_1,v_2,v_3,v_5,v_7\}$;
цикл   $c_{17} = \{e_{18},e_{19},e_{20}\} \leftrightarrow \{v_9,v_{10},v_{11}\}$.

Номер  <1,2,3,4,5,6,7,8,9,0,1,2,3,4,5,6,7,8,9,0>.
Вектор количества циклов по ребрам    $P_e$ = <2,2,1,1,1,1,1,1,1,0,0,1,0,0,0,1,0,1,1,1>
Номер  <1,2,3,4,5,6,7,8,9,0,1>.
Вектор количества циклов по вершинам $P_v$ = <4,3,3,2,2,2,2,2,4,2,2>.

**Комментарии**

В данной главе представлена генерация базисов подпространства циклов C(G) графа G, состоящая из изометрических циклов. Выделение подмножества независимых изометрических циклов (басиса подпространства циклов графа C) производится модифицированным алгоритмом Гаусса. Представлен метод удаления циклов из базиса (с соблюдением условия Эйлера) до получения подмножества циклов с нулевым значением кубического функционала Маклейна.

Для графа $G_2$, представленного на рис. 11.7, методом Монте-Карло выделено 9 вариантов базисов изометрических циклов. Плоский топологический рисунок для каждого базиса строится методом удаления циклов с соблюдением правила Эйлера. Процесс удаления циклов производится до достижения нулевого значения кубического функционала Маклейна. Для каждого варианта определен обод базиса изометрических циклов. Количество нулевых элементов в векторе $P_e$ характеризует количество удаленных ребер в процессе планаризации. Количество удаленных ребер будем обозначать двумя латинскими буквами Nu.

Следует заметить, что во всех топологических рисунках плоской конфигурации количество вершин остается неизменным и равным количеству вершин графа *n*.

Построим следующую таблицу:

| Вариант № | Значение функционала базиса $F(C_b)$ | Количество циклов в плоском суграфе $card(G')$ | Количество удаленных ребер Nu | Длина циклов | Сумма длин |
|---|---|---|---|---|---|
| 1 | 48 | 6 | 4 | $3\times3, 4\times4, 3\times5$ | 40 |
| 2 | 12 | 5 | 5 | $4\times3, 5\times4, 1\times5$ | 37 |
| 3 | 42 | 8 | 2 | $4\times3, 5\times4, 1\times5$ | 37 |
| 4 | 12 | 9 | 1 | $5\times3, 4\times4, 1\times5$ | 36 |
| 5 | 24 | 8 | 2 | $5\times3, 4\times4, 1\times5$ | 36 |
| 6 | 54 | 6 | 4 | $2\times3, 5\times4, 3\times5$ | 41 |
| 7 | 66 | 7 | 3 | $3\times3, 6\times4, 1\times5$ | 38 |
| 8 | 0 | 10 | 0 | $5\times3, 4\times4, 1\times5$ | 36 |
| 9 | 42 | 4 | 6 | $4\times3, 5\times4, 1\times5$ | 37 |

По результатам можно сделать следующий вывод: максимально плоский суграф может быть построен для базиса, имеющего минимальное значение кубического функционала Маклейна. С другой стороны, сумма длин циклов должна быть минимальной.

Таким образом, задача выделения максимально плоского суграфа сводится к задаче выделения базиса, состоящего из изометрических циклов графа, имеющего минимальное значение функционала Маклейна, и минимальной суммы длин циклов. Это можно записать в следующем виде:

$$C_b \in C_\tau, c_i \in C_b, i = [1,k], k(C_b) = m - n + 1, F(C_b) \to \min, \sum_{i=1}^{k} [c_i[ \to \min \qquad (11.13)$$

# Глава 12. МЕТОД НАИСКОРЕЙШЕГО СПУСКА
## 12.1. Метод наискорейшего спуска и выделение базиса

Рассмотренный ранее метод позволяет получать множество различных плоских конфигураций. Метод основан на случайном выборе последовательности расположения элементов в кортеже циклов, определяемой перестановкой. Такой подход характеризуется большими вычислительными затратами, поэтому желательно уменьшить количество изометрических циклов в системе. Метод наискорейшего спуска, применяемый в дискретной математике, позволяет выделить систему изометрических циклов с мощностью, соизмеримой с цикломатическим числом графа. Поэтому предлагается рассмотреть следующий метод последовательного анализа и отсеивания вариантов. Он заключается в последовательном исключении элементов из системы до получения подмножества с мощностью $card(m-n+1+k)$, где k – число от 1 до 4. В этом случае перебор вариантов связан с меньшими вычислительными затратами. При этом алгоритм поиска плоской конфигурации остается прежним.

Рассмотрим метод последовательного анализа и отсеивания вариантов на примере графа $G_4$.

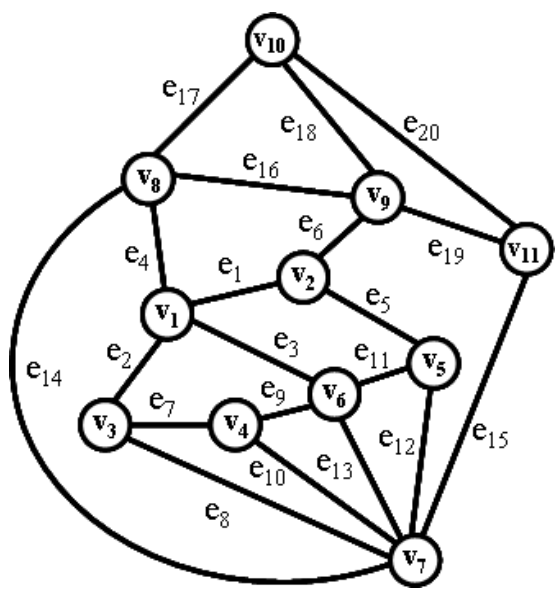

Рис. 12.1. Граф $G_2$.

Множество изометрических циклов графа:
цикл  $c_1 = \{e_1,e_3,e_5,e_{11}\} \leftrightarrow \{v_1,v_2,v_5,v_6\}$;
цикл  $c_2 = \{e_1,e_4,e_6,e_{16}\} \leftrightarrow \{v_1,v_2,v_8,v_9\}$;
цикл  $c_3 = \{e_2,e_3,e_7,e_9\} \leftrightarrow \{v_1,v_3,v_4,v_6\}$;
цикл  $c_4 = \{e_2,e_3,e_8,e_{13}\} \leftrightarrow \{v_1,v_3,v_6,v_7\}$;
цикл  $c_5 = \{e_2,e_4,e_8,e_{14}\} \leftrightarrow \{v_1,v_3,v_7,v_8\}$;
цикл  $c_6 = \{e_3,e_4,e_{13},e_{14}\} \leftrightarrow \{v_1,v_6,v_7,v_8\}$;
цикл  $c_7 = \{e_5,e_6,e_{12},e_{14},e_{16}\} \leftrightarrow \{v_2,v_5,v_7,v_8,v_9\}$;
цикл  $c_8 = \{e_5,e_6,e_{12},e_{15},e_{19}\} \leftrightarrow \{v_2,v_5,v_7,v_9,v_{11}\}$;
цикл  $c_9 = \{e_7,e_8,e_{10}\} \leftrightarrow \{v_3,v_4,v_7\}$;
цикл  $c_{10} = \{e_9,e_{10},e_{13}\} \leftrightarrow \{v_4,v_6,v_7\}$;
цикл  $c_{11} = \{e_{11},e_{12},e_{13}\} \leftrightarrow \{v_5,v_6,v_7\}$;

цикл  $c_{12} = \{e_1, e_2, e_5, e_8, e_{12}\} \leftrightarrow \{v_1, v_2, v_3, v_5, v_7\}$;
цикл  $c_{13} = \{e_1, e_4, e_5, e_{12}, e_{14}\} \leftrightarrow \{v_1, v_2, v_5, v_7, v_8\}$;
цикл  $c_{14} = \{e_{14}, e_{15}, e_{16}, e_{19}\} \leftrightarrow \{v_7, v_8, v_9, v_{11}\}$;
цикл  $c_{15} = \{e_{14}, e_{15}, e_{17}, e_{20}\} \leftrightarrow \{v_7, v_8, v_{10}, v_{11}\}$;
цикл  $c_{16} = \{e_{16}, e_{17}, e_{18}\} \leftrightarrow \{v_8, v_9, v_{10}\}$;
цикл  $c_{17} = \{e_{18}, e_{19}, e_{20}\} \leftrightarrow \{v_9, v_{10}, v_{11}\}$.

Цикломатическое число графа = 10

Номер  <1,2,3,4,5,6,7,8,9,0,1,2,3,4,5,6,7,8,9,0>
Вектор количества циклов по рёбрам     $P_e$ = <4,4,4,4,5,3,2,4,2,2,2,5,4,6,3,4,2,2,3,2>

Номер <1,2,3,4,5,6, 7 ,8,9,0,1>
Вектор количества циклов по вершинам $P_v$ = <8,6,5,3,6,6,12,7,5,3,4>.

Значение кубического функционала Маклейна для множества изометрических циклов графа равно 426.

**Определение 12.1.** *Алгебраической обратной производной структурного числа* называется структурное число $\delta A / \delta a$:

$$\frac{\delta A}{\delta a} = A \big| \text{столбцы, содержащие элемент } a \text{ опущены.} \tag{12.1}$$

Запись в виде обратной производной структурного числа $F(\frac{\delta C_\tau}{\delta c_3}) = 390$ говорит о том, что значение кубического функционала Маклейна для подмножества циклов с удаленным циклом $c_3$ из множества изометрических циклов $C_\tau$ равно 390. Будем последовательно удалять по одному циклу:

$F(\frac{\delta C_\tau}{\delta c_1}) = 354$;  $F(\frac{\delta C_\tau}{\delta c_2}) = 366$;  $F(\frac{\delta C_\tau}{\delta c_3}) = 390$;  $F(\frac{\delta C_\tau}{\delta c_4}) = 354$;  $F(\frac{\delta C_\tau}{\delta c_5}) = 312$;

$F(\frac{\delta C_\tau}{\delta c_6}) = 312$;  $F(\frac{\delta C_\tau}{\delta c_7}) = 270$;  $F(\frac{\delta C_\tau}{\delta c_8}) = 336$;  $F(\frac{\delta C_\tau}{\delta c_9}) = 408$;  $F(\frac{\delta C_\tau}{\delta c_{10}}) = 408$;

$F(\frac{\delta C_\tau}{\delta c_{11}}) = 372$;  $F(\frac{\delta C_\tau}{\delta c_{12}}) = 300$;  $F(\frac{\delta C_\tau}{\delta c_{13}}) = 258$;  $F(\frac{\delta C_\tau}{\delta c_{14}}) = 336$;  $F(\frac{\delta C_\tau}{\delta c_{15}}) = 360$;

$F(\frac{\delta C_\tau}{\delta c_{16}}) = 408$;  $F(\frac{\delta C_\tau}{\delta c_{17}}) = 420$.

Максимальное изменение значения функционала Маклейна получается после удаления цикла $c_{14}$. Удаляем цикл $c_{14}$. Продолжая процесс, получим:

$F(\frac{\delta^2 C_\tau}{\delta c_{13} \delta c_1}) = 216$;  $F(\frac{\delta^2 C_\tau}{\delta c_{13} \delta c_2}) = 222$;  $F(\frac{\delta^2 C_\tau}{\delta c_{13} \delta c_3}) = 222$;  $F(\frac{\delta^2 C_\tau}{\delta c_{13} \delta c_4}) = 186$;

$F(\frac{\delta^2 C_\tau}{\delta c_{13} \delta c_5}) = 180$;  $F(\frac{\delta^2 C_\tau}{\delta c_{13} \delta c_6}) = 180$;  $F(\frac{\delta^2 C_\tau}{\delta c_{13} \delta c_7}) = 162$;  $F(\frac{\delta^2 C_\tau}{\delta c_{13} \delta c_8}) = 204$;

$F(\frac{\delta^2 C_\tau}{\delta c_{13} \delta c_9}) = 240$;  $F(\frac{\delta^2 C_\tau}{\delta c_{13} \delta c_{10}}) = 240$;  $F(\frac{\delta^2 C_\tau}{\delta c_{13} \delta c_{11}}) = 222$;  $F(\frac{\delta^2 C_\tau}{\delta c_{13} \delta c_{12}}) = 180$;

$$F(\frac{\delta^2 C_\tau}{\delta c_{13}\delta c_{14}}) = 192; \quad F(\frac{\delta^2 C_\tau}{\delta c_{13}\delta c_{15}}) = 216; \quad F(\frac{\delta^2 C_\tau}{\delta c_{13}\delta c_{16}}) = 240; \quad F(\frac{\delta^2 C_\tau}{\delta c_{13}\delta c_{17}}) = 252.$$

Максимальное изменение значения функционала Маклейна получается после удаления цикла $c_7$. Удаляем $c_7$. Продолжая процесс, получим:

$$F(\frac{\delta^3 C_\tau}{\delta c_{13}\delta c_7\delta c_1}) = 132; \quad F(\frac{\delta^3 C_\tau}{\delta c_{13}\delta c_7\delta c_2}) = 144; \quad F(\frac{\delta^3 C_\tau}{\delta c_{13}\delta c_7\delta c_3}) = 126;$$

$$F(\frac{\delta^3 C_\tau}{\delta c_{13}\delta c_7\delta c_4}) = 90; \quad F(\frac{\delta^3 C_\tau}{\delta c_{13}\delta c_7\delta c_5}) = 102; \quad F(\frac{\delta^3 C_\tau}{\delta c_{13}\delta c_7\delta c_6}) = 102;$$

$$F(\frac{\delta^3 C_\tau}{\delta c_{13}\delta c_7\delta c_8}) = 138; \quad F(\frac{\delta^3 C_\tau}{\delta c_{13}\delta c_7\delta c_9}) = 144; \quad F(\frac{\delta^3 C_\tau}{\delta c_{13}\delta c_7\delta c_{10}}) = 144;$$

$$F(\frac{\delta^3 C_\tau}{\delta c_{13}\delta c_7\delta c_{11}}) = 138; \quad F(\frac{\delta^3 C_\tau}{\delta c_{13}\delta c_7\delta c_{12}}) = 108; \quad F(\frac{\delta^3 C_\tau}{\delta c_{13}\delta c_7\delta c_{14}}) = 126;$$

$$F(\frac{\delta^3 C_\tau}{\delta c_{13}\delta c_7\delta c_{15}}) = 138; \quad F(\frac{\delta^3 C_\tau}{\delta c_{13}\delta c_7\delta c_{16}}) = 156; \quad F(\frac{\delta^3 C_\tau}{\delta c_{13}\delta c_7\delta c_{17}}) = 156.$$

Максимальное изменение значения функционала Маклнйна получается после удаления цикла $c_4$. Удаляем цикл $c_4$. Продолжая процесс, получим:

$$F(\frac{\delta^4 C_\tau}{\delta c_{13}\delta c_7\delta c_4\delta c_1}) = 72; \quad F(\frac{\delta^4 C_\tau}{\delta c_{13}\delta c_7\delta c_4\delta c_2}) = 72; \quad F(\frac{\delta^4 C_\tau}{\delta c_{13}\delta c_7\delta c_4\delta c_3}) = 78;$$

$$F(\frac{\delta^4 C_\tau}{\delta c_{13}\delta c_7\delta c_4\delta c_5}) = 54; \quad F(\frac{\delta^4 C_\tau}{\delta c_{13}\delta c_7\delta c_4\delta c_6}) = 54; \quad F(\frac{\delta^4 C_\tau}{\delta c_{13}\delta c_7\delta c_4\delta c_8}) = 66;$$

$$F(\frac{\delta^4 C_\tau}{\delta c_{13}\delta c_7\delta c_4\delta c_9}) = 84; \quad F(\frac{\delta^4 C_\tau}{\delta c_{13}\delta c_7\delta c_4\delta c_{10}}) = 84; \quad F(\frac{\delta^4 C_\tau}{\delta c_{13}\delta c_7\delta c_4\delta c_{11}}) = 78;$$

$$F(\frac{\delta^4 C_\tau}{\delta c_{13}\delta c_7\delta c_4\delta c_{12}}) = 60; \quad F(\frac{\delta^4 C_\tau}{\delta c_{13}\delta c_7\delta c_4\delta c_{14}}) = 54; \quad F(\frac{\delta^4 C_\tau}{\delta c_{13}\delta c_7\delta c_4\delta c_{15}}) = 66;$$

$$F(\frac{\delta^4 C_\tau}{\delta c_{13}\delta c_7\delta c_4\delta c_{16}}) = 84; \quad F(\frac{\delta^4 C_\tau}{\delta c_{13}\delta c_7\delta c_4\delta c_{17}}) = 84.$$

Максимальное изменение значения функционала Маклейна получается после удаления циклов $c_5$, $c_6$ и $c_{14}$. Рассмотрим скорость изменения функционала:

$$F(\frac{\delta^3 C_\tau}{\delta c_{13}\delta c_7\delta c_5}) - F(\frac{\delta^4 C_\tau}{\delta c_{13}\delta c_7\delta c_4\delta c_5}) = 102-54 = 48;$$

$$F(\frac{\delta^3 C_\tau}{\delta c_{13}\delta c_7\delta c_6}) - F(\frac{\delta^4 C_\tau}{\delta c_{13}\delta c_7\delta c_4\delta c_6}) = 102-54 = 48;$$

$$F(\frac{\delta^3 C_\tau}{\delta c_{13}\delta c_7\delta c_{14}}) - F(\frac{\delta^4 C_\tau}{\delta c_{13}\delta c_7\delta c_4\delta c_{14}}) = 126-54 = 72.$$

Для удаления выбираем цикл $c_{14}$ с максимальной скоростью удаления цикла.

$$F(\frac{\delta^5 C_\tau}{\delta c_{13}\delta c_7\delta c_4\delta c_{14}\delta c_1}) = 36; F(\frac{\delta^5 C_\tau}{\delta c_{13}\delta c_7\delta c_4\delta c_{14}\delta c_2}) = 42;$$

$$F(\frac{\delta^5 C_\tau}{\delta c_{13}\delta c_7\delta c_4\delta c_{14}\delta c_3}) = 42; F(\frac{\delta^5 C_\tau}{\delta c_{13}\delta c_7\delta c_4\delta c_{14}\delta c_5}) = 30;$$

$$F(\frac{\delta^5 C_\tau}{\delta c_{13}\delta c_7\delta c_4\delta c_{14}\delta c_6}) = 30; F(\frac{\delta^5 C_\tau}{\delta c_{13}\delta c_7\delta c_4\delta c_{14}\delta c_8}) = 42; F(\frac{\delta^5 C_\tau}{\delta c_{13}\delta c_7\delta c_4\delta c_{14}\delta c_9}) = 48;$$

$$F(\frac{\delta^5 C_\tau}{\delta c_{13}\delta c_7\delta c_4\delta c_{14}\delta c_{10}}) = 48; F(\frac{\delta^5 C_\tau}{\delta c_{13}\delta c_7\delta c_4\delta c_{14}\delta c_{11}}) = 42; F(\frac{\delta^5 C_\tau}{\delta c_{13}\delta c_7\delta c_4\delta c_{14}\delta c_{12}}) = 24;$$

$$F(\frac{\delta^5 C_\tau}{\delta c_{13}\delta c_7\delta c_4\delta c_{14}\delta c_{15}}) = 48; F(\frac{\delta^5 C_\tau}{\delta c_{13}\delta c_7\delta c_4\delta c_{14}\delta c_{16}}) = 54; F(\frac{\delta^5 C_\tau}{\delta c_{13}\delta c_7\delta c_4\delta c_{14}\delta c_{17}}) = 54.$$

Максимальное изменение значения функционала Маклейна получается после удаления цикла $c_{12}$. Удаляем цикл $c_{12}$. Продолжаем процесс:

$$F(\frac{\delta^6 C_\tau}{\delta c_{13}\delta c_7\delta c_4\delta c_{14}\delta c_{12}\delta c_1}) = 18; F(\frac{\delta^6 C_\tau}{\delta c_{13}\delta c_7\delta c_4\delta c_{14}\delta c_{12}\delta c_2}) = 18;$$

$$F(\frac{\delta^6 C_\tau}{\delta c_{13}\delta c_7\delta c_4\delta c_{14}\delta c_{12}\delta c_3}) = 18; F(\frac{\delta^6 C_\tau}{\delta c_{13}\delta c_7\delta c_4\delta c_{14}\delta c_{12}\delta c_5}) = 12;$$

$$F(\frac{\delta^6 C_\tau}{\delta c_{13}\delta c_7\delta c_4\delta c_{14}\delta c_{12}\delta c_6}) = 0; F(\frac{\delta^6 C_\tau}{\delta c_{13}\delta c_7\delta c_4\delta c_{14}\delta c_{12}\delta c_8}) = 24;$$

$$F(\frac{\delta^6 C_\tau}{\delta c_{13}\delta c_7\delta c_4\delta c_{14}\delta c_{12}\delta c_9}) = 24; F(\frac{\delta^6 C_\tau}{\delta c_{13}\delta c_7\delta c_4\delta c_{14}\delta c_{12}\delta c_{10}}) = 18;$$

$$F(\frac{\delta^6 C_\tau}{\delta c_{13}\delta c_7\delta c_4\delta c_{14}\delta c_{12}\delta c_{11}}) = 18; F(\frac{\delta^6 C_\tau}{\delta c_{13}\delta c_7\delta c_4\delta c_{14}\delta c_{12}\delta c_{15}}) = 18;$$

$$F(\frac{\delta^6 C_\tau}{\delta c_{13}\delta c_7\delta c_4\delta c_{14}\delta c_{12}\delta c_{16}}) = 24; F(\frac{\delta^6 C_\tau}{\delta c_{13}\delta c_7\delta c_4\delta c_{14}\delta c_{12}\delta c_{17}}) = 24.$$

Максимальное изменение значения функционала Маклейна получается после удаления цикла $c_6$. Удаляем цикл $c_6$. Останавливаем процесс удаления циклов.

Исключим из множества изометрических циклов циклы $c_{13}, c_{17}, c_4, c_{14}, c_{12}, c_6$. В результате получим следующее подмножество циклов:

цикл $c_1 = \{e_1, e_3, e_5, e_{11}\} \leftrightarrow \{v_1, v_2, v_5, v_6\}$;
цикл $c_2 = \{e_1, e_4, e_6, e_{16}\} \leftrightarrow \{v_1, v_2, v_8, v_9\}$;
цикл $c_3 = \{e_2, e_3, e_7, e_9\} \leftrightarrow \{v_1, v_3, v_4, v_6\}$;
цикл $c_5 = \{e_2, e_4, e_8, e_{14}\} \leftrightarrow \{v_1, v_3, v_7, v_8\}$;
цикл $c_8 = \{e_5, e_6, e_{12}, e_{15}, e_{19}\} \leftrightarrow \{v_2, v_5, v_7, v_9, v_{11}\}$;
цикл $c_9 = \{e_7, e_8, e_{10}\} \leftrightarrow \{v_3, v_4, v_7\}$;
цикл $c_{10} = \{e_9, e_{10}, e_{13}\} \leftrightarrow \{v_4, v_6, v_7\}$;
цикл $c_{11} = \{e_{11}, e_{12}, e_{13}\} \leftrightarrow \{v_5, v_6, v_7\}$;
цикл $c_{15} = \{e_{14}, e_{15}, e_{17}, e_{20}\} \leftrightarrow \{v_7, v_8, v_{10}, v_{11}\}$;
цикл $c_{16} = \{e_{16}, e_{17}, e_{18}\} \leftrightarrow \{v_8, v_9, v_{10}\}$;
цикл $c_{17} = \{e_{18}, e_{19}, e_{20}\} \leftrightarrow \{v_9, v_{10}, v_{11}\}$.

Из случайной последовательности циклов <$c_9,c_{11},c_2,c_1,c_{17},c_8,c_{16},c_{15},c_5,c_{10},c_3$> выберем базис циклов:

цикл $c_1 = \{e_1,e_3,e_5,e_{11}\} \leftrightarrow \{v_1,v_2,v_5,v_6\}$;
цикл $c_2 = \{e_1,e_4,e_6,e_{16}\} \leftrightarrow \{v_1,v_2,v_8,v_9\}$;
цикл $c_3 = \{e_2,e_3,e_7,e_9\} \leftrightarrow \{v_1,v_3,v_4,v_6\}$;
цикл $c_5 = \{e_2,e_4,e_8,e_{14}\} \leftrightarrow \{v_1,v_3,v_7,v_8\}$;
цикл $c_8 = \{e_5,e_6,e_{12},e_{15},e_{19}\} \leftrightarrow \{v_2,v_5,v_7,v_9,v_{11}\}$;
цикл $c_9 = \{e_7,e_8,e_{10}\} \leftrightarrow \{v_3,v_4,v_7\}$;
цикл $c_{10} = \{e_9,e_{10},e_{13}\} \leftrightarrow \{v_4,v_6,v_7\}$;
цикл $c_{11} = \{e_{11},e_{12},e_{13}\} \leftrightarrow \{v_5,v_6,v_7\}$;
цикл $c_{15} = \{e_{14},e_{15},e_{17},e_{20}\} \leftrightarrow \{v_7,v_8,v_{10},v_{11}\}$;
цикл $c_{16} = \{e_{16},e_{17},e_{18}\} \leftrightarrow \{v_8,v_9,v_{10}\}$;
цикл $c_{17} = \{e_{18},e_{19},e_{20}\} \leftrightarrow \{v_9,v_{10},v_{11}\}$.

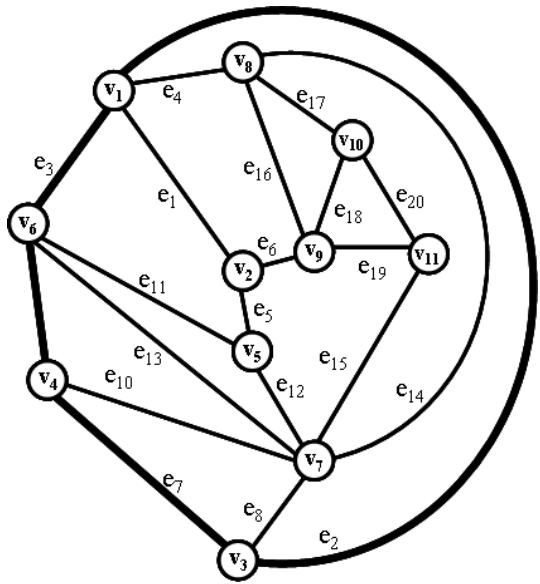

Рис. 12.2. Топологический рисунок плоской части графа

Топологический рисунок плоского суграфа (рис. 12.2):

цикл $c_1 = \{e_1,e_3,e_5,e_{11}\} \leftrightarrow \{v_1,v_2,v_5,v_6\}$;
цикл $c_2 = \{e_1,e_4,e_6,e_{16}\} \leftrightarrow \{v_1,v_2,v_8,v_9\}$;
цикл $c_5 = \{e_2,e_4,e_8,e_{14}\} \leftrightarrow \{v_1,v_3,v_7,v_8\}$;
цикл $c_8 = \{e_5,e_6,e_{12},e_{15},e_{19}\} \leftrightarrow \{v_2,v_5,v_7,v_9,v_{11}\}$;
цикл $c_9 = \{e_7,e_8,e_{10}\} \leftrightarrow \{v_3,v_4,v_7\}$;
цикл $c_{10} = \{e_9,e_{10},e_{13}\} \leftrightarrow \{v_4,v_6,v_7\}$;
цикл $c_{11} = \{e_{11},e_{12},e_{13}\} \leftrightarrow \{v_5,v_6,v_7\}$;
цикл $c_{15} = \{e_{14},e_{15},e_{17},e_{20}\} \leftrightarrow \{v_7,v_8,v_{10},v_{11}\}$;
цикл $c_{16} = \{e_{16},e_{17},e_{18}\} \leftrightarrow \{v_8,v_9,v_{10}\}$;
цикл $c_{17} = \{e_{18},e_{19},e_{20}\} \leftrightarrow \{v_9,v_{10},v_{11}\}$.

Таким образом, метод наискорейшего спуска позволяет определить базис подпространства циклов C(G), состоящий из изометрических циклов с минимальным значением кубического функционала Маклейна.

Причем базис циклов с минимальным значением функционала Маклейна, выбранный методом наискорейшего спуска, совпадает с базисом циклов, отобранным методом Монте-Карло.

Вычислительную сложность метода наискорейшего спуска можно оценить как квадрат количества циклов в графе:

$$\left(\frac{n(n-1)(n-2)}{6} \approx \frac{n^3}{6}\right)^2 \approx \frac{n^6}{36} = .$$

Следовательно, вычислительная сложность равна $o(n^6)$.

### 12.2. Непланарные графы

Будем рассматривать непланарный граф $G_1$ предсталенный на рис. 11.2.

Выделим дерево $T(G_1) = \{e_1,e_2,e_6,e_8,e_{10},e_{11},e_{16},e_{17},e_{18}\}$ в графе $G_1$. Выделим хорды $H(G_1) = \{e_3,e_4,e_5,e_7,e_9,e_{12},e_{13},e_{14},e_{15},e_{19},e_{20}\}$.

Множество изометрических циклов графа $G_1$:
цикл $c_1 = \{e_1,e_2,e_5\} \leftrightarrow \{v_1,v_2,v_3\}$;
цикл $c_2 = \{e_1,e_3,e_6,e_{10}\} \leftrightarrow \{v_1,v_2,v_4,v_5\}$;
цикл $c_3 = \{e_1,e_4,e_7\} \leftrightarrow \{v_1,v_2,v_7\}$;
цикл $c_4 = \{e_2,e_4,e_8\} \leftrightarrow \{v_1,v_3,v_7\}$;
цикл $c_5 = \{e_2,e_3,e_9,e_{12},e_{20}\} \leftrightarrow \{v_1,v_3,v_4,v_9,v_{10}\}$;
цикл $c_6 = \{e_3,e_4,e_{12},e_{18},e_{20}\} \leftrightarrow \{v_1,v_4,v_7,v_{10},v_9\}$;
цикл $c_7 = \{e_5,e_7,e_8\} \leftrightarrow \{v_2,v_3,v_7\}$;
цикл $c_8 = \{e_5,e_6,e_9,e_{13},e_{16}\} \leftrightarrow \{v_2,v_3,v_5,v_9,v_6\}$;
цикл $c_9 = \{e_6,e_7,e_{13},e_{16},e_{18}\} \leftrightarrow \{v_2,v_5,v_7,v_6,v_9\}$;
цикл $c_{10} = \{e_8,e_9,e_{18}\} \leftrightarrow \{v_3,v_7,v_9\}$;
цикл $c_{11} = \{e_{10},e_{11},e_{14}\} \leftrightarrow \{v_4,v_5,v_8\}$;
цикл $c_{12} = \{e_{10},e_{12},e_{13},e_{17}\} \leftrightarrow \{v_4,v_5,v_{10},v_6\}$;
цикл $c_{13} = \{e_{11},e_{12},e_{19}\} \leftrightarrow \{v_4,v_8,v_{10}\}$;
цикл $c_{14} = \{e_{13},e_{14},e_{15}\} \leftrightarrow \{v_5,v_6,v_8\}$;
цикл $c_{15} = \{e_{15},e_{17},e_{19}\} \leftrightarrow \{v_6,v_8,v_{10}\}$;
цикл $c_{16} = \{e_{16},e_{17},e_{20}\} \leftrightarrow \{v_6,v_9,v_{10}\}$.

Запишем матроид изометрических циклов в виде структурного числа $W(C_\tau)$ [2].

| | | | | |
|---|---|---|---|---|
| $e_3$ | $c_2$ | $c_5$ | $c_6$ | |
| $e_4$ | $c_3$ | $c_4$ | $c_6$ | |
| $e_5$ | $c_1$ | $c_7$ | $c_8$ | |
| $e_7$ | $c_3$ | $c_7$ | $c_9$ | |
| $e_9$ | $c_5$ | $c_8$ | $c_{10}$ | |
| $e_{12}$ | $c_5$ | $c_6$ | $c_{12}$ | $c_{13}$ |
| $e_{13}$ | $c_8$ | $c_9$ | $c_{12}$ | $c_{14}$ |
| $e_{14}$ | $c_{11}$ | $c_{14}$ | | |
| $e_{15}$ | $c_{14}$ | $c_{15}$ | | |
| $e_{19}$ | $c_{13}$ | $c_{15}$ | | |

| | | | |
|---|---|---|---|
| $e_{20}$ | $c_5$ | $c_6$ | $c_{16}$ |

Для выбора базиса, сформируем обратные производные первого порядка:

$F(\frac{\delta C_\tau}{\delta c_1}) = 114$; $F(\frac{\delta C_\tau}{\delta c_2}) = 108$; $F(\frac{\delta C_\tau}{\delta c_3}) = 114$; $F(\frac{\delta C_\tau}{\delta c_4}) = 114$; $F(\frac{\delta C_\tau}{\delta c_5}) = 90$;

$F(\frac{\delta C_\tau}{\delta c_6}) = 90$; $F(\frac{\delta C_\tau}{\delta c_7}) = 114$; $F(\frac{\delta C_\tau}{\delta c_8}) = 90$; $F(\frac{\delta C_\tau}{\delta c_9}) = 90$; $F(\frac{\delta C_\tau}{\delta c_{10}}) = 114$;

$F(\frac{\delta C_\tau}{\delta c_{11}}) = 126$; $F(\frac{\delta C_\tau}{\delta c_{12}}) = 84$; $F(\frac{\delta C_\tau}{\delta c_{13}}) = 114$; $F(\frac{\delta C_\tau}{\delta c_{14}}) = 114$; $F(\frac{\delta C_\tau}{\delta c_{15}}) = 126$;

$F(\frac{\delta C_\tau}{\delta c_{16}}) = 114$.

Максимальное изменение значения функционала Маклейна получается после удаления цикла $c_{12}$. Удаляем цикл $c_{12}$.

| | | | |
|---|---|---|---|
| $e_3$ | $c_2$ | $c_5$ | $c_6$ |
| $e_4$ | $c_3$ | $c_4$ | $c_6$ |
| $e_5$ | $c_1$ | $c_7$ | $c_8$ |
| $e_7$ | $c_3$ | $c_7$ | $c_9$ |
| $e_9$ | $c_5$ | $c_8$ | $c_{10}$ |
| $e_{12}$ | $c_5$ | $c_6$ | $c_{13}$ |
| $e_{13}$ | $c_8$ | $c_9$ | $c_{14}$ |
| $e_{14}$ | $c_{11}$ | $c_{14}$ | |
| $e_{15}$ | $c_{14}$ | $c_{15}$ | |
| $e_{19}$ | $c_{13}$ | $c_{15}$ | |
| $e_{20}$ | $c_5$ | $c_6$ | $c_{16}$ |

Продолжая процесс, получим:

$F(\frac{\delta^2 C_\tau}{\delta c_{12} \delta c_1}) = 66$; $F(\frac{\delta^2 C_\tau}{\delta c_{12} \delta c_2}) = 66$; $F(\frac{\delta^2 C_\tau}{\delta c_{12} \delta c_3}) = 66$; $F(\frac{\delta^2 C_\tau}{\delta c_{12} \delta c_4}) = 66$;

$F(\frac{\delta^2 C_\tau}{\delta c_{12} \delta c_5}) = 54$; $F(\frac{\delta^2 C_\tau}{\delta c_{12} \delta c_6}) = 54$; $F(\frac{\delta^2 C_\tau}{\delta c_{12} \delta c_7}) = 66$; $F(\frac{\delta^2 C_\tau}{\delta c_{12} \delta c_8}) = 54$;

$F(\frac{\delta^2 C_\tau}{\delta c_{12} \delta c_9}) = 54$; $F(\frac{\delta^2 C_\tau}{\delta c_{12} \delta c_{10}}) = 66$; $F(\frac{\delta^2 C_\tau}{\delta c_{12} \delta c_{11}}) = 84$; $F(\frac{\delta^2 C_\tau}{\delta c_{12} \delta c_{13}}) = 78$;

$F(\frac{\delta^2 C_\tau}{\delta c_{12} \delta c_{14}}) = 78$; $F(\frac{\delta^2 C_\tau}{\delta c_{12} \delta c_{15}}) = 84$; $F(\frac{\delta^2 C_\tau}{\delta c_{12} \delta c_{16}}) = 72$.

Здесь имеем четыре минимальных значения функционала Маклейна с одинаковой скоростью изменения

$$F(\frac{\delta C_\tau}{\delta c_5}) - F(\frac{\delta^2 C_\tau}{\delta c_{12}\delta c_5}) = 90 - 54 = 36;$$

$$F(\frac{\delta C_\tau}{\delta c_6}) - F(\frac{\delta^2 C_\tau}{\delta c_{12}\delta c_6}) = 90 - 54 = 36;$$

$$F(\frac{\delta C_\tau}{\delta c_8}) - F(\frac{\delta^2 C_\tau}{\delta c_{12}\delta c_8}) = 90 - 54 = 36;$$

$$F(\frac{\delta C_\tau}{\delta c_9}) - F(\frac{\delta^2 C_\tau}{\delta c_{12}\delta c_9}) = 90 - 54 = 36.$$

Удаляем цикл $c_9$.

| | | | |
|---|---|---|---|
| $e_3$ | $c_2$ | $c_5$ | $c_6$ |
| $e_4$ | $c_3$ | $c_4$ | $c_6$ |
| $e_5$ | $c_1$ | $c_7$ | $c_8$ |
| $e_7$ | $c_3$ | $c_7$ | |
| $e_9$ | $c_5$ | $c_8$ | $c_{10}$ |
| $e_{12}$ | $c_5$ | $c_6$ | $c_{13}$ |
| $e_{13}$ | $c_8$ | $c_{14}$ | |
| $e_{14}$ | $c_{11}$ | $c_{14}$ | |
| $e_{15}$ | $c_{14}$ | $c_{15}$ | |
| $e_{19}$ | $c_{13}$ | $c_{15}$ | |
| $e_{20}$ | $c_5$ | $c_6$ | $c_{16}$ |

Продолжая процесс, получим:

$$F(\frac{\delta^3 C_\tau}{\delta c_{12}\delta c_9\delta c_1}) = 36;\ F(\frac{\delta^3 C_\tau}{\delta c_{12}\delta c_9\delta c_2}) = 42;\ F(\frac{\delta^3 C_\tau}{\delta c_{12}\delta c_9\delta c_3}) = 42;$$

$$F(\frac{\delta^3 C_\tau}{\delta c_{12}\delta c_9\delta c_4}) = 36;\ F(\frac{\delta^3 C_\tau}{\delta c_{12}\delta c_9\delta c_5}) = 24;\ F(\frac{\delta^3 C_\tau}{\delta c_{12}\delta c_9\delta c_6}) = 30;$$

$$F(\frac{\delta^3 C_\tau}{\delta c_{12}\delta c_9\delta c_8}) = 42;\ F(\frac{\delta^3 C_\tau}{\delta c_{12}\delta c_9\delta c_{10}}) = 42;\ F(\frac{\delta^3 C_\tau}{\delta c_{12}\delta c_9\delta c_{11}}) = 42;$$

$$F(\frac{\delta^3 C_\tau}{\delta c_{12}\delta c_9\delta c_{12}}) = 54;\ F(\frac{\delta^3 C_\tau}{\delta c_{12}\delta c_9\delta c_{13}}) = 48;\ F(\frac{\delta^3 C_\tau}{\delta c_{12}\delta c_9\delta c_{14}}) = 54;$$

$$F(\frac{\delta^3 C_\tau}{\delta c_{12}\delta c_9\delta c_{15}}) = 48;\ F(\frac{\delta^3 C_\tau}{\delta c_{12}\delta c_9\delta c_{16}}) = 54.$$

Максимальное изменение значения функционала Маклейна получается после удаления цикла $c_5$. Удаляем цикл $c_5$.

| | | | |
|---|---|---|---|
| $e_3$ | $c_2$ | $c_6$ | |
| $e_4$ | $c_3$ | $c_4$ | $c_6$ |
| $e_5$ | $c_1$ | $c_7$ | $c_8$ |

| | | |
|---|---|---|
| $e_7$ | $c_3$ | $c_7$ |
| $e_9$ | $c_8$ | $c_{10}$ |
| $e_{12}$ | $c_6$ | $c_{13}$ |
| $e_{13}$ | $c_8$ | $c_{14}$ |
| $e_{14}$ | $c_{11}$ | $c_{14}$ |
| $e_{15}$ | $c_{14}$ | $c_{15}$ |
| $e_{19}$ | $c_{13}$ | $c_{15}$ |
| $e_{20}$ | $c_6$ | $c_{16}$ |

Продолжая процесс, получим:

$F(\dfrac{\delta^4 C_\tau}{\delta c_{12} \delta c_9 \delta c_5 \delta c_1}) = 12;\ F(\dfrac{\delta^4 C_\tau}{\delta c_{12} \delta c_9 \delta c_5 \delta c_2}) = 18;\ F(\dfrac{\delta^4 C_\tau}{\delta c_{12} \delta c_9 \delta c_5 \delta c_3}) = 12;$

$F(\dfrac{\delta^4 C_\tau}{\delta c_{12} \delta c_9 \delta c_5 \delta c_4}) = 12;\ F(\dfrac{\delta^4 C_\tau}{\delta c_{12} \delta c_9 \delta c_5 \delta c_6}) = 18;\ F(\dfrac{\delta^4 C_\tau}{\delta c_{12} \delta c_9 \delta c_5 \delta c_7}) = 12;$

$F(\dfrac{\delta^4 C_\tau}{\delta c_{12} \delta c_9 \delta c_5 \delta c_8}) = 18;\ F(\dfrac{\delta^4 C_\tau}{\delta c_{12} \delta c_9 \delta c_5 \delta c_{10}}) = 18;\ F(\dfrac{\delta^4 C_\tau}{\delta c_{12} \delta c_9 \delta c_5 \delta c_{11}}) = 24;$

$F(\dfrac{\delta^4 C_\tau}{\delta c_{12} \delta c_9 \delta c_5 \delta c_{13}}) = 24;\ F(\dfrac{\delta^4 C_\tau}{\delta c_{12} \delta c_9 \delta c_5 \delta c_{14}}) = 24;\ F(\dfrac{\delta^4 C_\tau}{\delta c_{12} \delta c_9 \delta c_5 \delta c_{15}}) = 24;$

$F(\dfrac{\delta^4 C_\tau}{\delta c_{12} \delta c_9 \delta c_5 \delta c_{16}}) = 24.$

В результате имеем четыре минимальных значения функционала Маклейна.

$F(\dfrac{\delta^3 C_\tau}{\delta c_{12} \delta c_9 \delta c_1}) - F(\dfrac{\delta^4 C_\tau}{\delta c_{12} \delta c_9 \delta c_5 \delta c_1}) = 36 - 12 = 24$ ;

$F(\dfrac{\delta^3 C_\tau}{\delta c_{12} \delta c_9 \delta c_3}) - F(\dfrac{\delta^4 C_\tau}{\delta c_{12} \delta c_9 \delta c_5 \delta c_3}) = 42 - 12 = 30$ ;

$F(\dfrac{\delta^3 C_\tau}{\delta c_{12} \delta c_9 \delta c_4}) - F(\dfrac{\delta^4 C_\tau}{\delta c_{12} \delta c_9 \delta c_5 \delta c_4}) = 36 - 12 = 24$ ;

$F(\dfrac{\delta^3 C_\tau}{\delta c_{12} \delta c_9 \delta c_7}) - F(\dfrac{\delta^4 C_\tau}{\delta c_{12} \delta c_9 \delta c_5 \delta c_7}) = 42 - 12 = 30$ ;

Удаляем цикл $c_7$.

| | | | |
|---|---|---|---|
| $e_3$ | $c_2$ | $c_6$ | |
| $e_4$ | $c_3$ | $c_4$ | $c_6$ |
| $e_5$ | $c_1$ | $c_8$ | |
| $e_7$ | $c_3$ | | |
| $e_9$ | $c_8$ | $c_{10}$ | |
| $e_{12}$ | $c_6$ | $c_{13}$ | |

| | | |
|---|---|---|
| $e_{13}$ | $c_8$ | $c_{14}$ |
| $e_{14}$ | $c_{11}$ | $c_{14}$ |
| $e_{15}$ | $c_{14}$ | $c_{15}$ |
| $e_{19}$ | $c_{13}$ | $c_{15}$ |
| $e_{20}$ | $c_6$ | $c_{16}$ |

Продолжая процесс, получим:

$$F(\frac{\delta^5 C_\tau}{\delta c_{12} \delta c_9 \delta c_5 \delta c_7 \delta c_1}) = 6; \quad F(\frac{\delta^5 C_\tau}{\delta c_{12} \delta c_9 \delta c_5 \delta c_7 \delta c_2}) = 6;$$

$$F(\frac{\delta^5 C_\tau}{\delta c_{12} \delta c_9 \delta c_5 \delta c_7 \delta c_3}) = 0; \quad F(\frac{\delta^5 C_\tau}{\delta c_{12} \delta c_9 \delta c_5 \delta c_7 \delta c_4}) = 6;$$

$$F(\frac{\delta^5 C_\tau}{\delta c_{12} \delta c_9 \delta c_5 \delta c_7 \delta c_6}) = 6; \quad F(\frac{\delta^5 C_\tau}{\delta c_{12} \delta c_9 \delta c_5 \delta c_7 \delta c_8}) = 12; \quad F(\frac{\delta^5 C_\tau}{\delta c_{12} \delta c_9 \delta c_5 \delta c_7 \delta c_{10}}) = 12;$$

$$F(\frac{\delta^5 C_\tau}{\delta c_{12} \delta c_9 \delta c_5 \delta c_7 \delta c_{11}}) = 12; \quad F(\frac{\delta^5 C_\tau}{\delta c_{12} \delta c_9 \delta c_5 \delta c_7 \delta c_{13}}) = 12; \quad F(\frac{\delta^5 C_\tau}{\delta c_{12} \delta c_9 \delta c_5 \delta c_7 \delta c_{14}}) = 12;$$

$$F(\frac{\delta^5 C_\tau}{\delta c_{12} \delta c_9 \delta c_5 \delta c_7 \delta c_{15}}) = 12; \quad F(\frac{\delta^5 C_\tau}{\delta c_{12} \delta c_9 \delta c_5 \delta c_7 \delta c_{16}}) = 12.$$

Удаление цикла $c_3$ приводит к удалению ребра $e_7$, что не допустимо. В базисе должны быть все ребра и все вершины. Тогда имеем четыре минимальных значения функционала Маклейна.

$$F(\frac{\delta^4 C_\tau}{\delta c_{12} \delta c_9 \delta c_5 \delta c_1}) - F(\frac{\delta^5 C_\tau}{\delta c_{12} \delta c_9 \delta c_5 \delta c_7 \delta c_1}) = 12 - 6 = 6;$$

$$F(\frac{\delta^4 C_\tau}{\delta c_{12} \delta c_9 \delta c_5 \delta c_2}) - F(\frac{\delta^5 C_\tau}{\delta c_{12} \delta c_9 \delta c_5 \delta c_7 \delta c_2}) = 18 - 6 = 12;$$

$$F(\frac{\delta^4 C_\tau}{\delta c_{12} \delta c_9 \delta c_5 \delta c_4}) - F(\frac{\delta^5 C_\tau}{\delta c_{12} \delta c_9 \delta c_5 \delta c_7 \delta c_4}) = 12 - 6 = 6;$$

$$F(\frac{\delta^4 C_\tau}{\delta c_{12} \delta c_9 \delta c_5 \delta c_6}) - F(\frac{\delta^5 C_\tau}{\delta c_{12} \delta c_9 \delta c_5 \delta c_7 \delta c_6}) = 18 - 6 = 12.$$

Удаляем цикл $c_6$.

| | | |
|---|---|---|
| $e_3$ | $c_2$ | |
| $e_4$ | $c_3$ | $c_4$ |
| $e_5$ | $c_1$ | $c_8$ |
| $e_7$ | $c_3$ | |
| $e_9$ | $c_8$ | $c_{10}$ |
| $e_{12}$ | $c_{13}$ | |
| $e_{13}$ | $c_8$ | $c_{14}$ |
| $e_{14}$ | $c_{11}$ | $c_{14}$ |
| $e_{15}$ | $c_{14}$ | $c_{15}$ |
| $e_{19}$ | $c_{13}$ | $c_{15}$ |
| $e_{20}$ | $c_{16}$ | |

Количество оставшихся циклов равно цикломатическому числу. Следовательно, это базис. Методом «прыгающая строка» проверим систему циклов на независимость.

| | |
|---|---|
| $e_3$ | $c_2$ |
| $e_4$ | $c_4$ |
| $e_5$ | $c_1$ |
| $e_7$ | $c_3$ |
| $e_9$ | $c_{10}$ |
| $e_{12}$ | $c_{13}$ |
| $e_{13}$ | $c_8$ |
| $e_{14}$ | $c_{11}$ |
| $e_{15}$ | $c_{14}$ |
| $e_{19}$ | $c_{15}$ |
| $e_{20}$ | $c_{16}$ |

Количество столбцов не кратно 2, следовательно, система независима.

$$F\left(\frac{\delta^5 C_\tau}{\delta c_{12} \delta c_9 \delta c_5 \delta c_7 \delta c_6}\right) = F(C_b) = 6.$$

Далее, будем удалять циклы с соблюдением условия Эйлера.

$$\frac{\partial C_b}{\partial c_2} = 0; \quad \frac{\partial C_b}{\partial c_3} = 0.$$

При удалении цикла $c_2$ удаляется ребро $e_3$ (рис. 12.3), а при удалении цикла $c_3$ – ребро $e_7$ (рис. 12.4).

Соответственно, в структурном числе $W(C_b)$ следует удалить либо цикл $c_2$ совместно с ребром $e_3$, либо удалить цикл $c_3$ совместно с ребром $e_7$.

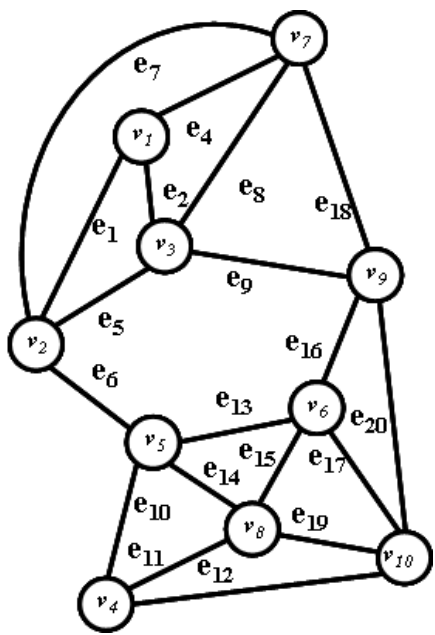
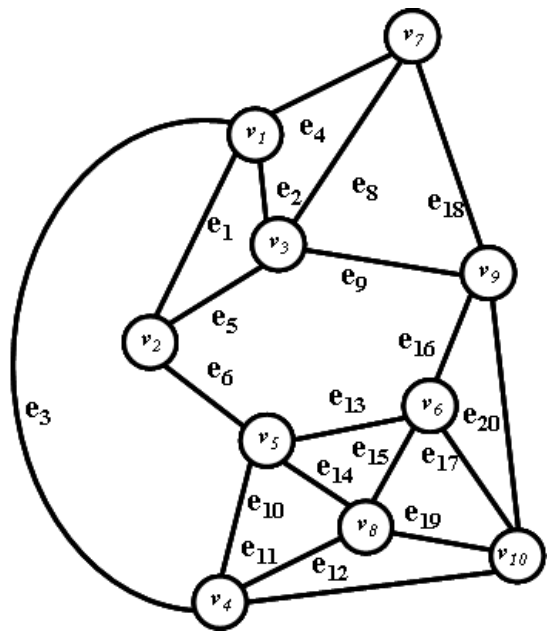

Рис. 12.3. Граф $G_1$ без ребра $e_3$.      Рис. 12.4. Граф $G_1$ без ребра $e_7$.

В качестве примера рассмотрим еще олин пример выделения плоской части графа, состоящей из изометрических циклов. Пусть задан граф $G_3$:

количество вершин графа = 12;
количество ребер графа = 24;
количество единичных циклов = 21;
цикломатическое число графа = 13.

Смежность графа:

вершина $v_1$:  $v_2$  $v_4$  $v_{11}$  $v_{12}$           вершина $v_2$:  $v_1$  $v_3$  $v_5$  $v_{12}$
вершина $v_3$:  $v_2$  $v_4$  $v_{10}$  $v_{12}$           вершина $v_4$:  $v_1$  $v_3$  $v_5$  $v_6$
вершина $v_5$:  $v_2$  $v_4$  $v_6$  $v_8$  $v_{10}$       вершина $v_6$:  $v_4$  $v_5$  $v_7$
вершина $v_7$:  $v_6$  $v_8$  $v_9$                         вершина $v_8$:  $v_5$  $v_7$  $v_9$  $v_{12}$
вершина $v_9$:  $v_7$  $v_8$  $v_{10}$  $v_{11}$           вершина $v_{10}$: $v_3$  $v_5$  $v_9$  $v_{11}$
вершина $v_{11}$: $v_1$  $v_9$  $v_{10}$  $v_{12}$         вершина $v_{12}$: $v_1$  $v_2$  $v_3$  $v_8$  $v_{11}$

Инцидентность графа:

ребро  $e_1$:  $(v_1,v_2)$ или $(v_2,v_1)$;                 ребро  $e_2$:  $(v_1,v_4)$ или $(v_4,v_1)$;
ребро  $e_3$:  $(v_1,v_{11})$ или $(v_{11},v_1)$;          ребро  $e_4$:  $(v_1,v_{12})$ или $(v_{12},v_1)$;
ребро  $e_5$:  $(v_2,v_3)$ или $(v_3,v_2)$;                 ребро  $e_6$:  $(v_2,v_5)$ или $(v_5,v_2)$;
ребро  $e_7$:  $(v_2,v_{12})$ или $(v_{12},v_2)$;          ребро  $e_8$:  $(v_3,v_4)$ или $(v_4,v_3)$;
ребро  $e_9$:  $(v_3,v_{10})$ или $(v_{10},v_3)$;          ребро  $e_{10}$: $(v_3,v_{12})$ или $(v_{12},v_3)$;
ребро  $e_{11}$: $(v_4,v_5)$ или $(v_5,v_4)$;              ребро  $e_{12}$: $(v_4,v_6)$ или $(v_6,v_4)$;
ребро  $e_{13}$: $(v_5,v_6)$ или $(v_6,v_5)$;              ребро  $e_{14}$: $(v_5,v_8)$ или $(v_8,v_5)$;
ребро  $e_{15}$: $(v_5,v_{10})$ или $(v_{10},v_5)$;        ребро  $e_{16}$: $(v_6,v_7)$ или $(v_7,v_6)$;
ребро  $e_{17}$: $(v_7,v_8)$ или $(v_8,v_7)$;              ребро  $e_{18}$: $(v_7,v_9)$ или $(v_9,v_7)$;
ребро  $e_{19}$: $(v_8,v_9)$ или $(v_9,v_8)$;              ребро  $e_{20}$: $(v_8,v_{12})$ или $(v_{12},v_8)$;
ребро  $e_{21}$: $(v_9,v_{10})$ или $(v_{10},v_9)$;        ребро  $e_{22}$: $(v_9,v_{11})$ или $(v_{11},v_9)$;
ребро  $e_{23}$: $(v_{10},v_{11})$ или $(v_{11},v_{10})$;  ребро  $e_{24}$: $(v_{11},v_{12})$ или $(v_{12},v_{11})$.

Множество изометрических циклов графа:

цикл   $c_1 = \{e_1,e_2,e_5,e_8\} \leftrightarrow \{v_1,v_2,v_3,v_4\}$;
цикл   $c_2 = \{e_1,e_2,e_6,e_{11}\} \leftrightarrow \{v_1,v_2,v_4,v_5\}$;
цикл   $c_3 = \{e_1,e_4,e_7\} \leftrightarrow \{v_1,v_2,v_{12}\}$;
цикл   $c_4 = \{e_2,e_4,e_8,e_{10}\} \leftrightarrow \{v_1,v_3,v_4,v_{12}\}$;
цикл   $c_5 = \{e_2,e_3,e_8,e_9,e_{23}\} \leftrightarrow \{v_1,v_3,v_4,v_{10},v_{11}\}$;
цикл   $c_6 = \{e_2,e_3,e_{11},e_{15},e_{23}\} \leftrightarrow \{v_1,v_4,v_5,v_{10},v_{11}\}$;
цикл   $c_7 = \{e_3,e_4,e_{24}\} \leftrightarrow \{v_1,v_{11},v_{12}\}$;
цикл   $c_8 = \{e_5,e_6,e_8,e_{11}\} \leftrightarrow \{v_2,v_3,v_4,v_5\}$;
цикл   $c_9 = \{e_5,e_6,e_9,e_{15}\} \leftrightarrow \{v_2,v_3,v_5,v_{10}\}$;
цикл   $c_{10} = \{e_5,e_7,e_{10}\} \leftrightarrow \{v_2,v_3,v_{12}\}$;
цикл   $c_{11} = \{e_6,e_7,e_{14},e_{20}\} \leftrightarrow \{v_2,v_5,v_8,v_{12}\}$;
цикл   $c_{12} = \{e_8,e_9,e_{11},e_{15}\} \leftrightarrow \{v_3,v_4,v_5,v_{10}\}$;
цикл   $c_{13} = \{e_9,e_{10},e_{23},e_{24}\} \leftrightarrow \{v_3,v_{10},v_{11},v_{12}\}$;
цикл   $c_{14} = \{e_{11},e_{12},e_{13}\} \leftrightarrow \{v_4,v_5,v_6\}$;
цикл   $c_{15} = \{e_{13},e_{14},e_{16},e_{17}\} \leftrightarrow \{v_5,v_6,v_7,v_8\}$;
цикл   $c_{16} = \{e_{14},e_{15},e_{19},e_{21}\} \leftrightarrow \{v_5,v_8,v_9,v_{10}\}$;
цикл   $c_{17} = \{e_{17},e_{18},e_{19}\} \leftrightarrow \{v_7,v_8,v_9\}$;
цикл   $c_{18} = \{e_{13},e_{15},e_{16},e_{18},e_{21}\} \leftrightarrow \{v_5,v_6,v_7,v_9,v_{10}\}$;
цикл   $c_{19} = \{e_{19},e_{20},e_{22},e_{24}\} \leftrightarrow \{v_8,v_9,v_{11},v_{12}\}$;
цикл   $c_{20} = \{e_{21},e_{22},e_{23}\} \leftrightarrow \{v_9,v_{10},v_{11}\}$;
цикл   $c_{21} = \{e_1,e_3,e_6,e_{15},e_{23}\} \leftrightarrow \{v_1,v_2,v_5,v_{10},v_{11}\}$.

Функциона Маклейна системы имеет значение равное 564.

$F(\frac{\delta C_\tau}{\delta c_1}) = 456$;  $F(\frac{\delta C_\tau}{\delta c_2}) = 438$;  $F(\frac{\delta C_\tau}{\delta c_3}) = 534$;  $F(\frac{\delta C_\tau}{\delta c_4}) = 480$;  $F(\frac{\delta C_\tau}{\delta c_5}) = 420$;

$F(\frac{\delta C_\tau}{\delta c_6}) = 378$;  $F(\frac{\delta C_\tau}{\delta c_7}) = 534$;  $F(\frac{\delta C_\tau}{\delta c_8}) = 438$;  $F(\frac{\delta C_\tau}{\delta c_9}) = 432$;  $F(\frac{\delta C_\tau}{\delta c_{10}}) = 534$;

$F(\frac{\delta C_\tau}{\delta c_{11}}) = 516$;  $F(\frac{\delta C_\tau}{\delta c_{12}}) = 414$;  $F(\frac{\delta C_\tau}{\delta c_{13}}) = 498$;  $F(\frac{\delta C_\tau}{\delta c_{14}}) = 522$;  $F(\frac{\delta C_\tau}{\delta c_{15}}) = 552$;

$F(\frac{\delta C_\tau}{\delta c_{16}}) = 486$;  $F(\frac{\delta C_\tau}{\delta c_{17}}) = 558$;  $F(\frac{\delta C_\tau}{\delta c_{18}}) = 492$;  $F(\frac{\delta C_\tau}{\delta c_{19}}) = 552$;  $F(\frac{\delta C_\tau}{\delta c_{20}}) = 522$;

$F(\frac{\delta C_\tau}{\delta c_{21}}) = 396$.

После удаления цикла $c_6$ функционал Маклейна принимает минимальное значение. Удаляем из системы цикл $c_6$ и продолжаем процесс удаления циклов далее.

$F(\frac{\delta^2 C_\tau}{\delta c_6 \delta c_1}) = 288$;  $F(\frac{\delta^2 C_\tau}{\delta c_6 \delta c_2}) = 288$;  $F(\frac{\delta^2 C_\tau}{\delta c_6 \delta c_3}) = 348$;  $F(\frac{\delta^2 C_\tau}{\delta c_6 \delta c_4}) = 312$;

$F(\frac{\delta^2 C_\tau}{\delta c_6 \delta c_5}) = 282$;  $F(\frac{\delta^2 C_\tau}{\delta c_6 \delta c_7}) = 360$;  $F(\frac{\delta^2 C_\tau}{\delta c_6 \delta c_8}) = 270$;  $F(\frac{\delta^2 C_\tau}{\delta c_6 \delta c_9}) = 270$;

$F(\frac{\delta^2 C_\tau}{\delta c_6 \delta c_{10}}) = 348$;  $F(\frac{\delta^2 C_\tau}{\delta c_6 \delta c_{11}}) = 330$;  $F(\frac{\delta^2 C_\tau}{\delta c_6 \delta c_{12}}) = 270$;  $F(\frac{\delta^2 C_\tau}{\delta c_6 \delta c_{13}}) = 330$;

$$F(\frac{\delta^2 C_\tau}{\delta c_6 \delta c_{14}}) = 354; \quad F(\frac{\delta^2 C_\tau}{\delta c_6 \delta c_{15}}) = 366; \quad F(\frac{\delta^2 C_\tau}{\delta c_6 \delta c_{16}}) = 324; \quad F(\frac{\delta^2 C_\tau}{\delta c_6 \delta c_{17}}) = 372;$$

$$F(\frac{\delta^2 C_\tau}{\delta c_6 \delta c_{18}}) = 330; \quad F(\frac{\delta^2 C_\tau}{\delta c_6 \delta c_{19}}) = 366; \quad F(\frac{\delta^2 C_\tau}{\delta c_6 \delta c_{20}}) = 354; \quad F(\frac{\delta^2 C_\tau}{\delta c_6 \delta c_{21}}) = 264.$$

После удаления цикла $c_{21}$ функционал Маклейна принимает минимальное значение. Удаляем из системы цикл $c_{21}$. Продолжаем процесс удаления циклов.

$$F(\frac{\delta^3 C_\tau}{\delta c_6 \delta c_{21} \delta c_1}) = 186; \quad F(\frac{\delta^3 C_\tau}{\delta c_6 \delta c_{21} \delta c_2}) = 204; \quad F(\frac{\delta^3 C_\tau}{\delta c_6 \delta c_{21} \delta c_3}) = 246; \quad F(\frac{\delta^3 C_\tau}{\delta c_6 \delta c_{21} \delta c_4}) = 198;$$

$$F(\frac{\delta^3 C_\tau}{\delta c_6 \delta c_{21} \delta c_5}) = 186; \quad F(\frac{\delta^3 C_\tau}{\delta c_6 \delta c_{21} \delta c_7}) = 252; \quad F(\frac{\delta^3 C_\tau}{\delta c_6 \delta c_{21} \delta c_8}) = 174; \quad F(\frac{\delta^3 C_\tau}{\delta c_6 \delta c_{21} \delta c_9}) = 192;$$

$$F(\frac{\delta^3 C_\tau}{\delta c_6 \delta c_{21} \delta c_{10}}) = 234; \quad F(\frac{\delta^3 C_\tau}{\delta c_6 \delta c_{21} \delta c_{11}}) = 234; \quad F(\frac{\delta^3 C_\tau}{\delta c_6 \delta c_{21} \delta c_{12}}) = 174;$$

$$F(\frac{\delta^3 C_\tau}{\delta c_6 \delta c_{21} \delta c_{13}}) = 228; \quad F(\frac{\delta^3 C_\tau}{\delta c_6 \delta c_{21} \delta c_{14}}) = 240; \quad F(\frac{\delta^3 C_\tau}{\delta c_6 \delta c_{21} \delta c_{15}}) = 252;$$

$$F(\frac{\delta^3 C_\tau}{\delta c_6 \delta c_{21} \delta c_{16}}) = 228; \quad F(\frac{\delta^3 C_\tau}{\delta c_6 \delta c_{21} \delta c_{17}}) = 258; \quad F(\frac{\delta^3 C_\tau}{\delta c_6 \delta c_{21} \delta c_{18}}) = 234;$$

$$F(\frac{\delta^3 C_\tau}{\delta c_6 \delta c_{21} \delta c_{19}}) = 252; \quad F(\frac{\delta^3 C_\tau}{\delta c_6 \delta c_{21} \delta c_{20}}) = 252.$$

Имеем два минимальных значения функционала Маклейна:

$$F(\frac{\delta^2 C_\tau}{\delta c_6 \delta c_8}) - F(\frac{\delta^3 C_\tau}{\delta c_6 \delta c_{21} \delta c_8}) = 270-174 = 96;$$

$$F(\frac{\delta^2 C_\tau}{\delta c_6 \delta c_{12}}) - F(\frac{\delta^3 C_\tau}{\delta c_6 \delta c_{21} \delta c_{12}}) = 270-174 = 96.$$

Удаляем цикл $c_{12}$. Продолжаем процесс удаления циклов.

$$F(\frac{\delta^4 C_\tau}{\delta c_6 \delta c_{21} \delta c_{12} \delta c_1}) = 114; \quad F(\frac{\delta^4 C_\tau}{\delta c_6 \delta c_{21} \delta c_{12} \delta c_2}) = 126; \quad F(\frac{\delta^4 C_\tau}{\delta c_6 \delta c_{21} \delta c_{12} \delta c_3}) = 156;$$

$$F(\frac{\delta^4 C_\tau}{\delta c_6 \delta c_{21} \delta c_{12} \delta c_4}) = 126; \quad F(\frac{\delta^4 C_\tau}{\delta c_6 \delta c_{21} \delta c_{12} \delta c_5}) = 126; \quad F(\frac{\delta^4 C_\tau}{\delta c_6 \delta c_{21} \delta c_{12} \delta c_7}) = 162;$$

$$F(\frac{\delta^4 C_\tau}{\delta c_6 \delta c_{21} \delta c_{12} \delta c_8}) = 114; \quad F(\frac{\delta^4 C_\tau}{\delta c_6 \delta c_{21} \delta c_{12} \delta c_9}) = 126; \quad F(\frac{\delta^4 C_\tau}{\delta c_6 \delta c_{21} \delta c_{12} \delta c_{10}}) = 144;$$

$$F(\frac{\delta^4 C_\tau}{\delta c_6 \delta c_{21} \delta c_{12} \delta c_{11}}) = 144; \quad F(\frac{\delta^4 C_\tau}{\delta c_6 \delta c_{21} \delta c_{12} \delta c_{13}}) = 150; \quad F(\frac{\delta^4 C_\tau}{\delta c_6 \delta c_{21} \delta c_{12} \delta c_{14}}) = 162;$$

$$F(\frac{\delta^4 C_\tau}{\delta c_6 \delta c_{21} \delta c_{12} \delta c_{15}}) = 162;\ F(\frac{\delta^4 C_\tau}{\delta c_6 \delta c_{21} \delta c_{12} \delta c_{16}}) = 150;\ F(\frac{\delta^4 C_\tau}{\delta c_6 \delta c_{21} \delta c_{12} \delta c_{17}}) = 168;$$

$$F(\frac{\delta^4 C_\tau}{\delta c_6 \delta c_{21} \delta c_{12} \delta c_{18}}) = 156;\ F(\frac{\delta^4 C_\tau}{\delta c_6 \delta c_{21} \delta c_{12} \delta c_{19}}) = 162;\ F(\frac{\delta^4 C_\tau}{\delta c_6 \delta c_{21} \delta c_{12} \delta c_{20}}) = 162.$$

Снова имеем два минимальных значения функционала Маклейна:

$$F(\frac{\delta^3 C_\tau}{\delta c_6 \delta c_{21} \delta c_1}) - F(\frac{\delta^4 C_\tau}{\delta c_6 \delta c_{21} \delta c_{12} \delta c_1}) = 186 - 114 = 72;$$

$$F(\frac{\delta^3 C_\tau}{\delta c_6 \delta c_{21} \delta c_8}) - F(\frac{\delta^4 C_\tau}{\delta c_6 \delta c_{21} \delta c_{12} \delta c_8}) = 174 - 114 = 60.$$

Удаляем цикл $c_1$. Продолжаем процесс удаления циклов.

$$F(\frac{\delta^5 C_\tau}{\delta c_6 \delta c_{21} \delta c_{12} \delta c_1 \delta c_2}) = 84;\ F(\frac{\delta^5 C_\tau}{\delta c_6 \delta c_{21} \delta c_{12} \delta c_1 \delta c_3}) = 102;\ F(\frac{\delta^5 C_\tau}{\delta c_6 \delta c_{21} \delta c_{12} \delta c_1 \delta c_4}) = 90;$$

$$F(\frac{\delta^5 C_\tau}{\delta c_6 \delta c_{21} \delta c_{12} \delta c_1 \delta c_5}) = 90;\ F(\frac{\delta^5 C_\tau}{\delta c_6 \delta c_{21} \delta c_{12} \delta c_1 \delta c_7}) = 102;\ F(\frac{\delta^5 C_\tau}{\delta c_6 \delta c_{21} \delta c_{12} \delta c_1 \delta c_8}) = 78;$$

$$F(\frac{\delta^5 C_\tau}{\delta c_6 \delta c_{21} \delta c_{12} \delta c_1 \delta c_9}) = 78;\ F(\frac{\delta^5 C_\tau}{\delta c_6 \delta c_{21} \delta c_{12} \delta c_1 \delta c_{10}}) = 96;\ F(\frac{\delta^5 C_\tau}{\delta c_6 \delta c_{21} \delta c_{12} \delta c_1 \delta c_{11}}) = 84;$$

$$F(\frac{\delta^5 C_\tau}{\delta c_6 \delta c_{21} \delta c_{12} \delta c_1 \delta c_{13}}) = 90;\ F(\frac{\delta^5 C_\tau}{\delta c_6 \delta c_{21} \delta c_{12} \delta c_1 \delta c_{14}}) = 102;\ F(\frac{\delta^5 C_\tau}{\delta c_6 \delta c_{21} \delta c_{12} \delta c_1 \delta c_{15}}) = 102;$$

$$F(\frac{\delta^5 C_\tau}{\delta c_6 \delta c_{21} \delta c_{12} \delta c_1 \delta c_{16}}) = 90;\ F(\frac{\delta^5 C_\tau}{\delta c_6 \delta c_{21} \delta c_{12} \delta c_1 \delta c_{17}}) = 108;\ F(\frac{\delta^5 C_\tau}{\delta c_6 \delta c_{21} \delta c_{12} \delta c_1 \delta c_{18}}) = 96.$$

$$F(\frac{\delta^5 C_\tau}{\delta c_6 \delta c_{21} \delta c_{12} \delta c_1 \delta c_{19}}) = 102;\ F(\frac{\delta^5 C_\tau}{\delta c_6 \delta c_{21} \delta c_{12} \delta c_1 \delta c_{20}}) = 102.$$

Вновь имеем два минимальных значения функционала Маклейна:

$$F(\frac{\delta^4 C_\tau}{\delta c_6 \delta c_{21} \delta c_{12} \delta c_8}) - F(\frac{\delta^5 C_\tau}{\delta c_6 \delta c_{21} \delta c_{12} \delta c_1 \delta c_8}) = 114 - 78 = 36;$$

$$F(\frac{\delta^4 C_\tau}{\delta c_6 \delta c_{21} \delta c_{12} \delta c_9}) - F(\frac{\delta^5 C_\tau}{\delta c_6 \delta c_{21} \delta c_{12} \delta c_1 \delta c_9}) = 126 - 78 = 48.$$

Удаляем цикл $c_9$. Продолжаем процесс удаления циклов.

$$F(\frac{\delta^6 C_\tau}{\delta c_6 \delta c_{21} \delta c_{12} \delta c_1 \delta c_9 \delta c_2}) = 60;\ F(\frac{\delta^6 C_\tau}{\delta c_6 \delta c_{21} \delta c_{12} \delta c_1 \delta c_9 \delta c_3}) = 66;$$

$$F(\frac{\delta^6 C_\tau}{\delta c_6 \delta c_{21} \delta c_{12} \delta c_1 \delta c_9 \delta c_4}) = 54;\ F(\frac{\delta^6 C_\tau}{\delta c_6 \delta c_{21} \delta c_{12} \delta c_1 \delta c_9 \delta c_5}) = 60;$$

$$F(\frac{\delta^6 C_\tau}{\delta c_6 \delta c_{21} \delta c_{12} \delta c_1 \delta c_9 \delta c_7}) = 66; \quad F(\frac{\delta^6 C_\tau}{\delta c_6 \delta c_{21} \delta c_{12} \delta c_1 \delta c_9 \delta c_8}) = 60;$$

$$F(\frac{\delta^6 C_\tau}{\delta c_6 \delta c_{21} \delta c_{12} \delta c_1 \delta c_9 \delta c_{10}}) = 66; \quad F(\frac{\delta^6 C_\tau}{\delta c_6 \delta c_{21} \delta c_{12} \delta c_1 \delta c_9 \delta c_{11}}) = 60;$$

$$F(\frac{\delta^6 C_\tau}{\delta c_6 \delta c_{21} \delta c_{12} \delta c_1 \delta c_9 \delta c_{13}}) = 60; \quad F(\frac{\delta^6 C_\tau}{\delta c_6 \delta c_{21} \delta c_{12} \delta c_1 \delta c_9 \delta c_{14}}) = 66;$$

$$F(\frac{\delta^6 C_\tau}{\delta c_6 \delta c_{21} \delta c_{12} \delta c_1 \delta c_9 \delta c_{15}}) = 66; \quad F(\frac{\delta^6 C_\tau}{\delta c_6 \delta c_{21} \delta c_{12} \delta c_1 \delta c_9 \delta c_{16}}) = 60;$$

$$F(\frac{\delta^6 C_\tau}{\delta c_6 \delta c_{21} \delta c_{12} \delta c_1 \delta c_9 \delta c_{17}}) = 72; \quad F(\frac{\delta^6 C_\tau}{\delta c_6 \delta c_{21} \delta c_{12} \delta c_1 \delta c_9 \delta c_{18}}) = 66;$$

$$F(\frac{\delta^6 C_\tau}{\delta c_6 \delta c_{21} \delta c_{12} \delta c_1 \delta c_9 \delta c_{19}}) = 66; \quad F(\frac{\delta^6 C_\tau}{\delta c_6 \delta c_{21} \delta c_{12} \delta c_1 \delta c_9 \delta c_{20}}) = 66.$$

Минимальное значение функционал Маклейна получает при удалении цикла $c_4$.

Удаляем цикл $c_4$ и продолжаем процесс.

$$F(\frac{\delta^7 C_\tau}{\delta c_6 \delta c_{21} \delta c_{12} \delta c_1 \delta c_9 \delta c_4 \delta c_2}) = 42; \quad F(\frac{\delta^7 C_\tau}{\delta c_6 \delta c_{21} \delta c_{12} \delta c_1 \delta c_9 \delta c_4 \delta c_3}) = 48;$$

$$F(\frac{\delta^7 C_\tau}{\delta c_6 \delta c_{21} \delta c_{12} \delta c_1 \delta c_9 \delta c_4 \delta c_5}) = 48; \quad F(\frac{\delta^7 C_\tau}{\delta c_6 \delta c_{21} \delta c_{12} \delta c_1 \delta c_9 \delta c_4 \delta c_7}) = 48;$$

$$F(\frac{\delta^7 C_\tau}{\delta c_6 \delta c_{21} \delta c_{12} \delta c_1 \delta c_9 \delta c_4 \delta c_8}) = 42; \quad F(\frac{\delta^7 C_\tau}{\delta c_6 \delta c_{21} \delta c_{12} \delta c_1 \delta c_9 \delta c_4 \delta c_{10}}) = 48;$$

$$F(\frac{\delta^7 C_\tau}{\delta c_6 \delta c_{21} \delta c_{12} \delta c_1 \delta c_9 \delta c_4 \delta c_{11}}) = 36; \quad F(\frac{\delta^7 C_\tau}{\delta c_6 \delta c_{21} \delta c_{12} \delta c_1 \delta c_9 \delta c_4 \delta c_{13}}) = 42;$$

$$F(\frac{\delta^7 C_\tau}{\delta c_6 \delta c_{21} \delta c_{12} \delta c_1 \delta c_9 \delta c_4 \delta c_{14}}) = 42; \quad F(\frac{\delta^7 C_\tau}{\delta c_6 \delta c_{21} \delta c_{12} \delta c_1 \delta c_9 \delta c_4 \delta c_{15}}) = 42;$$

$$F(\frac{\delta^7 C_\tau}{\delta c_6 \delta c_{21} \delta c_{12} \delta c_1 \delta c_9 \delta c_4 \delta c_{16}}) = 36; \quad F(\frac{\delta^7 C_\tau}{\delta c_6 \delta c_{21} \delta c_{12} \delta c_1 \delta c_9 \delta c_4 \delta c_{17}}) = 48;$$

$$FF(\frac{\delta^7 C_\tau}{\delta c_6 \delta c_{21} \delta c_{12} \delta c_1 \delta c_9 \delta c_4 \delta c_{18}}) = 42; \quad F(\frac{\delta^7 C_\tau}{\delta c_6 \delta c_{21} \delta c_{12} \delta c_1 \delta c_9 \delta c_4 \delta c_{19}}) = 42;$$

$$F(\frac{\delta^7 C_\tau}{\delta c_6 \delta c_{21} \delta c_{12} \delta c_1 \delta c_9 \delta c_4 \delta c_{20}}) = 42.$$

Имеется два минимальных значения функционала Маклейна:

$$F(\frac{\delta^6 C_\tau}{\delta c_6 \delta c_{21} \delta c_{12} \delta c_1 \delta c_9 \delta c_{11}}) - F(\frac{\delta^7 C_\tau}{\delta c_6 \delta c_{21} \delta c_{12} \delta c_1 \delta c_9 \delta c_4 \delta c_{11}}) = 60-36 = 24;$$

$$F(\frac{\delta^6 C_\tau}{\delta c_6 \delta c_{21} \delta c_{12} \delta c_1 \delta c_9 \delta c_{16}}) - F(\frac{\delta^7 C_\tau}{\delta c_6 \delta c_{21} \delta c_{12} \delta c_1 \delta c_9 \delta c_4 \delta c_{16}}) = 60-36 = 24.$$

Удаляем цикл $c_{16}$ и продолжаем процесс удаления циклов.

$$F(\frac{\delta^8 C_\tau}{\delta c_6 \delta c_{21} \delta c_{12} \delta c_1 \delta c_9 \delta c_4 \delta c_{16} \delta c_2}) = 30; \quad F(\frac{\delta^8 C_\tau}{\delta c_6 \delta c_{21} \delta c_{12} \delta c_1 \delta c_9 \delta c_4 \delta c_{16} \delta c_3}) = 36;$$

$$F(\frac{\delta^8 C_\tau}{\delta c_6 \delta c_{21} \delta c_{12} \delta c_1 \delta c_9 \delta c_4 \delta c_{16} \delta c_5}) = 30; \quad F(\frac{\delta^8 C_\tau}{\delta c_6 \delta c_{21} \delta c_{12} \delta c_1 \delta c_9 \delta c_4 \delta c_{16} \delta c_7}) = 30;$$

$$F(\frac{\delta^8 C_\tau}{\delta c_6 \delta c_{21} \delta c_{12} \delta c_1 \delta c_9 \delta c_4 \delta c_{16} \delta c_8}) = 36; \quad F(\frac{\delta^8 C_\tau}{\delta c_6 \delta c_{21} \delta c_{12} \delta c_1 \delta c_9 \delta c_4 \delta c_{16} \delta c_{10}}) = 30;$$

$$F(\frac{\delta^8 C_\tau}{\delta c_6 \delta c_{21} \delta c_{12} \delta c_1 \delta c_9 \delta c_4 \delta c_{16} \delta c_{11}}) = 24; \quad F(\frac{\delta^8 C_\tau}{\delta c_6 \delta c_{21} \delta c_{12} \delta c_1 \delta c_9 \delta c_4 \delta c_{16} \delta c_{13}}) = 24;$$

$$F(\frac{\delta^8 C_\tau}{\delta c_6 \delta c_{21} \delta c_{12} \delta c_1 \delta c_9 \delta c_4 \delta c_{16} \delta c_{14}}) = 36; \quad F(\frac{\delta^8 C_\tau}{\delta c_6 \delta c_{21} \delta c_{12} \delta c_1 \delta c_9 \delta c_4 \delta c_{16} \delta c_{15}}) = 30;$$

$$F(\frac{\delta^8 C_\tau}{\delta c_6 \delta c_{21} \delta c_{12} \delta c_1 \delta c_9 \delta c_4 \delta c_{16} \delta c_{17}}) = 30; \quad F(\frac{\delta^8 C_\tau}{\delta c_6 \delta c_{21} \delta c_{12} \delta c_1 \delta c_9 \delta c_4 \delta c_{16} \delta c_{18}}) = 24;$$

$$F(\frac{\delta^8 C_\tau}{\delta c_6 \delta c_{21} \delta c_{12} \delta c_1 \delta c_9 \delta c_4 \delta c_{16} \delta c_{19}}) = 24; \quad F(\frac{\delta^8 C_\tau}{\delta c_6 \delta c_{21} \delta c_{12} \delta c_1 \delta c_9 \delta c_4 \delta c_{16} \delta c_{20}}) = 24.$$

Имеется пять минимальных значений функционала Маклейна.

$$F(\frac{\delta^7 C_\tau}{\delta c_6 \delta c_{21} \delta c_{12} \delta c_1 \delta c_9 \delta c_4 \delta c_{13}}) - F(\frac{\delta^8 C_\tau}{\delta c_6 \delta c_{21} \delta c_{12} \delta c_1 \delta c_9 \delta c_4 \delta c_{16} \delta c_{13}}) = 42-24 = 18;$$

$$F(\frac{\delta^7 C_\tau}{\delta c_6 \delta c_{21} \delta c_{12} \delta c_1 \delta c_9 \delta c_4 \delta c_{14}}) - F(\frac{\delta^8 C_\tau}{\delta c_6 \delta c_{21} \delta c_{12} \delta c_1 \delta c_9 \delta c_4 \delta c_{16} \delta c_{14}}) = 42-24 = 18;$$

$$F(\frac{\delta^7 C_\tau}{\delta c_6 \delta c_{21} \delta c_{12} \delta c_1 \delta c_9 \delta c_4 \delta c_{18}}) - F(\frac{\delta^8 C_\tau}{\delta c_6 \delta c_{21} \delta c_{12} \delta c_1 \delta c_9 \delta c_4 \delta c_{16} \delta c_{18}}) = 42-24 = 18;$$

$$F(\frac{\delta^7 C_\tau}{\delta c_6 \delta c_{21} \delta c_{12} \delta c_1 \delta c_9 \delta c_4 \delta c_{19}}) - F(\frac{\delta^8 C_\tau}{\delta c_6 \delta c_{21} \delta c_{12} \delta c_1 \delta c_9 \delta c_4 \delta c_{16} \delta c_{19}}) = 42-24 = 18;$$

$$F(\frac{\delta^7 C_\tau}{\delta c_6 \delta c_{21} \delta c_{12} \delta c_1 \delta c_9 \delta c_4 \delta c_{20}}) - F(\frac{\delta^8 C_\tau}{\delta c_6 \delta c_{21} \delta c_{12} \delta c_1 \delta c_9 \delta c_4 \delta c_{16} \delta c_{20}}) = 42-24 = 18.$$

Выбрать для удаления цикл $c_{14}$ нельзя, так как исключение ребра $e_{12}$ приводит к нарушению условия Эйлера. Выбираем для удаления цикл $c_{13}$. Оставшаяся система циклов есть базис подпространства квазициклов графа. Количество циклов соответствует цикломатическому числу графа.

цикл $c_2 = \{e_1,e_2,e_6,e_{11}\} \leftrightarrow \{v_1,v_2,v_4,v_5\}$;
цикл $c_3 = \{e_1,e_4,e_7\} \leftrightarrow \{v_1,v_2,v_{12}\}$;
цикл $c_5 = \{e_2,e_3,e_8,e_9,e_{23}\} \leftrightarrow \{v_1,v_3,v_4,v_{10},v_{11}\}$;
цикл $c_7 = \{e_3,e_4,e_{24}\} \leftrightarrow \{v_1,v_{11},v_{12}\}$;
цикл $c_8 = \{e_5,e_6,e_8,e_{11}\} \leftrightarrow \{v_2,v_3,v_4,v_5\}$;
цикл $c_{10} = \{e_5,e_7,e_{10}\} \leftrightarrow \{v_2,v_3,v_{12}\}$;
цикл $c_{11} = \{e_6,e_7,e_{14},e_{20}\} \leftrightarrow \{v_2,v_5,v_8,v_{12}\}$;
цикл $c_{14} = \{e_{11},e_{12},e_{13}\} \leftrightarrow \{v_4,v_5,v_6\}$;
цикл $c_{15} = \{e_{13},e_{14},e_{16},e_{17}\} \leftrightarrow \{v_5,v_6,v_7,v_8\}$;
цикл $c_{17} = \{e_{17},e_{18},e_{19}\} \leftrightarrow \{v_7,v_8,v_9\}$;
цикл $c_{18} = \{e_{13},e_{15},e_{16},e_{18},e_{21}\} \leftrightarrow \{v_5,v_6,v_7,v_9,v_{10}\}$;
цикл $c_{19} = \{e_{19},e_{20},e_{22},e_{24}\} \leftrightarrow \{v_8,v_9,v_{11},v_{12}\}$;
цикл $c_{20} = \{e_{21},e_{22},e_{23}\} \leftrightarrow \{v_9,v_{10},v_{11}\}$;

Структурное число системы $W(Y)$ имеет вид:

| | | | |
|---|---|---|---|
| $e_2$ | $c_2$ | $c_5$ | |
| $e_3$ | $c_5$ | $c_7$ | |
| $e_4$ | $c_3$ | $c_7$ | |
| $e_6$ | $c_2$ | $c_8$ | $c_{11}$ |
| $e_7$ | $c_3$ | $c_{10}$ | $c_{11}$ |
| $e_9$ | $c_5$ | | |
| $e_{10}$ | $c_{10}$ | | |
| $e_{12}$ | $c_{14}$ | | |
| $e_{14}$ | $c_{11}$ | $c_{15}$ | |
| $e_{15}$ | $c_{18}$ | | |
| $e_{18}$ | $c_{17}$ | $c_{18}$ | |
| $e_{20}$ | $c_{11}$ | $c_{19}$ | |
| $e_{22}$ | $c_{19}$ | $c_{20}$ | |

Сжатое структурное число системы $W(Y)$ имеет следующий вид (удаленные циклы окрашены в красный цвет):

| | | | |
|---|---|---|---|
| $e_2$ | $c_2$ | **$c_5$** | |
| $e_3$ | **$c_5$** | $c_7$ | |
| $e_4$ | $c_3$ | **$c_7$** | |
| $e_6$ | **$c_2$** | $c_8$ | **$c_{11}$** |
| $e_7$ | **$c_3$** | **$c_{10}$** | $c_{11}$ |
| $e_9$ | $c_5$ | | |
| $e_{10}$ | $c_{10}$ | | |

| | | |
|---|---|---|
| $e_{12}$ | $c_{14}$ | |
| $e_{14}$ | $c_{11}$ | $c_{15}$ |
| $e_{15}$ | $c_{18}$ | |
| $e_{18}$ | $c_{17}$ | $c_{18}$ |
| $e_{20}$ | $c_{11}$ | $c_{19}$ |
| $e_{22}$ | $c_{19}$ | $c_{20}$ |

На втором этапе будем удалять циклы из системы, не нарушая правила Эйлера.

$$\frac{\partial C_b}{\partial c_5}=24; \quad \frac{\partial C_b}{\partial c_{10}}=18; \quad \frac{\partial C_b}{\partial c_{14}}=12; \quad \frac{\partial C_b}{\partial c_{18}}=18.$$

Для удаления выбираем цикл $c_{14}$, исключая ребро $e_{12}$. Продолжаем процесс.

$$\frac{\partial^2 C_b}{\partial c_{14} \partial c_5}=12; \quad ; \quad \frac{\partial^2 C_b}{\partial c_{14} \partial c_{10}}=6; \quad \frac{\partial^2 C_b}{\partial c_{14} \partial c_{18}}=12.$$

Для удаления выбираем цикл $c_{10}$, исключая ребро $e_{10}$. Продолжаем процесс.

$$\frac{\partial^3 C_b}{\partial c_{14} \partial c_{10} \partial c_5}=6; \quad \frac{\partial^3 C_b}{\partial c_{14} \partial c_{10} \partial c_8}=0; \quad \frac{\partial^3 C_b}{\partial c_{14} \partial c_{10} \partial c_{18}}=6.$$

Для удаления выбираем цикл $c_8$. При удалении исключается ребро $e_5$. Оставшаяся система циклов имеет значение функционала Маклейна равное нулю. Таким образом, мы выделили следующую плоскую часть непланарного графа (рис. 12.5).

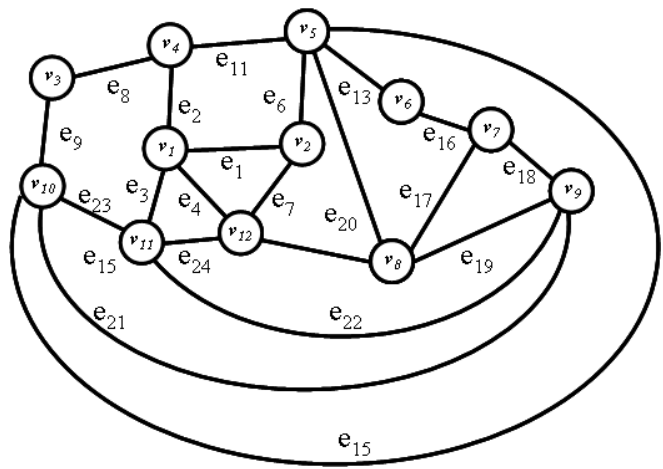

Рис. 12.5. Плоская часть графа $G_3$ с тремя удаленными ребрами.

Обратим внимание на следующий факт – во втором этапе процессе исключения циклов учавствуют циклы $c_5, c_8, c_{10}, c_{14}, c_{18}$. Но по условию удаления циклов исключены только циклы $c_8, c_{10}, c_{14}$.

Система циклов, характеризующая плоскую часть графа (рис. 12.5):

цикл   $c_2 = \{e_1, e_2, e_6, e_{11}\} \leftrightarrow \{v_1, v_2, v_4, v_5\}$;

цикл   $c_3 = \{e_1,e_4,e_7\} \leftrightarrow \{v_1,v_2,v_{12}\}$;
цикл   $c_5 = \{e_2,e_3,e_8,e_9,e_{23}\} \leftrightarrow \{v_1,v_3,v_4,v_{10},v_{11}\}$;
цикл   $c_7 = \{e_3,e_4,e_{24}\} \leftrightarrow \{v_1,v_{11},v_{12}\}$;
цикл   $c_{11} = \{e_6,e_7,e_{14},e_{20}\} \leftrightarrow \{v_2,v_5,v_8,v_{12}\}$;
цикл   $c_{15} = \{e_{13},e_{14},e_{16},e_{17}\} \leftrightarrow \{v_5,v_6,v_7,v_8\}$;
цикл   $c_{17} = \{e_{17},e_{18},e_{19}\} \leftrightarrow \{v_7,v_8,v_9\}$;
цикл   $c_{18} = \{e_{13},e_{15},e_{16},e_{18},e_{21}\} \leftrightarrow \{v_5,v_6,v_7,v_9,v_{10}\}$;
цикл   $c_{19} = \{e_{19},e_{20},e_{22},e_{24}\} \leftrightarrow \{v_8,v_9,v_{11},v_{12}\}$;
цикл   $c_{20} = \{e_{21},e_{22},e_{23}\} \leftrightarrow \{v_9,v_{10},v_{11}\}$;

Мы выделили плоскую часть графа $G_3$ с тремя удаленными ребрами. Попробуем построить плоскую часть графа с меньшим числом удаленных ребер.

Рассмотрим вопрос исключения оставшейся части циклов $c_5, c_{18}$, на втором этапе. Выделим в графе следующие плоские конфигурации:

цикл   $c_2 = \{e_1,e_2,e_6,e_{11}\} \leftrightarrow \{v_1,v_2,v_4,v_5\}$;
цикл   $c_3 = \{e_1,e_4,e_7\} \leftrightarrow \{v_1,v_2,v_{12}\}$;
цикл   $c_5 = \{e_2,e_3,e_8,e_9,e_{23}\} \leftrightarrow \{v_1,v_3,v_4,v_{10},v_{11}\}$;
цикл   $c_7 = \{e_3,e_4,e_{24}\} \leftrightarrow \{v_1,v_{11},v_{12}\}$;
цикл   $c_8 = \{e_5,e_6,e_8,e_{11}\} \leftrightarrow \{v_2,v_3,v_4,v_5\}$;
цикл   $c_{10} = \{e_5,e_7,e_{10}\} \leftrightarrow \{v_2,v_3,v_{12}\}$;
цикл   $c_{13} = \{e_9,e_{10},e_{23},e_{24}\} \leftrightarrow \{v_3,v_{10},v_{11},v_{12}\}$;
$c_2 \oplus c_3 \oplus c_5 \oplus c_7 \oplus c_8 \oplus c_{10} \oplus c_{13} = \varnothing$.

цикл   $c_{15} = \{e_{13},e_{14},e_{16},e_{17}\} \leftrightarrow \{v_5,v_6,v_7,v_8\}$;
цикл   $c_{16} = \{e_{14},e_{15},e_{19},e_{21}\} \leftrightarrow \{v_5,v_8,v_9,v_{10}\}$;
цикл   $c_{17} = \{e_{17},e_{18},e_{19}\} \leftrightarrow \{v_7,v_8,v_9\}$;
цикл   $c_{18} = \{e_{13},e_{15},e_{16},e_{18},e_{21}\} \leftrightarrow \{v_5,v_6,v_7,v_9,v_{10}\}$;
$c_{15} \oplus c_{16} \oplus c_{17} \oplus c_{18} = \varnothing$.

Чтобы избежать появления таких зависимых подсистем удалим из множества изометрических циклов самые длинные циклы – $c_5$ и $c_{18}$.

$$F\left(\frac{\delta^2 C_\tau}{\delta c_5 \delta c_{18}}\right) = 348.$$

Рассмотрим процесс построения другой системы циклов с заранее исключенными циклами $c_5$ и $c_{18}$.

Продолжим процесс удаления циклов из множества изометрических циклов:

$F\left(\frac{\delta^3 C_\tau}{\delta c_5 \delta c_{18} \delta c_1}\right)=276$; $F\left(\frac{\delta^3 C_\tau}{\delta c_5 \delta c_{18} \delta c_2}\right)=240$; $F\left(\frac{\delta^3 C_\tau}{\delta c_5 \delta c_{18} \delta c_3}\right)=318$; $F\left(\frac{\delta^3 C_\tau}{\delta c_5 \delta c_{18} \delta c_4}\right)=300$;

$F\left(\frac{\delta^3 C_\tau}{\delta c_5 \delta c_{18} \delta c_6}\right)=234$; $F\left(\frac{\delta^3 C_\tau}{\delta c_5 \delta c_{18} \delta c_7}\right)=330$; $F\left(\frac{\delta^3 C_\tau}{\delta c_5 \delta c_{18} \delta c_8}\right)=240$; $F\left(\frac{\delta^3 C_\tau}{\delta c_5 \delta c_{18} \delta c_9}\right)=252$;

$$F(\frac{\delta^3 C_\tau}{\delta c_5 \delta c_{18} \delta c_{10}})=318; \quad F(\frac{\delta^3 C_\tau}{\delta c_5 \delta c_{18} \delta c_{11}})=300; \quad F(\frac{\delta^3 C_\tau}{\delta c_5 \delta c_{18} \delta c_{12}})=252; \quad F(\frac{\delta^3 C_\tau}{\delta c_5 \delta c_{18} \delta c_{13}})=312;$$

$$F(\frac{\delta^3 C_\tau}{\delta c_5 \delta c_{18} \delta c_{14}})=312; \quad F(\frac{\delta^3 C_\tau}{\delta c_5 \delta c_{18} \delta c_{15}})=342; \quad F(\frac{\delta^3 C_\tau}{\delta c_5 \delta c_{18} \delta c_{16}})=300; \quad F(\frac{\delta^3 C_\tau}{\delta c_5 \delta c_{18} \delta c_{17}})=342;$$

$$F(\frac{\delta^3 C_\tau}{\delta c_5 \delta c_{18} \delta c_{19}})=336; \quad F(\frac{\delta^3 C_\tau}{\delta c_5 \delta c_{18} \delta c_{20}})=330; \quad F(\frac{\delta^3 C_\tau}{\delta c_5 \delta c_{18} \delta c_{21}})=234.$$

Имеем два элемента с минимальным значением функционала Маклейна. Оценим скорость изменения.

$$F(\frac{\delta^3 C_\tau}{\delta c_5 \delta c_{18} \delta c_6}) = 234;$$

$$F(\frac{\delta^3 C_\tau}{\delta c_5 \delta c_{18} \delta c_{21}}) = 234.$$

Для удаления нужно выбрать цикл $c_{21}$.

$$F(\frac{\delta^4 C_\tau}{\delta c_5 \delta c_{18} \delta c_{21} \delta c_1}) = 174; \quad F(\frac{\delta^4 C_\tau}{\delta c_5 \delta c_{18} \delta c_{21} \delta c_2}) = 156; \quad F(\frac{\delta^4 C_\tau}{\delta c_5 \delta c_{18} \delta c_{21} \delta c_3}) = 216;$$

$$F(\frac{\delta^4 C_\tau}{\delta c_5 \delta c_{18} \delta c_{21} \delta c_4}) = 186; \quad F(\frac{\delta^4 C_\tau}{\delta c_5 \delta c_{18} \delta c_{21} \delta c_6}) = 156; \quad F(\frac{\delta^4 C_\tau}{\delta c_5 \delta c_{18} \delta c_{21} \delta c_7}) = 222;$$

$$F(\frac{\delta^4 C_\tau}{\delta c_5 \delta c_{18} \delta c_{21} \delta c_8}) = 144; \quad F(\frac{\delta^4 C_\tau}{\delta c_5 \delta c_{18} \delta c_{21} \delta c_9}) = 174; \quad F(\frac{\delta^4 C_\tau}{\delta c_5 \delta c_{18} \delta c_{21} \delta c_{10}}) = 204;$$

$$F(\frac{\delta^4 C_\tau}{\delta c_5 \delta c_{18} \delta c_{21} \delta c_{11}}) = 204; \quad F(\frac{\delta^4 C_\tau}{\delta c_5 \delta c_{18} \delta c_{21} \delta c_{12}}) = 156; \quad F(\frac{\delta^4 C_\tau}{\delta c_5 \delta c_{18} \delta c_{21} \delta c_{13}}) = 210;$$

$$F(\frac{\delta^4 C_\tau}{\delta c_5 \delta c_{18} \delta c_{21} \delta c_{14}}) = 198; \quad F(\frac{\delta^4 C_\tau}{\delta c_5 \delta c_{18} \delta c_{21} \delta c_{15}}) = 228; \quad F(\frac{\delta^4 C_\tau}{\delta c_5 \delta c_{18} \delta c_{21} \delta c_{16}}) = 204;$$

$$F(\frac{\delta^4 C_\tau}{\delta c_5 \delta c_{18} \delta c_{21} \delta c_{17}}) = 228; \quad F(\frac{\delta^4 C_\tau}{\delta c_5 \delta c_{18} \delta c_{21} \delta c_{19}}) = 222; \quad F(\frac{\delta^4 C_\tau}{\delta c_5 \delta c_{18} \delta c_{21} \delta c_{20}}) = 228.$$

Минимальное значение функционал Маклейна принимает при удалении цикла $c_8$.

Удаляем цикл $c_8$ и продолжаем процесс исключения циклов.

$$F(\frac{\delta^5 C_\tau}{\delta c_5 \delta c_{18} \delta c_{21} \delta c_8 \delta c_1}) = 108; \quad F(\frac{\delta^5 C_\tau}{\delta c_5 \delta c_{18} \delta c_{21} \delta c_8 \delta c_2}) = 96;$$

$$F(\frac{\delta^5 C_\tau}{\delta c_5 \delta c_{18} \delta c_{21} \delta c_8 \delta c_3}) = 126; \quad F(\frac{\delta^5 C_\tau}{\delta c_5 \delta c_{18} \delta c_{21} \delta c_8 \delta c_4}) = 108;$$

$$F(\frac{\delta^5 C_\tau}{\delta c_5 \delta c_{18} \delta c_{21} \delta c_8 \delta c_6}) = 84; \quad F(\frac{\delta^5 C_\tau}{\delta c_5 \delta c_{18} \delta c_{21} \delta c_8 \delta c_7}) = 132;$$

$$F(\frac{\delta^5 C_\tau}{\delta c_5 \delta c_{18} \delta c_{21} \delta c_8 \delta c_9}) = 108; \quad F(\frac{\delta^5 C_\tau}{\delta c_5 \delta c_{18} \delta c_{21} \delta c_8 \delta c_{10}}) = 126;$$

$$F(\frac{\delta^5 C_\tau}{\delta c_5 \delta c_{18} \delta c_{21} \delta c_8 \delta c_{11}}) = 108; \quad F(\frac{\delta^5 C_\tau}{\delta c_5 \delta c_{18} \delta c_{21} \delta c_8 \delta c_{12}}) = 96;$$

$$F(\frac{\delta^5 C_\tau}{\delta c_5 \delta c_{18} \delta c_{21} \delta c_8 \delta c_{13}}) = 120; \quad F(\frac{\delta^5 C_\tau}{\delta c_5 \delta c_{18} \delta c_{21} \delta c_8 \delta c_{14}}) = 126;$$

$$F(\frac{\delta^5 C_\tau}{\delta c_5 \delta c_{18} \delta c_{21} \delta c_8 \delta c_{15}}) = 138; \quad F(\frac{\delta^5 C_\tau}{\delta c_5 \delta c_{18} \delta c_{21} \delta c_8 \delta c_{16}}) = 114;$$

$$F(\frac{\delta^5 C_\tau}{\delta c_5 \delta c_{18} \delta c_{21} \delta c_8 \delta c_{17}}) = 138; \quad F(\frac{\delta^5 C_\tau}{\delta c_5 \delta c_{18} \delta c_{21} \delta c_8 \delta c_{19}}) = 132;$$

$$F(\frac{\delta^5 C_\tau}{\delta c_5 \delta c_{18} \delta c_{21} \delta c_8 \delta c_{20}}) = 138.$$

Минимальное значение функционал Маклейна принимает при удалении цикла $c_6$.

Удаляем цикл $c_6$ и продолжаем процесс исключения циклов:

$$F(\frac{\delta^6 C_\tau}{\delta c_5 \delta c_{18} \delta c_{21} \delta c_8 \delta c_6 \delta c_1}) = 60; \quad F(\frac{\delta^6 C_\tau}{\delta c_5 \delta c_{18} \delta c_{21} \delta c_8 \delta c_6 \delta c_2}) = 60;$$

$$F(\frac{\delta^6 C_\tau}{\delta c_5 \delta c_{18} \delta c_{21} \delta c_8 \delta c_6 \delta c_3}) = 66; \quad F(\frac{\delta^6 C_\tau}{\delta c_5 \delta c_{18} \delta c_{21} \delta c_8 \delta c_6 \delta c_4}) = 60;$$

$$F(\frac{\delta^6 C_\tau}{\delta c_5 \delta c_{18} \delta c_{21} \delta c_8 \delta c_6 \delta c_7}) = 72; \quad F(\frac{\delta^6 C_\tau}{\delta c_5 \delta c_{18} \delta c_{21} \delta c_8 \delta c_6 \delta c_9}) = 60;$$

$$F(\frac{\delta^6 C_\tau}{\delta c_5 \delta c_{18} \delta c_{21} \delta c_8 \delta c_6 \delta c_{10}}) = 66; \quad F(\frac{\delta^6 C_\tau}{\delta c_5 \delta c_{18} \delta c_{21} \delta c_8 \delta c_6 \delta c_{11}}) = 66;$$

$$F(\frac{\delta^6 C_\tau}{\delta c_5 \delta c_{18} \delta c_{21} \delta c_8 \delta c_6 \delta c_{12}}) = 60; \quad F(\frac{\delta^6 C_\tau}{\delta c_5 \delta c_{18} \delta c_{21} \delta c_8 \delta c_6 \delta c_{13}}) = 66;$$

$$F(\frac{\delta^6 C_\tau}{\delta c_5 \delta c_{18} \delta c_{21} \delta c_8 \delta c_6 \delta c_{14}}) = 78; \quad F(\frac{\delta^6 C_\tau}{\delta c_5 \delta c_{18} \delta c_{21} \delta c_8 \delta c_6 \delta c_{15}}) = 78;$$

$$F(\frac{\delta^6 C_\tau}{\delta c_5 \delta c_{18} \delta c_{21} \delta c_8 \delta c_6 \delta c_{16}}) = 66; \quad F(\frac{\delta^6 C_\tau}{\delta c_5 \delta c_{18} \delta c_{21} \delta c_8 \delta c_6 \delta c_{17}}) = 78;$$

$$F(\frac{\delta^6 C_\tau}{\delta c_5 \delta c_{18} \delta c_{21} \delta c_8 \delta c_6 \delta c_{10}}) = 72; \quad F(\frac{\delta^6 C_\tau}{\delta c_5 \delta c_{18} \delta c_{21} \delta c_8 \delta c_6 \delta c_{20}}) = 84.$$

Имеем пять элементов с минимальным значением функционала Маклейна:

$$F(\frac{\delta^5 C_\tau}{\delta c_5 \delta c_{18} \delta c_{21} \delta c_8 \delta c_1}) - F(\frac{\delta^6 C_\tau}{\delta c_5 \delta c_{18} \delta c_{21} \delta c_8 \delta c_6 \delta c_1}) = 108-60 = 48;$$

$$F(\frac{\delta^5 C_\tau}{\delta c_5 \delta c_{18} \delta c_{21} \delta c_8 \delta c_2}) - F(\frac{\delta^6 C_\tau}{\delta c_5 \delta c_{18} \delta c_{21} \delta c_8 \delta c_6 \delta c_2}) = 96-60 = 36;$$

$$F(\frac{\delta^5 C_\tau}{\delta c_5 \delta c_{18} \delta c_{21} \delta c_8 \delta c_4}) - F(\frac{\delta^6 C_\tau}{\delta c_5 \delta c_{18} \delta c_{21} \delta c_8 \delta c_6 \delta c_4}) = 108-60 = 48;$$

$$F(\frac{\delta^5 C_\tau}{\delta c_5 \delta c_{18} \delta c_{21} \delta c_8 \delta c_9}) - F(\frac{\delta^6 C_\tau}{\delta c_5 \delta c_{18} \delta c_{21} \delta c_8 \delta c_6 \delta c_9}) = 108-60 = 48;$$

$$F(\frac{\delta^5 C_\tau}{\delta c_5 \delta c_{18} \delta c_{21} \delta c_8 \delta c_{12}}) - F(\frac{\delta^6 C_\tau}{\delta c_5 \delta c_{18} \delta c_{21} \delta c_8 \delta c_6 \delta c_{12}}) = 96-60 = 36.$$

Для исключения выбираем цикл $c_9$. Удаляем из системы циклов цикл $c_9$ и продолжаем процесс:

$$F(\frac{\delta^7 C_\tau}{\delta c_5 \delta c_{18} \delta c_{21} \delta c_8 \delta c_6 \delta c_9 \delta c_1}) = 42; \quad F(\frac{\delta^7 C_\tau}{\delta c_5 \delta c_{18} \delta c_{21} \delta c_8 \delta c_6 \delta c_9 \delta c_2}) = 42;$$

$$(\frac{\delta^7 C_\tau}{\delta c_5 \delta c_{18} \delta c_{21} \delta c_8 \delta c_6 \delta c_9 \delta c_3}) = 42; \quad F(\frac{\delta^7 C_\tau}{\delta c_5 \delta c_{18} \delta c_{21} \delta c_8 \delta c_6 \delta c_9 \delta c_4}) = 36;$$

$$(\frac{\delta^7 C_\tau}{\delta c_5 \delta c_{18} \delta c_{21} \delta c_8 \delta c_6 \delta c_9 \delta c_7}) = 48; \quad F(\frac{\delta^7 C_\tau}{\delta c_5 \delta c_{18} \delta c_{21} \delta c_8 \delta c_6 \delta c_9 \delta c_{10}}) = 48;$$

$$\left(\frac{\delta^7 C_\tau}{\delta c_5 \delta c_{18} \delta c_{21} \delta c_8 \delta c_6 \delta c_9 \delta c_{11}}\right) = 48; \quad F\left(\frac{\delta^7 C_\tau}{\delta c_5 \delta c_{18} \delta c_{21} \delta c_8 \delta c_6 \delta c_9 \delta c_{12}}\right) = 48;$$

$$\left(\frac{\delta^7 C_\tau}{\delta c_5 \delta c_{18} \delta c_{21} \delta c_8 \delta c_6 \delta c_9 \delta c_{13}}\right) = 54; \quad F\left(\frac{\delta^7 C_\tau}{\delta c_5 \delta c_{18} \delta c_{21} \delta c_8 \delta c_6 \delta c_9 \delta c_{14}}\right) = 54;$$

$$\left(\frac{\delta^7 C_\tau}{\delta c_5 \delta c_{18} \delta c_{21} \delta c_8 \delta c_6 \delta c_9 \delta c_{15}}\right) = 48; \quad F\left(\frac{\delta^7 C_\tau}{\delta c_5 \delta c_{18} \delta c_{21} \delta c_8 \delta c_6 \delta c_9 \delta c_{16}}\right) = 54;$$

$$\left(\frac{\delta^7 C_\tau}{\delta c_5 \delta c_{18} \delta c_{21} \delta c_8 \delta c_6 \delta c_9 \delta c_{17}}\right) = 54; \quad F\left(\frac{\delta^7 C_\tau}{\delta c_5 \delta c_{18} \delta c_{21} \delta c_8 \delta c_6 \delta c_9 \delta c_{19}}\right) = 48;$$

$$\left(\frac{\delta^7 C_\tau}{\delta c_5 \delta c_{18} \delta c_{21} \delta c_8 \delta c_6 \delta c_9 \delta c_{20}}\right) = 54.$$

Минимальное значение функционал Маклейна приобретает при исключении цикла $c_4$.

Удаляем цикл $c_4$ и продолжаем процесс:

$$F\left(\frac{\delta^8 C_\tau}{\delta c_5 \delta c_{18} \delta c_{21} \delta c_8 \delta c_6 \delta c_9 \delta c_4 \delta c_1}\right) = 30; \quad F\left(\frac{\delta^8 C_\tau}{\delta c_5 \delta c_{18} \delta c_{21} \delta c_8 \delta c_6 \delta c_9 \delta c_4 \delta c_2}\right) = 24;$$

$$F\left(\frac{\delta^8 C_\tau}{\delta c_5 \delta c_{18} \delta c_{21} \delta c_8 \delta c_6 \delta c_9 \delta c_4 \delta c_3}\right) = 24; \quad F\left(\frac{\delta^8 C_\tau}{\delta c_5 \delta c_{18} \delta c_{21} \delta c_8 \delta c_6 \delta c_9 \delta c_4 \delta c_7}\right) = 30;$$

$$F\left(\frac{\delta^8 C_\tau}{\delta c_5 \delta c_{18} \delta c_{21} \delta c_8 \delta c_6 \delta c_9 \delta c_4 \delta c_{10}}\right) = 30; \quad F\left(\frac{\delta^8 C_\tau}{\delta c_5 \delta c_{18} \delta c_{21} \delta c_8 \delta c_6 \delta c_9 \delta c_4 \delta c_{11}}\right) = 24;$$

$$F\left(\frac{\delta^8 C_\tau}{\delta c_5 \delta c_{18} \delta c_{21} \delta c_8 \delta c_6 \delta c_9 \delta c_4 \delta c_{12}}\right) = 30; \quad F\left(\frac{\delta^8 C_\tau}{\delta c_5 \delta c_{18} \delta c_{21} \delta c_8 \delta c_6 \delta c_9 \delta c_4 \delta c_{13}}\right) = 30;$$

$$F\left(\frac{\delta^8 C_\tau}{\delta c_5 \delta c_{18} \delta c_{21} \delta c_8 \delta c_6 \delta c_9 \delta c_4 \delta c_{14}}\right) = 30; \quad F\left(\frac{\delta^8 C_\tau}{\delta c_5 \delta c_{18} \delta c_{21} \delta c_8 \delta c_6 \delta c_9 \delta c_4 \delta c_{15}}\right) = 30;$$

$$F\left(\frac{\delta^8 C_\tau}{\delta c_5 \delta c_{18} \delta c_{21} \delta c_8 \delta c_6 \delta c_9 \delta c_4 \delta c_{16}}\right) = 24; \quad F\left(\frac{\delta^8 C_\tau}{\delta c_5 \delta c_{18} \delta c_{21} \delta c_8 \delta c_6 \delta c_9 \delta c_4 \delta c_{17}}\right) = 30;$$

$$F\left(\frac{\delta^8 C_\tau}{\delta c_5 \delta c_{18} \delta c_{21} \delta c_8 \delta c_6 \delta c_9 \delta c_4 \delta c_{19}}\right) = 24; \quad F\left(\frac{\delta^8 C_\tau}{\delta c_5 \delta c_{18} \delta c_{21} \delta c_8 \delta c_6 \delta c_9 \delta c_4 \delta c_{20}}\right) = 36.$$

Имеем пять элементов с минимальным значением функционала Маклейна:

$$F(\frac{\delta^7 C_\tau}{\delta c_5 \delta c_{18} \delta c_{21} \delta c_8 \delta c_6 \delta c_9 \delta c_2}) - F(\frac{\delta^8 C_\tau}{\delta c_5 \delta c_{18} \delta c_{21} \delta c_8 \delta c_6 \delta c_9 \delta c_4 \delta c_2}) = 42 - 24 = 18;$$

$$F(\frac{\delta^7 C_\tau}{\delta c_5 \delta c_{18} \delta c_{21} \delta c_8 \delta c_6 \delta c_9 \delta c_3}) - F(\frac{\delta^8 C_\tau}{\delta c_5 \delta c_{18} \delta c_{21} \delta c_8 \delta c_6 \delta c_9 \delta c_4 \delta c_3}) = 42 - 24 = 18;$$

$$F(\frac{\delta^7 C_\tau}{\delta c_5 \delta c_{18} \delta c_{21} \delta c_8 \delta c_6 \delta c_9 \delta c_{11}}) - F(\frac{\delta^8 C_\tau}{\delta c_5 \delta c_{18} \delta c_{21} \delta c_8 \delta c_6 \delta c_9 \delta c_4 \delta c_{11}}) = 48 - 24 = 24;$$

$$F(\frac{\delta^7 C_\tau}{\delta c_5 \delta c_{18} \delta c_{21} \delta c_8 \delta c_6 \delta c_9 \delta c_{16}}) - F(\frac{\delta^8 C_\tau}{\delta c_5 \delta c_{18} \delta c_{21} \delta c_8 \delta c_6 \delta c_9 \delta c_4 \delta c_{16}}) = 48 - 24 = 24;$$

$$F(\frac{\delta^7 C_\tau}{\delta c_5 \delta c_{18} \delta c_{21} \delta c_8 \delta c_6 \delta c_9 \delta c_{19}}) - F(\frac{\delta^8 C_\tau}{\delta c_5 \delta c_{18} \delta c_{21} \delta c_8 \delta c_6 \delta c_9 \delta c_4 \delta c_{19}}) = 48 - 24 = 24;$$

Для исключения выбираем цикл $c_{11}$. Количество оставшихся циклов равно цикломатическому числу графа. Прекращаем ислючение циклов. Система имеет вид:

цикл $c_1 = \{e_1, e_2, e_5, e_8\} \leftrightarrow \{v_1, v_2, v_3, v_4\}$;
цикл $c_2 = \{e_1, e_2, e_6, e_{11}\} \leftrightarrow \{v_1, v_2, v_4, v_5\}$;
цикл $c_3 = \{e_1, e_4, e_7\} \leftrightarrow \{v_1, v_2, v_{12}\}$;
цикл $c_7 = \{e_3, e_4, e_{24}\} \leftrightarrow \{v_1, v_{11}, v_{12}\}$;
цикл $c_{10} = \{e_5, e_7, e_{10}\} \leftrightarrow \{v_2, v_3, v_{12}\}$;
цикл $c_{12} = \{e_8, e_9, e_{11}, e_{15}\} \leftrightarrow \{v_3, v_4, v_5, v_{10}\}$;
цикл $c_{13} = \{e_9, e_{10}, e_{23}, e_{24}\} \leftrightarrow \{v_3, v_{10}, v_{11}, v_{12}\}$;
цикл $c_{14} = \{e_{11}, e_{12}, e_{13}\} \leftrightarrow \{v_4, v_5, v_6\}$;
цикл $c_{15} = \{e_{13}, e_{14}, e_{16}, e_{17}\} \leftrightarrow \{v_5, v_6, v_7, v_8\}$;
цикл $c_{16} = \{e_{14}, e_{15}, e_{19}, e_{21}\} \leftrightarrow \{v_5, v_8, v_9, v_{10}\}$;
цикл $c_{17} = \{e_{17}, e_{18}, e_{19}\} \leftrightarrow \{v_7, v_8, v_9\}$;
цикл $c_{19} = \{e_{19}, e_{20}, e_{22}, e_{24}\} \leftrightarrow \{v_8, v_9, v_{11}, v_{12}\}$;
цикл $c_{20} = \{e_{21}, e_{22}, e_{23}\} \leftrightarrow \{v_9, v_{10}, v_{11}\}$.

Структурное число системы $W(Y)$ имеет вид:

| | | |
|---|---|---|
| $e_2$ | $c_1$ | $c_2$ |
| $e_3$ | $c_7$ | |
| $e_4$ | $c_3$ | $c_7$ |
| $e_6$ | $c_2$ | |
| $e_7$ | $c_3$ | $c_{10}$ |
| $e_9$ | $c_{12}$ | |
| $e_{10}$ | $c_{10}$ | $c_{13}$ |
| $e_{12}$ | $c_{14}$ | |
| $e_{14}$ | $c_{15}$ | $c_{16}$ |
| $e_{15}$ | $c_{12}$ | $c_{16}$ |
| $e_{18}$ | $c_{17}$ | |
| $e_{20}$ | $c_{19}$ | |
| $e_{22}$ | $c_{19}$ | $c_{20}$ |

Сжатое структурное число системы W(Y) имеет один столбец (удаленные циклы в процессе сжатия окрашены красным цветом):

| | | |
|---|---|---|
| $e_2$ | $c_1$ | **$c_2$** |
| $e_3$ | $c_7$ | |
| $e_4$ | $c_3$ | **$c_7$** |
| $e_6$ | $c_2$ | |
| $e_7$ | **$c_3$** | $c_{10}$ |
| $e_9$ | $c_{12}$ | |
| $e_{10}$ | **$c_{10}$** | $c_{13}$ |
| $e_{12}$ | $c_{14}$ | |
| $e_{14}$ | $c_{15}$ | **$c_{16}$** |
| $e_{15}$ | **$c_{12}$** | $c_{16}$ |
| $e_{18}$ | $c_{17}$ | |
| $e_{20}$ | $c_{19}$ | |
| $e_{22}$ | **$c_{19}$** | $c_{20}$ |

Следовательно, это базис.

Будем исключать циклы из базиса с учетом выполнения правила Эйлера:

$$\frac{\partial C_b}{\partial c_2}=12;\quad \frac{\partial C_b}{\partial c_7}=18;\quad \frac{\partial C_b}{\partial c_{14}}=18;\quad \frac{\partial C_b}{\partial c_{15}}=24;\quad \frac{\partial C_b}{\partial c_{17}}=18;\quad \frac{\partial C_b}{\partial c_{19}}=12.$$

Для исключения выбираем цикл $c_2$. Удаляем цикл $c_2$ из базиса и продолжаем процесс:

$$\frac{\partial^2 C_b}{\partial c_2 \partial c_7}=6;\quad \frac{\partial^2 C_b}{\partial c_2 \partial c_{17}}=6;\quad \frac{\partial^2 C_b}{\partial c_2 \partial c_{19}}=0.$$

Удаляем цикл $c_{19}$ из системы. В результате построена система для описания топологического рисунка плоской части графа $G_3$ с двумя удаленными ребрами.

цикл $c_1 = \{e_1,e_2,e_5,e_8\} \leftrightarrow \{v_1,v_2,v_3,v_4\}$;
цикл $c_3 = \{e_1,e_4,e_7\} \leftrightarrow \{v_1,v_2,v_{12}\}$;
цикл $c_7 = \{e_3,e_4,e_{24}\} \leftrightarrow \{v_1,v_{11},v_{12}\}$;
цикл $c_{10} = \{e_5,e_7,e_{10}\} \leftrightarrow \{v_2,v_3,v_{12}\}$;
цикл $c_{12} = \{e_8,e_9,e_{11},e_{15}\} \leftrightarrow \{v_3,v_4,v_5,v_{10}\}$;
цикл $c_{13} = \{e_9,e_{10},e_{23},e_{24}\} \leftrightarrow \{v_3,v_{10},v_{11},v_{12}\}$;
цикл $c_{14} = \{e_{11},e_{12},e_{13}\} \leftrightarrow \{v_4,v_5,v_6\}$;
цикл $c_{15} = \{e_{13},e_{14},e_{16},e_{17}\} \leftrightarrow \{v_5,v_6,v_7,v_8\}$;
цикл $c_{16} = \{e_{14},e_{15},e_{19},e_{21}\} \leftrightarrow \{v_5,v_8,v_9,v_{10}\}$;
цикл $c_{17} = \{e_{17},e_{18},e_{19}\} \leftrightarrow \{v_7,v_8,v_9\}$;
цикл $c_{20} = \{e_{21},e_{22},e_{23}\} \leftrightarrow \{v_9,v_{10},v_{11}\}$.

Плоская часть графа представлена на рис. 12.6.

Рисунок графа показан на рис. 12.7.

Рис. 12.6. Плоская часть графа с двумя удаленными ребрами.

Рис. 12.7. Граф $G_3$.

## 12.3. Плоские конфигурации и метод Гаусса

**Определение 12.2.** Подмножество изометрических циклов с мощностью меньшей цикломатического числа графа, кольцевая сумма элементов которого есть пустое множество, называется *плоской конфигурацией*.

Плоская конфигурация – это зависимая система изометрических циклов. Она порождается в результате работы алгоритма Гаусса.

Рассмотрим свойства плоских конфигураций на примере графа $G_4$:

количество вершин графа = 10;

количество ребер графа = 23;

количество изометрических циклов = 32.

Матрица смежностей графа:

```
вершина  v1:   v2   v6   v7   v10
вершина  v2:   v1   v3   v5   v7
вершина  v3:   v2   v4   v5   v9
вершина  v4:   v3   v5   v6   v7   v9
вершина  v5:   v2   v3   v4   v6   v8   v10
вершина  v6:   v1   v4   v5   v7   v9
вершина  v7:   v1   v2   v4   v6   v8
вершина  v8:   v5   v7   v9   v10
вершина  v9:   v3   v4   v6   v8   v10
вершина  v10:  v1   v5   v8   v9
```

Элементы матрицы инциденций:

ребро $e_1$: $(v_1,v_2)$ или $(v_2,v_1)$     ребро $e_2$: $(v_1,v_6)$ или $(v_6,v_1)$
ребро $e_3$: $(v_1,v_7)$ или $(v_7,v_1)$     ребро $e_4$: $(v_1,v_{10})$ или $(v_{10},v_1)$
ребро $e_5$: $(v_2,v_3)$ или $(v_3,v_2)$     ребро $e_6$: $(v_2,v_5)$ или $(v_5,v_2)$
ребро $e_7$: $(v_2,v_7)$ или $(v_7,v_2)$     ребро $e_8$: $(v_3,v_4)$ или $(v_4,v_3)$
ребро $e_9$: $(v_3,v_5)$ или $(v_5,v_3)$     ребро $e_{10}$: $(v_3,v_9)$ или $(v_9,v_3)$
ребро $e_{11}$: $(v_4,v_5)$ или $(v_5,v_4)$     ребро $e_{12}$: $(v_4,v_6)$ или $(v_6,v_4)$
ребро $e_{13}$: $(v_4,v_7)$ или $(v_7,v_4)$     ребро $e_{14}$: $(v_4,v_9)$ или $(v_9,v_4)$
ребро $e_{15}$: $(v_5,v_6)$ или $(v_6,v_5)$     ребро $e_{16}$: $(v_5,v_8)$ или $(v_8,v_5)$
ребро $e_{17}$: $(v_5,v_{10})$ или $(v_{10},v_5)$     ребро $e_{18}$: $(v_6,v_7)$ или $(v_7,v_6)$
ребро $e_{19}$: $(v_6,v_9)$ или $(v_9,v_6)$     ребро $e_{20}$: $(v_7,v_8)$ или $(v_8,v_7)$
ребро $e_{21}$: $(v_8,v_9)$ или $(v_9,v_8)$     ребро $e_{22}$: $(v_8,v_{10})$ или $(v_{10},v_8)$
ребро $e_{23}$: $(v_9,v_{10})$ или $(v_{10},v_9)$

Множество изометрических циклов графа:

$c_1 = \{e_1,e_2,e_6,e_{15}\} \leftrightarrow \{v_1,v_2,v_5,v_6\}$;
$c_2 = \{e_1,e_3,e_7\} \leftrightarrow \{v_1,v_2,v_7\}$;
$c_3 = \{e_1,e_4,e_6,e_{17}\} \leftrightarrow \{v_1,v_2,v_5,v_{10}\}$;
$c_4 = \{e_2,e_3,e_{18}\} \leftrightarrow \{v_1,v_6,v_7\}$;
$c_5 = \{e_2,e_4,e_{15},e_{17}\} \leftrightarrow \{v_1,v_5,v_6,v_{10}\}$;
$c_6 = \{e_2,e_4,e_{19},e_{23}\} \leftrightarrow \{v_1,v_6,v_9,v_{10}\}$;
$c_7 = \{e_3,e_4,e_{20},e_{22}\} \leftrightarrow \{v_1,v_7,v_8,v_{10}\}$;
$c_8 = \{e_1,e_2,e_5,e_{10},e_{19}\} \leftrightarrow \{v_1,v_2,v_3,v_6,v_9\}$;
$c_9 = \{e_1,e_4,e_5,e_{10},e_{23}\} \leftrightarrow \{v_1,v_2,v_3,v_9,v_{10}\}$;
$c_{10} = \{e_5,e_6,e_9\} \leftrightarrow \{v_2,v_3,v_5\}$;
$c_{11} = \{e_5,e_7,e_8,e_{13}\} \leftrightarrow \{v_2,v_3,v_4,v_7\}$;
$c_{12} = \{e_6,e_7,e_{11},e_{13}\} \leftrightarrow \{v_2,v_4,v_5,v_7\}$;
$c_{13} = \{e_6,e_7,e_{15},e_{18}\} \leftrightarrow \{v_2,v_5,v_6,v_7\}$;
$c_{14} = \{e_6,e_7,e_{16},e_{20}\} \leftrightarrow \{v_2,v_5,v_7,v_8\}$;
$c_{15} = \{e_8,e_9,e_{11}\} \leftrightarrow \{v_3,v_4,v_5\}$;
$c_{16} = \{e_8,e_{10},e_{14}\} \leftrightarrow \{v_3,v_4,v_9\}$;
$c_{17} = \{e_9,e_{10},e_{15},e_{19}\} \leftrightarrow \{v_3,v_5,v_6,v_9\}$;
$c_{18} = \{e_9,e_{10},e_{16},e_{21}\} \leftrightarrow \{v_3,v_5,v_8,v_9\}$;
$c_{19} = \{e_9,e_{10},e_{17},e_{23}\} \leftrightarrow \{v_3,v_5,v_9,v_{10}\}$;
$c_{20} = \{e_{11},e_{12},e_{15}\} \leftrightarrow \{v_4,v_5,v_6\}$;
$c_{21} = \{e_{11},e_{13},e_{16},e_{20}\} \leftrightarrow \{v_4,v_5,v_7,v_8\}$;
$c_{22} = \{e_{11},e_{14},e_{16},e_{21}\} \leftrightarrow \{v_4,v_5,v_8,v_9\}$;
$c_{23} = \{e_{11},e_{14},e_{17},e_{23}\} \leftrightarrow \{v_4,v_5,v_9,v_{10}\}$;
$c_{24} = \{e_{12},e_{13},e_{18}\} \leftrightarrow \{v_4,v_6,v_7\}$;
$c_{25} = \{e_{12},e_{14},e_{19}\} \leftrightarrow \{v_4,v_6,v_9\}$;
$c_{26} = \{e_{13},e_{14},e_{20},e_{21}\} \leftrightarrow \{v_4,v_7,v_8,v_9\}$;
$c_{27} = \{e_{15},e_{16},e_{18},e_{20}\} \leftrightarrow \{v_5,v_6,v_7,v_8\}$;
$c_{28} = \{e_{15},e_{16},e_{19},e_{21}\} \leftrightarrow \{v_5,v_6,v_8,v_9\}$;
$c_{29} = \{e_{15},e_{17},e_{19},e_{23}\} \leftrightarrow \{v_5,v_6,v_9,v_{10}\}$;
$c_{30} = \{e_{16},e_{17},e_{22}\} \leftrightarrow \{v_5,v_8,v_{10}\}$;
$c_{31} = \{e_{18},e_{19},e_{20},e_{21}\} \leftrightarrow \{v_6,v_7,v_8,v_9\}$;
$c_{32} = \{e_{21},e_{22},e_{23}\} \leftrightarrow \{v_8,v_9,v_{10}\}$.

Для данной последовательности изометрических циклов применим модифицированный алгоритм Гаусса. Алгоритм выделит следующую систему плоских конфигураций. Количество плоских конфигураций определим из выражения $card\, C_\tau - m + n - 1$:

Количество 0-систем циклов = 18
система 1: $\{c_5,c_1,c_3\}$;
система 2: $\{c_9,c_6,c_8\}$;
система 3: $\{c_{13},c_4,c_1,c_2\}$;

система 4: $\{c_{15}, c_{11}, c_{12}, c_{10}\}$;
система 5: $\{c_{17}, c_8, c_{10}, c_1\}$;
система 6: $\{c_{19}, c_6, c_8, c_{10}, c_3\}$;
система 7: $\{c_{21}, c_{12}, c_{14}\}$;
система 8: $\{c_{22}, c_{18}, c_{11}, c_{12}, c_{16}, c_{10}\}$;
система 9: $\{c_{23}, c_6, c_8, c_{11}, c_{12}, c_{16}, c_3\}$;
система 10: $\{c_{24}, c_{12}, c_{20}, c_4, c_1, c_2\}$;
система 11: $\{c_{25}, c_8, c_{11}, c_{12}, c_{16}, c_{20}, c_1\}$;
система 12: $\{c_{26}, c_{18}, c_{11}, c_{14}, c_{16}, c_{10}\}$;
система 13: $\{c_{27}, c_{14}, c_4, c_1, c_2\}$;
система 14: $\{c_{28}, c_8, c_{10}, c_{18}, c_1\}$;
система 15: $\{c_{29}, c_6, c_1, c_3\}$;
система 16: $\{c_{30}, c_7, c_{14}, c_2, c_3\}$;
система 17: $\{c_{31}, c_8, c_{10}, c_{18}, c_{14}, c_4, c_2\}$;
система 18: $\{c_{32}, c_8, c_{10}, c_{18}, c_{14}, c_7, c_2, c_6\}$.

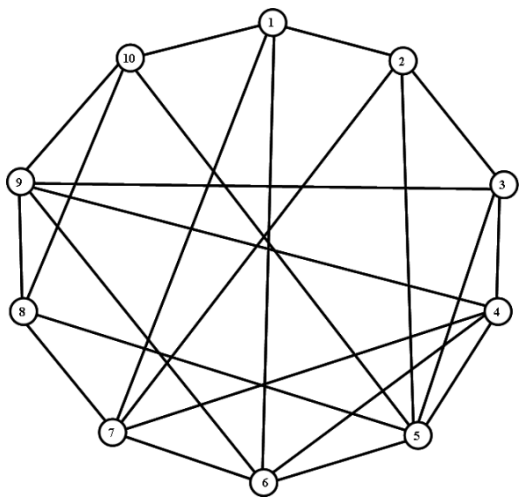 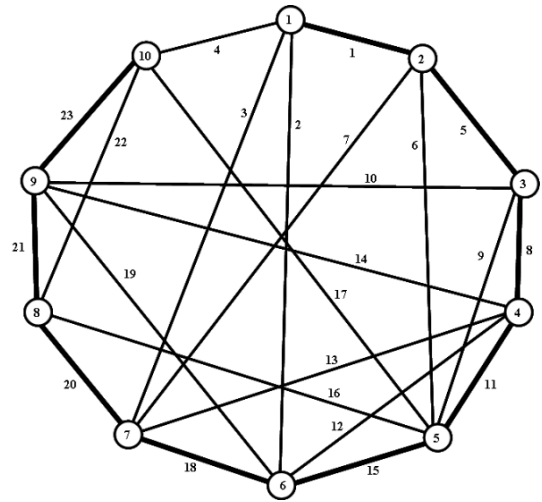

Рис. 12.8. Граф $G_4$.    Рис. 12.9. Выделение дерева в $G_4$.

| Номер цикла | 1 | 2 | 3 | 4 | 5 | 6 | 7 | 8 | 9 | 10 | 11 | 12 | 13 | 14 | 15 | 16 |
|---|---|---|---|---|---|---|---|---|---|---|---|---|---|---|---|---|
| количество систем с данным циклом | 8 | 6 | 5 | 4 | 1 | 5 | 2 | 8 | 1 | 8 | 5 | 6 | 1 | 6 | 1 | 4 |
| Номер цикла | 17 | 18 | 19 | 20 | 21 | 22 | 23 | 24 | 25 | 26 | 27 | 28 | 29 | 30 | 31 | 32 |
| количество систем с данным циклом | 1 | 4 | 1 | 2 | 1 | 1 | 1 | 1 | 1 | 1 | 1 | 1 | 1 | 1 | 1 | 1 |

Максимальное количество систем с циклами равно 8 – это три цикла: $c_1, c_8, c_{10}$. Для исключения выбираем циклы $c_8$ и $c_9$ как имеющие максимальную длину. Значение функционала Маклейна для множества изометрических циклов графа равно 2124.

Исключаем из системы циклы $c_8$ и $c_9$.

$$F\left(\frac{\delta^2 C_\tau}{\delta c_8 \delta c_9}\right) = 1704;$$

Последовательно исключаем циклы из системы с выбором минимального значения функционала Маклейна.

$$F(\frac{\delta^3 C_\tau}{\delta c_8 \delta c_9 \delta c_{28}}) = 1368;$$

$$F(\frac{\delta^4 C_\tau}{\delta c_8 \delta c_9 \delta c_{28} \delta c_{27}}) = 1122;$$

$$F(\frac{\delta^5 C_\tau}{\delta c_8 \delta c_9 \delta c_{28} \delta c_{27} \delta c_{23}}) = 930;$$

$$F(\frac{\delta^6 C_\tau}{\delta c_8 \delta c_9 \delta c_{28} \delta c_{27} \delta c_{23} \delta c_{13}}) = 780;$$

$$F(\frac{\delta^7 C_\tau}{\delta c_8 \delta c_9 \delta c_{28} \delta c_{27} \delta c_{23} \delta c_{13} \delta c_{21}}) = 636;$$

$$F(\frac{\delta^8 C_\tau}{\delta c_8 \delta c_9 \delta c_{28} \delta c_{27} \delta c_{23} \delta c_{13} \delta c_{21} \delta c_5}) = 510;$$

$$F(\frac{\delta^9 C_\tau}{\delta c_8 \delta c_9 \delta c_{28} \delta c_{27} \delta c_{23} \delta c_{13} \delta c_{21} \delta c_5 \delta c_{19}}) = 402;$$

$$F(\frac{\delta^{10} C_\tau}{\delta c_8 \delta c_9 \delta c_{28} \delta c_{27} \delta c_{23} \delta c_{13} \delta c_{21} \delta c_5 \delta c_{19} \delta c_{29}}) = 312;$$

$$F(\frac{\delta^{11} C_\tau}{\delta c_8 \delta c_9 \delta c_{28} \delta c_{27} \delta c_{23} \delta c_{13} \delta c_{21} \delta c_5 \delta c_{19} \delta c_{29} \delta c_{12}}) = 240;$$

$$F(\frac{\delta^{12} C_\tau}{\delta c_8 \delta c_9 \delta c_{28} \delta c_{27} \delta c_{23} \delta c_{13} \delta c_{21} \delta c_5 \delta c_{19} \delta c_{29} \delta c_{12} \delta c_{31}}) = 180;$$

$$F(\frac{\delta^{13} C_\tau}{\delta c_8 \delta c_9 \delta c_{28} \delta c_{27} \delta c_{23} \delta c_{13} \delta c_{21} \delta c_5 \delta c_{19} \delta c_{29} \delta c_{12} \delta c_{31} \delta c_3}) = 144;$$

$$F(\frac{\delta^{14} C_\tau}{\delta c_8 \delta c_9 \delta c_{28} \delta c_{27} \delta c_{23} \delta c_{13} \delta c_{21} \delta c_5 \delta c_{19} \delta c_{29} \delta c_{12} \delta c_{31} \delta c_3 \delta c_{17}}) = 108;$$

$$F(\frac{\delta^{15} C_\tau}{\delta c_8 \delta c_9 \delta c_{28} \delta c_{27} \delta c_{23} \delta c_{13} \delta c_{21} \delta c_5 \delta c_{19} \delta c_{29} \delta c_{12} \delta c_{31} \delta c_3 \delta c_{17} \delta c_{22}}) = 72;$$

$$F\left(\frac{\delta^{16}C_\tau}{\delta c_8 \delta c_9 \delta c_{28} \delta c_{27} \delta c_{23} \delta c_{13} \delta c_{21} \delta c_5 \delta c_{19} \delta c_{29} \delta c_{12} \delta c_{31} \delta c_3 \delta c_{17} \delta c_{22} \delta c_6}\right) = 54;$$

$$F\left(\frac{\delta^{17}C_\tau}{\delta c_8 \delta c_9 \delta c_{28} \delta c_{27} \delta c_{23} \delta c_{13} \delta c_{21} \delta c_5 \delta c_{19} \delta c_{29} \delta c_{12} \delta c_{31} \delta c_3 \delta c_{17} \delta c_{22} \delta c_6 \delta c_{11}}\right) = 36;$$

$$F\left(\frac{\delta^{18}C_\tau}{\delta c_8 \delta c_9 \delta c_{28} \delta c_{27} \delta c_{23} \delta c_{13} \delta c_{21} \delta c_5 \delta c_{19} \delta c_{29} \delta c_{12} \delta c_{31} \delta c_3 \delta c_{17} \delta c_{22} \delta c_6 \delta c_{11} \delta c_7}\right) = 18;$$

$P_e$ = <2,2,2,1,1,2,2,2,3,2,2,3,2,3,2,2,2,2,1,2,2,2,2>.

$$\frac{\partial C_b}{\partial c_{10}} = 12; \qquad \frac{\partial C_b}{\partial c_{25}} = 6.$$

Исключаем цикл $c_{25}$.

$P_e$ = <2,2,2,1,1,2,2,2,3,2,2,2,2,2,2,2,2,2,0,2,2,2,2>.

$$\frac{\partial^2 C_b}{\partial c_{25} \partial c_{10}} = 0.$$

Исключаем цикл $c_{10}$.

$P_e$ = <2,2,2,0,0,2,2,2,2,2,2,2,2,2,2,2,2,2,0,2,2,2,2>.

Может показаться, что мы нарушаем правило Эйлера, но на самом деле исключение производится с учетом подключения обода в состав системы циклов, характеризующих плоскую часть графа. Плоская часть графа имеет вид (рис. 12.10):

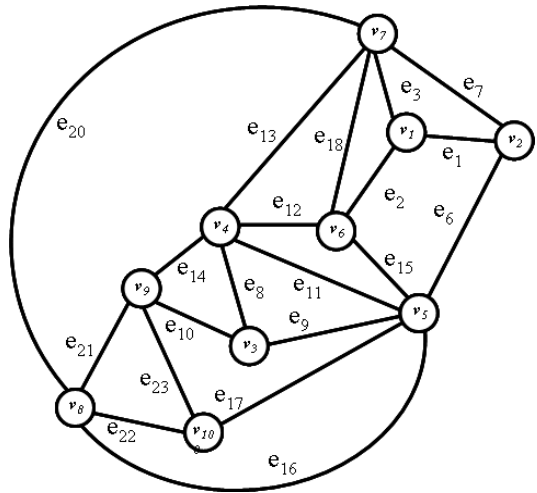

Рис. 12.10. Плоская часть графа $G_4$.

$c_1 = \{e_1, e_2, e_6, e_{15}\} \leftrightarrow \{v_1, v_2, v_5, v_6\}$;
$c_2 = \{e_1, e_3, e_7\} \leftrightarrow \{v_1, v_2, v_7\}$;
$c_4 = \{e_2, e_3, e_{18}\} \leftrightarrow \{v_1, v_6, v_7\}$;
$c_{14} = \{e_6, e_7, e_{16}, e_{20}\} \leftrightarrow \{v_2, v_5, v_7, v_8\}$;
$c_{15} = \{e_8, e_9, e_{11}\} \leftrightarrow \{v_3, v_4, v_5\}$;

$c_{16} = \{e_8, e_{10}, e_{14}\} \leftrightarrow \{v_3, v_4, v_9\}$;
$c_{19} = \{e_9, e_{10}, e_{17}, e_{23}\} \leftrightarrow \{v_3, v_5, v_9, v_{10}\}$;
$c_{20} = \{e_{11}, e_{12}, e_{15}\} \leftrightarrow \{v_4, v_5, v_6\}$;
$c_{24} = \{e_{12}, e_{13}, e_{18}\} \leftrightarrow \{v_4, v_6, v_7\}$;
$c_{26} = \{e_{13}, e_{14}, e_{20}, e_{21}\} \leftrightarrow \{v_4, v_7, v_8, v_9\}$;
$c_{30} = \{e_{16}, e_{17}, e_{22}\} \leftrightarrow \{v_5, v_8, v_{10}\}$;
$c_{32} = \{e_{21}, e_{22}, e_{23}\} \leftrightarrow \{v_8, v_9, v_{10}\}$.

Цикломатическое число графа без трех ребер равно 11.

## 12.4. Таблица циклов

Таблица изометрических циклов располагает циклы в последовательности увеличения их длин. Такое расположение, а также определение направления выбора, отражает многокритериальность задачи выделения плоского топологического рисунка.

Будем рассматривать несепарабельные графы (см. Определение 11.1).

Выделим максимально плоский суграф в графе $G_5$ (рис. 12.11):

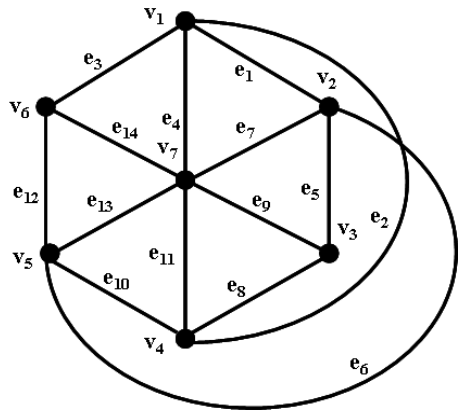

Рис. 12.11. Граф $G_5$.

количество вершин графа = 7;
количество ребер графа = 14;
количество изометрических циклов = 13;
цикломатическое число графа = 8.

Матрица смежностей графа:

вершина $v_1$: $v_2$ $v_4$ $v_6$ $v_7$
вершина $v_2$: $v_1$ $v_3$ $v_5$ $v_7$
вершина $v_3$: $v_2$ $v_4$ $v_7$
вершина $v_4$: $v_1$ $v_3$ $v_5$ $v_7$
вершина $v_5$: $v_2$ $v_4$ $v_6$ $v_7$
вершина $v_6$: $v_1$ $v_5$ $v_7$
вершина $v_7$: $v_1$ $v_2$ $v_3$ $v_4$ $v_5$ $v_6$

Элементы матрицы инциденций:

ребро $e_1$: $(v_1, v_2)$ или $(v_2, v_1)$;   ребро $e_2$: $(v_1, v_4)$ или $(v_4, v_1)$;
ребро $e_3$: $(v_1, v_6)$ или $(v_6, v_1)$;   ребро $e_4$: $(v_1, v_7)$ или $(v_7, v_1)$;
ребро $e_5$: $(v_2, v_3)$ или $(v_3, v_2)$;   ребро $e_6$: $(v_2, v_5)$ или $(v_5, v_2)$;
ребро $e_7$: $(v_2, v_7)$ или $(v_7, v_2)$;   ребро $e_8$: $(v_3, v_4)$ или $(v_4, v_3)$;
ребро $e_9$: $(v_3, v_7)$ или $(v_7, v_3)$;   ребро $e_{10}$: $(v_4, v_5)$ или $(v_5, v_4)$;

ребро  $e_{11}$: $(v_4,v_7)$ или $(v_7,v_4)$;   ребро  $e_{12}$: $(v_5,v_6)$ или $(v_6,v_5)$;
ребро  $e_{13}$: $(v_5,v_7)$ или $(v_7,v_5)$;   ребро  $e_{14}$: $(v_6,v_7)$ или $(v_7,v_6)$.

Множество изометрических циклов графа: $C_\tau$

цикл   $c_1 = \{e_1,e_2,e_5,e_8\} \leftrightarrow \{v_1,v_2,v_3,v_4\}$;
цикл   $c_2 = \{e_1,e_2,e_6,e_{10}\} \leftrightarrow \{v_1,v_2,v_4,v_5\}$;
цикл   $c_3 = \{e_1,e_3,e_6,e_{12}\} \leftrightarrow \{v_1,v_2,v_5,v_6\}$;
цикл   $c_4 = \{e_1,e_4,e_7\} \leftrightarrow \{v_1,v_2,v_7\}$;
цикл   $c_5 = \{e_2,e_3,e_{10},e_{12}\} \leftrightarrow \{v_1,v_4,v_5,v_6\}$;
цикл   $c_6 = \{e_2,e_4,e_{11}\} \leftrightarrow \{v_1,v_4,v_7\}$;
цикл   $c_7 = \{e_3,e_4,e_{14}\} \leftrightarrow \{v_1,v_6,v_7\}$;
цикл   $c_8 = \{e_5,e_6,e_8,e_{10}\} \leftrightarrow \{v_2,v_3,v_4,v_5\}$;
цикл   $c_9 = \{e_5,e_7,e_9\} \leftrightarrow \{v_2,v_3,v_7\}$;
цикл   $c_{10} = \{e_6,e_7,e_{13}\} \leftrightarrow \{v_2,v_5,v_7\}$;
цикл   $c_{11} = \{e_8,e_9,e_{11}\} \leftrightarrow \{v_3,v_4,v_7\}$;
цикл   $c_{12} = \{e_{10},e_{11},e_{13}\} \leftrightarrow \{v_4,v_5,v_7\}$;
цикл   $c_{13} = \{e_{12},e_{13},e_{14}\} \leftrightarrow \{v_5,v_6,v_7\}$.

Образуем вектор количества циклов, проходящих по ребрам графа
$P_e = <4,4,3,3,3,4,3,3,2,4,3,3,3,2>$.

Вектор количества циклов проходящих по ребру $P_e$ представляется в виде кортежа. Первый элемент кортежа, равный 4, характеризует количество циклов, проходящих по ребру $e_1$. Второй элемент кортежа, равный 4, характеризует количество циклов, проходящих по ребру $e_2$. Третий элемент кортежа, равный 3, характеризует количество циклов, проходящих по ребру $e_3$. И так далее.

Любому подмножеству из множества изометрических циклов графа можно поставить в соответствие кубический функционал Маклейна:

$$\sum_{i=1}^{m} p_i (p_i - 1)(p_i - 2) = \sum_{i=1}^{m} p_i^3 - 3\sum_{i=1}^{m} p_i^2 + 2\sum_{i=1}^{m} p_i$$

где $p_i$ – количество циклов, включающих ребро $e_i$.

Рассчитаем значение функционала Маклейна для графа $G_1$
$F(G_1) = 24+24+6+6+6+24+6+6+0+24+6+6+6+0=144$.

Составим таблицу циклов, располагая циклы по возрастанию их длин.

|  | 1 | 2 | 3 | 4 | 5 |
|---|---|---|---|---|---|
| цикл   $c_4 = \{e_1,e_4,e_7\}$ |  |  |  |  |  |
| цикл   $c_6 = \{e_2,e_4,e_{11}\}$ |  |  |  |  |  |
| цикл   $c_7 = \{e_3,e_4,e_{14}\}$ |  |  |  |  |  |
| цикл   $c_9 = \{e_5,e_7,e_9\}$ |  |  |  |  |  |
| цикл   $c_{10} = \{e_6,e_7,e_{13}\}$ |  |  |  |  |  |
| цикл   $c_{11} = \{e_8,e_9,e_{11}\}$ |  |  |  |  |  |
| цикл   $c_{12} = \{e_{10},e_{11},e_{13}\}$ |  |  |  |  |  |
| цикл   $c_{13} = \{e_{12},e_{13},e_{14}\}$ |  |  |  |  |  |
| цикл   $c_1 = \{e_1,e_2,e_5,e_8\}$ |  |  |  |  |  |

| 1 | 2 | 3 | 4 | 5 |
|---|---|---|---|---|
| цикл $c_2 = \{e_1,e_2,e_6,e_{10}\}$ | | | | |
| цикл $c_3 = \{e_1,e_3,e_6,e_{12}\}$ | | | | |
| цикл $c_5 = \{e_2,e_3,e_{10},e_{12}\}$ | | | | |
| цикл $c_8 = \{e_5,e_6,e_8,e_{10}\}$ | | | | |
| Значение функционала Маклейна | 144 | | | |

Столбец 1 характеризует запись цикла в ребрах. Столбец 2 характеризует запись значения функционала Маклейна при исключении цикла из системы. Столбец 3 характеризует предыдущую запись значения функционала Маклейна. Столбец 4 характеризует скорость изменения значения функционала Маклейна. В столбце 5 расспологается признак последующего исключения цикла.

Будем последовательно исключать циклы для получения топологического рисунка графа.

| 1 | 2 | 3 | 4 | 5 |
|---|---|---|---|---|
| цикл $c_4 = \{e_1,e_4,e_7\}$ | 114 | | | |
| цикл $c_6 = \{e_2,e_4,e_{11}\}$ | 114 | | | |
| цикл $c_7 = \{e_3,e_4,e_{14}\}$ | 132 | | | |
| цикл $c_9 = \{e_5,e_7,e_9\}$ | 132 | | | |
| цикл $c_{10} = \{e_6,e_7,e_{13}\}$ | 114 | | | |
| цикл $c_{11} = \{e_8,e_9,e_{11}\}$ | 132 | | | |
| цикл $c_{12} = \{e_{10},e_{11},e_{13}\}$ | 114 | | | |
| цикл $c_{13} = \{e_{12},e_{13},e_{14}\}$ | 132 | | | |
| цикл $c_1 = \{e_1,e_2,e_5,e_8\}$ | 96 | | | |
| цикл $c_2 = \{e_1,e_2,e_6,e_{10}\}$ | 72 | | | *min* |
| цикл $c_3 = \{e_1,e_3,e_6,e_{12}\}$ | 96 | | | |
| цикл $c_5 = \{e_2,e_3,e_{10},e_{12}\}$ | 96 | | | |
| цикл $c_8 = \{e_5,e_6,e_8,e_{10}\}$ | 96 | | | |
| Значение функционала Маклейна | 144 | | | |

Исключение из подмножества цикла $c_2$ определяет минимальное значение функционала Маклейна.

Вектор $P_e$ при исключении из множества $C_\tau$ цикла $c_2$ имеет вид:

$P_e = <3,3,3,3,3,3,3,3,2,3,3,3,3,2>$.

Продолжаем процесс исключения циклов

| 1 | 2 | 3 | 4 | 5 |
|---|---|---|---|---|
| цикл $c_4 = \{e_1,e_4,e_7\}$ | 54 | 114 | | |
| цикл $c_6 = \{e_2,e_4,e_{11}\}$ | 54 | 114 | | |
| цикл $c_7 = \{e_3,e_4,e_{14}\}$ | 60 | 132 | | |
| цикл $c_9 = \{e_5,e_7,e_9\}$ | 60 | 132 | | |
| цикл $c_{10} = \{e_6,e_7,e_{13}\}$ | 54 | 114 | | |
| цикл $c_{11} = \{e_8,e_9,e_{11}\}$ | 60 | 132 | | |
| цикл $c_{12} = \{e_{10},e_{11},e_{13}\}$ | 54 | 114 | | |
| цикл $c_{13} = \{e_{12},e_{13},e_{14}\}$ | 60 | 132 | | |
| цикл $c_1 = \{e_1,e_2,e_5,e_8\}$ | 48 | 96 | | |
| цикл $c_3 = \{e_1,e_3,e_6,e_{12}\}$ | 48 | 96 | | |
| цикл $c_5 = \{e_2,e_3,e_{10},e_{12}\}$ | 48 | 96 | | |

| | 1 | 2 | 3 | 4 | 5 |
|---|---|---|---|---|---|
| цикл $c_8 = \{e_5,e_6,e_8,e_{10}\}$ | 48 | 96 | *min* | | |
| Значение функционала Маклейна | 72 | | | | |

Среди кандидатов на удаление, имеющих минимальное значение и равные скорости измениения значения функционала Маклейна, выбираем первый с конца таблицы циклов элемент – $c_8$.

Вектор $P_e$ при исключении из множества $C_\tau$ циклов $c_2,c_8$ имеет вид:

$P_e = <3,3,3,3,2,2,3,2,2,2,3,3,3,2>$.

Продолжаем процесс исключения циклов

| | 1 | 2 | 3 | 4 | 5 |
|---|---|---|---|---|---|
| цикл $c_4 = \{e_1,e_4,e_7\}$ | 30 | 54 | 24 | | |
| цикл $c_6 = \{e_2,e_4,e_{11}\}$ | 30 | 54 | 24 | | |
| цикл $c_7 = \{e_3,e_4,e_{14}\}$ | 36 | 60 | | | |
| цикл $c_9 = \{e_5,e_7,e_9\}$ | 42 | 60 | | | |
| цикл $c_{10} = \{e_6,e_7,e_{13}\}$ | 36 | 54 | | | |
| цикл $c_{11} = \{e_8,e_9,e_{11}\}$ | 42 | 60 | | | |
| цикл $c_{12} = \{e_{10},e_{11},e_{13}\}$ | 36 | 54 | | | |
| цикл $c_{13} = \{e_{12},e_{13},e_{14}\}$ | 36 | 60 | | | |
| цикл $c_1 = \{e_1,e_2,e_5,e_8\}$ | 36 | 48 | | | |
| цикл $c_3 = \{e_1,e_3,e_6,e_{12}\}$ | 30 | 48 | 18 | | |
| цикл $c_5 = \{e_2,e_3,e_{10},e_{12}\}$ | 30 | 48 | 18 | *min* | |
| Значение функционала Маклейна | 48 | | | | |

Для исключения выбираем цикл $c_5$, имеющий при исключении минимальное значение функционала Маклейна и большую длину.

Вектор $P_e$ при исключении из множества $C_\tau$ циклов $c_2,c_8,c_5$ имеет вид:

$P_e = <3,2,2,3,2,2,3,2,2,1,3,2,3,2>$.

Продолжаем процесс исключения циклов

| | 1 | 2 | 3 | 4 | 5 |
|---|---|---|---|---|---|
| цикл $c_4 = \{e_1,e_4,e_7\}$ | 12 | 30 | | *min* | |
| цикл $c_6 = \{e_2,e_4,e_{11}\}$ | 18 | 30 | | | |
| цикл $c_7 = \{e_3,e_4,e_{14}\}$ | 24 | 36 | | | |
| цикл $c_9 = \{e_5,e_7,e_9\}$ | 24 | 42 | | | |
| цикл $c_{10} = \{e_6,e_7,e_{13}\}$ | 18 | 36 | | | |
| цикл $c_{11} = \{e_8,e_9,e_{11}\}$ | 24 | 42 | | | |
| цикл $c_{12} = \{e_{10},e_{11},e_{13}\}$ | 18 | 36 | | | |
| цикл $c_{13} = \{e_{12},e_{13},e_{14}\}$ | 24 | 36 | | | |
| цикл $c_1 = \{e_1,e_2,e_5,e_8\}$ | 24 | 36 | | | |
| цикл $c_3 = \{e_1,e_3,e_6,e_{12}\}$ | 24 | 30 | | | |
| Значение функционала Маклейна | 30 | | | | |

Для исключения выбираем цикл $c_4$, имеющий при исключении минимальное значение функционала Маклейна.

Вектор $P_e$ при исключении из множества $C_\tau$ циклов $c_2,c_8,c_5,c_4$ имеет вид:

$P_e = <2,2,2,2,2,2,2,2,2,1,3,2,3,2>$.

Продолжаем процесс исключения циклов

| 1 | 2 | 3 | 4 | 5 |
|---|---|---|---|---|
| цикл $c_6 = \{e_2,e_4,e_{11}\}$ | 6 | 18 | 12 | |
| цикл $c_7 = \{e_3,e_4,e_{14}\}$ | 12 | 24 | | |
| цикл $c_9 = \{e_5,e_7,e_9\}$ | 12 | 24 | | |
| цикл $c_{10} = \{e_6,e_7,e_{13}\}$ | 6 | 18 | 12 | |
| цикл $c_{11} = \{e_8,e_9,e_{11}\}$ | 6 | 24 | 18 | |
| цикл $c_{12} = \{e_{10},e_{11},e_{13}\}$ | 0 | 18 | | |
| цикл $c_{13} = \{e_{12},e_{13},e_{14}\}$ | 6 | 24 | 18 | *min* |
| цикл $c_1 = \{e_1,e_2,e_5,e_8\}$ | 12 | 24 | | |
| цикл $c_3 = \{e_1,e_3,e_6,e_{12}\}$ | 12 | 24 | | |
| Значение функционала Маклейна | 12 | | | |

Минимальное значение функционал Маклейна принимает после исключения цикла $c_{12}$. Однако в этом случае значение координаты для десятого ребра в векторе $P_e$ равно нулю, что не допустимо для подмножества циклов, характеризующих базис подпространства циклов.

Поэтому, для исключения выбираем цикл $c_{13}$, имеющий при исключении следующее минимальное значение функционала Маклейна и большую скорость изменения значений функционала Маклейна.

Вектор $P_e$ при исключении из множества $C_\tau$ циклов $c_2,c_8,c_5,c_4,c_{13}$ имеет вид:

$P_e = <2,2,2,2,2,2,2,2,2,1,3,1,2,1>$.

Количество оставшихся циклов в подмножестве равно цикломатическому числу графа и образует базис подпространства циклов C.

Продолжим процесс исключения циклов для получения нулевого значения функционала Маклейна.

| 1 | 2 | 3 | 4 | 5 |
|---|---|---|---|---|
| цикл $c_6 = \{e_2,e_4,e_{11}\}$ | 0 | | | нарушение |
| цикл $c_7 = \{e_3,e_4,e_{14}\}$ | 6 | | | |
| цикл $c_9 = \{e_5,e_7,e_9\}$ | 6 | | | |
| цикл $c_{10} = \{e_6,e_7,e_{13}\}$ | 6 | | | |
| цикл $c_{11} = \{e_8,e_9,e_{11}\}$ | 0 | | | нарушение |
| цикл $c_{12} = \{e_{10},e_{11},e_{13}\}$ | 0 | | | *min* |
| цикл $c_1 = \{e_1,e_2,e_5,e_8\}$ | 6 | | | |
| цикл $c_3 = \{e_1,e_3,e_6,e_{12}\}$ | 6 | | | |
| Значение функционала Маклейна | 6 | | | |

Исключение циклов $c_6,c_{11}$ приводит к нарушению условия Эйлера (удаление ребра должно приводить к исключению одного и только одного цикла).

Вектор $P_e$ при исключении из множества $C_\tau$ циклов $c_2,c_8,c_5,c_4,c_{13},c_{12}$ имеет вид:

$P_e = <2,2,2,2,2,2,2,2,2,0,2,1,1,1>$.

Значение функционала Маклейна для выбранного подмножества изометрических циклов равно нулю и характеризует плоский топологический рисунок графа (рис. 12.12).

цикл $c_1 = \{e_1,e_2,e_5,e_8\} \leftrightarrow \{v_1,v_2,v_3,v_4\}$;

цикл $c_3 = \{e_1, e_3, e_6, e_{12}\} \leftrightarrow \{v_1, v_2, v_5, v_6\}$;
цикл $c_6 = \{e_2, e_4, e_{11}\} \leftrightarrow \{v_1, v_4, v_7\}$;
цикл $c_7 = \{e_3, e_4, e_{14}\} \leftrightarrow \{v_1, v_6, v_7\}$;
цикл $c_9 = \{e_5, e_7, e_9\} \leftrightarrow \{v_2, v_3, v_7\}$;
цикл $c_{10} = \{e_6, e_7, e_{13}\} \leftrightarrow \{v_2, v_5, v_7\}$;
цикл $c_{11} = \{e_8, e_9, e_{11}\} \leftrightarrow \{v_3, v_4, v_7\}$.

Обод выбранной системы циклов равен кольцевой сумме циклов:

$c_0 = c_{13} = c_1 \oplus c_3 \oplus c_6 \oplus c_7 \oplus c_9 \oplus c_{10} \oplus c_{11} = \{e_{12}, e_{13}, e_{14}\}$.

В результате можно поменять базис подпространства циклов:

$c_0 = c_3 = c_1 \oplus c_{13} \oplus c_6 \oplus c_7 \oplus c_9 \oplus c_{10} \oplus c_{11} = \{e_1, e_3, e_6, e_{12}\}$.

Независимая система циклов характеризует другой плоский топологический рисунок (рис. 12.13).

цикл $c_1 = \{e_1, e_2, e_5, e_8\} \leftrightarrow \{v_1, v_2, v_3, v_4\}$;
цикл $c_3 = \{e_1, e_3, e_6, e_{12}\} \leftrightarrow \{v_1, v_2, v_5, v_6\}$;
цикл $c_6 = \{e_2, e_4, e_{11}\} \leftrightarrow \{v_1, v_4, v_7\}$;
цикл $c_7 = \{e_3, e_4, e_{14}\} \leftrightarrow \{v_1, v_6, v_7\}$;
цикл $c_9 = \{e_5, e_7, e_9\} \leftrightarrow \{v_2, v_3, v_7\}$;
цикл $c_{10} = \{e_6, e_7, e_{13}\} \leftrightarrow \{v_2, v_5, v_7\}$;
цикл $c_{11} = \{e_8, e_9, e_{11}\} \leftrightarrow \{v_3, v_4, v_7\}$.

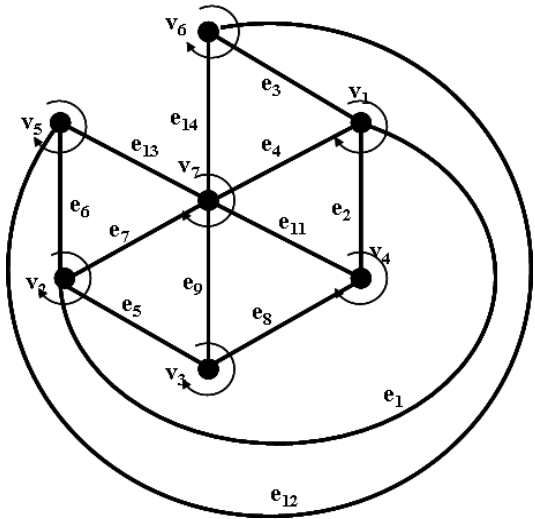 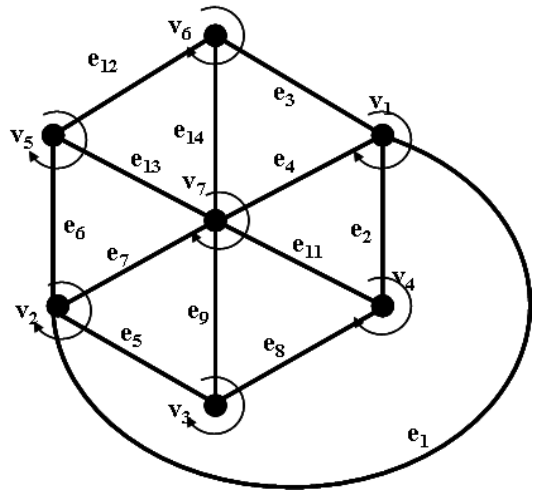

Рис. 12.12. Топологический рисунок графа.   Рис. 12.13. Топологический рисунок графа.

### Комментарии

В главе рассмотрен метод наискорейшего спуска для выделения базиса изометрических циклов подпространства квазициклов C(G) с минимальным значением функционала Маклейна. Показано, что алгоритм, построенный на основе метода наискорейшего спуска, носит полиномиальный характер и имеет вычислительную сложность $o(n^6)$.

Описан процесс построения топологического рисунка плоской части непланарного графа. Процесс построения разбит на два этапа:

1) этап выделения базиса изометрических циклов подпространства квазициклов $C(G)$ с минимальным значением кубического функционала Маклейна;

2) этап удаления циклов из базиса с обязательным выполнением условия Эйлера до достижения нулевого значения функционала Маклейна.

В разделе 12.3 рассмотрен вопрос влияния длин изометрических циклов на построение плоской конфигурации. Показано, что выделение максимально плоского суграфа следует производить из базиса, состоящего из подмножества изометрических циклов как можно меньшей длины. Причем подмножество первоначального удаления циклов определяется условием Эйлера. Например, исключение трех циклов $c_8, c_{10}, c_{14}$ из базиса графа $G_3$ приводит к удалению трех ребер и одновременно порождает подмножество циклов $c_5, c_{18}$, исключение которых приводит к построению плоского топологического рисунка графа $G_3$ с двумя удаленными ребрами.

## Глава 13. Дополнение топологического рисунка
## 13.1. Дополнение плоских конфигураций

Ранее мы рассмотрели процесс построения топологического рисунка плоского суграфа, состоящего из двух этапов. Однако существует и третий этап, наглядно проявляющийся на графах с большим числом ребер и вершин.

Далее будем рассматривать свойства подмножеств изометрических циклов, имеющих нулевое значение функционала Маклейна, на примере случайного графа $G_3$, состоящего из 20 вершин и 61 ребер:

Количество вершин графа = 20.
Количество ребер графа = 61.
Количество единичных циклов = 117.

Смежность графа:

```
вершина  v₁:  v₂  v₃  v₄  v₅  v₆  v₈  v₁₅ v₂₀
вершина  v₂:  v₁  v₇  v₉  v₁₈
вершина  v₃:  v₁  v₅  v₁₁ v₁₄ v₁₅ v₁₇ v₁₈ v₁₉
вершина  v₄:  v₁  v₅  v₁₁ v₁₄ v₁₅ v₁₇ v₂₀
вершина  v₅:  v₁  v₃  v₄  v₁₁ v₁₄ v₁₇ v₁₈
вершина  v₆:  v₁  v₁₂ v₁₃ v₁₄ v₁₅ v₁₆
вершина  v₇:  v₂  v₁₁ v₁₂ v₁₃ v₁₆ v₁₈
вершина  v₈:  v₁  v₁₂ v₁₉
вершина  v₉:  v₂  v₁₁ v₁₃ v₁₅ v₁₆ v₁₇
вершина  v₁₀: v₁₃ v₁₄ v₁₈
вершина  v₁₁: v₃  v₄  v₅  v₇  v₉  v₁₂ v₁₄ v₁₇
вершина  v₁₂: v₆  v₇  v₈  v₁₁ v₁₄ v₁₇ v₁₈ v₁₉
вершина  v₁₃: v₆  v₇  v₉  v₁₀
вершина  v₁₄: v₃  v₄  v₅  v₆  v₁₀ v₁₁ v₁₂
вершина  v₁₅: v₁  v₃  v₄  v₆  v₉  v₁₆ v₁₇
вершина  v₁₆: v₆  v₇  v₉  v₁₅
вершина  v₁₇: v₃  v₄  v₅  v₉  v₁₁ v₁₂ v₁₅ v₁₈ v₂₀
вершина  v₁₈: v₂  v₃  v₅  v₇  v₁₀ v₁₂ v₁₇ v₂₀
вершина  v₁₉: v₃  v₈  v₁₂ v₂₀
вершина  v₂₀: v₁  v₄  v₁₇ v₁₈ v₁₉
```

Элементы матрицы инциденций:

ребро e₁:  (v₁,v₂) или (v₂,v₁);           ребро e₂:  (v₁,v₃) или (v₃,v₁);
ребро e₃:  (v₁,v₄) или (v₄,v₁);           ребро e₄:  (v₁,v₅) или (v₅,v₁);
ребро e₅:  (v₁,v₆) или (v₆,v₁);           ребро e₆:  (v₁,v₈) или (v₈,v₁);
ребро e₇:  (v₁,v₁₅) или (v₁₅,v₁);         ребро e₈:  (v₁,v₂₀) или (v₂₀,v₁);
ребро e₉:  (v₂,v₇) или (v₇,v₂);           ребро e₁₀: (v₂,v₉) или (v₉,v₂);
ребро e₁₁: (v₂,v₁₈) или (v₁₈,v₂);         ребро e₁₂: (v₃,v₅) или (v₅,v₃);
ребро e₁₃: (v₃,v₁₁) или (v₁₁,v₃);         ребро e₁₄: (v₃,v₁₄) или (v₁₄,v₃);
ребро e₁₅: (v₃,v₁₅) или (v₁₅,v₃);         ребро e₁₆: (v₃,v₁₇) или (v₁₇,v₃);
ребро e₁₇: (v₃,v₁₈) или (v₁₈,v₃);         ребро e₁₈: (v₃,v₁₉) или (v₁₉,v₃);
ребро e₁₉: (v₄,v₅) или (v₅,v₄);           ребро e₂₀: (v₄,v₁₁) или (v₁₁,v₄);
ребро e₂₁: (v₄,v₁₄) или (v₁₄,v₄);         ребро e₂₂: (v₄,v₁₅) или (v₁₅;v₄);
ребро e₂₃: (v₄,v₁₇) или (v₁₇,v₄);         ребро e₂₄: (v₄,v₂₀) или (v₂₀,v₄);

ребро $e_{25}$: $(v_5,v_{11})$ или $(v_{11},v_5)$;           ребро $e_{26}$: $(v_5,v_{14})$ или $(v_{14},v_5)$;
ребро $e_{27}$: $(v_5,v_{17})$ или $(v_{17},v_5)$;           ребро $e_{28}$: $(v_5,v_{18})$ или $(v_{18},v_5)$;
ребро $e_{29}$: $(v_6,v_{12})$ или $(v_{12},v_6)$;           ребро $e_{30}$: $(v_6,v_{13})$ или $(v_{13},v_6)$;
ребро $e_{31}$: $(v_6,v_{14})$ или $(v_{14},v_6)$;           ребро $e_{32}$: $(v_6,v_{15})$ или $(v_{15},v_6)$;
ребро $e_{33}$: $(v_6,v_{16})$ или $(v_{16},v_6)$;           ребро $e_{34}$: $(v_7,v_{11})$ или $(v_{11},v_7)$;
ребро $e_{35}$: $(v_7,v_{12})$ или $(v_{12},v_7)$;           ребро $e_{36}$: $(v_7,v_{13})$ или $(v_{13},v_7)$;
ребро $e_{37}$: $(v_7,v_{16})$ или $(v_{16},v_7)$;           ребро $e_{38}$: $(v_7,v_{18})$ или $(v_{18},v_7)$;
ребро $e_{39}$: $(v_8,v_{12})$ или $(v_{12},v_8)$;           ребро $e_{40}$: $(v_8,v_{19})$ или $(v_{19},v_8)$
ребро $e_{41}$: $(v_9,v_{11})$ или $(v_{11},v_9)$;           ребро $e_{42}$: $(v_9,v_{13})$ или $(v_{13},v_9)$;
ребро $e_{43}$: $(v_9,v_{15})$ или $(v_{15},v_9)$;           ребро $e_{44}$: $(v_9,v_{16})$ или $(v_{16},v_9)$;
ребро $e_{45}$: $(v_9,v_{17})$ или $(v_{17},v_9)$;           ребро $e_{46}$: $(v_{10},v_{13})$ или $(v_{13},v_{10})$;
ребро $e_{47}$: $(v_{10},v_{14})$ или $(v_{14},v_{10})$;     ребро $e_{48}$: $(v_{10},v_{18})$ или $(v_{18},v_{10})$;
ребро $e_{49}$: $(v_{11},v_{12})$ или $(v_{12},v_{11})$;     ребро $e_{50}$: $(v_{11},v_{14})$ или $(v_{14},v_{11})$;
ребро $e_{51}$: $(v_{11},v_{17})$ или $(v_{17},v_{11})$;     ребро $e_{52}$: $(v_{12},v_{14})$ или $(v_{14},v_{12})$;
ребро $e_{53}$: $(v_{12},v_{17})$ или $(v_{17},v_{12})$;     ребро $e_{54}$: $(v_{12}.v_{18})$ или $(v_{18},v_{12})$;
ребро $e_{55}$: $(v_{12},v_{19})$ или $(v_{19},v_{12})$;     ребро $e_{56}$: $(v_{15},v_{16})$ или $(v_{16},v_{15})$;
ребро $e_{57}$: $(v_{15},v_{17})$ или $(v_{17},v_{15})$;     ребро $e_{58}$: $(v_{17,},v_{18})$ или $(v_{18},v_{17})$;
ребро $e_{59}$: $(v_{17},v_{20})$ или $(v_{20},v_{17})$;     ребро $e_{60}$: $(v_{18},v_{20})$ или $(v_{20},v_{18})$;
ребро $e_{61}$: $(v_{19},v_{20})$ или $(v_{20},v_{19})$.

Множество изометрических циклов графа:

цикл $c_1 = \{e_1,e_2,e_{11},e_{17}\} \leftrightarrow \{v_1,v_2,v_3,v_{18}\}$;
цикл $c_2 = \{e_1,e_3,e_9,e_{20},e_{34}\} \leftrightarrow \{v_1,v_2,v_4,v_7,v_{11}\}$;
цикл $c_3 = \{e_1,e_4,e_{11},e_{28}\} \leftrightarrow \{v_1,v_2,v_5,v_{18}\}$;
цикл $c_4 = \{e_1,e_5,e_9,e_{29},e_{35}\} \leftrightarrow \{v_1,v_2,v_6,v_7,v_{12}\}$;
цикл $c_5 = \{e_1,e_5,e_9,e_{30},e_{36}\} \leftrightarrow \{v_1,v_2,v_6,v_7,v_{13}\}$;
цикл $c_6 = \{e_1, e_5,e_9,e_{33},e_{37}\} \leftrightarrow \{v_1,v_2,v_6,v_7,v_{16}\}$;
цикл $c_7 = \{e_1,e_6,e_9,e_{35},e_{39}\} \leftrightarrow \{v_1,v_2,v_7,v_8,v_{12}\}$;
цикл $c_8 = \{e_1,e_7,e_{10},e_{43}\} \leftrightarrow \{v_1,v_2,v_9,v_{15}\}$;
цикл $c_9 = \{e_1,e_8,e_{11},e_{60}\} \leftrightarrow \{v_1,v_2,v_{18},v_{20}\}$;
цикл $c_{10} = \{e_2,e_3,e_{13},e_{20}\} \leftrightarrow \{v_1,v_3,v_4,v_{11}\}$;
цикл $c_{11} = \{e_2,e_3,e_{14},e_{21}\} \leftrightarrow \{v_1,v_3,v_4,v_{14}\}$;
цикл $c_{12} = \{e_2,e_3,e_{16},e_{23}\} \leftrightarrow \{v_1,v_3,v_4,v_{17}\}$;
цикл $c_{13} = \{e_2,e_4,e_{12}\} \leftrightarrow \{v_1,v_3,v_5\}$;
цикл $c_{14} = \{e_2,e_5,e_{14},e_{31}\} \leftrightarrow \{v_1,v_3,v_6,v_{14}\}$;
цикл $c_{15} = \{e_2,e_6,e_{18},e_{40}\} \leftrightarrow \{v_1,v_3,v_8,v_{19}\}$;
цикл $c_{16} = \{e_2,e_7,e_{15}\} \leftrightarrow \{v_1,v_3,v_{15}\}$;
цикл $c_{17} = \{e_2,e_8,e_{16},e_{59}\} \leftrightarrow \{v_1,v_3,v_{17},v_{20}\}$;
цикл $c_{18} = \{e_2,e_8,e_{17},e_{60}\} \leftrightarrow \{v_1,v_3,v_{18},v_{20}\}$;
цикл $c_{19} = \{e_2,e_8,e_{18},e_{61}\} \leftrightarrow \{v_1,v_3,v_{19},v_{20}\}$;
цикл $c_{20} = \{e_3,e_4,e_{19}\} \leftrightarrow \{v_1,v_4,v_5\}$;
цикл $c_{21} = \{e_3,e_5,e_{21},e_{31}\} \leftrightarrow \{v_1,v_4,v_6,v_{14}\}$;
цикл $c_{22} = \{e_3,e_7,e_{22}\} \leftrightarrow \{v_1,v_4,v_{15}\}$;
цикл $c_{23} = \{e_3,e_8,e_{24}\} \leftrightarrow \{v_1,v_4,v_{20}\}$;
цикл $c_{24} = \{e_4,e_5,e_{26},e_{31}\} \leftrightarrow \{v_1,v_5,v_6,v_{14}\}$;
цикл $c_{25} = \{e_4,e_7,e_{27},e_{57}\} \leftrightarrow \{v_1,v_5,v_{15},v_{17}\}$;
цикл $c_{26} = \{e_4,e_8,e_{27},e_{59}\} \leftrightarrow \{v_1,v_5,v_{17},v_{20}\}$;
цикл $c_{27} = \{e_4,e_8,e_{28},e_{60}\} \leftrightarrow \{v_1,v_5,v_{18},v_{20}\}$;
цикл $c_{28} = \{e_1,e_5,e_{10},e_{30},e_{42}\} \leftrightarrow \{v_1,v_2,v_6,v_9,v_{13}\}$;
цикл $c_{29} = \{e_5, e_6,e_{29},e_{39}\} \leftrightarrow \{v_1,v_6,v_8,v_{12}\}$;
цикл $c_{30} = \{e_5,e_7,e_{32}\} \leftrightarrow \{v_1,v_6,v_{15}\}$;

цикл $c_{31} = \{e_6, e_8, e_{40}, e_{61}\} \leftrightarrow \{v_1, v_8, v_{19}, v_{20}\}$;
цикл $c_{32} = \{e_7, e_8, e_{57}, e_{59}\} \leftrightarrow \{v_1, v_{15}, v_{17}, v_{20}\}$;
цикл $c_{33} = \{e_9, e_{10}, e_{34}, e_{41}\} \leftrightarrow \{v_2, v_7, v_9, v_{11}\}$;
цикл $c_{34} = \{e_9, e_{10}, e_{36}, e_{42}\} \leftrightarrow \{v_2, v_7, v_9, v_{13}\}$;
цикл $c_{35} = \{e_9, e_{10}, e_{37}, e_{44}\} \leftrightarrow \{v_2, v_7, v_9, v_{16}\}$;
цикл $c_{36} = \{e_9, e_{11}, e_{38}\} \leftrightarrow \{v_2, v_7, v_{18}\}$;
цикл $c_{37} = \{e_{10}, e_{11}, e_{45}, e_{58}\} \leftrightarrow \{v_2, v_9, v_{17}, v_{18}\}$;
цикл $c_{38} = \{e_{12}, e_{13}, e_{25}\} \leftrightarrow \{v_3, v_5, v_{11}\}$;
цикл $c_{39} = \{e_{12}, e_{14}, e_{26}\} \leftrightarrow \{v_3, v_5, v_{14}\}$;
цикл $c_{40} = \{e_{12}, e_{15}, e_{19}, e_{22}\} \leftrightarrow \{v_3, v_4, v_5, v_{15}\}$;
цикл $c_{41} = \{e_{12}, e_{16}, e_{27}\} \leftrightarrow \{v_3, v_5, v_{17}\}$;
цикл $c_{42} = \{e_{12}, e_{17}, e_{28}\} \leftrightarrow \{v_3, v_5, v_{18}\}$;
цикл $c_{43} = \{e_{13}, e_{14}, e_{50}\} \leftrightarrow \{v_3, v_{11}, v_{14}\}$;
цикл $c_{44} = \{e_{13}, e_{15}, e_{20}, e_{22}\} \leftrightarrow \{v_3, v_4, v_{11}, v_{15}\}$;
цикл $c_{45} = \{e_{13}, e_{15}, e_{41}, e_{43}\} \leftrightarrow \{v_3, v_9, v_{11}, v_{15}\}$;
цикл $c_{46} = \{e_{13}, e_{16}, e_{51}\} \leftrightarrow \{v_3, v_{11}, v_{17}\}$;
цикл $c_{47} = \{e_{13}, e_{17}, e_{34}, e_{38}\} \leftrightarrow \{v_3, v_7, v_{11}, v_{18}\}$;
цикл $c_{48} = \{e_{13}, e_{17}, e_{49}, e_{54}\} \leftrightarrow \{v_3, v_{11}, v_{12}, v_{18}\}$;
цикл $c_{49} = \{e_{13}, e_{18}, e_{49}, e_{55}\} \leftrightarrow \{v_3, v_{11}, v_{12}, v_{19}\}$;
цикл $c_{50} = \{e_{14}, e_{15}, e_{21}, e_{22}\} \leftrightarrow \{v_3, v_4, v_{14}, v_{15}\}$;
цикл $c_{51} = \{e_{14}, e_{15}, e_{31}, e_{32}\} \leftrightarrow \{v_3, v_6, v_{14}, v_{15}\}$;
цикл $c_{52} = \{e_{14}, e_{16}, e_{21}, e_{23}\} \leftrightarrow \{v_3, v_4, v_{14}, v17\}$;
цикл $c_{53} = \{e_{14}, e_{16}, e_{52}, e_{53}\} \leftrightarrow \{v_3, v_{12}, v_{14}, v_{17}\}$;
цикл $c_{54} = \{e_{14}, e_{17}, e_{47}, e_{48}\} \leftrightarrow \{v_3, v_{10}, v_{14}, v_{18}\}$;
цикл $c_{55} = \{e_{14}, e_{17}, e_{52}, e_{54}\} \leftrightarrow \{v_3, v_{12}, v_{14}, v_{18}\}$;
цикл $c_{56} = \{e_{14}, e_{18}, e_{52}, e_{55}\} \leftrightarrow \{v_3, v_{12}, v_{14}, v_{19}\}$;
цикл $c_{57} = \{e_{15}, e_{16}, e_{57}\} \leftrightarrow \{v_3, v_{15}, v_{17}\}$;
цикл $c_{58} = \{e_{16}, e_{17}, e_{58}\} \leftrightarrow \{v_3, v_{17}, v_{18}\}$;
цикл $c_{59} = \{e_{16}, e_{18}, e_{53}, e_{55}\} \leftrightarrow \{v_3, v_{12}, v_{17}, v_{19}\}$;
цикл $c_{60} = \{e_{16}, e_{18}, e_{59}, e_{61}\} \leftrightarrow \{v_3, v_{17}, v_{19}, v_{20}\}$;
цикл $c_{61} = \{e_{17}, e_{18}, e_{54}, e_{55}\} \leftrightarrow \{v_3, v_{12}, v_{18}, v_{19}\}$;
цикл $c_{62} = \{e_{17}, e_{18}, e_{60}, e_{61}\} \leftrightarrow \{v_3, v_{18}, v_{19}, v_{20}\}$;
цикл $c_{63} = \{e_{19}, e_{20}, e_{25}\} \leftrightarrow \{v_4, v_5, v_{11}\}$;
цикл $c_{64} = \{e_{19}, e_{21}, e_{26}\} \leftrightarrow \{v_4, v_5, v_{14}\}$;
цикл $c_{65} = \{e_{19}, e_{23}, e_{27}\} \leftrightarrow \{v_4, v_5, v_{17}\}$;
цикл $c_{66} = \{e_{19}, e_{24}, e_{28}, e_{60}\} \leftrightarrow \{v_4, v_5, v_{18}, v_{20}\}$;
цикл $c_{67} = \{e_{20}, e_{21}, e_{50}\} \leftrightarrow \{v_4, v_{11}, v_{14}\}$;
цикл $c_{68} = \{e_{20}, e_{22}, e_{41}, e_{43}\} \leftrightarrow \{v_4, v_9, v_{11}, v_{15}\}$;
цикл $c_{69} = \{e_{20}, e_{23}, e_{51}\} \leftrightarrow \{v_4, v_{11}, v_{17}\}$;
цикл $c_{70} = \{e_{21}, e_{22}, e_{31}, e_{32}\} \leftrightarrow \{v_4, v_6, v_{14}, v_{15}\}$;
цикл $c_{71} = \{e_{21}, e_{23}, e_{52}, e_{53}\} \leftrightarrow \{v_4, v_{12}, v_{14}, v_{17}\}$;
цикл $c_{72} = \{e_{21}, e_{24}, e_{47}, e_{48}, e_{60}\} \leftrightarrow \{v_4, v_{10}, v_{14}, v_{18}, v_{20}\}$;
цикл $c_{73} = \{e_{22}, e_{23}, e_{57}\} \leftrightarrow \{v_4, v_{15}, v_{17}\}$;
цикл $c_{74} = \{e_{23}, e_{24}, e_{59}\} \leftrightarrow \{v_4, v_{17}, v_{20}\}$;
цикл $c_{75} = \{e_{14}, e_{18}, e_{21}, e_{24}, e_{61}\} \leftrightarrow \{v_3, v_4, v_{14}, v_{19}, v_{20}\}$;
цикл $c_{76} = \{e_{21}, e_{24}, e_{52}, e_{55}, e_{61}\} \leftrightarrow \{v_4, v_{12}, v_{14}, v_{19}, v_{20}\}$;
цикл $c_{77} = \{e_{25}, e_{26}, e_{50}\} \leftrightarrow \{v_5, v_{11}, v_{14}\}$;
цикл $c_{78} = \{e_{25}, e_{27}, e_{51}\} \leftrightarrow \{v_5, v_{11}, v_{17}\}$;
цикл $c_{79} = \{e_{25}, e_{28}, e_{34}, e_{38}\} \leftrightarrow \{v_5, v_7, v_{11}, v_{18}\}$;
цикл $c_{80} = \{e_{25}, e_{28}, e_{49}, e_{54}\} \leftrightarrow \{v_5, v_{11}, v_{12}, v_{18}\}$;
цикл $c_{81} = \{e_{26}, e_{27}, e_{52}, e_{53}\} \leftrightarrow \{v_5, v_{12}, v_{14}, v_{17}\}$;

цикл $c_{82} = \{e_{26}, e_{28}, e_{47}, e_{48}\} \leftrightarrow \{v_5, v_{10}, v_{14}, v_{18}\}$;
цикл $c_{83} = \{e_{26}, e_{28}, e_{52}, e_{54}\} \leftrightarrow \{v_5, v_{12}, v_{14}, v_{18}\}$;
цикл $c_{84} = \{e_{27}, e_{28}, e_{58}\} \leftrightarrow \{v_5, v_{17}, v_{18}\}$;
цикл $c_{85} = \{e_{29}, e_{30}, e_{35}, e_{36}\} \leftrightarrow \{v_6, v_7, v_{12}, v_{13}\}$;
цикл $c_{86} = \{e_{29}, e_{31}, e_{52}\} \leftrightarrow \{v_6, v_{12}, v_{14}\}$;
цикл $c_{87} = \{e_{29}, e_{32}, e_{53}, e_{57}\} \leftrightarrow \{v_6, v_{12}, v_{15}, v_{17}\}$;
цикл $c_{88} = \{e_{29}, e_{33}, e_{35}, e_{37}\} \leftrightarrow \{v_6, v_7, v_{12}, v_{16}\}$;
цикл $c_{89} = \{e_{30}, e_{31}, e_{46}, e_{47}\} \leftrightarrow \{v_6, v_{10}, v_{13}, v_{14}\}$;
цикл $c_{90} = \{e_{30}, e_{32}, e_{42}, e_{43}\} \leftrightarrow \{v_6, v_9, v_{13}, v_{15}\}$;
цикл $c_{91} = \{e_{30}, e_{33}, e_{36}, e_{37}\} \leftrightarrow \{v_6, v_7, v_{13}, v_{16}\}$;
цикл $c_{92} = \{e_{30}, e_{33}, e_{42}, e_{44}\} \leftrightarrow \{v_6, v_9, v_{13}, v_{16}\}$;
цикл $c_{93} = \{e_{32}, e_{33}, e_{56}\} \leftrightarrow \{v_6, v_{15}, v_{16}\}$;
цикл $c_{94} = \{e_{34}, e_{35}, e_{49}\} \leftrightarrow \{v_7, v_{11}, v_{12}\}$;
цикл $c_{95} = \{e_{34}, e_{36}, e_{41}, e_{42}\} \leftrightarrow \{v_7, v_9, v_{11}, v_{13}\}$;
цикл $c_{96} = \{e_{34}, e_{37}, e_{41}, e_{44}\} \leftrightarrow \{v_7, v_9, v_{11}, v_{16}\}$;
цикл $c_{97} = \{e_{34}, e_{38}, e_{51}, e_{58}\} \leftrightarrow \{v_7, v_{11}, v_{17}, v_{18}\}$;
цикл $c_{98} = \{e_{35}, e_{38}, e_{54}\} \leftrightarrow \{v_7, v_{12}, v_{18}\}$;
цикл $c_{99} = \{e_{36}, e_{37}, e_{42}, e_{44}\} \leftrightarrow \{v_7, v_9, v_{13}, v_{16}\}$;
цикл $c_{100} = \{e_{36}, e_{38}, e_{46}, e_{48}\} \leftrightarrow \{v_7, v_{10}, v_{13}, v_{18}\}$;
цикл $c_{101} = \{e_{15}, e_{17}, e_{37}, e_{38}, e_{56}\} \leftrightarrow \{v_3, v_7, v_{15}, v_{16}, v_{18}\}$;
цикл $c_{102} = \{e_{37}, e_{38}, e_{56}, e_{57}, e_{58}\} \leftrightarrow \{v_7, v_{15}, v_{16}, v_{17}, v_{18}\}$;
цикл $c_{103} = \{e_{39}, e_{40}, e_{55}\} \leftrightarrow \{v_8, v_{12}, v_{19}\}$;
цикл $c_{104} = \{e_{41}, e_{45}, e_{51}\} \leftrightarrow \{v_9, v_{11}, v_{17}\}$;
цикл $c_{105} = \{e_{42}, e_{45}, e_{46}, e_{48}, e_{58}\} \leftrightarrow \{v_9, v_{10}, v_{13}, v_{17}, v_{18}\}$;
цикл $c_{106} = \{e_{43}, e_{44}, e_{56}\} \leftrightarrow \{v_9, v_{15}, v_{16}\}$;
цикл $c_{107} = \{e_{43}, e_{45}, e_{57}\} \leftrightarrow \{v_9, v_{15}, v_{17}\}$;
цикл $c_{108} = \{e_{47}, e_{48}, e_{52}, e_{54}\} \leftrightarrow \{v_{10}, v_{12}, v_{14}, v_{18}\}$;
цикл $c_{109} = \{e_{49}, e_{50}, e_{52}\} \leftrightarrow \{v_{11}, v_{12}, v_{14}\}$;
цикл $c_{110} = \{e_{49}, e_{51}, e_{53}\} \leftrightarrow \{v_{11}, v_{12}, v_{17}\}$;
цикл $c_{111} = \{e_{41}, e_{42}, e_{46}, e_{47}, e_{50}\} \leftrightarrow \{v_9, v_{10}, v_{11}, v_{13}, v_{14}\}$;
цикл $c_{112} = \{e_{53}, e_{54}, e_{58}\} \leftrightarrow \{v_{12}, v_{17}, v_{18}\}$;
цикл $c_{113} = \{e_{53}, e_{55}, e_{59}, e_{61}\} \leftrightarrow \{v_{12}, v_{17}, v_{19}, v_{20}\}$;
цикл $c_{114} = \{e_{54}, e_{55}, e_{60}, e_{61}\} \leftrightarrow \{v_{12}, v_{18}, v_{19}, v_{20}\}$;
цикл $c_{115} = \{e_{13}, e_{15}, e_{34}, e_{37}, e_{56}\} \leftrightarrow \{v_3, v_7, v_{11}, v_{15}, v_{16}\}$;
цикл $c_{116} = \{e_{20}, e_{22}, e_{34}, e_{37}, e_{56}\} \leftrightarrow \{v_4, v_7, v_{11}, v_{15}, v_{16}\}$;
цикл $c_{117} = \{e_{58}, e_{59}, e_{60}\} \leftrightarrow \{v_{17}, v_{18}, v_{20}\}$.

Выделим подмножество изометрических циклов, описывающее плоскую конфигурацию, состоящую из 14 циклов. Так как суграф связан, не имеет точек сочленения и отсутствует разрыв вершин, то можно построить топологический рисунок суграфа.

Множество изометрических циклов пронумеровано согласно последовательности их расположения в подмножестве:

$c_1 = \{e_9, e_{11}, e_{38}\} \leftrightarrow \{v_7, v_{18}, v_2\}$;
$c_2 = \{e_{34}, e_{36}, e_{41}, e_{42}\} \leftrightarrow \{v_{11}, v_9, v_{13}, v_7\}$;
$c_3 = \{e_{12}, e_{14}, e_{26}\} \leftrightarrow \{v_5, v_{14}, v_3\}$;
$c_4 = \{e_2, e_8, e_{17}, e_{60}\} \leftrightarrow \{v_3, v_{18}, v_{20}, v_1\}$;
$c_5 = \{e_{25}, e_{26}, e_{50}\} \leftrightarrow \{v_{11}, v_{14}, v_5\}$;
$c_6 = \{e_{12}, e_{17}, e_{28}\} \leftrightarrow \{v_5, v_{18}, v_3\}$;
$c_7 = \{e_{32}, e_{33}, e_{56}\} \leftrightarrow \{v_{15}, v_{16}, v_6\}$;

$c_8 = \{e_{58}, e_{59}, e_{60}\} \leftrightarrow \{v_{18}, v_{20}, v_{17}\}$;
$c_9 = \{e_{21}, e_{22}, e_{31}, e_{32}\} \leftrightarrow \{v_{14}, v_6, v_{15}, v_4\}$;
$c_{10} = \{e_6, e_8, e_{40}, e_{61}\} \leftrightarrow \{v_8, v_{19}, v_{20}, v_1\}$;
$c_{11} = \{e_2, e_5, e_{14}, e_{31}\} \leftrightarrow \{v_3, v_{14}, v_6, v_1\}$;
$c_{12} = \{e_{25}, e_{28}, e_{34}, e_{38}\} \leftrightarrow \{v_{11}, v_7, v_{18}, v_5\}$;
$c_{13} = \{e_5, e_6, e_{29}, e_{39}\} \leftrightarrow \{v_6, v_{12}, v_8, v_1\}$;
$c_{14} = \{e_{41}, e_{42}, e_{46}, e_{47}, e_{50}\} \leftrightarrow \{v_{11}, v_{14}, v_{10}, v_{13}, v_9\}$.

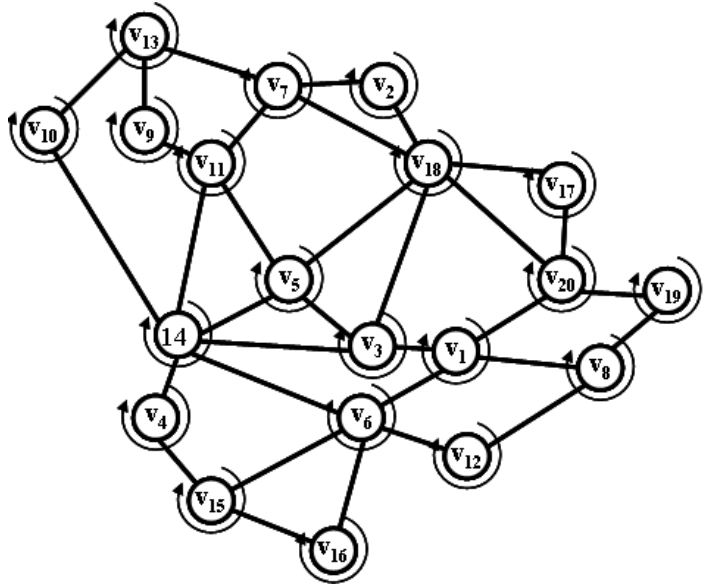

Рис. 13.1. Топологический рисунок суграфа $G_{31}$, состоящий из 14 циклов.

Для построения топологического рисунка плоского суграфа выделим следующее подмножество, состоящее из 15 изометрических циклов:

$c_1 = \{e_{43}, e_{44}, e_{56}\} \leftrightarrow \{v_{15}, v_{16}, v_9\}$;
$c_2 = \{e_{53}, e_{54}, e_{58}\} \leftrightarrow \{v_{17}, v_{18}, v_{12}\}$;
$c_3 = \{e_{34}, e_{36}, e_{41}, e_{42}\} \leftrightarrow \{v_{11}, v_9, v_{13}, v_7\}$;
$c_4 = \{e_9, e_{11}, e_{38}\} \leftrightarrow \{v_7, v_{18}, v_2\}$;
$c_5 = \{e_{30}, e_{33}, e_{36}, e_{37}\} \leftrightarrow \{v_{13}, v_7, v_{16}, v_6\}$;
$c_6 = \{e_2, e_8, e_{18}, e_{61}\} \leftrightarrow \{v_3, v_{19}, v_{20}, v_1\}$;
$c_7 = \{e_{20}, e_{22}, e_{41}, e_{43}\} \leftrightarrow \{v_{11}, v_9, v_{15}, v_4\}$;
$c_8 = \{e_{37}, e_{38}, e_{56}, e_{57}, e_{58}\} \leftrightarrow \{v_{16}, v_{15}, v_{17}, v_{18}, v_7\}$;
$c_9 = \{e_{12}, e_{14}, e_{26}\} \leftrightarrow \{v_5, v_{14}, v_3\}$;
$c_{10} = \{e_1, e_3, e_9, e_{20}, e_{34}\} \leftrightarrow \{v_2, v_7, v_{11}, v_4, v_1\}$;
$c_{11} = \{e_{14}, e_{17}, e_{47}, e_{48}\} \leftrightarrow \{v_{14}, v_{10}, v_{18}, v_3\}$;
$c_{12} = \{e_{30}, e_{33}, e_{42}, e_{44}\} \leftrightarrow \{v_{13}, v_9, v_{16}, v_6\}$;
$c_{13} = \{e_{22}, e_{23}, e_{57}\} \leftrightarrow \{v_{15}, v_{17}, v_4\}$;
$c_{14} = \{e_1, e_2, e_{11}, e_{17}\} \leftrightarrow \{v_2, v_{18}, v_3, v_1\}$;
$c_{15} = \{e_6, e_8, e_{40}, e_{61}\} \leftrightarrow \{v_8, v_{19}, v_{20}, v_1\}$.

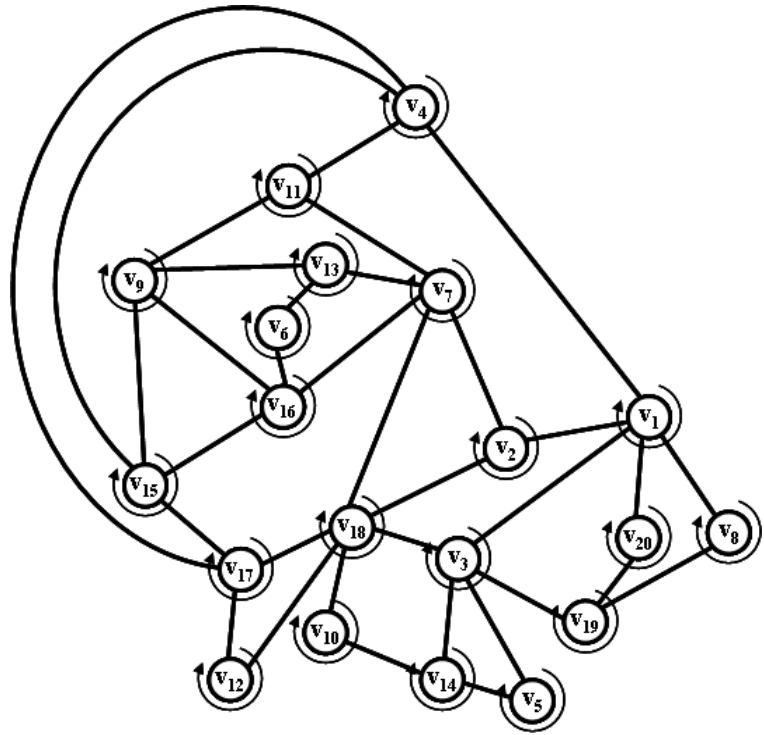

Рис. 13.2. Топологический рисунок суграфа $G_{32}$, состоящий из 15 циклов.

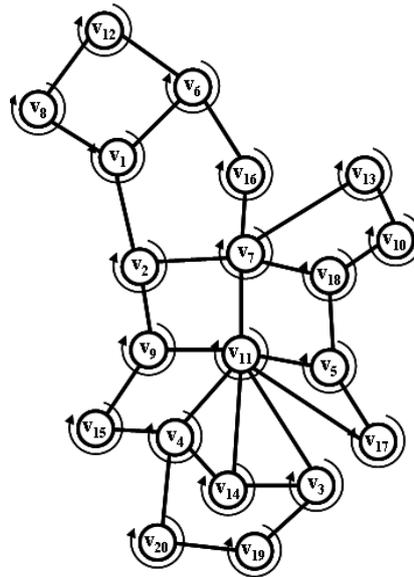

Рис. 13.3. Топологический рисунок плоского суграфа $G_{33}$, состоящий из 8 циклов.

Выделим следующее подмножество, состоящее из 8 изометрических циклов:

$c_1 = \{e_{14}, e_{18}, e_{21}, e_{24}, e_{61}\} \leftrightarrow \{v_{14}, v_4, v_{20}, v_{19}, v_3\}$;
$c_2 = \{e_{13}, e_{14}, e_{50}\} \leftrightarrow \{v_{11}, v_{14}, v_3\}$;
$c_3 = \{e_9, e_{10}, e_{34}, e_{41}\} \leftrightarrow \{v_7, v_{11}, v_9, v_2\}$;
$c_4 = \{e_{25}, e_{27}, e_{51}\} \leftrightarrow \{v_{11}, v_{17}, v_5\}$;
$c_5 = \{e_1, e_5, e_9, e_{33}, e_{37}\} \leftrightarrow \{v_2, v_7, v_{16}, v_6, v_1\}$;
$c_6 = \{e_{20}, e_{22}, e_{41}, e_{43}\} \leftrightarrow \{v_{11}, v_9, v_{15}, v_4\}$;
$c_7 = \{e_{36}, e_{38}, e_{46}, e_{48}\} \leftrightarrow \{v_{13}, v_{10}, v_{18}, v_7\}$;
$c_8 = \{e_{20}, e_{21}, e_{50}\} \leftrightarrow \{v_{11}, v_{14}, v_4\}$;
$c_9 = \{e_{25}, e_{28}, e_{34}, e_{38}\} \leftrightarrow \{v_{11}, v_7, v_{18}, v_5\}$;
$c_{10} = \{e_5, e_6, e_{29}, e_{39}\} \leftrightarrow \{v_6, v_{12}, v_8, v_1\}$.

Выделим следующее подмножество, состоящее из 17 изометрических циклов:

$c_1 = \{e_{27}, e_{28}, e_{58}\} \leftrightarrow \{v_{17}, v_{18}, v_5\}$

$c_2 = \{e_2, e_6, e_{18}, e_{40}\} \leftrightarrow \{v_3, v_{19}, v_8, v_1\}$

$c_3 = \{e_{19}, e_{24}, e_{28}, e_{60}\} \leftrightarrow \{v_5, v_{18}, v_{20}, v_4\}$

$c_4 = \{e_{20}, e_{21}, e_{50}\} \leftrightarrow \{v_{11}, v_{14}, v_4\}$

$c_5 = \{e_1, e_5, e_9, e_{29}, e_{35}\} \leftrightarrow \{v_2, v_7, v_{12}, v_6, v_1\}$

$c_6 = \{e_{34}, e_{37}, e_{41}, e_{44}\} \leftrightarrow \{v_{11}, v_9, v_{16}, v_7\}$

$c_7 = \{e_{15}, e_{16}, e_{57}\} \leftrightarrow \{v_{15}, v_{17}, v_3\}$

$c_8 = \{e_2, e_7, e_{15}\} \leftrightarrow \{v_3, v_{15}, v_1\}$

$c_9 = \{e_4, e_7, e_{27}, e_{57}\} \leftrightarrow \{v_5, v_{17}, v_{15}, v_1\}$

$c_{10} = \{e_3, e_4, e_{19}\} \leftrightarrow \{v_4, v_5, v_1\}$

$c_{11} = \{e_{16}, e_{17}, e_{58}\} \leftrightarrow \{v_{17}, v_{18}, v_3\}$

$c_{12} = \{e_{17}, e_{18}, e_{54}, e_{55}\} \leftrightarrow \{v_{18}, v_{12}, v_{19}, v_3\}$

$c_{13} = \{e_1, e_3, e_9, e_{20}, e_{34}\} \leftrightarrow \{v_2, v_7, v_{11}, v_4, v_1\}$

$c_{14} = \{e_{36}, e_{37}, e_{42}, e_{44}\} \leftrightarrow \{v_{13}, v_9, v_{16}, v_7\}$

$c_{15} = \{e_5, e_6, e_{29}, e_{39}\} \leftrightarrow \{v_6, v_{12}, v_8, v_1\}$

$c_{16} = \{e_{41}, e_{42}, e_{46}, e_{47}, e_{50}\} \leftrightarrow \{v_{11}, v_{14}, v_{10}, v_{13}, v_9\}$

$c_{17} = \{e_{39}, e_{40}, e_{55}\} \leftrightarrow \{v_{12}, v_{19}, v_8\}$

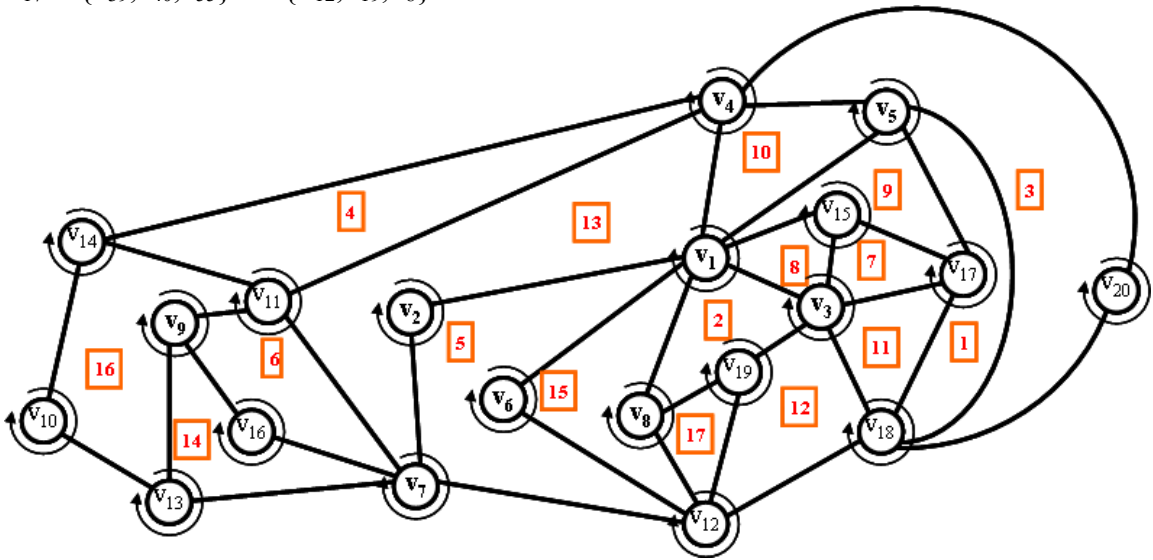

Рис. 13.4. Топологический рисунок суграфа $G_{34}$, состоящий из 17 циклов.

### 13.2. Топологический рисунок графа и векторная алгебра пересечений

Представим обод выделенного плоского суграфа в виде координатно-базисной системы векторов [27-29].

Например, обод плоского суграфа $G_{31}$, представленного топологическим рисунком 1.7, имеет вид:

$c_0 = \{e_{36}, e_9, e_{11}, e_{58}, e_{59}, e_{61}, e_{40}, e_{39}, e_{29}, e_{33}, e_{56}, e_{22}, e_{21}, e_{47}, e_{46}\} \leftrightarrow$
$\leftrightarrow <v_{13}, v_7, v_2, v_{18}, v_{17}, v_{20}, v_{19}, v_8, v_{12}, v_6, v_{16}, v_{15}, v_4, v_{14}, v_{10}>$.

Из оставшихся соединений построим рисунок соединений с вершинами координатно-базисной системы, расположенной на ободе плоского суграфа. Рассмотрим проекцию соединений на координатно-базисную систему:

пр$(v_{12}, v_{17}) =$ пр$\{e_{53}\} = \{e_{59}, e_{61}, e_{40}, e_{39}\} \wedge \{e_{29}, e_{33}, e_{56}, e_{22}, e_{21}, e_{47}, e_{46}, e_{36}, e_9, e_{11}, e_{58}\}$;

пр$(v_{12}, v_{18}) =$ пр$\{e_{54}\} = \{e_{58}, e_{59}, e_{61}, e_{40}, e_{39}\} \wedge \{e_{29}, e_{33}, e_{56}, e_{22}, e_{21}, e_{47}, e_{46}, e_{36}, e_9, e_{11}\}$;

пр($v_{12},v_7$) = пр{$e_{35}$} = {$e_9,e_{11},e_{58},e_{59},e_{61},e_{40},e_{39}$} ∧ {$e_{29},e_{33},e_{56},e_{22},e_{21},e_{47},e_{46},e_{36}$};
пр($v_{12},v_{14}$) = пр{$e_{52}$} = {$e_{47},e_{46},e_{36},e_9,e_{11},e_{58},e_{59},e_{61},e_{40},e_{39}$} ∧ {$e_{29},e_{33},e_{56},e_{22},e_{21}$};
пр($v_6,v_{13}$) = пр{$e_{30}$} = {$e_{36},e_9,e_{11},e_{58},e_{59},e_{61},e_{40},e_{39},e_{29}$} ∧ {$e_{33},e_{56},e_{22},e_{21},e_{47},e_{46}$};
пр($v_7,v_{16}$) = пр{$e_{37}$} = {$e_9,e_{11},e_{58},e_{59},e_{61},e_{40},e_{39},e_{29},e_{33}$} ∧ {$e_{56},e_{22},e_{21},e_{47},e_{46},e_{36}$};
пр($v_{10},v_{18}$) = пр{$e_{48}$} = {$e_{46},e_{36},e_9,e_{11}$} ∧ {$e_{58},e_{59},e_{61},e_{40},e_{39},e_{29},e_{33},e_{56},e_{22},e_{21},e_{47}$};
пр($v_4,v_{17}$) = пр{$e_{23}$} = {$e_{59},e_{61},e_{40},e_{39},e_{29},e_{33},e_{56},e_{22}$} ∧ {$e_{21},e_{47},e_{46},e_{36},e_9,e_{11},e_{58}$};
пр($v_4,v_{20}$) = пр{$e_{24}$} = {$e_{61},e_{40},e_{39},e_{29},e_{33},e_{56}, e_{22}$} ∧ { $e_{21},e_{47},e_{46},e_{36},e_9,e_{11},e_{58},e_{59}$}.
пр($v_{15},v_{17}$) = пр{$e_{57}$} = {$e_{59},e_{61},e_{40},e_{39},e_{29},e_{33},e_{56}$} ∧ {$e_{22},e_{21},e_{47},e_{46},e_{36},e_9,e_{11},e_{58}$}.

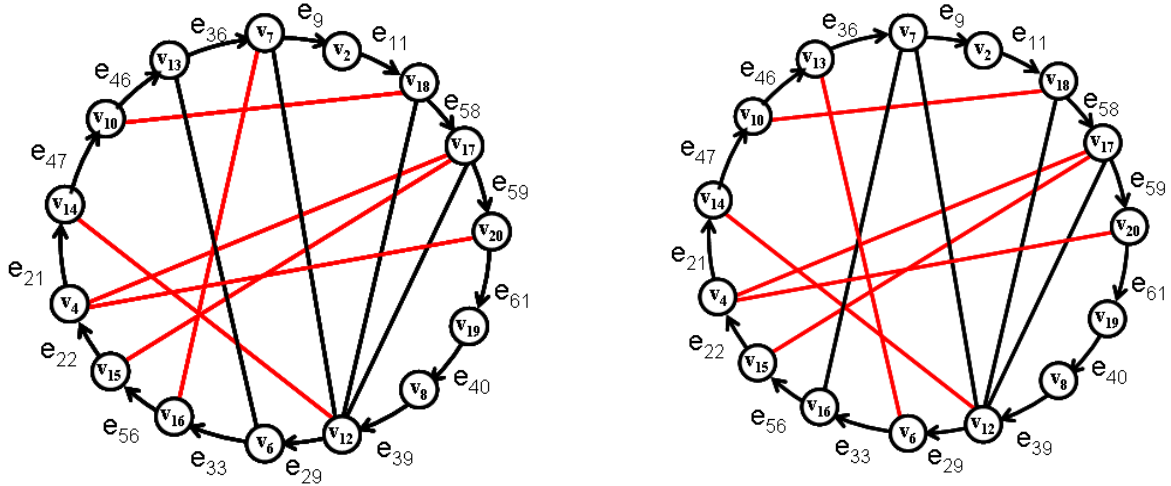

Рис. 13.5. Способы выделения непересекающихся соединений.

Определим пересечения соединений, рассматривая пересечение соединений как пересеченинhe их проекций на КБС. Рассмотрим пересечение соединений $e_{35}$ и $e_{23}$.

{$e_{35}$} ∩ {$e_{23}$} = пр{$e_{35}$} ∩ пр{$e_{23}$} =
{$e_9,e_{11},e_{58},e_{59},e_{61},e_{40},e_{39}$} ∩ {$e_{59},e_{61},e_{40},e_{39},e_{29},e_{33},e_{56},e_{22}$} = {$e_{59},e_{61},e_{40},e_{39}$};
{$e_9,e_{11},e_{58},e_{59},e_{61},e_{40},e_{39}$} ∩ {$e_{21},e_{47},e_{46},e_{36},e_9,e_{11},e_{58}$}={$e_9,e_{11},e_{58}$};
{$e_{29},e_{33},e_{56},e_{22},e_{21},e_{47},e_{46},e_{36}$} ∩ {$e_{59},e_{61},e_{40},e_{39},e_{29},e_{33},e_{56},e_{22}$}={$e_{29},e_{33},e_{56},e_{22}$};
{$e_{29},e_{33},e_{56},e_{22},e_{21},e_{47},e_{46},e_{36}$} ∩ {$e_{21},e_{47},e_{46},e_{36},e_9,e_{11},e_{58}$}={$e_{21},e_{47},e_{46},e_{36}$};
соединения {$e_{35}$} и{$e_{23}$} пересекаются, так как пересекаются их проекции.

Рассмотрим пересечение соединений $e_{35}$ и $e_{23}$.

{$e_{35}$} ∩ {$e_{30}$} = пр{$e_{35}$} ∩ пр{$e_{30}$} =
{$e_9,e_{11},e_{58},e_{59},e_{61},e_{40},e_{39}$} ∩ {$e_{36},e_9,e_{11},e_{58},e_{59},e_{51},e_{40},e_{39},e_{29}$}=
{$e_9,e_{11},e_{58},e_{59},e_{61},e_{40},e_{39}$} ⊆ {$e_{36},e_9,e_{11},e_{58},e_{59},e_{51},e_{40},e_{39},e_{29}$}=∅ ;
{$e_9,e_{11},e_{58},e_{59},e_{61},e_{40},e_{39}$} ∩ {$e_{33},e_{56},e_{22},e_{21},e_{47},e_{46}$}=∅ ;
{$e_{29},e_{33},e_{56},e_{22},e_{21},e_{47},e_{46},e_{36}$} ∩ {$e_{36},e_9,e_{11},e_{58},e_{59},e_{51},e_{40},e_{39},e_{29}$}=∅ ;
{$e_{29},e_{33},e_{56},e_{22},e_{21},e_{47},e_{46},e_{36}$} ∩ {$e_{33},e_{56},e_{22},e_{21},e_{47},e_{46}$} =
{$e_{33},e_{56},e_{22},e_{21},e_{47},e_{46}$} ⊆ {$e_{29},e_{33},e_{56},e_{22},e_{21},e_{47},e_{46},e_{36}$}=∅ .

Проверим все 45 пар соединений на пересечение и исключим наиболее пересекающиеся соединения. В результате остались четыре не пересекающихся соединения {$e_{53}$},{$e_{54}$},{$e_{35}$} и {$e_{30}$}, выделенные черным цветом (рис. 13.5, слева). Выбранные ребра и их проекции образуют дополнительное подмножество простых циклов для включения в топологический рисунок.

Система дополнительных простых циклов $G_{21}$ имеет вид:
$c_{d1}$ = {$e_{53},e_{59},e_{61},e_{40},e_{39}$} ↔ {$v_8,v_{12},v_{17},v_{19},v_{20}$};

$c_{d2} = \{e_{53}, e_{54}, e_{58}\} \leftrightarrow \{v_{12}, v_{17}, v_{18}\}$;
$c_{d3} = \{e_{35}, e_{54}, e_9, e_{11}\} \leftrightarrow \{v_2, v_7, v_{12}, v_{18}\}$;
$c_{d4} = \{e_{30}, e_{35}, e_{29}, e_{36}\} \leftrightarrow \{v_6, v_7, v_{12}, v_{13}\}$.

На рис. 13.5 справа, представлен альтернативный вариант выбора соединений.

Рис. 13.6. Топологический рисунок плоского суграфа $G_{31}$ с 36 рёбрами.

Определим обод плоского суграфа $G_{31}$

$c_0 = c_{КБС} \oplus c_{d1} \oplus c_{d2} \oplus c_{d3} \oplus c_{d4} = \{e_{30}, e_{33}, e_{56}, e_{22}, e_{21}, e_{47}, e_{46}\} \leftrightarrow$
$\leftrightarrow \{v_{10}, v_{13}, v_6, v_{16}, v_{15}, v_4, v_{14}\}$.

Рассмотрим построение дополнительных соединений для плоского суграфа $G_{32}$.

Рис. 13.7. Построение дополнительных соединений относительно обода плоского суграфа $G_{32}$.

Определим простые циклы, порождённые построением соединений относительно координатно-базисной системы.

$c_0 = c_{КБС} = \{e_3, e_6, e_{40}, e_{18}, e_{12}, e_{26}, e_{27}, e_{47}, e_{48}, e_{54}, e_{53}, e_{23}\} \leftrightarrow$
$\leftrightarrow \{v_4, v_1, v_8, v_{19}, v_3, v_5, v_{14}, v_{10}, v_{18}, v_{12}, v_{17}\}$;

Система дополнительных циклов имеет вид:

$c_{d1} = \{e_6, e_{40}, e_{18}, e_{12}, e_4\} \leftrightarrow \{v_1, v_8, v_{19}, v_3, v_5\}$;
$c_{d2} = \{e_3, e_4, e_{19}\} \leftrightarrow \{v_1, v_4, v_5\}$;
$c_{d3} = \{e_{19}, e_{21}, e_{26}\} \leftrightarrow \{v_4, v_5, v_{14}\}$;
$c_{d4} = \{e_{52}, e_{54}, e_{47}, e_{48}\} \leftrightarrow \{v_4, v_{14}, v_{12}, v_{17}\}$.

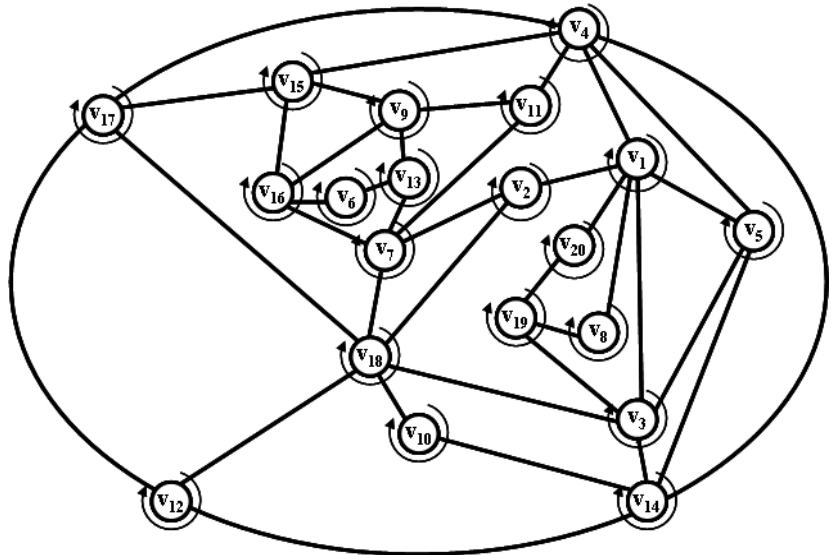

Рис. 13.8. Топологический рисуглк суграфа $G_{32}$, с дополнительными циклами.

Обод топологического рисунка суграфа $G_{32}$ с дополнительными циклами:

$c_0 = c_{КБС} \oplus c_{d1} \oplus c_{d2} \oplus c_{d3} \oplus c_{d4}$
$c_0 = c_{КБС} = \{e_{21}, e_{52}, e_{53}, e_{23}\} \leftrightarrow \{v_4, v_{14}, v_{12}, v_{17}\}$.

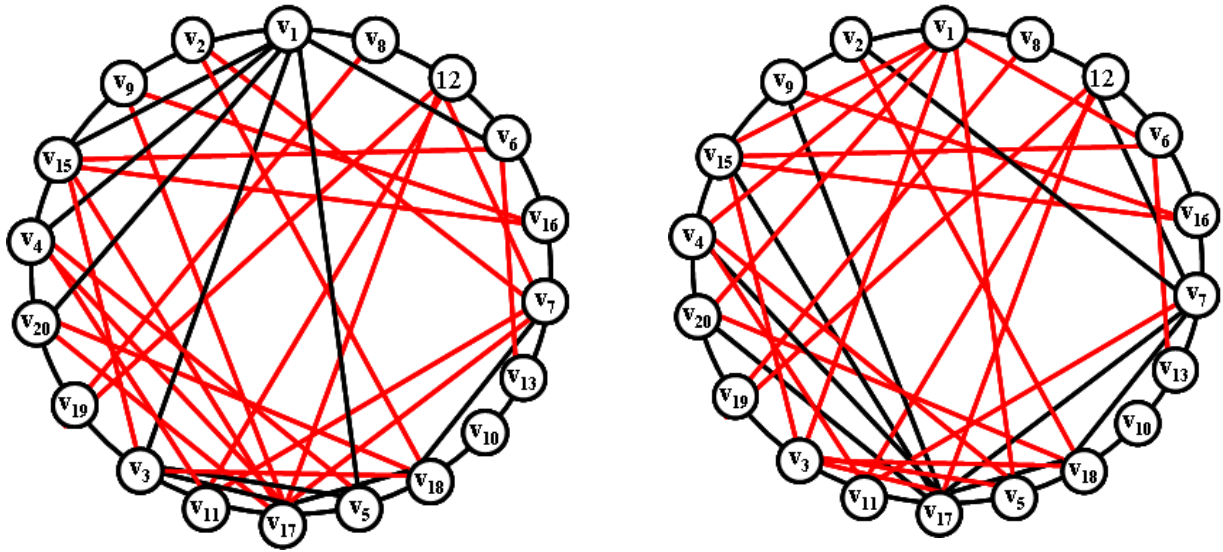

Рис. 13.9. Построение дополнительных соединений относительно обода плоского суграфа $G_{33}$.

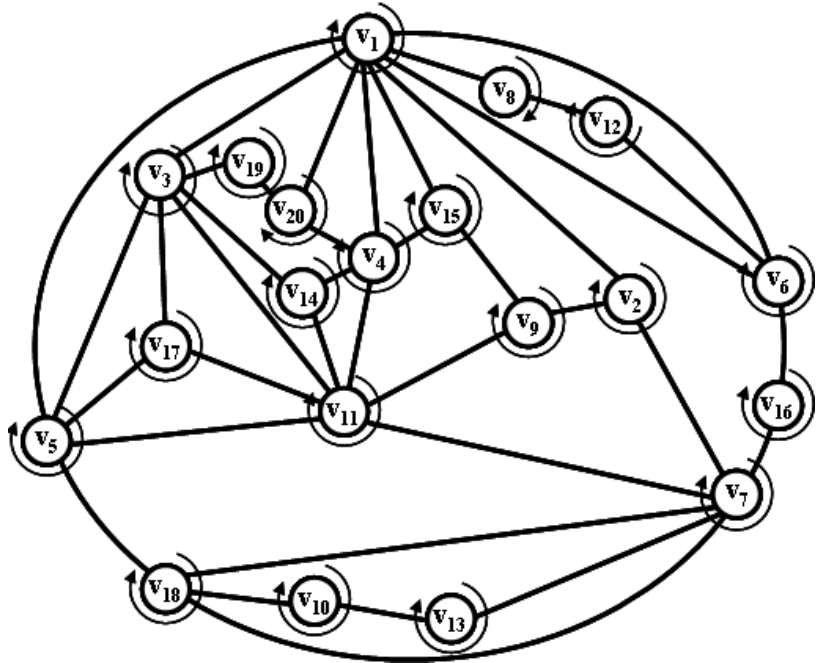

Рис. 13.10. Топологический рисунок суграфа $G_{33}$ с дополнительными циклами.

На рисунке 13.11 представлен процесс подключения новых ребер (обозначены пунктирными линиями) к ободу суграфа $G_{24}$. В результате образуются два простых цикла.

$c_{d1} = \{e_{35}, e_{38}, e_{54}\} \leftrightarrow \{v_7, v_{12}, v_{18}\}$;

$c_{d2} = \{e_{36}, e_{46}, e_{48}, e_{54}\} \leftrightarrow \{v_7, v_{10}, v_{13}, v_{18}\}$.

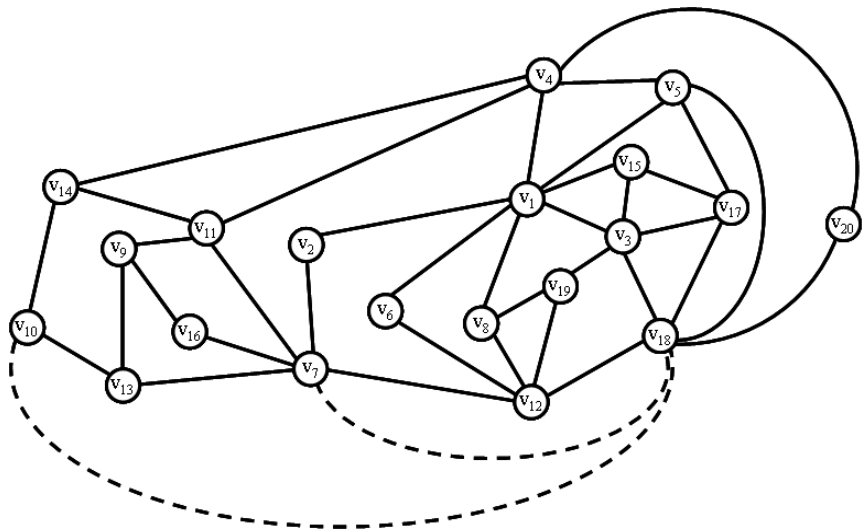

Рис. 13.11. Формирование простых циклов в суграфе $G_{34}$.

Для определения топологического рисунка плоского суграфа построим таблицу количества ребер.

Таблица количества ребер в плоском суграфе.

| суграф | $G_{21}$ | $G_{22}$ | $G_{23}$ | $G_{24}$ |
|---|---|---|---|---|
| количество ребер | 36 | 38 | 38 | 38 |

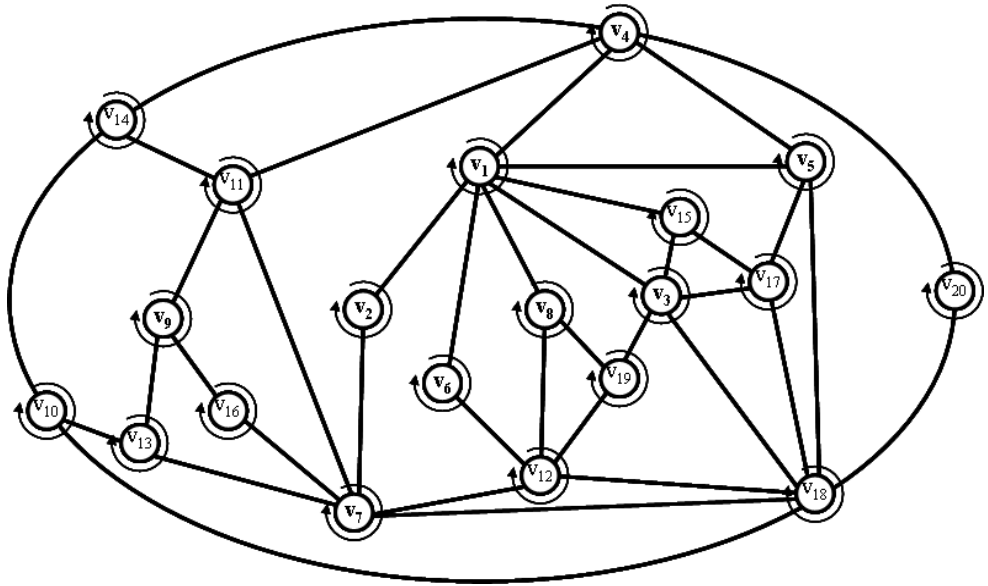

Рис. 13.12. Топологический рисунок суграфа $G_{34}$ с дополнительными циклами.

**Комментарии**

В главе рассматривается третий этап построения топологического рисунка плоского суграфа. Построение осуществляется на основе координатно-базисной системы векторов, образованной ободом плоского суграфа, вычисленного на втором этапе. Показано, что построение дополнительной системы циклов возможно несколькими способами.

## Глава 14. ГАМИЛЬТОНОВ ЦИКЛ ГРАФА
### 14.1. Удаление циклов и граф циклов

Представим, что выделена плоская часть графа с нулевым значением функционала Маклейна. Выделение плоской части является первым этапом решения задачи выделения гамильтонова цикла. На втором этапе построения гамильтонова цикла мы должны удалять цикл с одновременным удалением ребра, соблюдая правило: при удалении цикла удаляется одно и только одно ребро (для выполнения закона Эйлера).

$$k_c - (m - n + 1) = 0, \qquad (14.1)$$

$k_c$ – количество независимых циклов.

Поставим в соответствие каждому циклу графа G(V,E) вершину графа циклов H, а соприкасающиеся ребра циклов в графе G будем считать ребрами графа H.

**Определение 14.1.** Связанный неориентированный несепарабельный граф, вершины которого соответствуют простым циклам топологического рисунка плоского суграфа G, а ребра соответствуют ребрам плоского суграфа G, за исключением ребер обода, называется *графом циклов* H.

В качестве примера, для описания второго этапа построения гамильтонова цикла рассмотрим граф, представленный на рисунке 14.1. Важную роль в определении гамильтонова цикла путем удаления циклов в плоском графе является граф циклов H. На рисунке 14.1 он представлен пунктирными ребрами, соединяющими вершины, находящиеся в центре цикла.

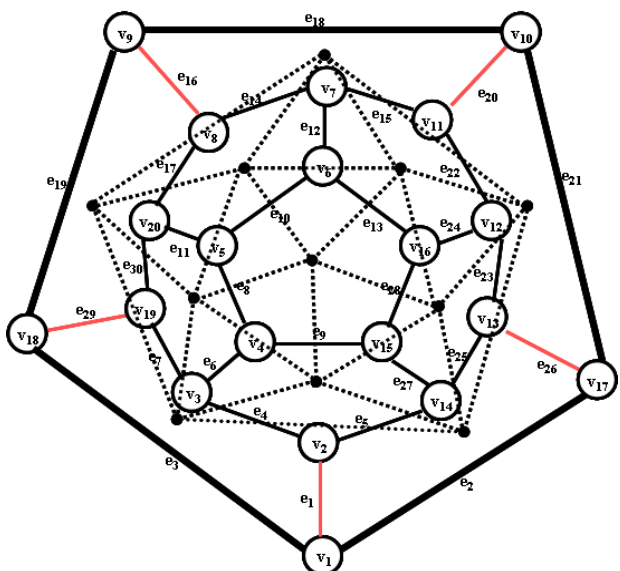

Рис. 14.1. Плоский граф G и его граф циклов H.

Множество изометрических циклов графа:

$c_1 = \{v_{12}, v_{11}, v_7, v_6, v_{16}, v_{12}\} \rightarrow \{e_{22}, e_{15}, e_{12}, e_{13}, e_{24}\}$;
$c_2 = \{v_{11}, v_{10}, v_9, v_8, v_7, v_{11}\} \rightarrow \{e_{20}, e_{18}, e_{16}, e_{14}, e_{15}\}$;
$c_3 = \{v_2, v_3, v_{19}, v_{18}, v_1, v_2\} \rightarrow \{e_4, e_7, e_{29}, e_3, e_1\}$;
$c_4 = \{v_3, v_4, v_5, v_{20}, v_{19}, v_3\} \rightarrow \{e_6, e_8, e_{11}, e_{30}, e_7\}$;
$c_5 = \{v_5, v_6, v_7, v_8, v_{20}, v_5\} \rightarrow \{e_{10}, e_{12}, e_{14}, e_{17}, e_{11}\}$;

$c_6 = \{v_{10},v_{11},v_{12},v_{13},v_{17},v_{10}\} \to \{e_{20},e_{22},e_{23},e_{26},e_{21}\};$
$c_7 = \{v_2,v_{14},v_{15},v_4,v_3,v_2\} \to \{e_5,e_{27},e_9,e_6,e_4\};$
$c_8 = \{v_{15},v_{16},v_6,v_5,v_4,v_{15}\} \to \{e_{28},e_{13},e_{10},e_8,e_9\};$
$c_9 = \{v_{14},v_{13},v_{12},v_{16},v_{15},v_{14}\} \to \{e_{25},e_{23},e_{24},e_{28},e_{27}\};$
$c_{10} = \{v_2,v_1,v_{17},v_{13},v_{14},v_2\} \to \{e_1,e_2,e_{26},e_{25},e_5\};$
$c_{11} = \{v_9,v_{18},v_{19},v_{20},v_8,v_9\} \to \{e_{19},e_{29},e_{30},e_{17},e_{16}\};$
$c_0 = c_{12} = \{v_{18},v_9,v_{10},v_{17},v_1,v_{18}\} \to \{e_{19},e_{18},e_{21},e_2,e_3\}.$

Вектор количества циклов, проходящих по рёбрам:

$P_e = <2,2,2,2,2,2,2,2,2,2,2,2,2,2,2,2,2,2,2,2,2,2,2,2,2,2,2,2,2,2>.$

Вектор количества циклов, проходящих по вершинам:

$P_v = <3,3,3,3,3,3,3,3,3,3,3,3,3,3,3,3,3,3,3,3>.$

Если выбрать в качестве обода цикл $c_{12}$, то вектора будут иметь вид –

Вектор количества циклов проходящих по рёбрам:

$P_e = <2,1,1,2,2,2,2,2,2,2,2,2,2,2,2,2,2,1,1,2,1,2,2,2,2,2,2,2,2,2>.$

Вектор количества циклов проходящих по вершинам:

$P_v = <2,3,3,3,3,3,3,3,2,2,3,3,3,3,3,3,2,2,3,3>.$

Граф циклов H для системы изометрических циклов имеет вид:

Рис. 14.2. Граф циклов H для подмножества циклов.

На рисунке 14.2 цифры вершин графа циклов соответствуют номерам изометрических циклов. В качестве обода выбран цикл $c_{12}$. При таком выборе обода графа, цикл $c_{12}$ не будем указывать в графе циклов H.

Количество рёбер в ободе всегда можно определить как количество единиц в векторе $P_e$.

**Определение 14.2**. Два простых цикла системы, имеющие одно общее ребро и отличающиеся друг от друга хотя бы одной вершиной, называются *соприкасающимися циклами*.

Все циклы, соприкасающиеся (имеющие общие рёбра) с ободом, должны быть помечены. Будем удалять только помеченные циклы, имеющие максимальную валентность в графе циклов H. Пометим вершины 2,3,6,10,11 в графе циклов H. В данном случае, все помеченные вершины имеют валентность равную 4. Поэтому для удаления можно выбрать любой цикл. Для удаления выбираем цикл $c_2$ и ребро $e_{18}$ (рис. 14.3). Цикл $c_2$ будем считать опорным циклом.

После удаления цикла $c_2$ множество изометрических циклов имеет вид:

$c_1 = \{v_{12},v_{11},v_7,v_6,v_{16},v_{12}\} \to \{e_{22},e_{15},e_{12},e_{13},e_{24}\}$;
$c_3 = \{v_2,v_3,v_{19},v_{18},v_1,v_2\} \to \{e_4,e_7,e_{29},e_3,e_1\}$;
$c_4 = \{v_3,v_4,v_5,v_{20},v_{19},v_3\} \to \{e_6,e_8,e_{11},e_{30},e_7\}$;
$c_5 = \{v_5,v_6,v_7,v_8,v_{20},v_5\} \to \{e_{10},e_{12},e_{14},e_{17},e_{11}\}$;
$c_6 = \{v_{10},v_{11},v_{12},v_{13},v_{17},v_{10}\} \to \{e_{20},e_{22},e_{23},e_{26},e_{21}\}$;
$c_7 = \{v_2,v_{14},v_{15},v_4,v_3,v_2\} \to \{e_5,e_{27},e_9,e_6,e_4\}$;
$c_8 = \{v_{15},v_{16},v_6,v_5,v_4,v_{15}\} \to \{e_{28},e_{13},e_{10},e_8,e_9\}$;
$c_9 = \{v_{14},v_{13},v_{12},v_{16},v_{15},v_{14}\} \to \{e_{25},e_{23},e_{24},e_{28},e_{27}\}$;
$c_{10} = \{v_2,v_1,v_{17},v_{13},v_{14},v_2\} \to \{e_1,e_2,e_{26},e_{25},e_5\}$;
$c_{11} = \{v_9,v_{18},v_{19},v_{20},v_8,v_9\} \to \{e_{19},e_{29},e_{30},e_{17},e_{16}\}$.

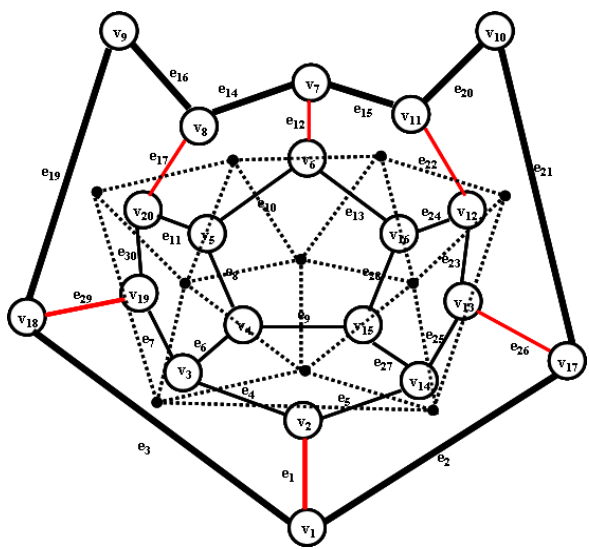

Рис. 14.3. Рисунок суграфа графа G с удаленным ребром $e_{18}$.

Вектор количества циклов, проходящих по рёбрам:

$P_e(G') = \langle 2,1,1,2,2,2,2,2,2,2,2,2,2,1,1,1,2,0,1,1,1,2,2,2,2,2,2,2,2,2 \rangle$.

Вектор количества циклов, проходящих по вершинам:

$P_v(G') = \langle 2,3,3,3,3,3,3,2,2,1,1,2,3,3,3,3,3,2,2,3,3 \rangle$.

Часть графа циклов H представлена на рис. 14.4. Имеем четыре цикла соприкасающихся с опорным циклом $c_2$ – $c_1,c_6,c_{11},c_5$. Пометим эти циклы для удаления. Циклы $c_{11}$ и $c_6$ нельзя удалять из системы циклов, так как их удаление приводит к удалению вершин из суграфа. Можно удалить цикл $c_1$ или $c_5$. Удаляем цикл $c_1$, имеющий максимальную валентность в графе циклов H равную 4, и ребро $e_{15}$ (рис. 14.5). Теперь опорным циклом будем считать цикл $c_1$.

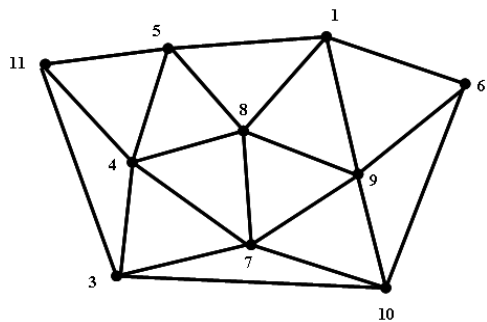

Рис. 14.4. Часть графа циклов H.

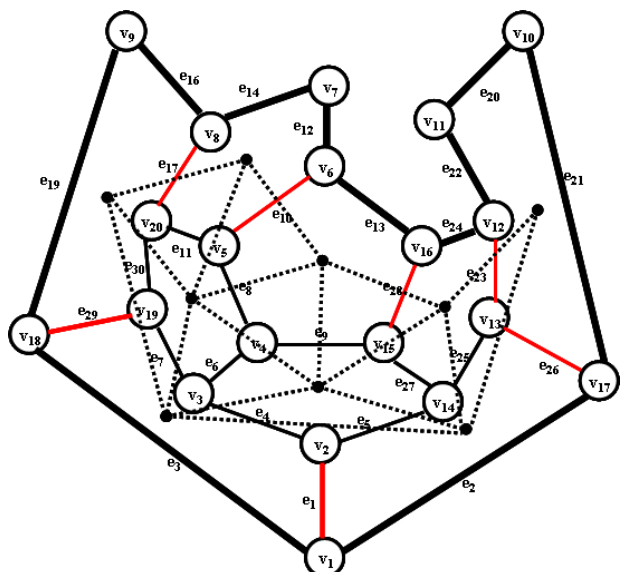

Рис. 14.5. Рисунок суграфа графа G с удаленными ребрами $e_{18}$, $e_{15}$.

Множество оставшихся изометрических циклов суграфа:

$c_3 = \{v_2,v_3,v_{19},v_{18},v_1,v_2\} \rightarrow \{e_4,e_7,e_{29},e_3,e_1\}$;
$c_4 = \{v_3,v_4,v_5,v_{20},v_{19},v_3\} \rightarrow \{e_6,e_8,e_{11},e_{30},e_7\}$;
$c_5 = \{v_5,v_6,v_7,v_8,v_{20},v_5\} \rightarrow \{e_{10},e_{12},e_{14},e_{17},e_{11}\}$;
$c_6 = \{v_{10},v_{11},v_{12},v_{13},v_{17},v_{10}\} \rightarrow \{e_{20},e_{22},e_{23},e_{26},e_{21}\}$;
$c_7 = \{v_2,v_{14},v_{15},v_4,v_3,v_2\} \rightarrow \{e_5,e_{27},e_9,e_6,e_4\}$;
$c_8 = \{v_{15},v_{16},v_6,v_5,v_4,v_{15}\} \rightarrow \{e_{28},e_{13},e_{10},e_8,e_9\}$;
$c_9 = \{v_{14},v_{13},v_{12},v_{16},v_{15},v_{14}\} \rightarrow \{e_{25},e_{23},e_{24},e_{28},e_{27}\}$;
$c_{10} = \{v_2,v_1,v_{17},v_{13},v_{14},v_2\} \rightarrow \{e_1,e_2,e_{26},e_{25},e_5\}$;
$c_{11} = \{v_9,v_{18},v_{19},v_{20},v_8,v_9\} \rightarrow \{e_{19},e_{29},e_{30},e_{17},e_{16}\}$.

Вектор количества циклов, проходящих по рёбрам:

$P_e(G') = <2,1,1,2,2,2,2,2,2,2,2,1,1,1,0,1,2,0,1,1,1,1,2,1,2,2,2,2,2,2>$.

Вектор количества циклов, проходящих по вершинам:

$P_v(G') = <2,3,3,3,3,2,1,2,1,1,1,2,3,3,3,2,1,2,3,3>$.

Часть графа циклов H представлена на рис. 14.6. Опорный цикл $c_1$ соприкасается с циклами $c_8$ и $c_9$, имеющими валентность в графе циклов H равную 4. Пометим эти циклы. Выбираем для удаления цикл $c_8$ и ребро $e_{13}$ (рис. 14.7). Цикл $c_8$ будем считать опорным циклом.

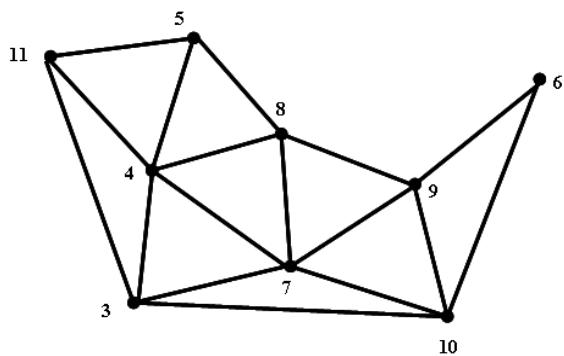

Рис. 14.6. Часть графа циклов H.

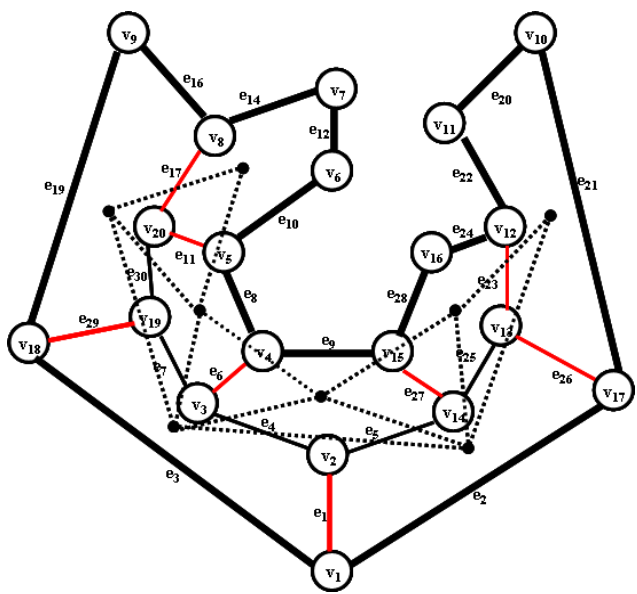

Рис. 14.7. Рисунок суграфа графа G с удаленными ребрами $e_{18}$, $e_{15}$, $e_{13}$.

Множество изометрических циклов суграфа после удаления цикла $c_8$:

$c_3 = \{v_2,v_3,v_{19},v_{18},v_1,v_2\} \to \{e_4,e_7,e_{29},e_3,e_1\}$;
$c_4 = \{v_3,v_4,v_5,v_{20},v_{19},v_3\} \to \{e_6,e_8,e_{11},e_{30},e_7\}$;
$c_5 = \{v_5,v_6,v_7,v_8,v_{20},v_5\} \to \{e_{10},e_{12},e_{14},e_{17},e_{11}\}$;
$c_6 = \{v_{10},v_{11},v_{12},v_{13},v_{17},v_{10}\} \to \{e_{20},e_{22},e_{23},e_{26},e_{21}\}$;
$c_7 = \{v_2,v_{14},v_{15},v_4,v_3,v_2\} \to \{e_5,e_{27},e_9,e_6,e_4\}$;
$c_9 = \{v_{14},v_{13},v_{12},v_{16},v_{15},v_{14}\} \to \{e_{25},e_{23},e_{24},e_{28},e_{27}\}$;
$c_{10} = \{v_2,v_1,v_{17},v_{13},v_{14},v_2\} \to \{e_1,e_2,e_{26},e_{25},e_5\}$;
$c_{11} = \{v_9,v_{18},v_{19},v_{20},v_8,v_9\} \to \{e_{19},e_{29},e_{30},e_{17},e_{16}\}$.

Часть графа циклов H изменилась, она представлена на рис. 14.8. Соприкасающие циклы с опорным циклом $c_8 - c_4$ и $c_7$. Пометим их. Для удаления выбираем цикл $c_4$ с максимальной валентностью в H и одновременно удаляем ребро $e_8$ (рис. 14.9).

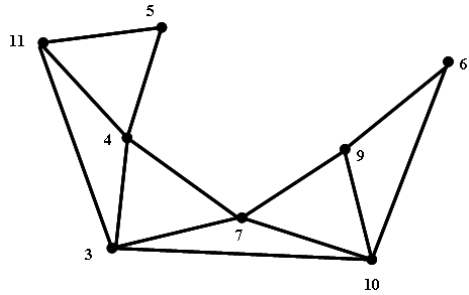

Рис. 14.8. Часть графа циклов H.

Вектор количества циклов, проходящих по рёбрам:

$P_e(G') = <2,1,1,2,2,2,2,1,1,1,2,1,0,1,0,1,2,0,1,1,1,1,2,1,2,2,2,1,2,2>$.

Вектор количества циклов, проходящих по вершинам:

$P_v(G') = <2,3,3,2,2,1,1,2,1,1,1,2,3,3,2,1,2,2,3,3>$.

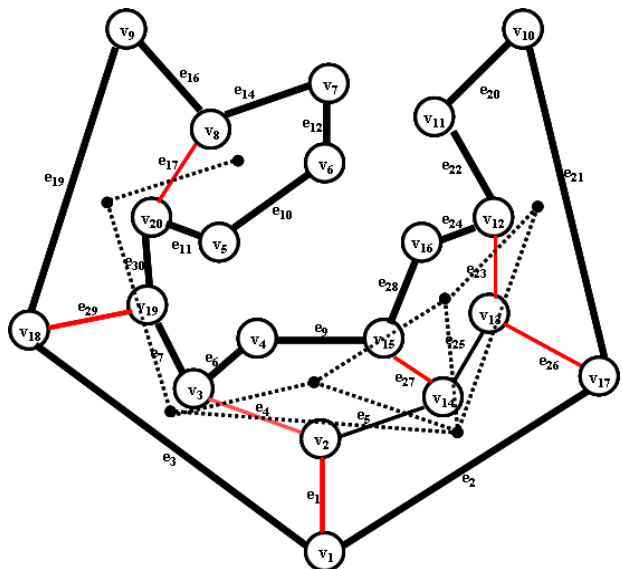

Рис. 14.9. Рисунок суграфа графа G с удалёнными рёбрами $e_{18}, e_{15}, e_{13}, e_8$.

Множество изометрических циклов суграфа после удаления цикла $c_4$:

$c_3 = \{v_2,v_3,v_{19},v_{18},v_1,v_2\} \to \{e_4,e_7,e_{29},e_3,e_1\}$;
$c_5 = \{v_5,v_6,v_7,v_8,v_{20},v_5\} \to \{e_{10},e_{12},e_{14},e_{17},e_{11}\}$;
$c_6 = \{v_{10},v_{11},v_{12},v_{13},v_{17},v_{10}\} \to \{e_{20},e_{22},e_{23},e_{26},e_{21}\}$;
$c_7 = \{v_2,v_{14},v_{15},v_4,v_3,v_2\} \to \{e_5,e_{27},e_9,e_6,e_4\}$;
$c_9 = \{v_{14},v_{13},v_{12},v_{16},v_{15},v_{14}\} \to \{e_{25},e_{23},e_{24},e_{28},e_{27}\}$;
$c_{10} = \{v_2,v_1,v_{17},v_{13},v_{14},v_2\} \to \{e_1,e_2,e_{26},e_{25},e_5\}$;
$c_{11} = \{v_9,v_{18},v_{19},v_{20},v_8,v_9\} \to \{e_{19},e_{29},e_{30},e_{17},e_{16}\}$.

Вектор количества циклов, проходящих по рёбрам:

$P_e(G') = <2,1,1,2,2,1,1,0,1,1,1,1,0,1,0,1,2,0,1,1,1,1,2,1,2,2,2,2,2,1>$.

Вектор количества циклов, проходящих по вершинам:

$P_v(G') = <2,3,2,1,1,1,1,2,1,1,1,2,3,3,2,1,2,2,2,2>$.

Часть графа циклов H представлена на рис. 14.10. Максимальная валентность в графе H у помеченных циклов с опорным циклом $c_4$ только у цикла $c_7$. Удаление цикла $c_7$ приводит к удалению из графа вершины $v_4$. Поэтому помечаем все циклы, соприкасающиеся с обручем выбранного подмножества изометрических циклов. Такой цикл существует, это $c_{10}$. Остается удалить цикл $c_{10}$ и ребро $e_2$ (рис. 14.11).

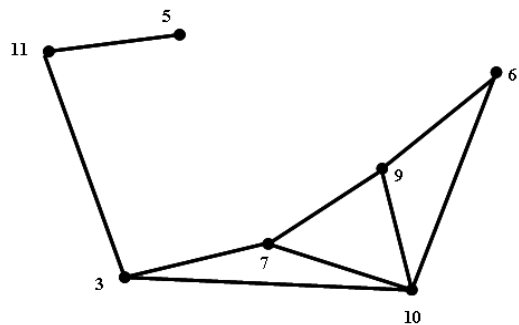

Рис. 14.10. Часть графа циклов H.

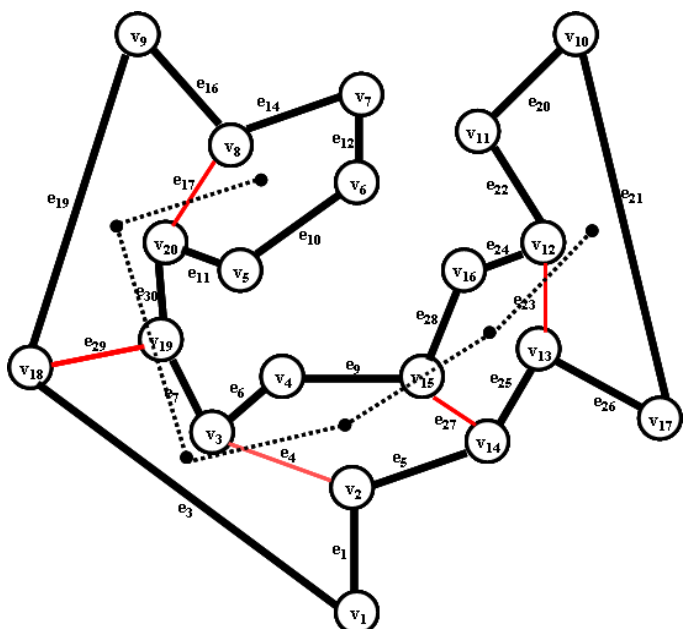

Рис. 14.11. Рисунок суграфа графа G с удаленными ребрами $e_{18}$, $c_{15}$, $e_{13}$, $e_8$, $e_2$.

Множество изометрических циклов суграфа:

$c_3 = \{v_2, v_3, v_{19}, v_{18}, v_1, v_2\} \rightarrow \{e_4, e_7, e_{29}, e_3, e_1\}$;
$c_5 = \{v_5, v_6, v_7, v_8, v_{20}, v_5\} \rightarrow \{e_{10}, e_{12}, e_{14}, e_{17}, e_{11}\}$;
$c_6 = \{v_{10}, v_{11}, v_{12}, v_{13}, v_{17}, v_{10}\} \rightarrow \{e_{20}, e_{22}, e_{23}, e_{26}, e_{21}\}$;
$c_7 = \{v_2, v_{14}, v_{15}, v_4, v_3, v_2\} \rightarrow \{e_5, e_{27}, e_9, e_6, e_4\}$;
$c_9 = \{v_{14}, v_{13}, v_{12}, v_{16}, v_{15}, v_{14}\} \rightarrow \{e_{25}, e_{23}, e_{24}, e_{28}, e_{27}\}$;
$c_{11} = \{v_9, v_{18}, v_{19}, v_{20}, v_8, v_9\} \rightarrow \{e_{19}, e_{29}, e_{30}, e_{17}, e_{16}\}$.

Вектор количества циклов, проходящих по ребрам:

$P_e(G') = <1,0,1,2,1,1,1,0,1,1,1,1,0,1,0,1,2,0,1,1,1,1,2,1,1,1,2,1,2,1>$.

Вектор количества циклов, проходящих по вершинам:

$P_v(G') = <1,2,2,1,1,1,1,2,1,1,1,2,2,2,2,1,1,2,2,2>$.

Часть графа циклов H представлена на рис. 14.12. Обод системы циклов представляет собой гамильтонов цикл (рис. 14.13).

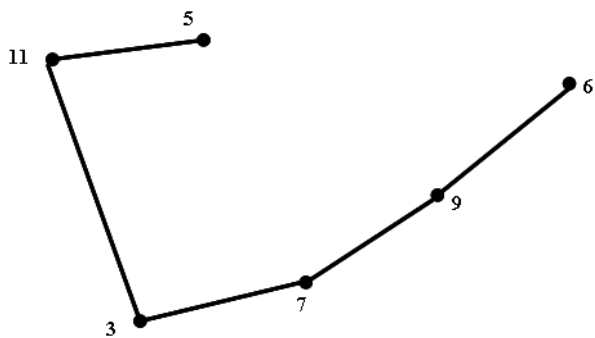

Рис. 14.12. Часть графа циклов H.

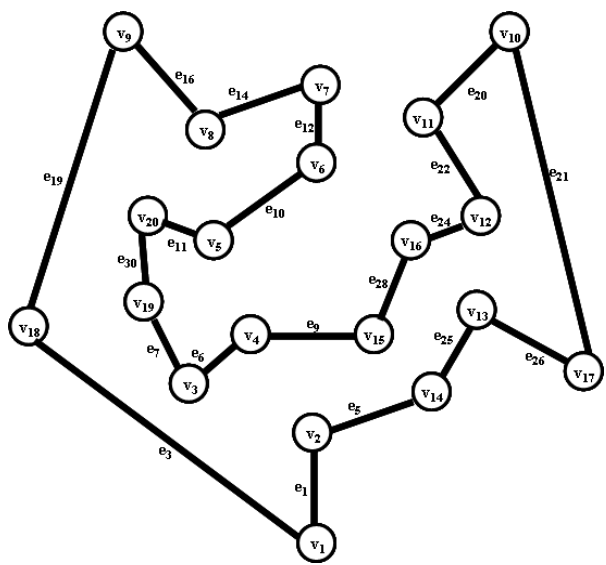

Рис. 14.13. Гамильтонов цикл.

Вектор $P_e$ имеет двадцать единиц – число, равное количеству вершин в графе. Следовательно, это гамильтонов цикл (рис.14.12).

Рассмотрим случай, когда после удаления цикла $c_8$ удаляется цикл $c_7$ (рис. 14.14).

В этом случае удаление цикла $c_7$ с максимальной валентностью приводит к графу циклов с мостом (рис. 14.15). Граф циклов с мостом опасен тем, что в будущем удаление любой конечной вершины моста приведет к разрыву графа циклов на большее число компонент связности. В данном случае образовался мост с петлей (правая сторона).

Будем рассматривать случаи выделения гамильтонова цикла для плоских суграфов. Исключим из рассмотрения плоские суграфы с точкой сочленения. Для плоского суграфа (рис. 14.16) граф циклов H является мультиграфом.

Рис. 14.14. Рисунок суграфа графа G с удаленными ребрами $e_{18}$, $e_{15}$, $e_{13}$, $e_9$.

Рис. 14.15. Часть графа циклов H.

Рис. 14.16. Суграф, граф циклов H которого представляется мультиграфом.

## 14.2. Гамильтонов цикл в непланарном графе

Мы рассмотрели способ выделения гамильтонова цикла в плоском суграфе. А теперь рассмотрим процесс построения гамильтонова цикла в непланарном графе. Следует заметить, что любой цикл в графе, в том числе и гамильтонов, может быть образован как линейная комбинация простых циклов [32].

Пусть имеется множество изометрических циклов графа $C_\tau$. Будем строить цепочку циклов, состоящую только из соприкасающихся циклов. Если в графе существует гамильтонов цикл, то среди множества цепочек, состоящих из соприкасающихся циклов, найдется цепочка, обод которой будет состоять из ребер в количестве равном числу вершин в графе.

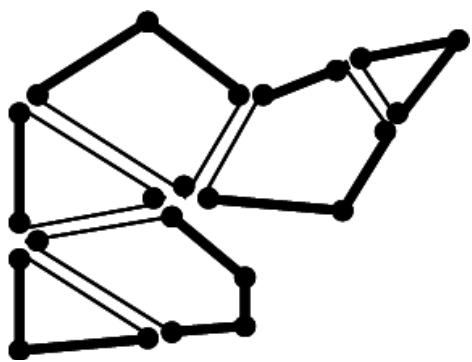 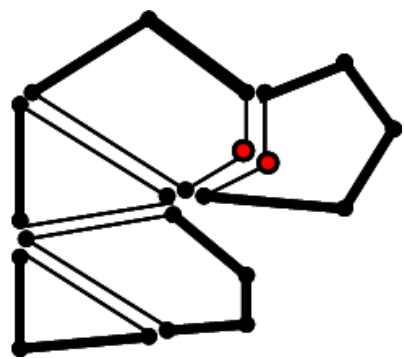

Рис. 14.17. Гамильтонов цикл построенный из соприкасающихся циклов.

Рис. 14.18. Образование не гамильтонова цикла.

На рисунке 14.17 представлен результат процесса построения цепочки соприкасающихся циклов, обод которой образует гамильтонов цикл. На рис. 14.18 показана невозможность построения гамильтонова цикла, так как в цепочку построения попали два несоприкасающихся цикла (циклы, у которых имеются два или более общих ребер).

Построение гамильтонова цикла цепочкой соприкасающихся циклов образует подмножество простых циклов с нулевым значением функционала Маклейна и ободом с количеством ребер равным числу вершин графа.

Таким образом, можно сказать, что гамильтонов цикл есть плоский суграф, который можно представить топологическим рисунком.

В качестве примера рассмотрим построение цепочек, состоящих из соприкасающихся циклов, применяя ранее описанный алгоритм для плоских суграфов графа (рис. 11.7).

Для плоского суграфа имеем следующий топологический рисунок плоской части и граф циклов H (вариант 1, рис. 14.19).

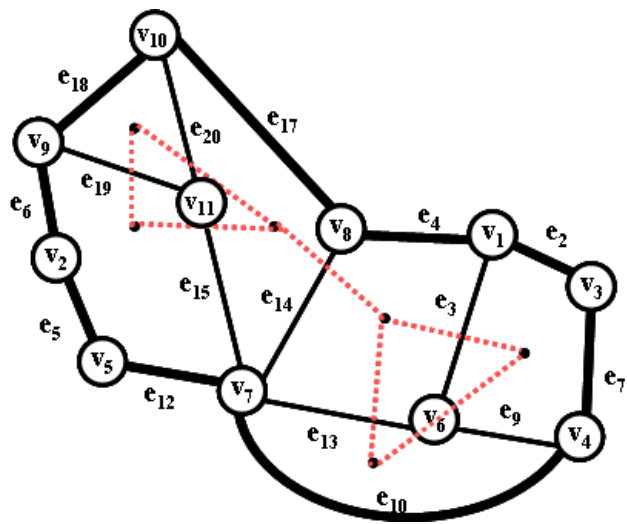

Рис. 14.19. Топологический рисунок плоского суграфа и граф циклов H.

Введем правило удаления циклов: удаляются циклы, имеющие максимальное значение валентности среди помеченных вершин в графе циклов H и это удаление происходит до тех пор, пока граф циклов H не превратится в дерево.

Топологический рисунок плоского графа:

цикл $c_3 = \{e_2,e_3,e_7,e_9\} \leftrightarrow \{v_1,v_3,v_4,v_6\}$;
цикл $c_6 = \{e_3,e_4,e_{13},e_{14}\} \leftrightarrow \{v_1,v_6,v_7,v_8\}$;
цикл $c_8 = \{e_5,e_6,e_{12},e_{15},e_{19}\} \leftrightarrow \{v_2,v_5,v_7,v_9,v_{11}\}$;
цикл $c_{10} = \{e_9,e_{10},e_{13}\} \leftrightarrow \{v_4,v_6,v_7\}$;
цикл $c_{15} = \{e_{14},e_{15},e_{17},e_{20}\} \leftrightarrow \{v_7,v_8,v_{10},v_{11}\}$;
цикл $c_{17} = \{e_{18},e_{19},e_{20}\} \leftrightarrow \{v_9,v_{10},v_{11}\}$.

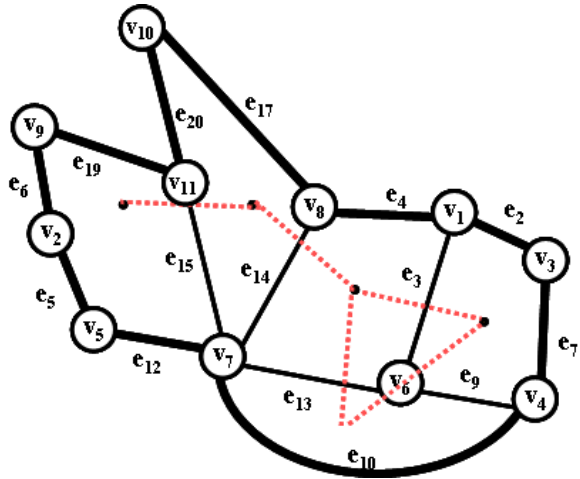

Рис. 14.20. Топологический рисунок после удаления цикла $c_{17}$.

Удаляем цикл $c_{17}$ и ребро $e_{18}$.

Топологический рисунок плоского суграфа:

цикл $c_3 = \{e_2,e_3,e_7,e_9\} \leftrightarrow \{v_1,v_3,v_4,v_6\}$;
цикл $c_6 = \{e_3,e_4,e_{13},e_{14}\} \leftrightarrow \{v_1,v_6,v_7,v_8\}$;
цикл $c_8 = \{e_5,e_6,e_{12},e_{15},e_{19}\} \leftrightarrow \{v_2,v_5,v_7,v_9,v_{11}\}$;
цикл $c_{10} = \{e_9,e_{10},e_{13}\} \leftrightarrow \{v_4,v_6,v_7\}$;
цикл $c_{15} = \{e_{14},e_{15},e_{17},e_{20}\} \leftrightarrow \{v_7,v_8,v_{10},v_{11}\}$.

Номер $<1,2,3,4,5,6,7,8,9,0,1,2,3,4,5,6,7,8,9,0>$.
Вектор количества циклов по рёбрам $P_e = <0,1,2,1,1,1,1,0,2,1,0,1,2,2,2,0,1,0,1,1>$.

Номер $<1,2,3,4,5,6,7,8,9,0,1>$.
Вектор количества циклов по вершинам $P_v = <2,1,1,2,1,3,4,2,1,1,2>$.

Топологический рисунок плоского графа:

цикл $c_3 = \{e_2,e_3,e_7,e_9\} \leftrightarrow \{v_1,v_3,v_4,v_6\}$;
цикл $c_6 = \{e_3,e_4,e_{13},e_{14}\} \leftrightarrow \{v_1,v_6,v_7,v_8\}$;
цикл $c_8 = \{e_5,e_6,e_{12},e_{15},e_{19}\} \leftrightarrow \{v_2,v_5,v_7,v_9,v_{11}\}$;
цикл $c_{15} = \{e_{14},e_{15},e_{17},e_{20}\} \leftrightarrow \{v_7,v_8,v_{10},v_{11}\}$;
цикл $c_{17} = \{e_{18},e_{19},e_{20}\} \leftrightarrow \{v_9,v_{10},v_{11}\}$.

Номер $<1,2,3,4,5,6,7,8,9,0,1,2,3,4,5,6,7,8,9,0>$.
Вектор количества циклов по рёбрам $P_e = <0,1,2,1,1,1,1,0,1,0,0,1,1,2,2,0,1,0,1,1>$.

Номер $<1,2,3,4,5,6,7,8,9,0,1>$.
Вектор количества циклов по вершинам $P_v = <2,1,1,1,1,2,3,2,1,1,2>$.

Количество рёбер в ободе равно количеству вершин *n* = 11.

Обод данного подмножества циклов есть гамильтонов цикл (рис. 14.21).

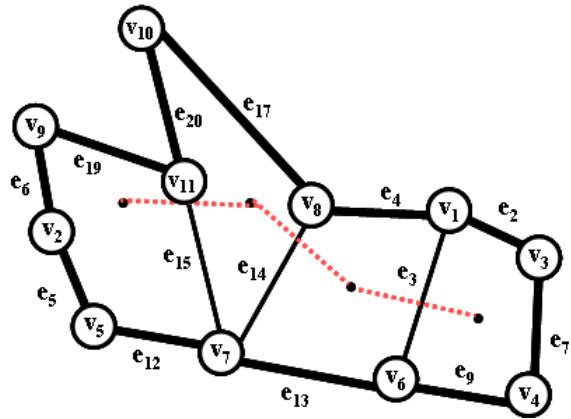

Рис. 14.21. Топологический рисунок после удаления цикла $c_{10}$.

Рассмотрим суграф (вариант 2, рис. 14.20).

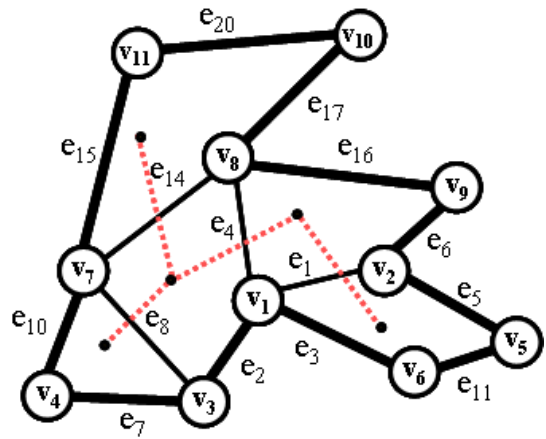

Рис. 14.22. Гамильтонов цикл плоского суграфа.

Топологический рисунок суграфа:

цикл $c_1 = \{e_1, e_3, e_5, e_{11}\} \leftrightarrow \{v_1, v_2, v_5, v_6\}$;
цикл $c_2 = \{e_1, e_4, e_6, e_{16}\} \leftrightarrow \{v_1, v_2, v_8, v_9\}$;
цикл $c_5 = \{e_2, e_4, e_8, e_{14}\} \leftrightarrow \{v_1, v_3, v_7, v_8\}$;
цикл $c_9 = \{e_7, e_8, e_{10}\} \leftrightarrow \{v_3, v_4, v_7\}$;
цикл $c_{15} = \{e_{14}, e_{15}, e_{17}, e_{20}\} \leftrightarrow \{v_7, v_8, v_{10}, v_{11}\}$;

            Номер   <1,2,3,4,5,6,7,8,9,0,1,2,3,4,5,6,7,8,9,0>.
Вектор количества циклов по рёбрам   $P_e$ = <2,1,1,2,1,1,1,2,0,1,1,0,0,2,1,1,1,0,0,1>.

            Номер   <1,2,3,4,5,6,7,8,9,0,1>.
Вектор количества циклов по вершинам $P_v$ = <3,2,2,1,1,1,3,3,1,1,1>.

Обод топологического рисунка представляет гамильтонов цикл (рис. 14.22).

Рассмотрим плоский суграф (вариант 3, рис. 14.23).

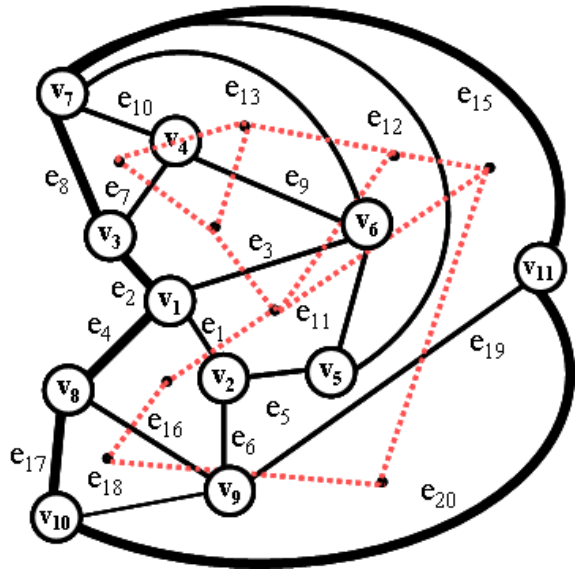

Рис. 14.23. Топологичекий рисунок плоского суграфа.

Топологический рисунок:

цикл $c_1 = \{e_1,e_3,e_5,e_{11}\} \leftrightarrow \{v_1,v_2,v_5,v_6\}$;
цикл $c_2 = \{e_1,e_4,e_6,e_{16}\} \leftrightarrow \{v_1,v_2,v_8,v_9\}$;
цикл $c_3 = \{e_2,e_3,e_7,e_9\} \leftrightarrow \{v_1,v_3,v_4,v_6\}$;
цикл $c_8 = \{e_5,e_6,e_{12},e_{15},e_{19}\} \leftrightarrow \{v_2,v_5,v_7,v_9,v_{11}\}$;
цикл $c_9 = \{e_7,e_8,e_{10}\} \leftrightarrow \{v_3,v_4,v_7\}$;
цикл $c_{10} = \{e_9,e_{10},e_{13}\} \leftrightarrow \{v_4,v_6,v_7\}$;
цикл $c_{11} = \{e_{11},e_{12},e_{13}\} \leftrightarrow \{v_5,v_6,v_7\}$;
цикл $c_{16} = \{e_{16},e_{17},e_{18}\} \leftrightarrow \{v_8,v_9,v_{10}\}$;
цикл $c_{17} = \{e_{18},e_{19},e_{20}\} \leftrightarrow \{v_9,v_{10},v_{11}\}$.

Номер $<1,2,3,4,5,6,7,8,9,0,1,2,3,4,5,6,7,8,9,0>$.
Вектор количества циклов по рёбрам $P_e = <2,1,2,1,2,2,2,1,2,2,2,2,2,0,1,2,1,2,2,1>$.

Номер $<1,2,3,4,5,6,7,8,9,0,1>$.
Вектор количества циклов по вершинам $P_v = <3,3,2,3,3,4,3,2,4,2,2>$.

Удаляем цикл $c_8$ как имеющий большую валентность в графе циклов.

Топологический рисунок:

цикл $c_1 = \{e_1,e_3,e_5,e_{11}\} \leftrightarrow \{v_1,v_2,v_5,v_6\}$;
цикл $c_2 = \{e_1,e_4,e_6,e_{16}\} \leftrightarrow \{v_1,v_2,v_8,v_9\}$;
цикл $c_3 = \{e_2,e_3,e_7,e_9\} \leftrightarrow \{v_1,v_3,v_4,v_6\}$;
цикл $c_9 = \{e_7,e_8,e_{10}\} \leftrightarrow \{v_3,v_4,v_7\}$;
цикл $c_{10} = \{e_9,e_{10},e_{13}\} \leftrightarrow \{v_4,v_6,v_7\}$;
цикл $c_{11} = \{e_{11},e_{12},e_{13}\} \leftrightarrow \{v_5,v_6,v_7\}$;
цикл $c_{16} = \{e_{16},e_{17},e_{18}\} \leftrightarrow \{v_8,v_9,v_{10}\}$;
цикл $c_{17} = \{e_{18},e_{19},e_{20}\} \leftrightarrow \{v_9,v_{10},v_{11}\}$.

Номер $<1,2,3,4,5,6,7,8,9,0,1,2,3,4,5,6,7,8,9,0>$.
Вектор количества циклов по рёбрам $P_e = <2,0,1,1,2,2,1,1,2,2,2,2,2,0,1,1,1,2,2,1>$.

Номер $<1,2,3,4,5,6,7,8,9,0,1>$.
Вектор количества циклов по вершинам $P_v = <3,2,2,3,2,4,3,2,3,2,1>$.

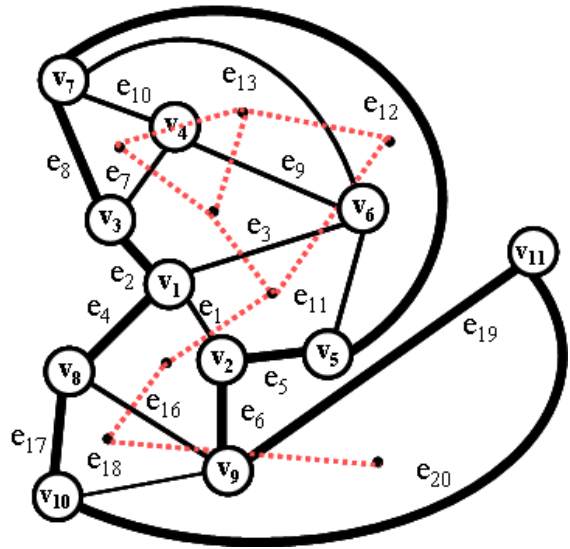

Рис. 14.24. Топологичекий рисунок плоской части графа
после удаления цикла $c_8$.

Удаляем цикл $c_3$ как имеющий большую валентность в графе циклов.

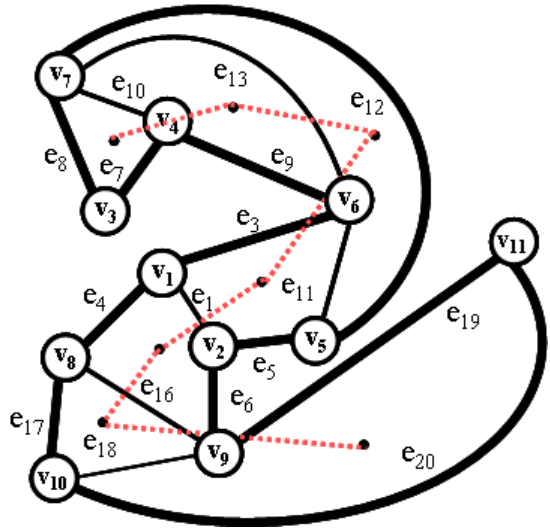

Рис. 14.25. Гамильтонов цикл плоского суграфа после удаления цикла $c_3$.

Топологический рисунок:

цикл $c_1 = \{e_1, e_3, e_5, e_{11}\} \leftrightarrow \{v_1, v_2, v_5, v_6\}$;
цикл $c_2 = \{e_1, e_4, e_6, e_{16}\} \leftrightarrow \{v_1, v_2, v_8, v_9\}$;
цикл $c_3 = \{e_2, e_3, e_7, e_9\} \leftrightarrow \{v_1, v_3, v_4, v_6\}$;
цикл $c_9 = \{e_7, e_8, e_{10}\} \leftrightarrow \{v_3, v_4, v_7\}$;
цикл $c_{10} = \{e_9, e_{10}, e_{13}\} \leftrightarrow \{v_4, v_6, v_7\}$;
цикл $c_{11} = \{e_{11}, e_{12}, e_{13}\} \leftrightarrow \{v_5, v_6, v_7\}$;
цикл $c_{16} = \{e_{16}, e_{17}, e_{18}\} \leftrightarrow \{v_8, v_9, v_{10}\}$;
цикл $c_{17} = \{e_{18}, e_{19}, e_{20}\} \leftrightarrow \{v_9, v_{10}, v_{11}\}$.

Номер    <1,2,3,4,5,6,7,8,9,0,1,2,3,4,5,6,7,8,9,0>.
Вектор количества циклов по ребрам    $P_e = <2,0,1,1,1,1,1,1,1,2,2,1,2,0,0,2,1,2,1,1>$.
Номер    <1,2,3,4,5,6,7,8,9,0,1>.
Вектор количества циклов по вершинам $P_v = <2,2,1,2,2,3,2,2,3,2,1>$.

Обод топологического рисунка представляет гамильтонов цикл (рис. 14.25).

Рассмотрим плоский суграф (вариант 4, рис. 14.24).

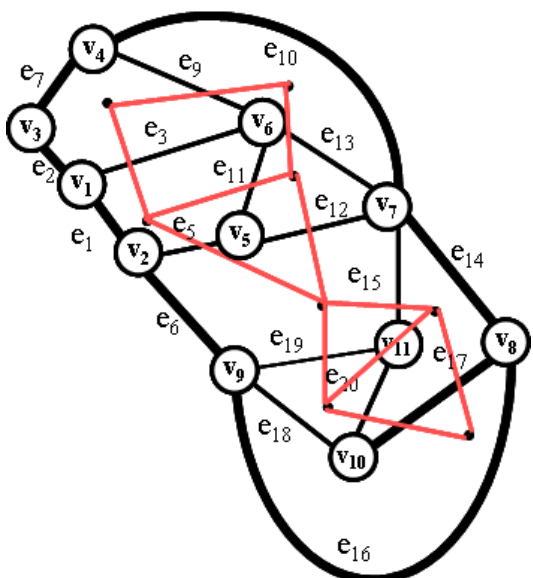

Рис. 14.26. Плоский суграф.

Топологический рисунок суграфа:

цикл $c_1 = \{e_1,e_3,e_5,e_{11}\} \leftrightarrow \{v_1,v_2,v_5,v_6\}$;
цикл $c_3 = \{e_2,e_3,e_7,e_9\} \leftrightarrow \{v_1,v_3,v_4,v_6\}$;
цикл $c_8 = \{e_5,e_6,e_{12},e_{15},e_{19}\} \leftrightarrow \{v_2,v_5,v_7,v_9,v_{11}\}$;
цикл $c_{10} = \{e_9,e_{10},e_{13}\} \leftrightarrow \{v_4,v_6,v_7\}$;
цикл $c_{11} = \{e_{11},e_{12},e_{13}\} \leftrightarrow \{v_5,v_6,v_7\}$;
цикл $c_{15} = \{e_{14},e_{15},e_{17},e_{20}\} \leftrightarrow \{v_7,v_8,v_{10},v_{11}\}$;
цикл $c_{16} = \{e_{16},e_{17},e_{18}\} \leftrightarrow \{v_8,v_9,v_{10}\}$;
цикл $c_{17} = \{e_{18},e_{19},e_{20}\} \leftrightarrow \{v_9,v_{10},v_{11}\}$.

$$\text{Номер } <1,2,3,4,5,6,7,8,9,0,1,2,3,4,5,6,7,8,9,0>.$$
Вектор количества циклов по рёбрам $P_e = <1,1,2,0,2,1,1,0,2,1,2,2,2,1,2,1,2,2,2,2>$.

$$\text{Номер } <1,2,3,4,5,6,7,8,9,0,1>.$$
Вектор количества циклов по вершинам $P_v = <2,2,1,2,3,4,4,2,3,3,3>$.

Удаляем цикл $c_1$.

Топологический рисунок суграфа:

цикл $c_3 = \{e_2,e_3,e_7,e_9\} \leftrightarrow \{v_1,v_3,v_4,v_6\}$;
цикл $c_8 = \{e_5,e_6,e_{12},e_{15},e_{19}\} \leftrightarrow \{v_2,v_5,v_7,v_9,v_{11}\}$;
цикл $c_{10} = \{e_9,e_{10},e_{13}\} \leftrightarrow \{v_4,v_6,v_7\}$;
цикл $c_{11} = \{e_{11},e_{12},e_{13}\} \leftrightarrow \{v_5,v_6,v_7\}$;
цикл $c_{15} = \{e_{14},e_{15},e_{17},e_{20}\} \leftrightarrow \{v_7,v_8,v_{10},v_{11}\}$;
цикл $c_{16} = \{e_{16},e_{17},e_{18}\} \leftrightarrow \{v_8,v_9,v_{10}\}$;
цикл $c_{17} = \{e_{18},e_{19},e_{20}\} \leftrightarrow \{v_9,v_{10},v_{11}\}$.

$$\text{Номер } <1,2,3,4,5,6,7,8,9,0,1,2,3,4,5,6,7,8,9,0>.$$
Вектор количества циклов по рёбрам $P_e = <0,1,1,0,1,1,1,0,2,1,1,2,2,1,2,1,2,2,2,2>$.

$$\text{Номер } <1,2,3,4,5,6,7,8,9,0,1>.$$
Вектор количества циклов по вершинам $P_v = <1,1,1,2,2,3,4,2,3,3,3>$.

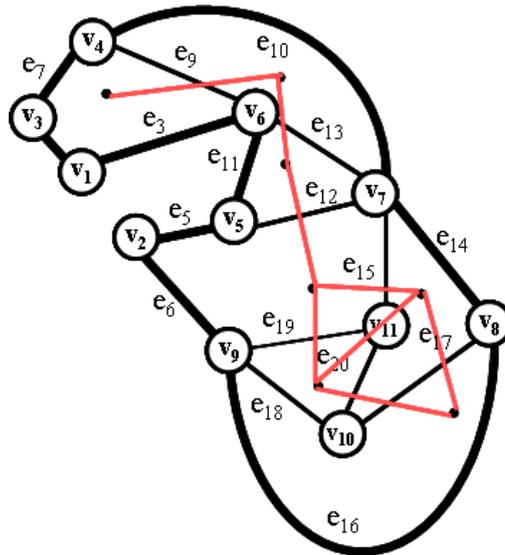

Рис. 14.27. Плоская часть графа после удаления цикла $c_1$.

Удаляем цикл $c_{15}$.

Топологический рисунок суграфа:

цикл $c_3 = \{e_2,e_3,e_7,e_9\} \leftrightarrow \{v_1,v_3,v_4,v_6\}$;
цикл $c_8 = \{e_5,e_6,e_{12},e_{15},e_{19}\} \leftrightarrow \{v_2,v_5,v_7,v_9,v_{11}\}$;
цикл $c_{10} = \{e_9,e_{10},e_{13}\} \leftrightarrow \{v_4,v_6,v_7\}$;
цикл $c_{11} = \{e_{11},e_{12},e_{13}\} \leftrightarrow \{v_5,v_6,v_7\}$;
цикл $c_{16} = \{e_{16},e_{17},e_{18}\} \leftrightarrow \{v_8,v_9,v_{10}\}$;
цикл $c_{17} = \{e_{18},e_{19},e_{20}\} \leftrightarrow \{v_9,v_{10},v_{11}\}$.

Номер <1,2,3,4,5,6,7,8,9,0,1,2,3,4,5,6,7,8,9,0>.
Вектор количества циклов по рёбрам $P_e = $ <0,1,1,0,1,1,1,0,2,1,1,2,2,0,1,1,1,2,2,1>.

Номер <1,2,3,4,5,6,7,8,9,0,1>.
Вектор количества циклов по вершинам $P_v = $ <1,1,1,2,2,3,2,2,3,3,2>.

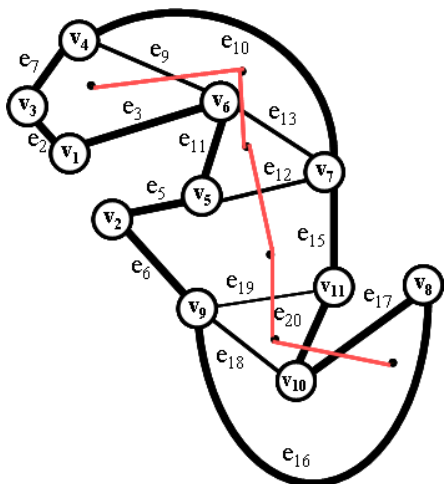

Рис. 14.28. Гамильтонов цикл после удаления цикла $c_{15}$.

Рассмотрим плоский суграф (вариант 5, рис. 14.29).

Топологический рисунок плоской части:

цикл $c_1 = \{e_1,e_3,e_5,e_{11}\} \leftrightarrow \{v_1,v_2,v_5,v_6\}$;
цикл $c_2 = \{e_1,e_4,e_6,e_{16}\} \leftrightarrow \{v_1,v_2,v_8,v_9\}$;

цикл $c_3 = \{e_2,e_3,e_7,e_9\} \leftrightarrow \{v_1,v_3,v_4,v_6\}$;
цикл $c_5 = \{e_2,e_4,e_8,e_{14}\} \leftrightarrow \{v_1,v_3,v_7,v_8\}$;
цикл $c_9 = \{e_7,e_8,e_{10}\} \leftrightarrow \{v_3,v_4,v_7\}$;
цикл $c_{16} = \{e_{16},e_{17},e_{18}\} \leftrightarrow \{v_8,v_9,v_{10}\}$;
цикл $c_{17} = \{e_{18},e_{19},e_{20}\} \leftrightarrow \{v_9,v_{10},v_{11}\}$.

Номер <1,2,3,4,5,6,7,8,9,0,1,2,3,4,5,6,7,8,9,0>.
Вектор количества циклов по рёбрам $P_e$ = <2,2,2,2,1,1,2,2,1,1,1,0,0,1,0,2,1,2,1,1>.

Номер <1,2,3,4,5,6,7,8,9,0,1>.
Вектор количества циклов по вершинам $P_v$ = <4,2,3,2,2,2,2,3,3,2,1>.

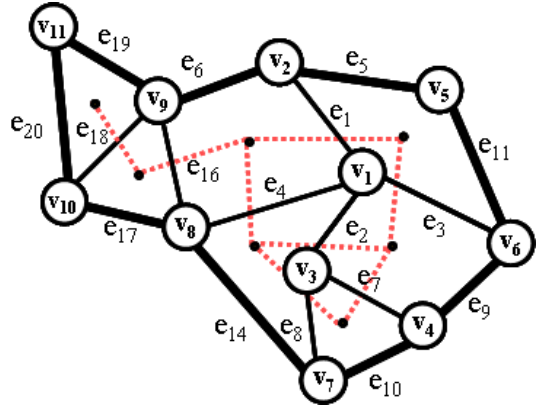

Рис. 14.29. Плоский суграф.

Гамильтонов цикл образуется путем удаления цикла $c_3$ с максимальной валентностью в графе циклов (рис. 14.30).

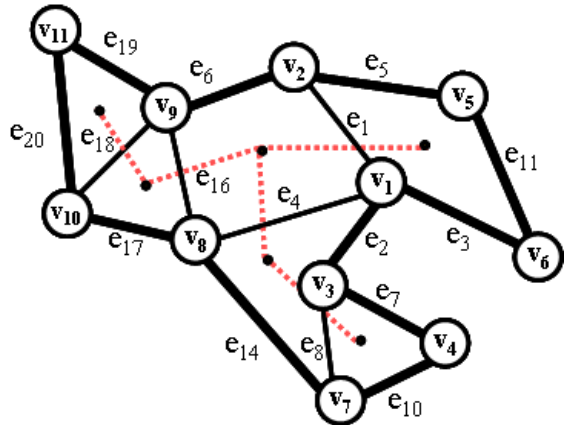

Рис. 14.30. Гамильтонов цикл после удаления цикла $c_3$.

Топологический рисунок плоской части:

цикл $c_1 = \{e_1,e_3,e_5,e_{11}\} \leftrightarrow \{v_1,v_2,v_5,v_6\}$;
цикл $c_2 = \{e_1,e_4,e_6,e_{16}\} \leftrightarrow \{v_1,v_2,v_8,v_9\}$;
цикл $c_3 = \{e_2,e_3,e_7,e_9\} \leftrightarrow \{v_1,v_3,v_4,v_6\}$;
цикл $c_5 = \{e_2,e_4,e_8,e_{14}\} \leftrightarrow \{v_1,v_3,v_7,v_8\}$;
цикл $c_9 = \{e_7,e_8,e_{10}\} \leftrightarrow \{v_3,v_4,v_7\}$;
цикл $c_{16} = \{e_{16},e_{17},e_{18}\} \leftrightarrow \{v_8,v_9,v_{10}\}$;
цикл $c_{17} = \{e_{18},e_{19},e_{20}\} \leftrightarrow \{v_9,v_{10},v_{11}\}$.

Номер <1,2,3,4,5,6,7,8,9,0,1,2,3,4,5,6,7,8,9,0>.
Вектор количества циклов по рёбрам $P_e$ = <2,1,1,2,1,1,1,2,0,1,1,0,0,1,0,2,1,2,1,1>.
Номер <1,2,3,4,5,6,7,8,9,0,1>.
Вектор количества циклов по вершинам $P_v$ = <3,2,2,1,2,1,3,3,2,1>.

Рассмотрим плоский суграф (вариант 6, рис. 14.31).

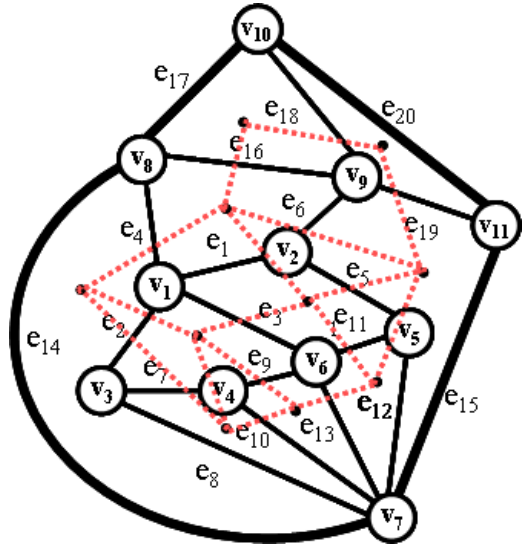

Рис. 14.31. Плоский суграф.

цикл $c_1 = \{e_1,e_3,e_5,e_{11}\} \leftrightarrow \{v_1,v_2,v_5,v_6\}$;
цикл $c_2 = \{e_1,e_4,e_6,e_{16}\} \leftrightarrow \{v_1,v_2,v_8,v_9\}$;
цикл $c_3 = \{e_2,e_3,e_7,e_9\} \leftrightarrow \{v_1,v_3,v_4,v_6\}$;
цикл $c_5 = \{e_2,e_4,e_8,e_{14}\} \leftrightarrow \{v_1,v_3,v_7,v_8\}$;
цикл $c_8 = \{e_5,e_6,e_{12},e_{15},e_{19}\} \leftrightarrow \{v_2,v_5,v_7,v_9,v_{11}\}$;
цикл $c_9 = \{e_7,e_8,e_{10}\} \leftrightarrow \{v_3,v_4,v_7\}$;
цикл $c_{10} = \{e_9,e_{10},e_{13}\} \leftrightarrow \{v_4,v_6,v_7\}$;
цикл $c_{11} = \{e_{11},e_{12},e_{13}\} \leftrightarrow \{v_5,v_6,v_7\}$;
цикл $c_{16} = \{e_{16},e_{17},e_{18}\} \leftrightarrow \{v_8,v_9,v_{10}\}$;
цикл $c_{17} = \{e_{18},e_{19},e_{20}\} \leftrightarrow \{v_9,v_{10},v_{11}\}$.

$\phantom{aaaaaaaaaaaaaaaaaaaaaaaaaaaaaaa}$ Номер $<1,2,3,4,5,6,7,8,9,0,1,2,3,4,5,6,7,8,9,0>$.
Вектор количества циклов по рёбрам $\phantom{aa} P_e = <2,2,2,2,2,2,2,2,2,2,2,2,2,1,1,2,1,1,1,0>$.

$\phantom{aaaaaaaaaaaaaaaaaaaaaaaaaaaaaaa}$ Номер $<1,2,3,4,5,6,7,8,9,0,1>$.
Вектор количества циклов по вершинам $P_v = <4,3,3,3,3,4,5,3,4,2,2>$.

Удаляем цикл $c_8$ имеющий максимальную валентность в графе циклов Н.

цикл $c_1 = \{e_1,e_3,e_5,e_{11}\} \leftrightarrow \{v_1,v_2,v_5,v_6\}$;
цикл $c_2 = \{e_1,e_4,e_6,e_{16}\} \leftrightarrow \{v_1,v_2,v_8,v_9\}$;
цикл $c_3 = \{e_2,e_3,e_7,e_9\} \leftrightarrow \{v_1,v_3,v_4,v_6\}$;
цикл $c_5 = \{e_2,e_4,e_8,e_{14}\} \leftrightarrow \{v_1,v_3,v_7,v_8\}$;
цикл $c_8 = \{e_5,e_6,e_{12},e_{15},e_{19}\} \leftrightarrow \{v_2,v_5,v_7,v_9,v_{11}\}$;
цикл $c_9 = \{e_7,e_8,e_{10}\} \leftrightarrow \{v_3,v_4,v_7\}$;
цикл $c_{10} = \{e_9,e_{10},e_{13}\} \leftrightarrow \{v_4,v_6,v_7\}$;
цикл $c_{11} = \{e_{11},e_{12},e_{13}\} \leftrightarrow \{v_5,v_6,v_7\}$;
цикл $c_{16} = \{e_{16},e_{17},e_{18}\} \leftrightarrow \{v_8,v_9,v_{10}\}$;
цикл $c_{17} = \{e_{18},e_{19},e_{20}\} \leftrightarrow \{v_9,v_{10},v_{11}\}$.

$\phantom{aaaaaaaaaaaaaaaaaaaaaaaaaaaaaaa}$ Номер $<1,2,3,4,5,6,7,8,9,0,1,2,3,4,5,6,7,8,9,0>$.
Вектор количества циклов по рёбрам $\phantom{aa} P_e = <2,2,2,2,1,1,2,2,2,2,2,2,1,2,1,0,2,1,1,1,0>$.

$\phantom{aaaaaaaaaaaaaaaaaaaaaaaaaaaaaaa}$ Номер $<1,2,3,4,5,6,7,8,9,0,1>$.
Вектор количества циклов по вершинам $P_v = <4,2,3,3,2,4,4,3,3,2,1>$.

Удаляем цикл $c_5$ имеющий максимальную валентность в графе циклов H (рис. 14.32).

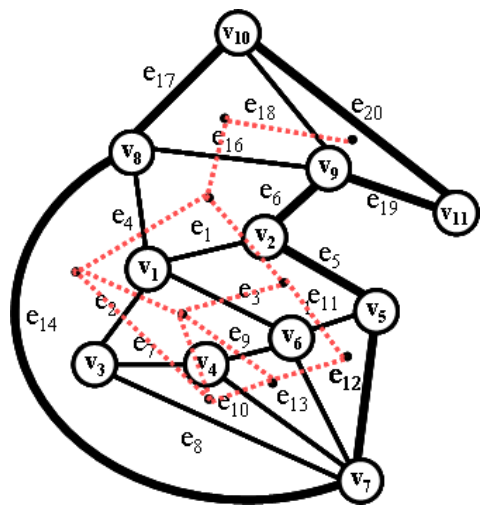

Рис. 14.32. Плоский суграф после удаления ребра $e_{15}$.

цикл $c_1 = \{e_1,e_3,e_5,e_{11}\} \leftrightarrow \{v_1,v_2,v_5,v_6\}$;
цикл $c_2 = \{e_1,e_4,e_6,e_{16}\} \leftrightarrow \{v_1,v_2,v_8,v_9\}$;
цикл $c_3 = \{e_2,e_3,e_7,e_9\} \leftrightarrow \{v_1,v_3,v_4,v_6\}$;
цикл $c_5 = \{e_2,e_4,e_8,e_{14}\} \leftrightarrow \{v_1,v_3,v_7,v_8\}$;
цикл $c_9 = \{e_7,e_8,e_{10}\} \leftrightarrow \{v_3,v_4,v_7\}$;
цикл $c_{10} = \{e_9,e_{10},e_{13}\} \leftrightarrow \{v_4,v_6,v_7\}$;
цикл $c_{11} = \{e_{11},e_{12},e_{13}\} \leftrightarrow \{v_5,v_6,v_7\}$;
цикл $c_{16} = \{e_{16},e_{17},e_{18}\} \leftrightarrow \{v_8,v_9,v_{10}\}$;
цикл $c_{17} = \{e_{18},e_{19},e_{20}\} \leftrightarrow \{v_9,v_{10},v_{11}\}$.

Номер <1,2,3,4,5,6,7,8,9,0,1,2,3,4,5,6,7,8,9,0>.
Вектор количества циклов по рёбрам $P_e = $ <2,1,2,1,1,1,2,1,2,2,2,1,2,0,0,2,1,1,1,0>.

Номер <1,2,3,4,5,6,7,8,9,0,1>.
Вектор количества циклов по вершинам $P_v = $ <3,2,2,3,2,4,3,2,3,2,1>.

Удаляем цикл $c_3$ имеющий максимальную валентность в графе циклов H (рис. 14.33).

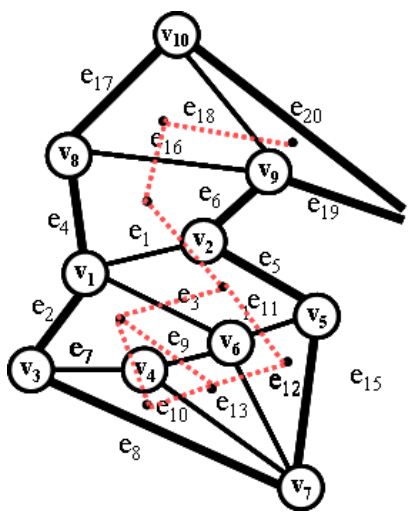

Рис. 14.33. Плоский суграф после удаления ребра $e_2$.

цикл $c_1 = \{e_1,e_3,e_5,e_{11}\} \leftrightarrow \{v_1,v_2,v_5,v_6\}$;
цикл $c_2 = \{e_1,e_4,e_6,e_{16}\} \leftrightarrow \{v_1,v_2,v_8,v_9\}$;
цикл $c_3 = \{e_2,e_3,e_7,e_9\} \leftrightarrow \{v_1,v_3,v_4,v_6\}$;
цикл $c_5 = \{e_2,e_4,e_8,e_{14}\} \leftrightarrow \{v_1,v_3,v_7,v_8\}$;

цикл $c_9 = \{e_7, e_8, e_{10}\} \leftrightarrow \{v_3, v_4, v_7\}$;
цикл $c_{10} = \{e_9, e_{10}, e_{13}\} \leftrightarrow \{v_4, v_6, v_7\}$;
цикл $c_{11} = \{e_{11}, e_{12}, e_{13}\} \leftrightarrow \{v_5, v_6, v_7\}$;
цикл $c_{16} = \{e_{16}, e_{17}, e_{18}\} \leftrightarrow \{v_8, v_9, v_{10}\}$;
цикл $c_{17} = \{e_{18}, e_{19}, e_{20}\} \leftrightarrow \{v_9, v_{10}, v_{11}\}$.

Номер <1,2,3,4,5,6,7,8,9,0,1,2,3,4,5,6,7,8,9,0>.
Вектор количества циклов по ребрам $P_e$ = <2,0,1,1,1,1,1,1,1,2,2,1,2,0,0,2,1,1,1,0>.

Номер <1,2,3,4,5,6,7,8,9,0,1>.
Вектор количества циклов по вершинам $P_v$ = <2,2,1,2,2,3,3,2,3,2,1>.

Выпишем гамильтоновы циклы:

<$e_2, e_4, e_5, e_6, e_7, e_9, e_{12}, e_{13}, e_{17}, e_{19}, e_{20}$>;
<$e_2, e_3, e_5, e_6, e_7, e_{10}, e_{11}, e_{15}, e_{16}, e_{17}, e_{20}$>;
<$e_3, e_4, e_5, e_6, e_7, e_8, e_9, e_{12}, e_{17}, e_{19}, e_{20}$>;
<$e_2, e_3, e_5, e_6, e_7, e_{10}, e_{11}, e_{14}, e_{17}, e_{19}, e_{20}$>.

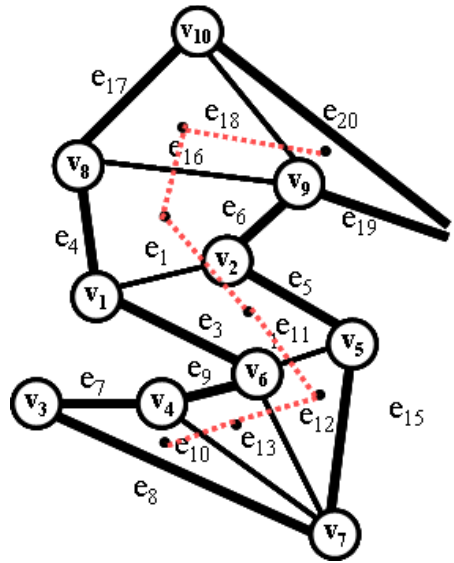

Рис. 14.34. Гамильтонов цикл.

### 14.3. Множество гамильтоновых циклов в плоском суграфе

Определим количество составляющих гамильтонова цикла для соответствующих плоских суграфов:

плоский суграф (вариант 1) – 4 цикла,
плоский суграф (вариант 2) – 5 циклов,
плоский суграф (вариант 3) – 7 циклов,
плоский суграф (вариант 4) – 6 циклов,
плоский суграф (вариант 5) – 6 циклов,
плоский суграф (вариант 6) – 7 циклов.

Следует заметить, что может быть построено несколько гамильтоновых циклов в одном выделенном плоском суграфе. Например, для плоского суграфа (вариант 2), различные гамильтоновы циклы представлены на рис. 14.35.

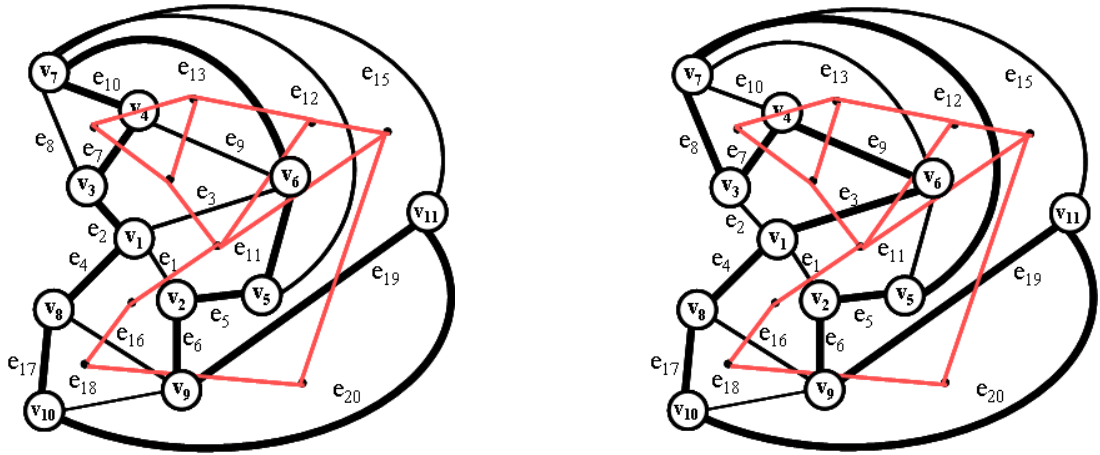

Рис. 14.35. Различные гамильтоновы циклы в суграфе.

Различные гамильтоновы циклы для плоского суграфа (вариант 4), представлены топологическим рисунком (рис. 14.36).

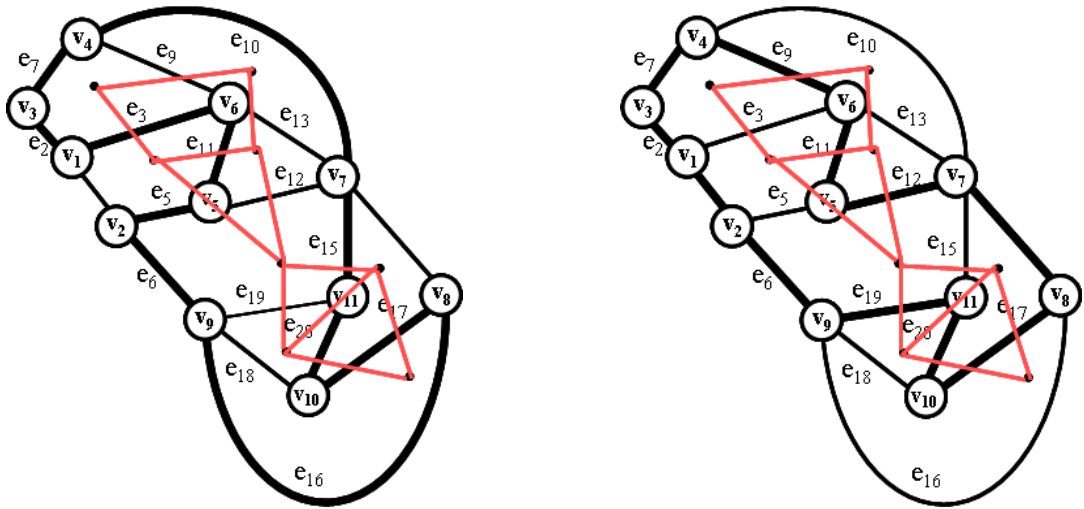

Рис. 14.36. Различные гамильтоновы циклы $G_5^*$.

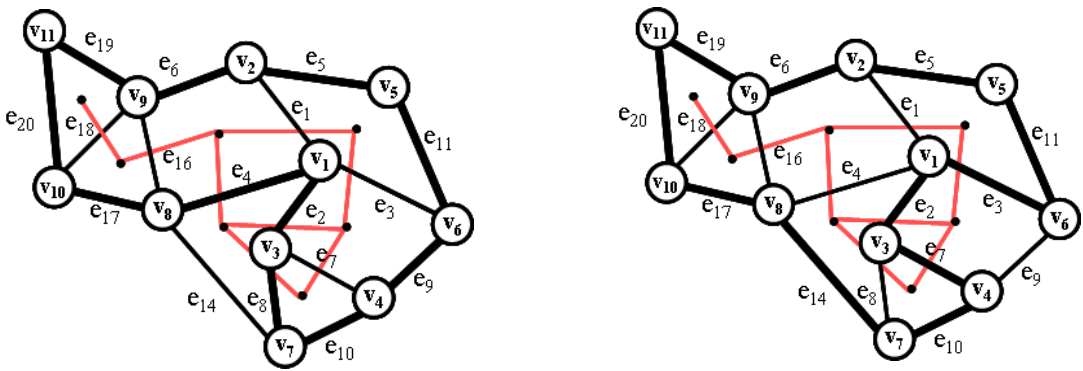

Рис. 14.37. Различные гамильтоновы циклы (вариант 5).

Два различных гамильтоновых цикла для плоского суграфа (вариант 5) представлены топологическим рисунком (рис. 14.37).

Четыре различных гамильтоновых цикла для плоского суграфа (вариант 6) представлены топологическим рисунком (рис. 14.38).

### 14.4. Алгоритм построения гамильтонова цикла в плоском суграфе

Опишем алгоритм построения гамильтонова цикла в плоском суграфе.

Будем рассматривать только несепарабельные суграфы.

**Определение 14.3**. Назовем циклы *соприкасающимися*, если существует ребро, одновременно принадлежащее двум циклам $c_i$ и $c_j$.

$$(e_k \in c_i) \& (c_k \in c_j) \qquad (14.1)$$

Соответственно, такое ребро будем называть *соприкасающимся*.

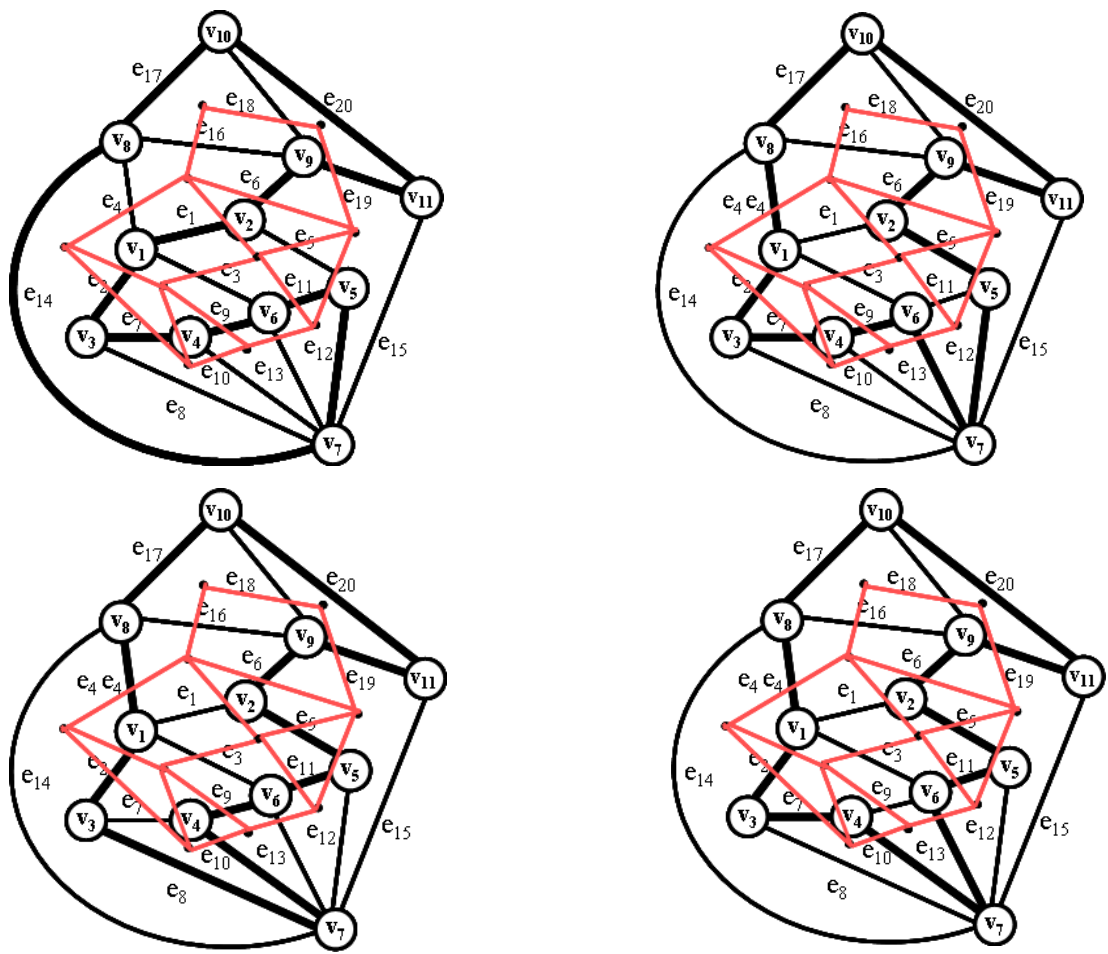

Рис. 14.38. Различные гамильтоновы циклы плоского суграфа (вариант 6).

**Шаг 1. [Построение графа циклов H].** Ставим в соответствие каждому изометрическому циклу суграфа $G'(V, E')$ вершину графа циклов H, а соприкасающиеся ребра циклов в суграфе будем считать ребрами графа H. Определяем валентность вершин графа циклов H.

**Шаг 2. [Отбраковка суграфов].** Если суграф циклов не несепарабельный, то исключаем суграф из рассмотрения.

**Шаг 3. [Выбор цикла].** Исключаем из системы изометрических циклов цикл, имеющий максимальную валентность в графе циклов. Выбранный цикл считаем опорным и помечаем все циклы, соприкасающиеся с опорным циклом.

**Шаг 4. [Выбор помеченного цикла].** Выбираем среди помеченных циклов цикл с максимальной валентностью вершины в графе циклов Н и удаляем его из системы изометрических циклов, затем объявляем его опорным циклом. Если удаление цикла приводит к образованию графа циклов с мостом, то производится выбор другого помеченного цикла. Если среди помеченных циклов не существует цикла для выбора, то есть все циклы приводят к удалению вершин из суграфа, идем на шаг 6.

**Шаг 5. [Построение обруча].** Строим обруч для оставшейся системы изометрических циклов. Идем на шаг 4.

**Шаг 6. [Образование дерева в графе циклов Н].** Если в процессе исключения циклов не образовалось дерево в графе Н, то идем на шаг 3.

**Шаг 7. [Оценка результата].** Если образовано дерево в графе циклов, то выделен гамильтовов цикл. Если в процессе исключения цикла получен суграф с несколькими компонентами связности, то гамильтонова цикла в суграфе не существует.

### 14.5. Структура гамильтонова цикла

Рассмотрим вопрос о составляющих циклов в гамильтоновом цикле. Определим состоит ли гамильтонов цикл только из изометрических циклов. Рассмотрим плоский суграф, представленный на рис. 14.39.

Количество вершин графа = 18
Количество ребер графа = 39
Количество единичных циклов = 30

Смежность графа:

вершина $v_1$: $v_2$ $v_3$ $v_5$ $v_{10}$ $v_{11}$ $v_{12}$ $v_{18}$
вершина $v_2$: $v_1$ $v_3$ $v_8$ $v_{10}$
вершина $v_3$: $v_1$ $v_2$ $v_4$ $v_5$
вершина $v_4$: $v_3$ $v_5$ $v_6$ $v_9$ $v_{10}$
вершина $v_5$: $v_1$ $v_3$ $v_4$ $v_{10}$ $v_{11}$ $v_{13}$ $v_{18}$
вершина $v_6$: $v_4$ $v_7$ $v_9$
вершина $v_7$: $v_6$ $v_8$ $v_9$
вершина $v_8$: $v_2$ $v_7$ $v_9$ $v_{10}$
вершина $v_9$: $v_4$ $v_6$ $v_7$ $v_8$
вершина $v_{10}$: $v_1$ $v_2$ $v_4$ $v_5$ $v_8$
вершина $v_{11}$: $v_1$ $v_5$ $v_{12}$ $v_{13}$
вершина $v_{12}$: $v_1$ $v_{11}$ $v_{15}$ $v_{18}$
вершина $v_{13}$: $v_5$ $v_{11}$ $v_{14}$ $v_{17}$ $v_{18}$
вершина $v_{14}$: $v_{13}$ $v_{16}$ $v_{17}$
вершина $v_{15}$: $v_{12}$ $v_{16}$ $v_{17}$ $v_{18}$
вершина $v_{16}$: $v_{14}$ $v_{15}$ $v_{17}$
вершина $v_{17}$: $v_{13}$ $v_{14}$ $v_{15}$ $v_{16}$

вершина $v_{18}$:  $v_1$  $v_5$  $v_{12}$  $v_{13}$  $v_{15}$

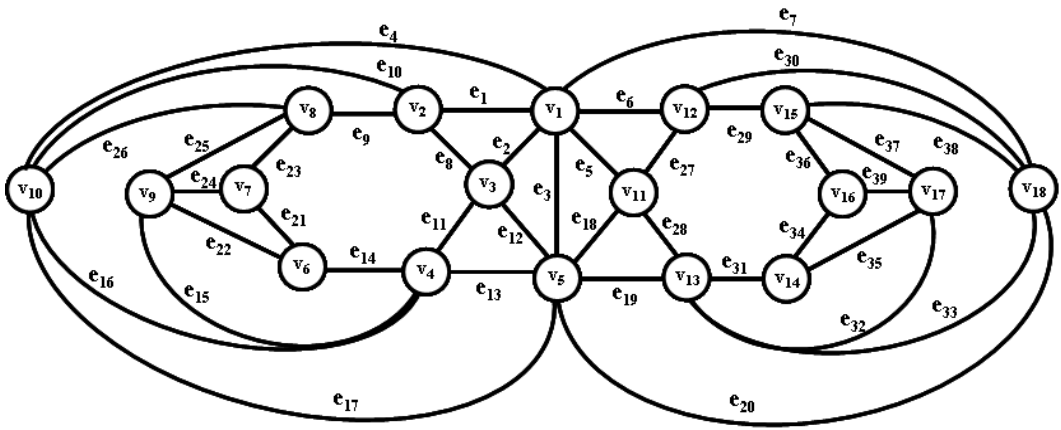

Рис. 14.39. Плоский суграф.

Инцидентность графа:

ребро $e_1$: $(v_1,v_2)$ или $(v_2,v_1)$;   ребро $e_2$: $(v_1,v_3)$ или $(v_3,v_1)$;
ребро $e_3$: $(v_1,v_5)$ или $(v_5,v_1)$;   ребро $e_4$: $(v_1,v_{10})$ или $(v_{10},v_1)$;
ребро $e_5$: $(v_1,v_{11})$ или $(v_{11},v_1)$;   ребро $e_6$: $(v_1,v_{12})$ или $(v_{12},v_1)$;
ребро $e_7$: $(v_1,v_{18})$ или $(v_{18},v_1)$;   ребро $e_8$: $(v_2,v_3)$ или $(v_3,v_2)$;
ребро $e_9$: $(v_2,v_8)$ или $(v_8,v_2)$;   ребро $e_{10}$: $(v_2,v_{10})$ или $(v_{10},v_2)$;
ребро $e_{11}$: $(v_3,v_4)$ или $(v_4,v_3)$;   ребро $e_{12}$: $(v_3,v_5)$ или $(v_5,v_3)$;
ребро $e_{13}$: $(v_4,v_5)$ или $(v_5,v_4)$;   ребро $e_{14}$: $(v_4,v_6)$ или $(v_6,v_4)$;
ребро $e_{15}$: $(v_4,v_9)$ или $(v_9,v_4)$;   ребро $e_{16}$: $(v_4,v_{10})$ или $(v_{10},v_4)$;
ребро $e_{17}$: $(v_5,v_{10})$ или $(v_{10},v_5)$;   ребро $e_{18}$: $(v_5,v_{11})$ или $(v_{11},v_5)$;
ребро $e_{19}$: $(v_5,v_{13})$ или $(v_{13},v_5)$;   ребро $e_{20}$: $(v_5,v_{18})$ или $(v_{18},v_5)$;
ребро $e_{21}$: $(v_6,v_7)$ или $(v_7,v_6)$;   ребро $e_{22}$: $(v_6,v_9)$ или $(v_9,v_6)$;
ребро $e_{23}$: $(v_7,v_8)$ или $(v_8,v_7)$;   ребро $e_{24}$: $(v_7,v_9)$ или $(v_9,v_7)$;
ребро $e_{25}$: $(v_8,v_9)$ или $(v_9,v_8)$;   ребро $e_{26}$: $(v_8,v_{10})$ или $(v_{10},v_8)$;
ребро $e_{27}$: $(v_{11},v_{12})$ или $(v_{12},v_{11})$;   ребро $e_{28}$: $(v_{11},v_{13})$ или $(v_{13},v_{11})$;
ребро $e_{29}$: $(v_{12},v_{15})$ или $(v_{15},v_{12})$;   ребро $e_{30}$: $(v_{12},v_{18})$ или $(v_{18},v_{12})$;
ребро $e_{31}$: $(v_{13},v_{14})$ или $(v_{14},v_{13})$;   ребро $e_{32}$: $(v_{13},v_{17})$ или $(v_{17},v_{13})$;
ребро $e_{33}$: $(v_{13},v_{18})$ или $(v_{18},v_{13})$;   ребро $e_{34}$: $(v_{14},v_{16})$ или $(v_{16},v_{14})$;
ребро $e_{35}$: $(v_{14},v_{17})$ или $(v_{17},v_{14})$;   ребро $e_{36}$: $(v_{15},v_{16})$ или $(v_{16},v_{15})$;
ребро $e_{37}$: $(v_{15},v_{17})$ или $(v_{17},v_{15})$;   ребро $e_{38}$: $(v_{15},v_{18})$ или $(v_{18},v_{15})$;
ребро $e_{39}$: $(v_{16},v_{17})$ или $(v_{17},v_{16})$.

Множество изометрических циклов графа:

$c_1 = \{e_1,e_2,e_8\} \leftrightarrow \{v_1,v_2,v_3\}$;
$c_2 = \{e_1,e_4,e_{10}\} \leftrightarrow \{v_1,v_2,v_{10}\}$;
$c_3 = \{e_2,e_3,e_{12}\} \leftrightarrow \{v_1,v_3,v_5\}$;
$c_5 = \{e_3,e_4,e_{17}\} \leftrightarrow \{v_1,v_5,v_{10}\}$;
$c_6 = \{e_3,e_5,e_{18}\} \leftrightarrow \{v_1,v_5,v_{11}\}$;
$c_7 = \{e_3,e_7,e_{20}\} \leftrightarrow \{v_1,v_5,v_{18}\}$;
$c_8 = \{e_5,e_6,e_{27}\} \leftrightarrow \{v_1,v_{11},v_{12}\}$;
$c_{10} = \{e_6,e_7,e_{30}\} \leftrightarrow \{v_1,v_{12},v_{18}\}$;
$c_{14} = \{e_9,e_{10},e_{26}\} \leftrightarrow \{v_2,v_8,v_{10}\}$;
$c_{15} = \{e_{11},e_{12},e_{13}\} \leftrightarrow \{v_3,v_4,v_5\}$;
$c_{16} = \{e_{13},e_{16},e_{17}\} \leftrightarrow \{v_4,v_5,v_{10}\}$;
$c_{17} = \{e_{14},e_{15},e_{22}\} \leftrightarrow \{v_4,v_6,v_9\}$;

$c_{18} = \{e_{15},e_{16},e_{25},e_{26}\} \leftrightarrow \{v_4,v_8,v_9,v_{10}\}$;
$c_{19} = \{e_{18},e_{19},e_{28}\} \leftrightarrow \{v_5,v_{11},v_{13}\}$;
$c_{21} = \{e_{19},e_{20},e_{33}\} \leftrightarrow \{v_5,v_{13},v_{18}\}$;
$c_{22} = \{e_{21},e_{22},e_{24}\} \leftrightarrow \{v_6,v_7,v_9\}$;
$c_{23} = \{e_{23},e_{24},e_{25}\} \leftrightarrow \{v_7,v_8,v_9\}$;
$c_{26} = \{e_{29},e_{30},e_{38}\} \leftrightarrow \{v_{12},v_{15},v_{18}\}$;
$c_{27} = \{e_{31},e_{32},e_{35}\} \leftrightarrow \{v_{13},v_{14},v_{17}\}$;
$c_{28} = \{e_{32},e_{33},e_{37},e_{38}\} \leftrightarrow \{v_{13},v_{15},v_{17},v_{18}\}$;
$c_{29} = \{e_{34},e_{35},e_{39}\} \leftrightarrow \{v_{14},v_{16},v_{17}\}$;
$c_{30} = \{e_{36},e_{37},e_{39}\} \leftrightarrow \{v_{15},v_{16},v_{17}\}$.

Топологический рисунок плоской части графа состоит из следующих изометрических циклов:

$c_1 = \{e_1,e_2,e_8\} \leftrightarrow \{v_1,v_2,v_3\}$;
$c_2 = \{e_1,e_4,e_{10}\} \leftrightarrow \{v_1,v_2,v_{10}\}$;
$c_3 = \{e_2,e_3,e_{12}\} \leftrightarrow \{v_1,v_3,v_5\}$;
$c_4 = \{e_2,e_4,e_{11},e_{16}\} \leftrightarrow \{v_1,v_3,v_4,v_{10}\}$;
$c_6 = \{e_3,e_5,e_{18}\} \leftrightarrow \{v_1,v_5,v_{11}\}$;
$c_8 = \{e_5,e_6,e_{27}\} \leftrightarrow \{v_1,v_{11},v_{12}\}$;
$c_9 = \{e_5,e_7,e_{28},e_{33}\} \leftrightarrow \{v_1,v_{11},v_{13},v_{18}\}$;
$c_{10} = \{e_6,e_7,e_{30}\} \leftrightarrow \{v_1,v_{12},v_{18}\}$;
$c_{11} = \{e_8,e_{10},e_{11},e_{16}\} \leftrightarrow \{v_2,v_3,v_4,v_{10}\}$;
$c_{12} = \{e_8,e_{10},e_{12},e_{17}\} \leftrightarrow \{v_2,v_3,v_5,v_{10}\}$;
$c_{13} = \{e_8,e_9,e_{11},e_{15},e_{25}\} \leftrightarrow \{v_2,v_3,v_4,v_8,v_9\}$;
$c_{14} = \{e_9,e_{10},e_{26}\} \leftrightarrow \{v_2,v_8,v_{10}\}$;
$c_{15} = \{e_{11},e_{12},e_{13}\} \leftrightarrow \{v_3,v_4,v_5\}$;
$c_{16} = \{e_{13},e_{16},e_{17}\} \leftrightarrow \{v_4,v_5,v_{10}\}$;
$c_{17} = \{e_{14},e_{15},e_{22}\} \leftrightarrow \{v_4,v_6,v_9\}$;
$c_{18} = \{e_{15},e_{16},e_{25},e_{26}\} \leftrightarrow \{v_4,v_8,v_9,v_{10}\}$;
$c_{19} = \{e_{18},e_{19},e_{28}\} \leftrightarrow \{v_5,v_{11},v_{13}\}$;
$c_{20} = \{e_{18},e_{20},e_{27},e_{30}\} \leftrightarrow \{v_5,v_{11},v_{12},v_{18}\}$;
$c_{21} = \{e_{19},e_{20},e_{33}\} \leftrightarrow \{v_5,v_{13},v_{18}\}$;
$c_{22} = \{e_{21},e_{22},e_{24}\} \leftrightarrow \{v_6,v_7,v_9\}$;
$c_{23} = \{e_{23},e_{24},e_{25}\} \leftrightarrow \{v_7,v_8,v_9\}$;
$c_{24} = \{e_{27},e_{28},e_{30},e_{33}\} \leftrightarrow \{v_{11},v_{12},v_{13},v_{18}\}$;
$c_{25} = \{e_{27},e_{28},e_{29},e_{32},e_{37}\} \leftrightarrow \{v_{11},v_{12},v_{13},v_{15},v_{17}\}$;
$c_{26} = \{e_{29},e_{30},e_{38}\} \leftrightarrow \{v_{12},v_{15},v_{18}\}$;
$c_{27} = \{e_{31},e_{32},e_{35}\} \leftrightarrow \{v_{13},v_{14},v_{17}\}$;
$c_{28} = \{e_{32},e_{33},e_{37},e_{38}\} \leftrightarrow \{v_{13},v_{15},v_{17},v_{18}\}$;
$c_{29} = \{e_{34},e_{35},e_{39}\} \leftrightarrow \{v_{14},v_{16},v_{17}\}$;
$c_{30} = \{e_{36},e_{37},e_{39}\} \leftrightarrow \{v_{15},v_{16},v_{17}\}$.

Номер  <1,2,3,4,5,6,7,8,9,0,1,2,3,4,5,6,7,8,9,0,1,2,3,4,5,6,7,8,9,0.1,2,3,4,5,6,7,8,9>.
Вектор  $P_e$ = <2,2,2,1,2,2,1,1,1,2,1,2,2,1,2,2,1,2,2,1,1,2,1,2,2,2,1,1,1,2,1,2,2,1,2,1,2,2,2>

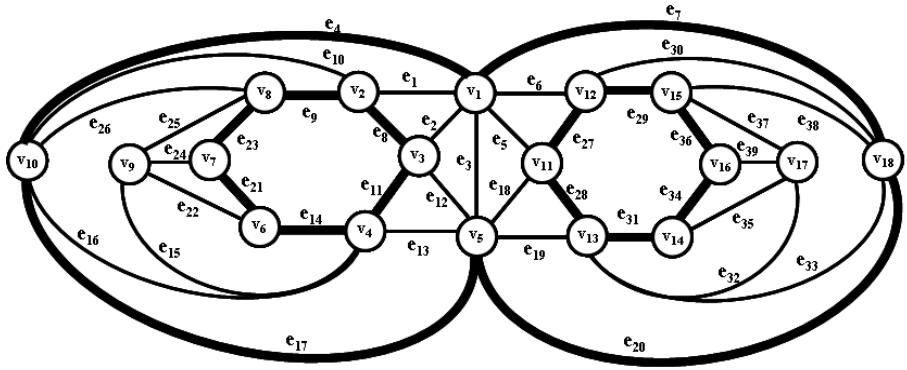

Рис. 14.40. Обод суграфа.

Заметим, что в обод суграфа входят два простых цикла:

$c_{31} = \{e_8, e_9, e_{11}, e_{14}, e_{21}, e_{23}\} \leftrightarrow \{v_2, v_3, v_4, v_6, v_7, v_8\}$;
$c_{32} = \{e_{27}, e_{28}, e_{29}, e_{31}, e_{34}, e_{36}\} \leftrightarrow \{v_{11}, v_{12}, v_{13}, v_{14}, v_{15}, v_{16}\}$.

Циклы $c_{31}$ и $c_{32}$ можно включить в состав топологического рисунка графа.

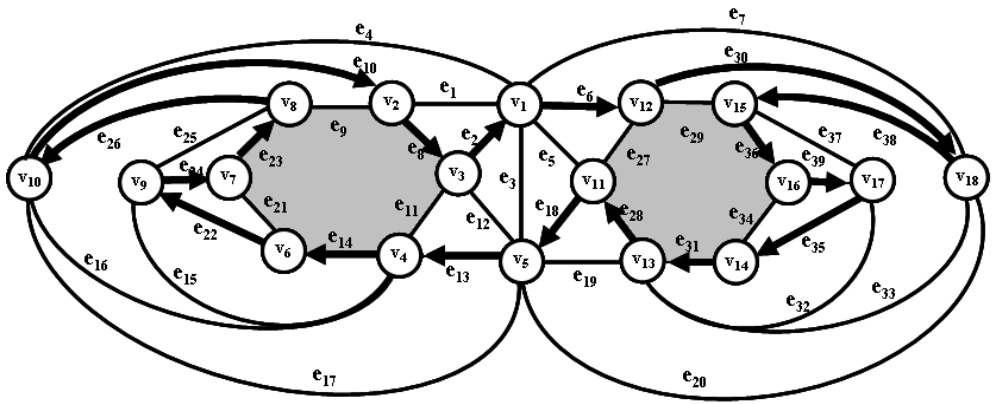

Рис. 14.41. Гамильтонов цикл.

Как следует из построения, в гамильтонов цикл наряду с изометрическими циклами могут входить и простые циклы (рис. 14.41).

**Комментарии**

Будем исходить из предположения, что гамильтонов цикл состоит из простых циклов, принадлежащих топологическому рисунку плоского суграфа.

В данной главе рассмотрены различные случаи выделения гамильтонова цикла для несеапарабельного плоского суграфа. Показано, что граф циклов для гамильтоного цикла представляет собой дерево. Перспективно использовать максимально плоский суграф для построения гамильтонова цикла. Приведен алгоритм для выделения гамильтонова цикла в несепарабельном плоском суграфе.

# Глава 15. НЕГАМИЛЬТОНОВЫ ГРАФЫ
## 15.1. Двухсвязные графы

Рассмотрим построение гамильтонова цикла в двусвязном графе.

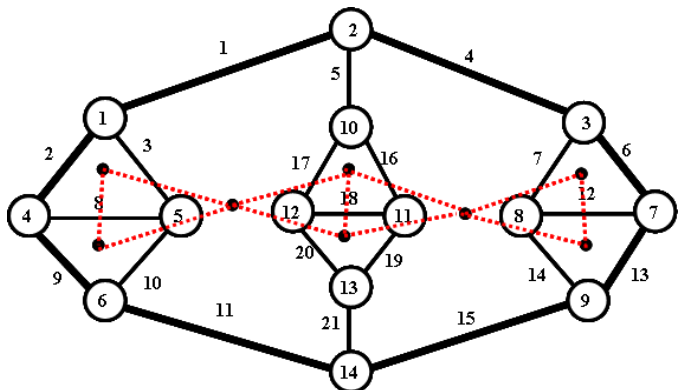

Рис. 15.1. Двухсвязный граф.

Строим граф циклов H. Выбираем цикл с максимальной валентностью вершины c = {$e_4,e_5,e_7,e_{14},e_{15},e_{16},e_{19},e_{21}$} и строим обод {$e_1,e_2,e_5,e_6,e_7,e_9,e_{11},e_{13},e_{14},e_{16},e_{19},e_{21}$} (рис. 15.2). Образуется две части обода.

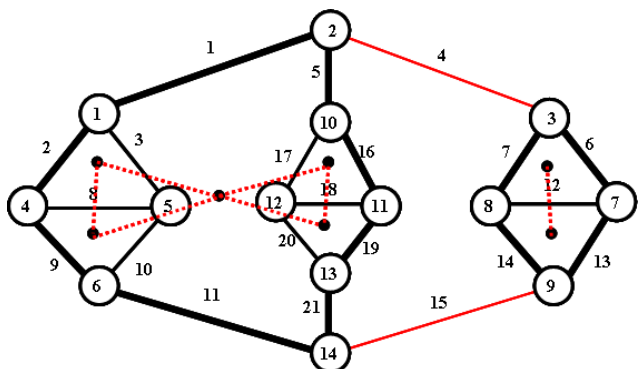

Рис. 15.2. Удаление цикла {$e_4,e_5,e_7,e_{14},e_{15},e_{16},e_{19},e_{21}$}

Удаляем цикл {$e_1,e_3,e_8$}, строим обруч {$e_1,e_3,e_5,e_6,e_7,e_8,e_9,e_{11},e_{13},e_{14},e_{16},e_{19},e_{21}$} (рис. 15.3).

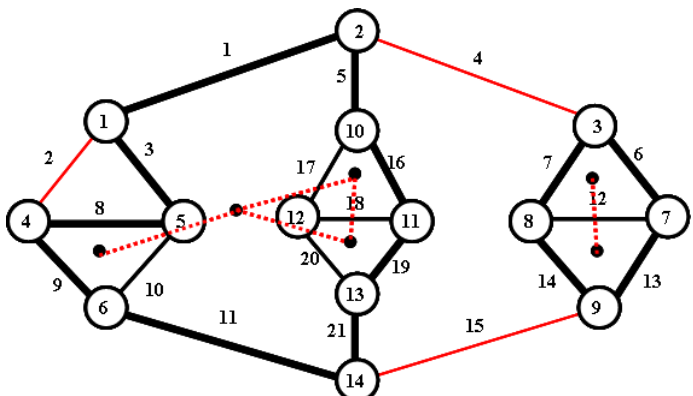

Рис. 15.3. Удаление цикла {$e_1,e_3,e_8$}.

Удаляем цикл {$e_{16},e_{17},e_{18}$} и получаем гамильтонов квазицикл (рис. 15.4).

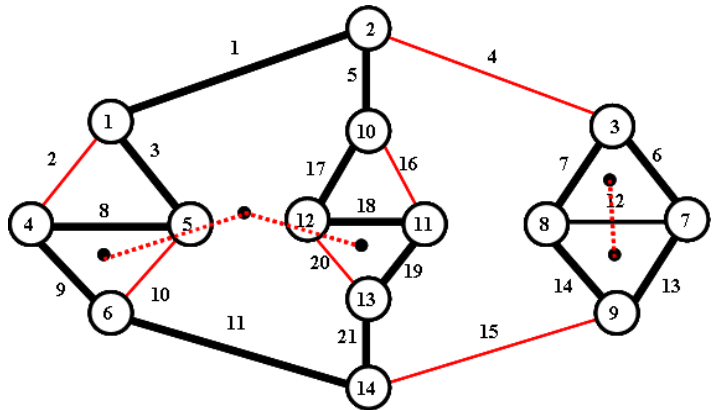

Рис. 15.4. Гамильтонов квазицикл.

### 15.2. Граф Татта

Рассмотрим построение гамильноногового цикла для известного графа Татта. Выделяем обод графа $\{e_1, e_3, e_4, e_7, e_{28}, e_{39}, e_{48}, e_{50}, e_{62}\}$ (рис. 15.5).

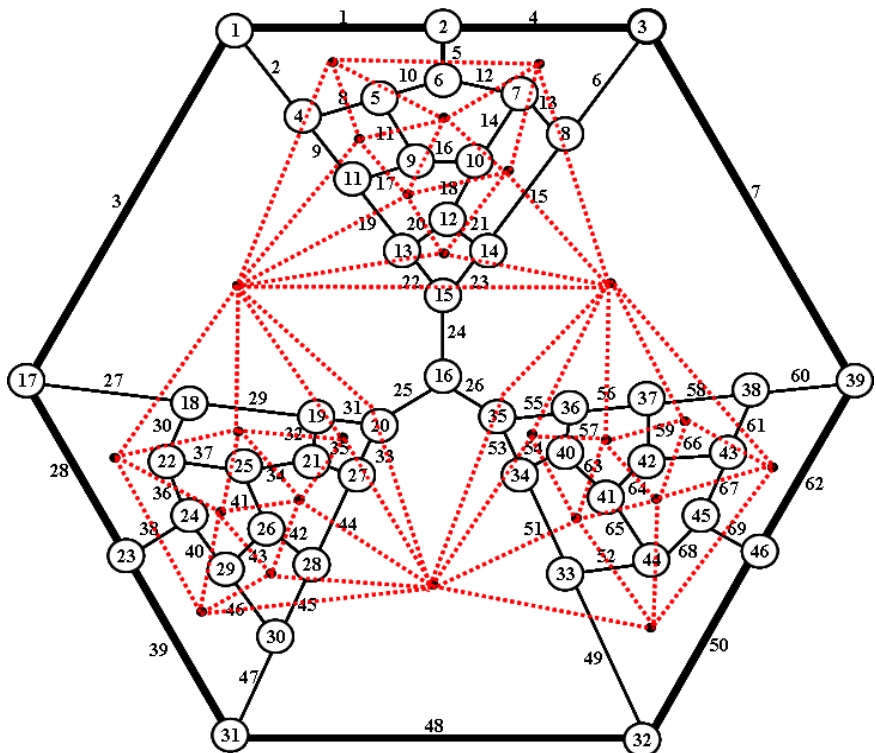

Рис. 15.5. Граф Татта и его обод.

Удаляем цикл, имеющий максимальную валентность вершины в графе циклов H – $\{e_6, e_7, e_{15}, e_{23}, e_{24}, e_{26}, e_{55}, e_{56}, e_{58}, e_{60}\}$ и строим обод для плоского суграфа (рис. 15.6).

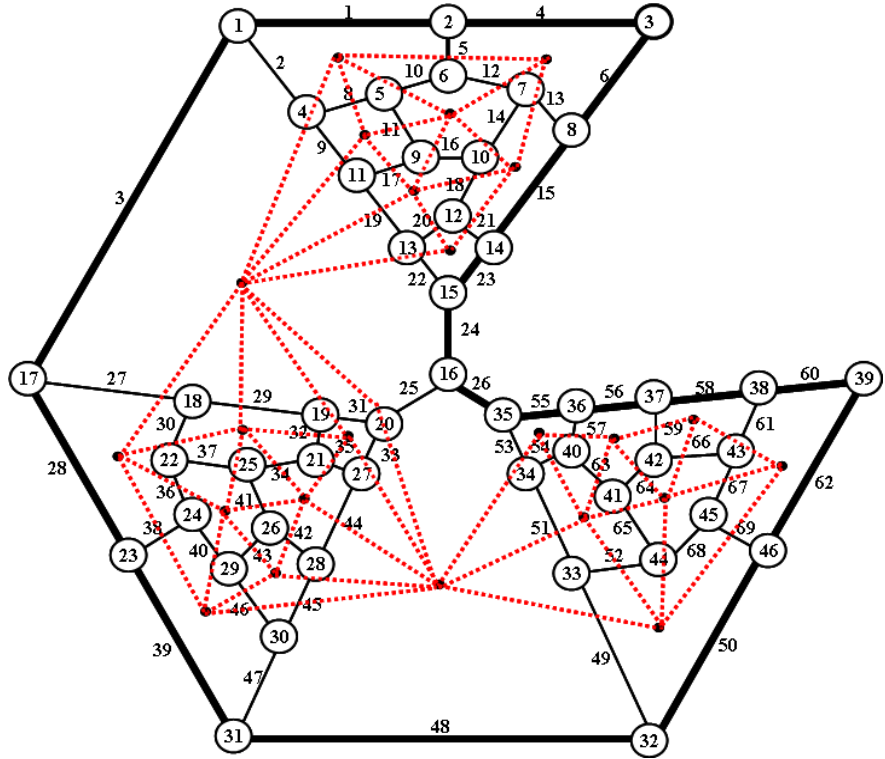

Рис. 15.6. Обод плоского суграфа.

Вновь удаляем цикл, имеющий максимальную валентность вершины в графе циклов H – $\{e_{25}, e_{26}, e_{33}, e_{44}, e_{45}, e_{46}, e_{47}, e_{48}, e_{49}, e_{51}, e_{53}\}$ и строим обод для плоского суграфа (рис. 15.7). В результате получаем обод, состоящий из двух замкнутых маршрутов.

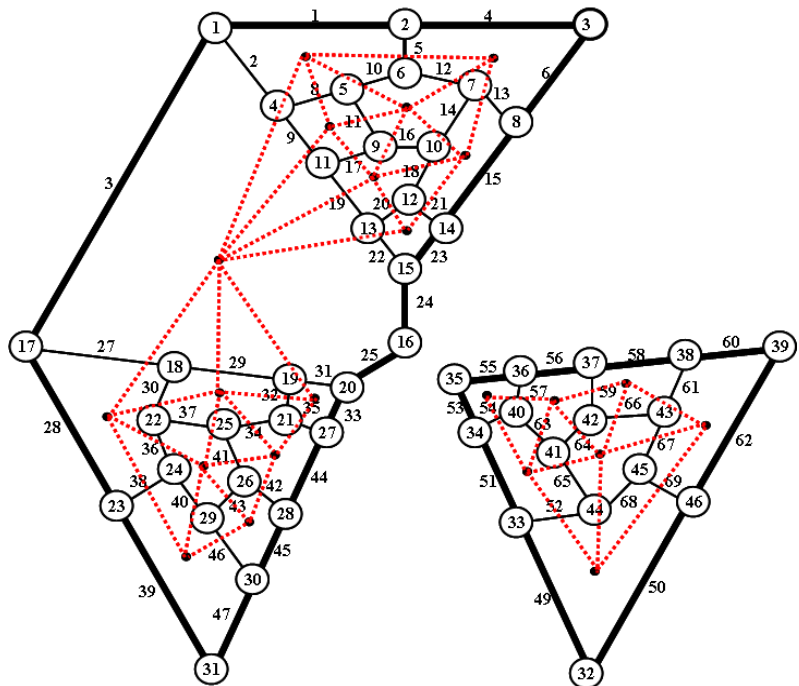

Рис. 15.7. Топологический рисунок плоского суграфа.

Удаление цикла $\{e_2, e_3, e_9, e_{19}, e_{22}, e_{24}, e_{25}, e_{27}, e_{29}, e_{31}\}$ приводит к удалению вершины $v_{16}$. Поэтому удаляем цикл $\{e_{56}, e_{57}, e_{59}, e_{63}, e_{64}\}$, имеющий соответствующую вершину с максимальной валентностью в графе циклов H. Строим обод суграфа (рис. 15.8).

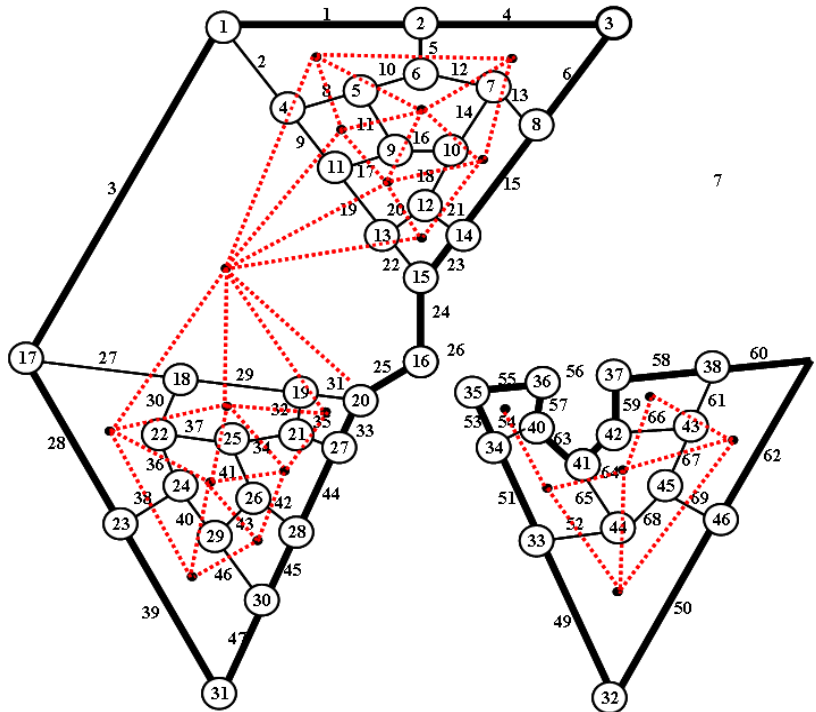

Рис. 15.8. Топологический рисунок плоского суграфа.

Удаляем цикл $\{e_{64}, e_{65}, e_{66}, e_{67}, e_{68}\}$, имеющий соответствующую вершину с максимальной валентностью в графе циклов H. Строим обод суграфа (рис. 15.9).

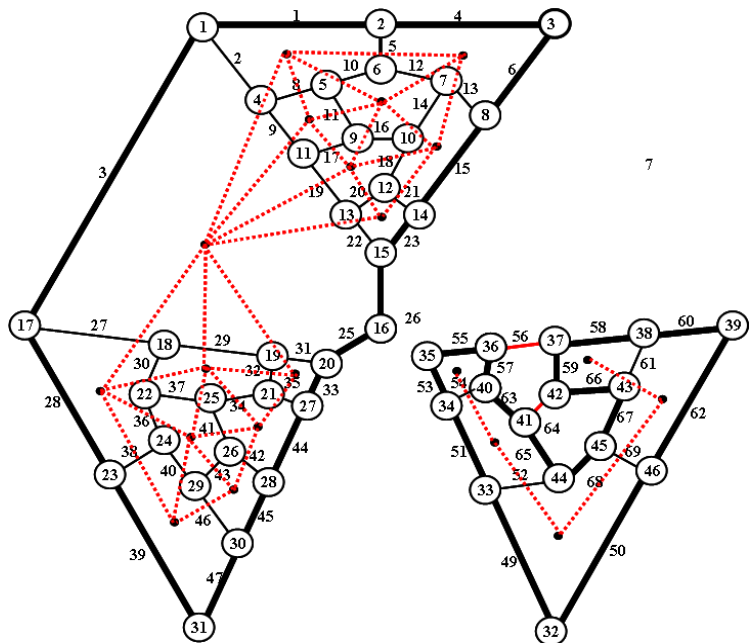

Рис. 15.9. Топологический рисунок плоского суграфа.

Удаляем цикл $\{e_{13}, e_{14}, e_{15}, e_{18}, e_{21}\}$, имеющий соответствующую вершину с максимальной валентностью в графе циклов H. Строим обод суграфа (рис. 15.10).

Удаляем цикл {$e_{10}, e_{11}, e_{12}, e_{14}, e_{16}$}, имеющий соответствующую вершину с максимальной валентностью в графе циклов H. Строим обод суграфа (рис. 15.11).

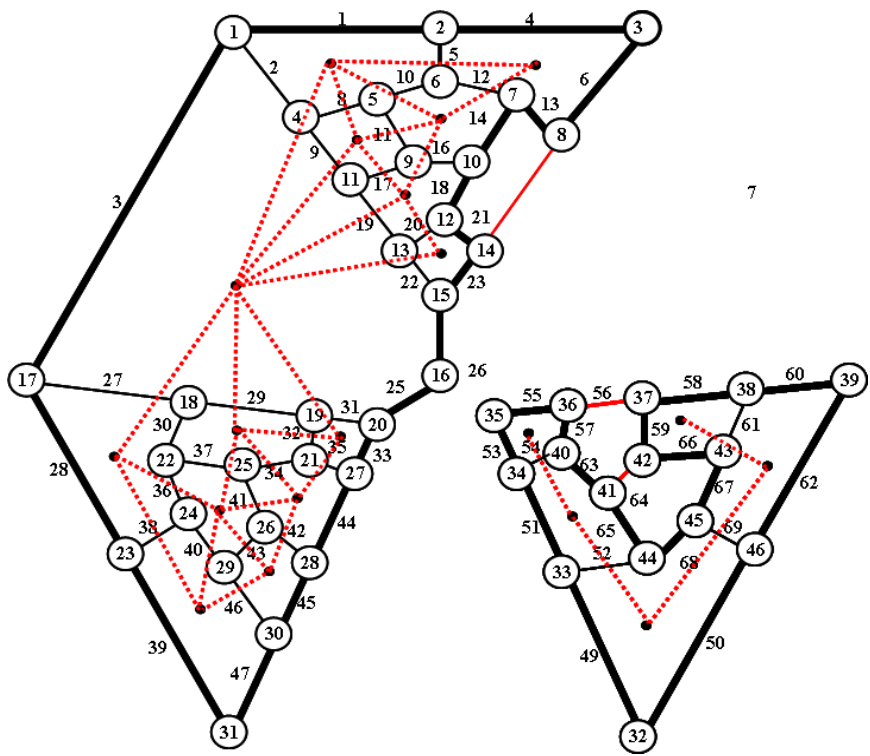

Рис. 15.10. Топологический рисунок плоского суграфа.

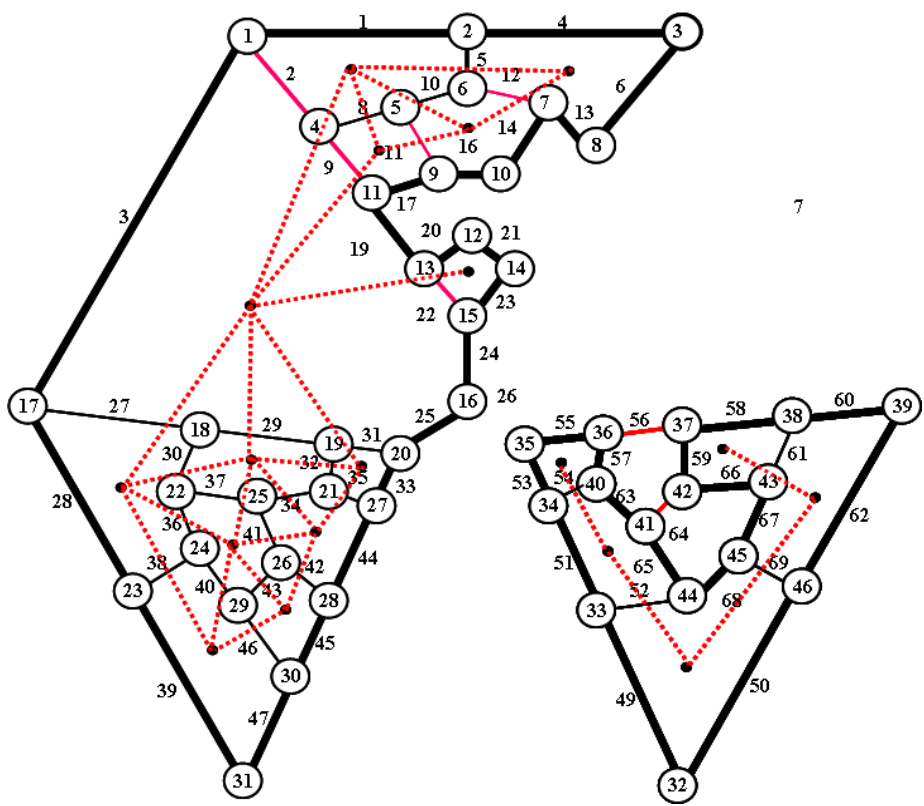

Рис. 15.11. Топологический рисунок плоского суграфа.

Удаляем цикл {$e_1, e_2, e_5, e_8, e_{10}$}, имеющий соответствующую вершину с максимальной валентностью в графе циклов H. Строим обод суграфа (рис. 15.12).

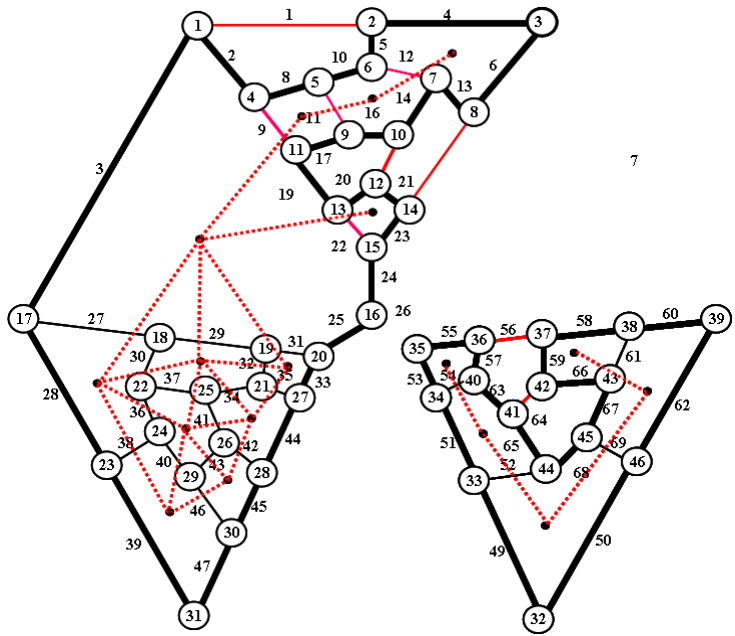

Рис. 15.12. Топологический рисунок плоского суграфа.

Удаляем цикл {$e_{27}, e_{28}, e_{30}, e_{36}, e_{38}$}, имеющий соответствующую вершину с максимальной валентностью в графе циклов H. Строим обод суграфа (рис. 15.13).

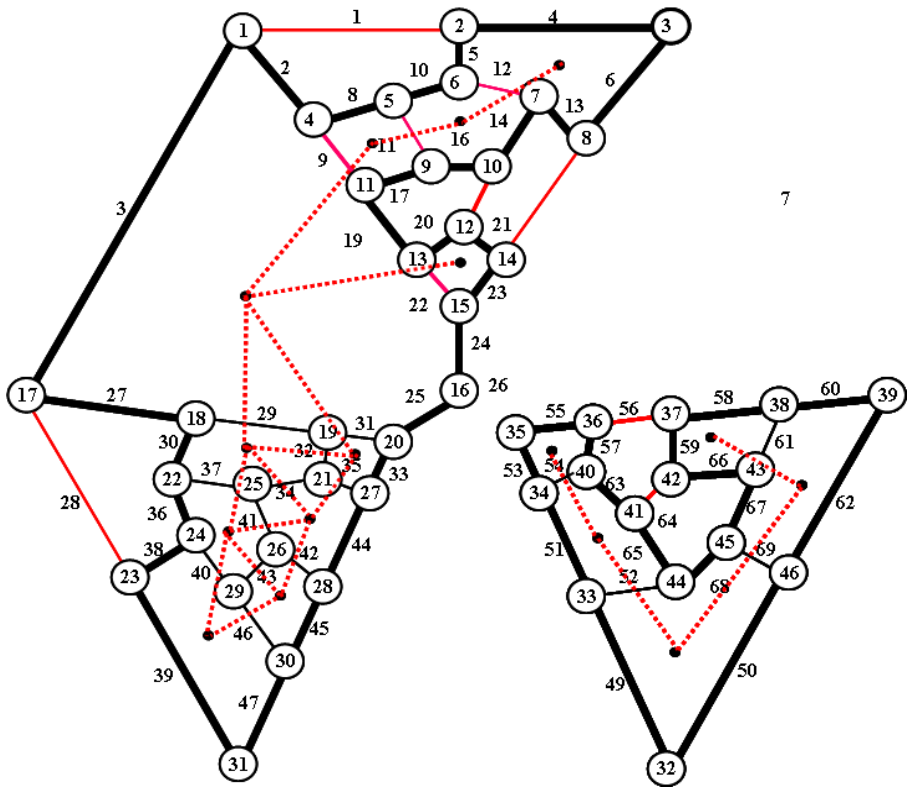

Рис. 15.13. Топологический рисунок плоского суграфа.

Удаляем цикл {$e_{36}, e_{37}, e_{40}, e_{41}, e_{43}$}, имеющий соответствующую вершину с максимальной валентностью в графе циклов H. Строим обод суграфа (рис. 15.14).

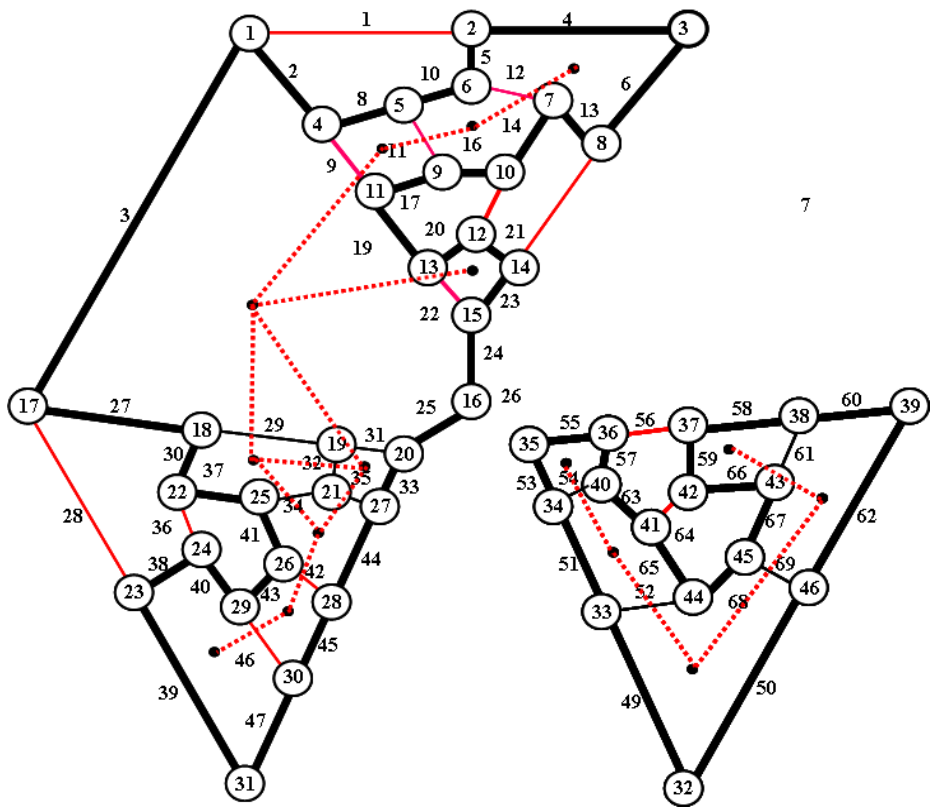

Рис. 15.14. Топологический рисунок плоского суграфа.

Удаляем цикл $\{e_{31}, e_{32}, e_{33}, e_{35}\}$, имеющий соответствующую вершину с максимальной валентностью в графе циклов H. Строим обод суграфа (рис. 15.15).

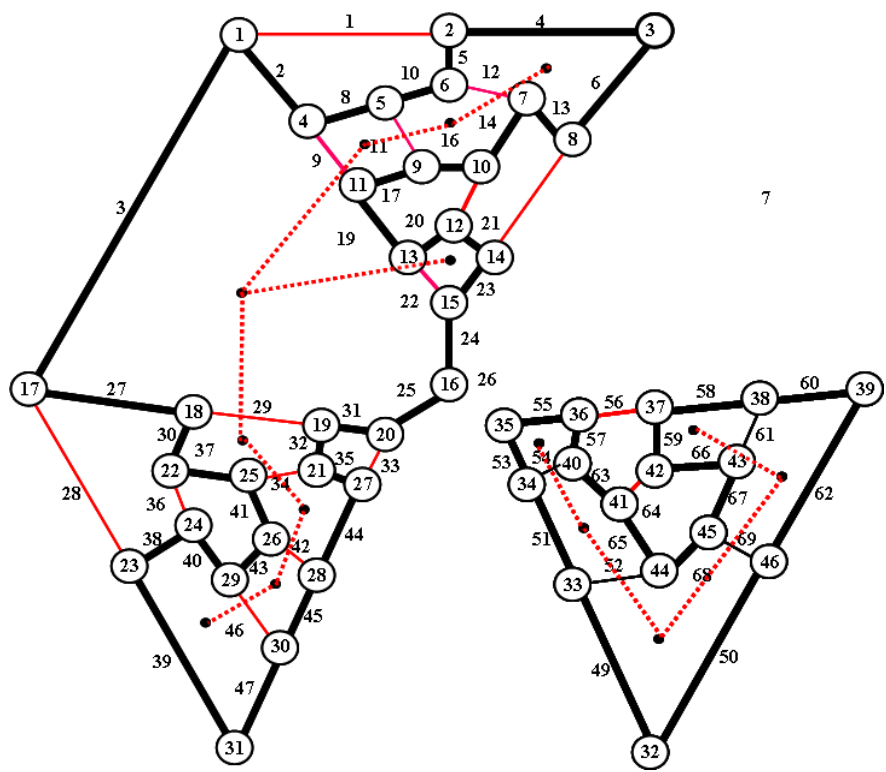

Рис. 15.15. Топологический рисунок плоского суграфа.

Для получения гамильтонова квазицикла удаляем избыточные ребра и получаем два диска нечетной длины (рис. 15.16).

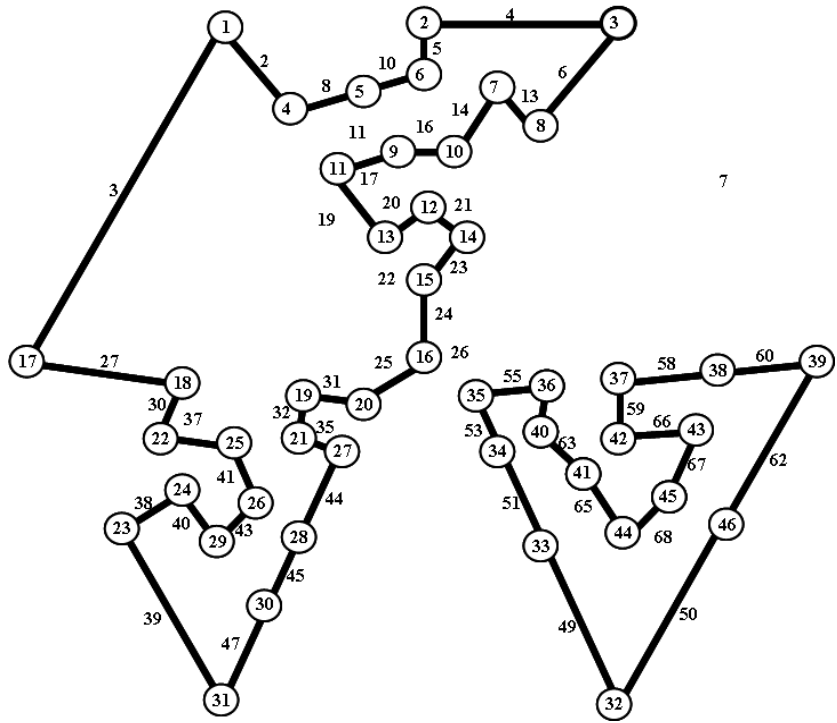

Рис. 15.16. Гамильтонов квазицикл с дисками нечетной длины.

Для получения квазицикла с дисками четной длины поступим следующим образом. Выделим цикл четной длины и будем считать его диском (рис. 3.17).

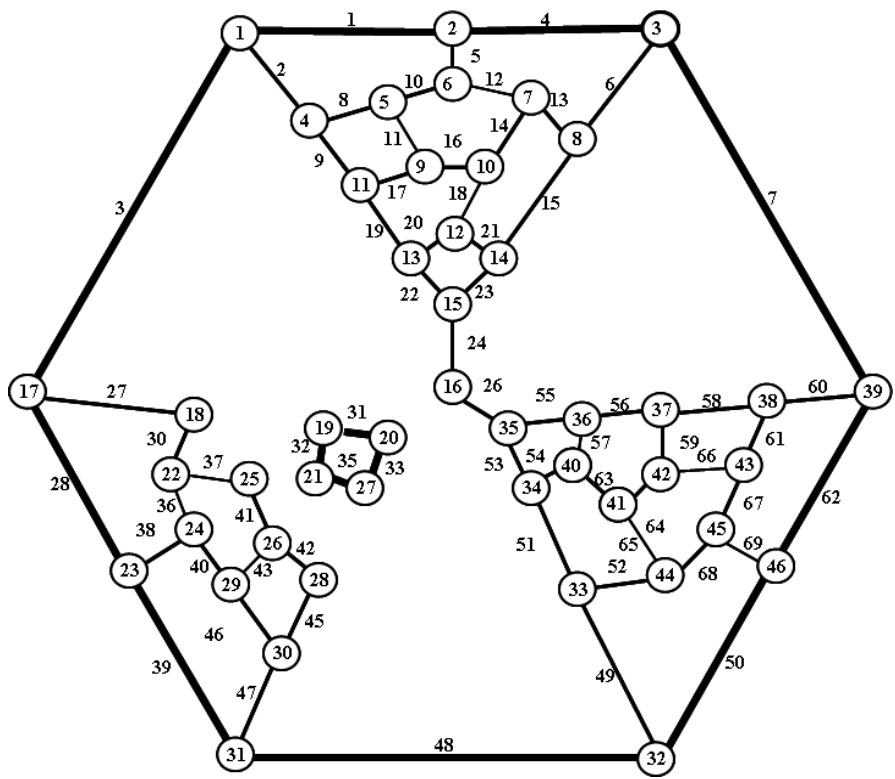

Рис. 15.17. Диск четной длины.

В оставшемся суграфе выделим гамильнов цикл четной длины. В результате получим гамильтонов квазицикл с дисками четной длины (рис. 3.18).

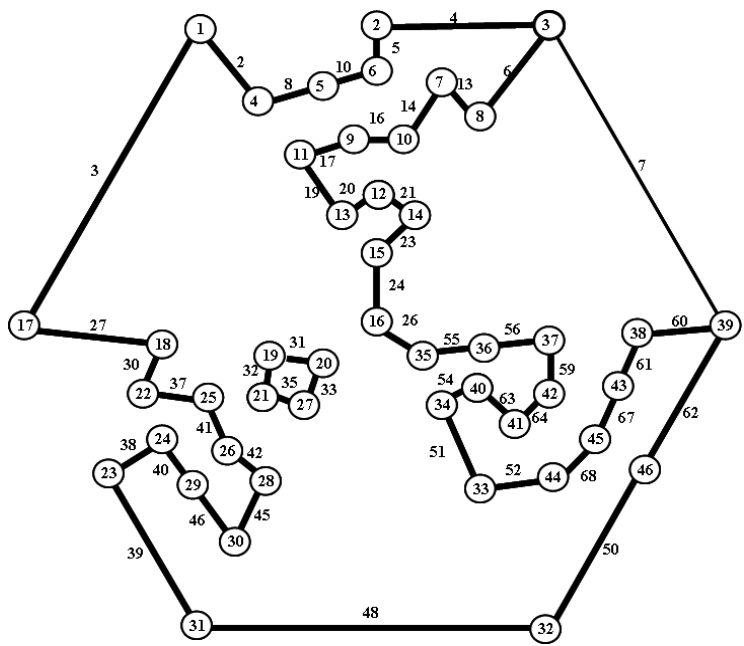

Рис. 15.18. Гамильтонов квазицикл с дисками четной длины.

### 15.3. Граф Петерсена

Рассмотрим граф Петерсена (рис. 15.19).

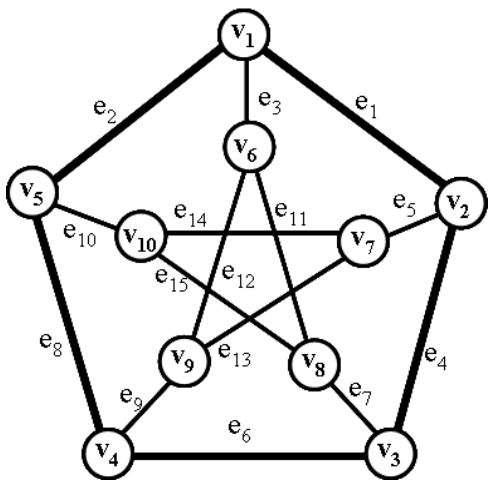

Рис. 15.19. Граф Петерсена и его обод.

Множество изометрических циклов графа:

$c_1 = \{e_1, e_2, e_4, e_6, e_8\} \leftrightarrow \{v_1, v_2, v_3, v_4, v_5\}$;
$c_2 = \{e_1, e_2, e_5, e_{10}, e_{14}\} \leftrightarrow \{v_1, v_2, v_5, v_7, v_{10}\}$;
$c_3 = \{e_1, e_3, e_4, e_7, e_{11}\} \leftrightarrow \{v_1, v_2, v_3, v_6, v_8$
$c_4 = \{e_1, e_3, e_5, e_{12}, e_{13}\} \leftrightarrow \{v_1, v_2, v_6, v_7, v_9\}$;
$c_5 = \{e_2, e_3, e_{10}, e_{11}, e_{15}\} \leftrightarrow \{v_1, v_5, v_6, v_8, v_{10}\}$;
$c_6 = \{e_2, e_3, e_8, e_9, e_{12}\} \leftrightarrow \{v_1, v_4, v_5, v_6, v_9\}$;
$c_7 = \{e_4, e_5, e_6, e_9, e_{13}\} \leftrightarrow \{v_2, v_3, v_4, v_7, v_9\}$;
$c_8 = \{e_4, e_5, e_7, e_{14}, e_{15}\} \leftrightarrow \{v_2, v_3, v_7, v_8, v_{10}\}$;
$c_9 = \{e_6, e_7, e_9, e_{11}, e_{12}\} \leftrightarrow \{v_3, v_4, v_6, v_8, v_9\}$;
$c_{10} = \{e_6, e_7, e_8, e_{10}, e_{15}\} \leftrightarrow \{v_3, v_4, v_5, v_8, v_{10}\}$;

$c_{11} = \{e_8,e_9,e_{10},e_{13},e_{14}\} \leftrightarrow \{v_4,v_5,v_7,v_9,v_{10}\}$;
$c_{12} = \{e_{11},e_{12},e_{13},e_{14},e_{15}\} \leftrightarrow \{v_6,v_7,v_8,v_9,v_{10}\}$.

Вектор валентности: $S(G) = <3,3,3,3,3,3,3,3,3,3>$.

Вектор $W(G) = <4,4,4,4,4,4,4,4,4,4,4,4,4,4,4>$.

Вектор $P_e(G) = <4,4,4,4,4,4,4,4,4,4,4,4,4,4,4>$.

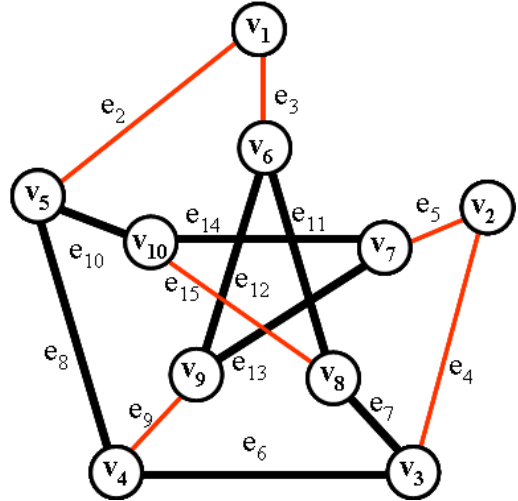

Рис. 15.20. Обод суграфа после удаления ребра $e_1$.

Исключим из системы цикл $c_1$ и получим обод графа (рис. 15.19).

Вектор валентности: $S(G') = <3,3,3,3,3,3,3,3,3,3>$.

Вектор $W(G') = <4,4,4,4,4,4,4,4,4,4,4,4,4,4,4>$.

Вектор $P_e(G') = <3,3,4,3,4,3,4,3,4,4,4,4,4,4,4>$.

Следующим шагом исключим из графа ребро $e_1$ и циклы, содержащие это ребро (рис. 15.20).

$c_5 = \{e_2,e_3,e_{10},e_{11},e_{15}\} \leftrightarrow \{v_1,v_5,v_6,v_8,v_{10}\}$;
$c_6 = \{e_2,e_3,e_8,e_9,e_{12}\} \leftrightarrow \{v_1,v_4,v_5,v_6,v_9\}$;
$c_7 = \{e_4,e_5,e_6,e_9,e_{13}\} \leftrightarrow \{v_2,v_3,v_4,v_7,v_9\}$;
$c_8 = \{e_4,e_5,e_7,e_{14},e_{15}\} \leftrightarrow \{v_2,v_3,v_7,v_8,v_{10}\}$;
$c_9 = \{e_6,e_7,e_9,e_{11},e_{12}\} \leftrightarrow \{v_3,v_4,v_6,v_8,v_9\}$;
$c_{10} = \{e_6,e_7,e_8,e_{10},e_{15}\} \leftrightarrow \{v_3,v_4,v_5,v_8,v_{10}\}$;
$c_{11} = \{e_8,e_9,e_{10},e_{13},e_{14}\} \leftrightarrow \{v_4,v_5,v_7,v_9,v_{10}\}$;
$c_{12} = \{e_{11},e_{12},e_{13},e_{14},e_{15}\} \leftrightarrow \{v_6,v_7,v_8,v_9,v_{10}\}$.

Вектор валентности: $S(G') = <2,2,3,3,3,3,3,3,3,3>$.

Вектор $W(G') = <0,3,3,3,3,4,4,4,4,4,4,4,4,4,4>$.

Вектор $P_e(G') = <0,2,2,2,2,3,3,3,4,3,3,3,3,3,4>$.

Для того чтобы присоединить вершины, не вошедшие в обод суграфа, можем для удаления выбрать ребра $e_9$ или $e_{15}$. Выбираем ребро $e_{15}$ (рис. 15.21).

$c_6 = \{e_2,e_3,e_8,e_9,e_{12}\} \leftrightarrow \{v_1,v_4,v_5,v_6,v_9\}$;
$c_7 = \{e_4,e_5,e_6,e_9,e_{13}\} \leftrightarrow \{v_2,v_3,v_4,v_7,v_9\}$;
$c_9 = \{e_6,e_7,e_9,e_{11},e_{12}\} \leftrightarrow \{v_3,v_4,v_6,v_8,v_9\}$;

$c_{11} = \{e_8, e_9, e_{10}, e_{13}, e_{14}\} \leftrightarrow \{v_4, v_5, v_7, v_9, v_{10}\}$;

Вектор валентности: $S(G') = <2,2,3,3,3,3,3,3,2,2>$.

Вектор $W(G') = <0,3,3,3,3,4,3,4,4,3,3,4,4,3,0>$.

Вектор $P_e(G') = <0,1,1,1,1,2,1,2,4,1,1,2,2,1,0>$.

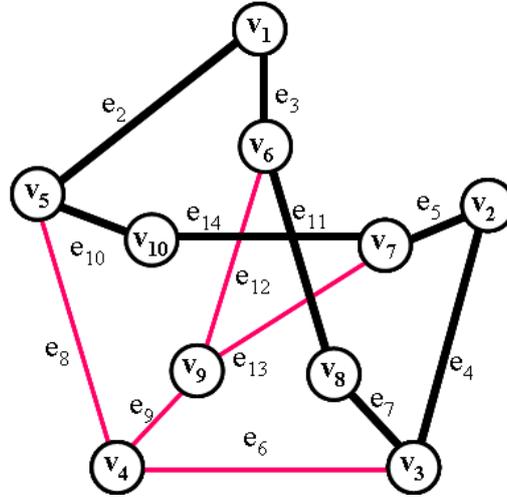

Рис. 15.21. Обод суграфа после удаления ребра $e_{15}$.

Как видим, удаление ребра $e_{15}$ не привело к подключению вершин к ободу. Следовательно мы не смогли построить плоскую часть графа методом удаления ребер. Отсюда следует вывод об отсутствии у графа гамильтонова цикла.

Рассмотрим способ построения плоской части графа методом удаления ребер.

$c_1 = \{e_1, e_2, e_4, e_6, e_8\} \leftrightarrow \{v_1, v_2, v_3, v_4, v_5\}$;
$c_2 = \{e_1, e_2, e_5, e_{10}, e_{14}\} \leftrightarrow \{v_1, v_2, v_5, v_7, v_{10}\}$;
$c_3 = \{e_1, e_3, e_4, e_7, e_{11}\} \leftrightarrow \{v_1, v_2, v_3, v_6, v_8$
$c_4 = \{e_1, e_3, e_5, e_{12}, e_{13}\} \leftrightarrow \{v_1, v_2, v_6, v_7, v_9\}$;
$c_5 = \{e_2, e_3, e_{10}, e_{11}, e_{15}\} \leftrightarrow \{v_1, v_5, v_6, v_8, v_{10}\}$;
$c_6 = \{e_2, e_3, e_8, e_9, e_{12}\} \leftrightarrow \{v_1, v_4, v_5, v_6, v_9\}$;
$c_7 = \{e_4, e_5, e_6, e_9, e_{13}\} \leftrightarrow \{v_2, v_3, v_4, v_7, v_9\}$;
$c_8 = \{e_4, e_5, e_7, e_{14}, e_{15}\} \leftrightarrow \{v_2, v_3, v_7, v_8, v_{10}\}$;
$c_9 = \{e_6, e_7, e_9, e_{11}, e_{12}\} \leftrightarrow \{v_3, v_4, v_6, v_8, v_9\}$;
$c_{10} = \{e_6, e_7, e_8, e_{10}, e_{15}\} \leftrightarrow \{v_3, v_4, v_5, v_8, v_{10}\}$;
$c_{11} = \{e_8, e_9, e_{10}, e_{13}, e_{14}\} \leftrightarrow \{v_4, v_5, v_7, v_9, v_{10}\}$;
$c_{12} = \{e_{11}, e_{12}, e_{13}, e_{14}, e_{15}\} \leftrightarrow \{v_6, v_7, v_8, v_9, v_{10}\}$.

Вектор $P_e(G) = <4,4,4,4,4,4,4,4,4,4,4,4,4,4,4>$.

Удаляем цикл $c_1$, вектор $P_e(G) = <3,3,4,3,4,3,4,3,4,4,4,4,4,4,4>$.

Удаляем цикл $c_{12}$, вектор $P_e(G) = <3,3,4,3,4,3,4,3,4,4,3,3,3,3,3>$.

Удаляем цикл $c_8$, вектор $P_e(G) = <3,3,4,2,3,3,3,3,4,4,3,3,3,2,2>$.

Удаляем цикл $c_{11}$, вектор $P_e(G) = <3,3,4,2,3,3,3,2,3,3,3,3,2,1,2>$.

Удаляем цикл $c_6$, вектор $P_e(G) = <3,2,3,2,3,3,3,1,2,3,3,2,1,1,2>$.

Удаляем цикл $c_4$, вектор $P_e(G) = <2,2,2,2,2,3,3,1,2,3,3,1,1,2,3>$.

Удаляем цикл $c_9$, вектор $P_e(G) = <2,2,2,2,2,2,2,1,1,3,2,0,1,2,3>$.

Удаляем цикл $c_2$, вектор $P_e(G) = <1,1,2,2,1,2,2,1,1,2,2,0,1,0,2>$.

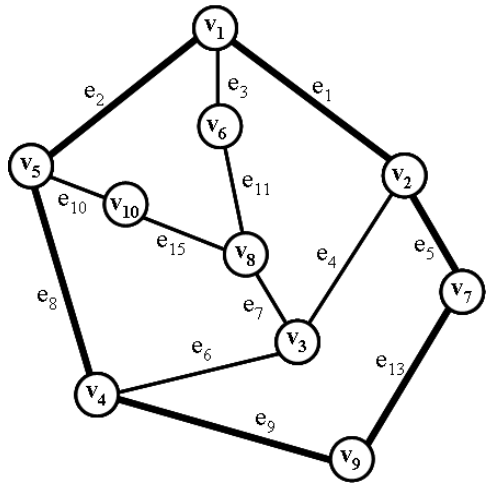

Рис. 15.22. Максимально плоский суграф.

Граф циклов H для максимально плоского суграфа сепарабелен, так как имеет в своем составе мультиребра. Следовательно, граф Петерсена не гамильтонов.

**Комментарии**

В данной главе рассматриваются негамильтоновы графы. Показано, что гамильтонов цикл может состоять из нескольких компонент связности, если в процессе удаления вершин с максимальной валентностью в графе циклов H появляются двусвязные суграфы. Либо граф циклов H является сепарабельным.

# ЗАКЛЮЧЕНИЕ

Произведен сравнительный анализ двух методов выделения базиса подпространства циклов C(G) с минимальным значением кубического функционала Маклейна из множества изометрических циклов графа. Показано, что метод наискорейшего спуска обладает полиномиальной вычислительной сложностью в отличие от случайного перебора вариантов метода генерации подмножества независимых циклов. Вычислительная сложность метода случайного неуправляемого выделения базиса близка к полному перебору всевозможных систем, состоящих из изометрических циклов с мощностью равной цикломатическому числу графа. И, таким образом, вычислительная сложность этого метода может быть оценена как экспоненциальная.

Приближенный процесс построения топологического рисунка максимально плоского суграфа разбит на три этапа:

- выделение из множества изометрических циклов базиса подпространства циклов C(G) с минимальным значением функционала Маклейна;

- выделение из базиса циклов подмножества циклов, удовлетворяющего нулевому значению кубического функционала Маклейна, методом исключения изометрических циклов, удовлетворяющих условию Эйлера;

- построение дополнительных простых циклов относительно координатно-базисной системы векторов, построенной на ободе плоской конфигурации.

Рассмотрена связь между максимальным плоским суграфом и построением гамильтонова цикла графа.

Представлен метод построения гамильтонова цикла.